\input amstex

\define\scrO{\Cal O}
\define\Pee{{\Bbb P}}
\define\Zee{{\Bbb Z}}
\define\Cee{{\Bbb C}}

\define\Spec{\operatorname{Spec}}
\define\Pic{\operatorname{Pic}}
\define\Ker{\operatorname{Ker}}

\define\Sym{\operatorname{Sym}}
\define\Hom{\operatorname{Hom}}

\define\Ext{\operatorname{Ext}}
\define\Id{\operatorname{Id}}
\define\ch{\operatorname{ch}}
\define\Todd{\operatorname{Todd}}
\define\ddual{\spcheck{}\spcheck}

\define\proof{\demo{Proof}}
\define\endproof{\qed\enddemo}
\define\endstatement{\endproclaim}
\define\theorem#1{\proclaim{Theorem #1}}
\define\lemma#1{\proclaim{Lemma #1}}
\define\proposition#1{\proclaim{Proposition #1}}
\define\corollary#1{\proclaim{Corollary #1}}
\define\claim#1{\proclaim{Claim #1}}

\define\section#1{\specialhead #1 \endspecialhead}
\define\ssection#1{\medskip\noindent{\it  #1}\medskip\nobreak}

\documentstyle{amsppt}

\topmatter
\title Vector bundles over elliptic fibrations
\endtitle
\author {Robert Friedman,  John W. Morgan and Edward Witten}
\endauthor
\address Department of Mathematics, Columbia University, New York, NY
10027, USA\endaddress
\email rf\@math.columbia.edu, jm\@math.columbia.edu  \endemail
\address School of Natural Sciences, Institute for Advanced Study,
Princeton NJ 08540 USA\endaddress
\email  witten\@ias.edu  \endemail
\thanks The first author was partially supported by NSF grant DMS-96-22681.
The second author was partially supported by NSF grant DMS-94-02988. The
third author was partially supported by NSF grant PHY-95-13835.
\endthanks

\endtopmatter

\document

\section{Introduction}

Let $\pi \: Z \to B$ be an elliptic fibration with a section. The goal of
this paper is to study holomorphic vector bundles over $Z$. We are mainly
concerned with  vector bundles $V$ with trivial determinant, or more
generally such that $\det V$ has trivial restriction to each fiber, so that
$\det V$ is the pullback of a line bundle on
$B$. (The case where $\det V$ has nonzero degree on every fiber is in a
certain sense simpler, since it usually reduces to the case considered here
for a bundle of smaller rank.) We give two constructions of vector bundles,
one based on the idea of a spectral cover of
$B$ and the other based on the idea of extensions of certain fixed bundles
over the elliptic manifold $Z$.  Each of these constructions has advantages
and the combination of the two seems to give the most comprehensive
information.

Vector bundles over a single elliptic curve were first classified by Atiyah
\cite{1}; however, he did not attempt to construct universal bundles or
work in families. The case of rank two bundles over an elliptic surface was
studied in
\cite{3}, \cite{6}, \cite{4} with a view toward making computations in
Donaldson theory. The motivation for this paper and the more general study
of the moduli of principal $G$-bundles over families of elliptic curves
(which will be treated in another paper) grew out  of questions arising in
the recent study of
$F$-theory by physicists.  The explanation of these connections was given in
\cite{7}. For these applications $Z$ is assumed to be a Calabi-Yau manifold,
usually of dimension two or three. However, most of the results on vector
bundles and more generally
$G$-bundles are true with no assumptions on $Z$. The case of a general
simple and simply connected complex Lie group $G$ involves a fair amount of
algebraic group theory and will be treated elsewhere, but the case
$G=SL_n(\Cee)$ can be done in a quite explicit and concrete way, and that
is the subject of this paper.

For both mathematical and physical reasons, we shall be primarily
interested in constructing stable vector bundles on $Z$. Of course,
stability must be defined with respect to a suitable ample divisor.
Following well-established principles, the natural ample divisors to work
with are those of the form
$H_0 + N\pi^*H$ for $N\gg 0$, where $H_0$ is some fixed ample divisor on
$Z$ and
$H$ is an ample divisor on $B$. If $V$ is stable with respect to such a
divisor, then
$V|f$ is semistable with respect to almost all fibers $f$. (However the
converse is not necessarily true.) One special feature of vector bundles
$V$ with trivial determinant on an elliptic curve is that, if the rank of
$V$ is at least two, then
$V$ is never properly stable. Moreover, if $V$ has rank $n>1$, then
$V$ is never simple; in fact, the endomorphism algebra  of $V$ has
dimension at least
$n$. But there is still a relative coarse moduli space
$\Cal M_{Z/B}$, which turns out to be a $\Pee^{n-1}$-bundle over $B$. A
stable vector bundle on $Z$ defines a rational section of $\Cal M_{Z/B}$.
Conversely, a regular section of $\Cal M_{Z/B}$ defines a vector bundle
over $Z$, and in fact it defines many such bundles. Our goal will be to
describe all such bundles, to see how the properties of the section are
reflected in the properties of these bundles, and to find sufficient
conditions for the bundles in question to be stable.

In the first three sections we consider a single (generalized) elliptic
curve $E$.  In Section 1 we construct a coarse moduli space for
$S$-equivalence classes of semistable $SL_n(\Cee)$-bundles over $E$. It is a
projective space $\Pee^{n-1}$, in fact it is the projective space of the
complete linear system $|np_0|$ where $p_0\in E$ is the origin of the group
law.  It turns out that each $S$-equivalence class of semistable bundles has
a ``best" representative, the so-called regular representative. The
defining property of these bundles, at least when $E$ is smooth, is that 
their automorphism groups are of the smallest possible dimension, namely
$n$.  We view them as analogues of regular elements in the group
$SL_n(\Cee)$.  The moduli space we construct  is also the coarse moduli
space for isomorphism classes of  regular semistable $SL_n(\Cee)$-bundles
over $E$. As we shall see, the regular bundles are the bundles which arise
if we try to fit together the
$S$-equivalence classes in order to find universal holomorphic bundles over
$\Pee^{n-1}\times E$.

In Section 2, assuming that $E$ is smooth, we construct  a tautological
bundle $U$ over
$\Pee^{n-1}\times E$ which is  regular semistable and with trivial
determinant on each slice $\{x\}\times E$ and such that
$U|\{x\}\times E$ corresponds to the regular bundle over $E$ whose
$S$-equivalence class is $x$.   There is not a unique such bundle over
$\Pee^{n-1}\times E$, and we proceed to construct all such.  The idea is
that there is an
$n$-sheeted covering $T\to \Pee^{n-1}$ called the {\sl spectral cover\/},
such that $U$ is obtained by pushing down a Poincar\'e line bundle 
${\Cal P}\to T\times E$ under the covering map.  It turns out that  every
bundle over $\Pee^{n-1}\times E$ which is  of the correct isomorphism class
on each slice $\{x\}\times E$ is obtained by pushing down ${\Cal P}\otimes
p_1^*M$ for some line bundle $M$ on $T$. There is a generalization of this
result to cover the case of families of regular semistable bundles on $E$
parameterized by arbitrary spaces $S$.

In Section 3 we turn to a different construction of ``universal" bundles
over $\Pee^{n-1}\times E$. Here we consider the space of extensions of two
fixed bundles with determinants $\scrO_E(\pm p_0)$. For a fixed rank $d$,
there is a unique stable bundle $W_d$ of rank $d$ such that $\det W_d \cong
\scrO_E(p_0)$. For
$1\leq d\leq n-1$, we consider the space of all nonsplit extensions $V$ of
the form
$$0 \to W_d\spcheck \to V \to W_{n-d} \to 0.$$ The moduli space of all such
extensions is simply $\Pee H^1(W_{n-d}\spcheck
\otimes W_d\spcheck) \cong \Pee^{n-1}$.  Over $\Pee^{n-1} \times E$ there
is a universal extension whose restriction to each fiber is regular
semistable.  There is thus an induced map from the $\Pee^{n-1}$ of
extensions to the coarse moduli space defined in Section 1, which is
$|np_0|\cong
\Pee^{n-1}$.  By a direct analysis we show that this map is an isomorphism. 
Actually, there are
$n-1$ different versions of this construction, depending on the choice of
the integer $d$, but the projective spaces that they produce are all
canonically identified. On the other hand, the universal extensions
associated with different versions of the construction are non-isomorphic
universal bundles. Finally, we relate these families of bundles  to the
ones arising from the spectral cover construction, which we can then extend
to the case where
$E$ is singular. We remark here that we can interpret the construction of
Section 3 as parametrizing those bundles whose structure group can be
reduced to a maximal parabolic subgroup
$P$ of $SL_n$, such that the induced bundle on the Levi factor is required
to be
$W_d\spcheck \oplus W_{n-d}$ in the obvious sense.  This interpretation can
then be generalized to other complex simple groups \cite{8}.

In Section 4 we generalize the results of the first three sections to a
family $\pi \: Z\to B$ of elliptic curves with a section $\sigma$. By taking
cohomology along the fibers of $\pi$, we produce a vector bundle over the
base, namely $\pi _*\scrO_Z(n\sigma)=\Cal V_n$,  which globalizes
$H^0(E,\scrO_E(np_0))$.  The associated projective bundle
$\Pee\pi _*\scrO_Z(n\sigma)=\Pee\Cal V_n$ then becomes the appropriate
relative coarse moduli space.  We show that
$\pi _*\scrO_Z(n\sigma)$ has a natural splitting as a direct sum of line
bundles. This decomposition is closely related to the fact that the
coefficients of the characteristic polynomial of an element in
$\frak{sl}_n$ are a polynomial basis for the algebra of  polynomial
functions on $\frak{sl}_n$ invariant under the adjoint action. Having
constructed  the relative coarse moduli space, we give a relative version
of the constructions of Sections 2 and 3 to produce  bundles over $\Pee\Cal
V_n\times _BZ$. The extension construction generalizes easily.  The bundles
we used over a single elliptic curve have natural extensions to any elliptic
fibration.  We form the relative extension bundle and the universal
relative extension  in direct analogy with the case of a single elliptic
curve. Relative versions of results from Section 3 show that the relative
extension space is identified with $\Pee\Cal V_n$. Following the pattern of
Section 3, we use the extension picture to define a universal spectral
cover of $\Pee\Cal V_n$, and in turn use this spectral cover to construct
new universal vector bundles. Finally, we calculate the Chern classes of
the universal bundles we have constructed.

In Section 5, using the theory developed in the first four sections, we
study vector bundles
$V$ over an elliptic fibration $\pi\: Z\to B$ such that the restriction of
$V$ to every fiber is regular and semistable. To such a bundle $V$, we
associate a section $A(V)$ of $\Pee\Cal V_n$ and a cover $C_A\to B$ of
degree $n$, the {\sl spectral cover\/} of $B$ determined by $V$.
Conversely, $V$ is determined by $A$ and by the choice of a line bundle on
$C_A$. After computing some determinants and Chern classes, we discuss the
possible line bundles which can exist on the spectral cover. Then we turn
to specific types of bundles. After describing symmetric bundles, which are
interesting from the point of view of $F$-theory, we turn to bundles
corresponding to a degenerate section. First we consider the most
degenerate case, and then we consider reducible sections where the
restriction of
$V$ to every fiber has a section. Finally, we relate reducible sections to
the existence of certain subbundles of $V$. 

In Section 6, we consider  bundles $V$ whose restriction to a generic fiber
is regular and semistable, but such that there exist fibers $E_b$ where
$V|E_b$ is either unstable or it is semistable but not regular. If $V$ fails
to be regular or semistable in codimension one, it can be improved by
elementary modifications to a reflexive sheaf whose restriction to every
fiber outside a codimension two set is regular and semistable. We describe
this process and, as an illustration, analyze the tangent bundle to an
elliptic surface. On the other hand, if the locus of bad fibers has
codimension at least two, no procedure exists for improving $V$, and we
must analyze it directly. The case of instability in codimension two or
higher corresponds to the case where the rational section $A$ determined by
$V$ does not actually define a regular section (this case can also lead to
reflexive but non-locally free sheaves). The case where $V$ has irregular
restriction to certain fibers in codimension at least two corresponds to
singular spectral covers. We give some examples of such behavior, without
trying to be definitive. Our construction can be viewed as a generalization
of the method of Section 3 to certain non-maximal parabolic subgroups of
$SL_n$.

Finally, we turn in Section 7 to the problem of deciding when the bundles
$V$ constructed by our methods are stable. This is the most interesting
case for both mathematical and physical reasons. While we do not try to
give necessary and sufficient conditions, we show that, in case the
spectral cover $C_A$ of $B$ determined by $V$ is irreducible, then $V$ is
stable with respect to all ample divisors of the form $H_0 + N\pi ^*H$,
where $H_0$ is an ample divisor on $Z$ and
$H$ is an ample divisor on $B$, and $N\gg 0$. We are only able to give an
effective bound on $N$ in case $\dim B =1$, i\.e\. $Z$ is an elliptic
surface, but it seems likely that such an effective bound exists in general.

We will have to deal systematically with singular fibers of $Z\to B$, and
the price that must be paid for analyzing this case is a heavy dose of
commutative algebra. In an attempt to make the paper more readable, we have
tried to isolate these arguments where possible. We  collect here some
preliminary definitions and technical results. While these results are
well-known, we could not find an adequate reference for many of them.

\ssection{Notation and conventions.}

All schemes are assumed to be separated and of finite type over $\Cee$. A
sheaf is always a coherent sheaf. We will identify a vector bundle with its
locally free sheaf of sections, covariantly. If $V$ is a vector bundle,
then $\Pee V$ is the projective space bundle whose associated sheaf of
graded algebras is $\bigoplus _{k\geq 0}\Sym ^kV\spcheck$; thus these
conventions are opposite to those of EGA or
\cite{10}. Given sheaves $\Cal S, \Cal S'$, we denote by $Hom(\Cal S, \Cal
S')$ the sheaf of homomorphisms from $\Cal S$ to $\Cal S'$ and by
$\Hom(\Cal S, \Cal S') = H^0(Hom(\Cal S, \Cal S'))$ the group of all such
homomorphisms. Likewise
$Ext ^k(\Cal S, \Cal S')$ is the Ext sheaf and $\Ext ^k(\Cal S, \Cal S')$
is the global Ext group (related to the local Ext groups by the local to
global spectral sequence).

\ssection{0.1. Elliptic curves and elliptic fibrations.}

Recall that a {\sl  Weierstrass equation\/} is a homogeneous cubic equation
of the form
$$Y^2Z=4X^3-g_2XZ^2-g_3Z^3,\tag{0.1}$$ with $g_2,g_3$  constants. We will
refer to the curve $E$ in $\Pee^2$ defined by such an equation, together
with the marked point $p_0=[0,1,0]$ at infinity, as a {\sl  Weierstrass
cubic}.  Setting  
$$\Delta(g_2,g_3)=g_2^3- 27g_3^2,$$  if $\Delta(g_2,g_3)\neq 0$, then (0.1)
defines a smooth cubic curve in $\Pee^2$ with the marked point
$[0,1,0]$, i\.e\., defines the structure of an elliptic curve. If
$\Delta(g_2,g_3) =0$, then the corresponding plane cubic $E$ is a singular
curve with arithmetic genus $p_a(E) =1$. If
$(g_2, g_3)$ is a smooth point of the locus $\Delta(g_2,g_3) =0$, then the
corresponding plane cubic curve is a rational curve with a single node.  The
smooth points of such a curve form a group isomorphic to $\Cee^*$ with
identity element $p_0$. The point $g_2= g_3=0$ is the unique singular point
of
$\Delta(g_2,g_3) =0$ and the corresponding plane curve  is a rational curve
with a single cusp. Once again its smooth points form a group, isomorphic
to $\Cee$, with identity element $p_0$.   These are all possible reduced and
irreducible curves of arithmetic genus one.

Next we consider the relative version of a Weierstrass equation. Let $\pi
\: Z \to B$ be a flat morphism of relative dimension one, such that the
general fiber is a smooth elliptic curve and all fibers are isomorphic to
reduced irreducible plane cubics. Here we will assume that
$B$ is a smooth variety (although the case of a complex manifold is
similar).  We shall always suppose that
$\pi$ has a section $\sigma$, i\.e\. there exists a divisor $\sigma$
contained in the smooth points of $Z$ such that $\pi|\sigma$ is an
isomorphism. Let $L = R^1\pi_*\scrO_Z \cong
\scrO_Z(-\sigma)|\sigma$, viewed as a line bundle on $B$. Then there are
sections
$G_2 \in H^0(B; L^{\otimes 4})$ and $G_3 \in H^0(B; L^{\otimes 6})$ such
that
$\Delta (G_2, G_3)\neq 0$ as a section of $L^{\otimes 12}$, and $Z$ is
isomorphic to the subvariety of $\Pee(\scrO_B\oplus L^2\oplus L^3)$ defined
by the Weierstrass equation $Y^2Z=4X^3-G_2XZ^2-G_3Z^3$. Conversely, given
the line bundle $L$ on $B$ and sections $G_2 \in H^0(B; L^{\otimes 4})$,
$G_3 \in H^0(B; L^{\otimes 6})$ such that
$\Delta (G_2, G_3)\neq 0$, the equation $Y^2Z=4X^3-G_2XZ^2-G_3Z^3$ defines a
hypersurface $Z$ in $\Pee(\scrO_B\oplus L^2\oplus L^3)$, such that the
projection to $B$ is a flat morphism whose fibers are reduced irreducible
plane curves, generically smooth. We will not need to assume that $Z$ is
smooth; it is always Gorenstein and the relative dualizing sheaf
$\omega_{Z/B}$ is isomorphic to $ L$. Thus, the dualizing sheaf
$\omega _Z$ is isomorphic to $\pi^*K_B\otimes L$.

Let us describe explicitly the case where the divisors associated to $G_2$
and
$G_3$ are smooth and meet transversally. This means in particular that if
$G_2$ and
$G_3$ are chosen generically, then $G_2^3 - 27G_3^2$ defines a section of
$L^{12}$. We shall denote by $\overline{\Gamma}$ the zero set of this
section. Then  $\overline{\Gamma}$ is smooth except where $G_2 = G_3 = 0$,
where it has singularities which are locally trivial families of cusps. The
fiber of $\pi$ over a smooth point of
$\overline{\Gamma}$ is a nodal plane cubic, and over a point where
$G_2 = G_3 = 0$ the fibers of $\pi$ are cusps. Let $\Gamma$ be the locus of
points where
$\pi$ is singular. Thus $\Gamma$ maps bijectively onto
$\overline{\Gamma}$. There are local analytic coordinates on $B$ so that,
near a cuspidal fiber $Z$ has the local equation $y^2 = x^3 + sx+t$. Here
$x,y$ are a set of fiber coordinates for $\Pee(\scrO_B\oplus L^2\oplus
L^3)$ away from the line at infinity and
$x,y,s,t$ form part of a set of local coordinates for $\Pee(\scrO_B\oplus
L^2\oplus L^3)$. Thus $x,y,s$ are coordinates for $Z$. The local equation
for
$\overline{\Gamma}$ is $4s^3 + 27 t^2 = 0$. The equations for the singular
point of the fiber over $\overline{\Gamma}$ are as follows:
$y=0, s = -3x^2, t= 2x^3$. In particular, $\Gamma$ is  smooth, and is the
normalization of $\overline{\Gamma}$.  The morphism from $Z$ to
$B$ is given locally by $(s,t)$, where $t= y^2-x^3-sx$.

\ssection{0.2. Rank one torsion free sheaves.}

Let $E$ be a singular Weierstrass cubic and let $E_{\text{reg}}$ be the set
of smooth points of $E$. The arithmetic genus $p_a(E)$ is one. We let $n\:
\tilde E \to E$ be the normalization map. The generalized Jacobian
$J(E)$ is the group of line bundles of degree zero on $E$, and (as in the
smooth case) is isomorphic to
$E_{\text{reg}}$ via the map $e\in E_{\text{reg}}\mapsto \scrO_E(e-p_0)$.
Just as we can compactify $E_{\text{reg}}$ to $E$ by adding the singular
point, we can compactify $J(E)$ to the {\sl compactified generalized
Jacobian\/} $\bar J(E)$, by adding the unique rank one torsion free sheaf
which is not locally free. Here a sheaf
$\Cal S$ over
$E$ is  {\sl torsion-free\/} if it has no nonzero sections which are
supported on a proper closed subset (i\.e\. a finite set). In particular,
the restriction of
$\Cal S$ to the smooth points of
$E$ is a vector bundle, and so has a well-defined rank, which we also call
the rank of
$\Cal S$. If
$\Cal S$ is a torsion-free sheaf on $E$ we let $\deg \Cal S = \chi (\Cal S)
+ (p_a(E)-1)(\operatorname{rank} \Cal S) = \chi (\Cal S)$. (This agrees
with the usual Riemann-Roch formula in case $E$ is smooth.) Thus the degree
of such sheaves is additive in  exact sequences, and if
$\Cal S'\subseteq \Cal  S$ such that the quotient is supported at a  finite
set of points, then $\deg \Cal S' \leq 
\deg \Cal S$ with equality if and only if $\Cal S' = \Cal S$. If $\Cal S$ is
torsion free and $V$ is locally free, then $\deg (V\otimes \Cal S)= (\deg
V)(\operatorname{rank}\Cal S) + (\deg \Cal S)(\operatorname{rank} V)$. To
see this, first use the fact that there is a filtration of $V$ by
subbundles whose successive quotients are line bundles, so by the
additivity of degree we can reduce to the case where  $V$ is a line bundle.
In this case, we may write $V= \scrO_E(d_1- d_2)$, where $d_1$ and $d_2$
are effective divisors supported on the smooth points of $E$, and then use 
the exact sequences
$$0 \to \Cal S \otimes \scrO_E(-d_2) \to \Cal S \to \Cal S\otimes
\scrO_{d_2} \to 0$$ and 
$$0 \to \Cal S \otimes \scrO_E(-d_2) \to \Cal S \otimes \scrO_E(d_1-d_2)\to
\Cal S\otimes \scrO_{d_1} \to 0,$$ together with the usual properties, to
conclude that $\deg (V\otimes \Cal S)= (\deg V)(\operatorname{rank}\Cal S)
+ (\deg \Cal S)$ in case $V$ is a line bundle. Thus we have established the
formula in general.

Next let us show that there is a unique torsion free rank one sheaf $\Cal
F$ which compactifies the generalized Jacobian.

\lemma{0.2} There is a unique rank one torsion free sheaf $\Cal F$ on $E$ of
degree zero which is not locally free. It satisfies:
\roster
\item"{(i)}"  $Hom(\Cal F, \Cal F) =n_*\scrO_{\tilde E}$.
\item"{(ii)}" $\Cal F\spcheck \cong \Cal F$.
\item"{(iii)}" For all line bundles $\lambda$ of degree zero, $Hom(\lambda,
\Cal F) = Hom (\Cal F, \lambda) = \Cal F$ and $\Hom(\lambda, \Cal F) = \Hom
(\Cal F, \lambda) = 0$. Likewise $\Ext^1(\Cal F, \lambda) =  \Ext
^1(\lambda, \Cal F) = 0$.
\endroster
\endstatement
\proof The first statement is essentially a local result. Let $R$ be the
local ring of
$E$ at the singular point and let $\tilde R$ be the normaliztion of $R$. If
locally
$\Cal F$ corresponds to the $R$-module $M$, let $\tilde M$ be the $\tilde
R$-module
$M\otimes _R\tilde R$ modulo torsion. Then by construction $\tilde M$ is a
torsion free rank one $\tilde R$-module, so that we may choose a $\tilde
R$-module isomorphism from $\tilde M$ to $\tilde R$. Since $M$ is torsion
free, the natural map from $M$ to $\tilde M \cong \tilde R$ is injective,
identifying $M$ as an
$R$-submodule of $\tilde R$ which generates $\tilde R$ as a $\tilde
R$-module. Thus $M$ contains a unit of $\tilde R$, which after a change of
basis we may assume to be $1$, and furthermore $M$ contains $R\cdot
1=R\subseteq \tilde R$. But since the singularity of $E$ is a node or a
cusp, $\ell(\tilde R/R) =1$, and so either $M=R$ or $M=\tilde R$. Note that
there are two isomorphic non-locally free $R$-modules of rank one: $\tilde
R$ and $\frak m$, where $\frak m$ is the maximal ideal of $R$. The ideal
$\frak m$ is the conductor of the extension
$\tilde R$ of $R$, and $\Hom_R(\tilde R, R) \cong \frak m$, where the
isomorphism is canonical.

By the above, every rank one torsion free sheaf on
$E$ is either a line bundle or of the form $n_*L$, where $L$ is a line
bundle on
$\tilde E$. Now $\tilde E\cong \Pee^1$, and $\deg n_*\scrO_{\Pee^1}(a) =
a+1$. Thus
$n_*\scrO_{\Pee^1}(-1)$ is the unique rank one torsion free sheaf on $E$ of
degree zero which is not locally free. Note that, if $\frak m_x$ is the
ideal sheaf of the singular point $x\in E$, then  $\deg
\frak m_x = -1$, by using the exact sequence
$$0 \to \frak m_x \to \scrO_E \to \Cee _x \to 0.$$ Thus $\frak
m_x=n_*\scrO_{\tilde E}(-2)$.

To see (i), note that 
$$Hom (n_*\scrO_{\Pee^1}(-1), n_*\scrO_{\Pee^1}(-1)) =
n_*Hom(n^*n_*\scrO_{\Pee^1}(-1), \scrO_{\Pee^1}(-1)).$$ Since $n$ is
finite, the natural map $n^*n_*\scrO_{\Pee^1}(-1) \to
\scrO_{\Pee^1}(-1)$ is surjective, and its kernel is torsion. Thus 
$$Hom(n^*n_*\scrO_{\Pee^1}(-1), \scrO_{\Pee^1}(-1)) \cong
Hom(\scrO_{\Pee^1}(-1), \scrO_{\Pee^1}(-1)) = \scrO_{\Pee^1}= \scrO_{\tilde
E},$$ proving (i). To see (ii), we have invariantly that 
$$Hom(n_*\scrO_{\tilde E},
\scrO_E) = \frak m_x=n_*\scrO_{\tilde E}(-2).$$ Thus tensoring with
$\scrO_E(-p_0)$ and using $n_*\scrO_{\tilde E}\otimes
\scrO_E(-p_0) = n_*n^*\scrO_E(-p_0) = n_*\scrO_{\tilde E}(-1)$ gives
$$\gather Hom(n_*\scrO_{\tilde E}(-1),
\scrO_E)= Hom(n_*\scrO_{\tilde E}\otimes
\scrO_E(-p_0),\scrO_E)= \\ =n_*\scrO_{\tilde E}(-2)\otimes
\scrO_E(p_0) = n_*\scrO_{\tilde E}(-1),
\endgather$$ which is the statement that $\Cal F\spcheck
\cong \Cal F$.  To see (iii), if $\lambda$ is a line bundle of degree zero,
then
$Hom(\lambda,
\Cal F) =\lambda^{-1}\otimes \Cal F$ is a nonlocally free sheaf of degree
zero, and hence it is isomorphic to $\Cal F$ by uniqueness. Likewise 
$$Hom (\Cal F,
\lambda) \cong \lambda \otimes \Cal F\spcheck \cong \lambda \otimes \Cal F
\cong
\Cal F.$$ Moreover, $\Hom(\lambda, \Cal F) = H^0(\Cal F) = 0$, since by
degree considerations a nonzero map $\lambda ^{-1}\to \Cal F$ would have to
be an isomorphism, contradicting the fact that $\Cal F$ is not locally
free. The proof that $\Hom (\Cal F,
\lambda)=0$ is similar. Now $\Ext^1(\Cal F, \lambda)$ is Serre dual to
$\Hom(\lambda,\Cal F)=0$, since $\lambda$ is locally free. Also, $\Ext
^1(\lambda,
\Cal F) = H^1(\lambda^{-1}\otimes \Cal F) = H^1(\Cal F) =0$, since
$h^0(\Cal F) =
\deg \Cal F =0$.
\endproof

\remark{Remark} In case $E$ is nodal, $\Ext ^1(\Cal F, \Cal F)$ is not
Serre dual to
$\Hom (\Cal F, \Cal F)$, and in fact $\Ext ^1(\Cal F, \Cal F)\cong H^0(Ext
^1(\Cal F, \Cal F))$ has dimension two. In this case $\Pee\Ext ^1(\Cal F,
\Cal F)\cong \Pee^1$ can be identified with the normalization of $E$. The
preimages
$\{x_1, x_2\}$ of the singular point give two different non-locally free
extensions, and the remaining locally free extensions $V$ of $\Cal F$ by
$\Cal F$ are parametrized by
$\Pee^1-\{x_1, x_2\}\cong \Cee^*$. The set of such $V$ is in
$1-1$ correspondence with $J(E) \cong E$ via the determinant.
\endremark
\medskip

Next we define the compactified generalized Jacobian of $E$. Let $\Delta
_0$ be the diagonal in $E\times E$ and let $I_{\Delta _0}$ be its ideal
sheaf. We let
$\scrO_{E\times E}(\Delta _0) = I_{\Delta _0}\spcheck$ and 
$$\Cal P_0 = \scrO_{E\times E}(\Delta _0-(E\times \{p_0\})) = I_{\Delta
_0}\spcheck\otimes \pi _2^*\scrO_E(-p_0).$$

\lemma{0.3} In the above notation,
\roster
\item"{(i)}" $\Cal P_0$ is flat over both factors of $E\times E$, and $\Cal
P_0\spcheck$ is locally isomorphic to $I_{\Delta _0}$.
\item"{(ii)}" If $e$ is a smooth point of $E$, the restriction of $\Cal
P_0$ to the slice $\{e\}\times E$ is $\scrO_E(e-p_0)$. If $x$ is the
singular point of $E$, the  restriction of $\Cal P_0$ to the slice
$\{x\}\times E$ is $\Cal F$.
\item"{(iii)}" Suppose that $S$ is a scheme and that $\Cal L$ is a coherent
sheaf on $S\times E$, flat over $S$, such that for every slice $\{s\}\times
E$, the restriction of $\Cal L$ to $\{s\}\times E$ is a rank one torsion
free sheaf on $E$ of degree zero. Then there exists a unique morphism $f\:
S \to E$ and a line bundle $M$ on $S$ such that $\Cal L = (f\times
\Id)^*\Cal P_0\otimes \pi _1^*M$.
\endroster
\endstatement 
\proof We shall just outline the proof of this essentially standard result.
The proofs of (i) and (ii) in case
$\Cal P_0$ is replaced by
$I_{\Delta _0}$, with the necessary changes in (ii), are easy: From the
exact sequence 
$$0 \to I_{\Delta _0} \to \scrO_{E\times E} \to \scrO_{\Delta _0} \to 0,$$
and the fact that both $E\times E$ and $\Delta _0$ are flat over each
factor, we see that $I_{\Delta _0}$ is flat over both factors, and the
restriction of
$I_{\Delta _0}$ to the slice $\{e\}\times E$ is $\scrO_E(-e)$, if $e\neq
x$, and is $\frak m_x$ in case $e=x$. To handle the case of $\Cal P_0$, the
main point is to check that $\scrO_{E\times E}(\Delta _0) = I_{\Delta
_0}\spcheck$ is locally isomorphic to $I_{\Delta _0}$, and that the
inclusion $I_{\Delta _0} \to
\scrO_{E\times E}$ dualizes to give an exact sequence
$$0 \to \scrO_{E\times E} \to \scrO_{E\times E}(\Delta _0) \to
\scrO_{\Delta _0}
\to 0.$$ This may be checked by hand, by working out a local resolution of
$I_{\Delta _0}$. We omit the details.

Another, less concrete, proof which generalizes to a flat family $\pi\: Z
\to B$ is as follows. Dualizing the inclusion of $I_{\Delta _0} \to
\scrO_{E\times E}$ gives an exact sequence
$$0 \to \scrO_{E\times E} \to I_{\Delta _0}\spcheck \to
Ext^1(\scrO_{\Delta_0},\scrO_{E\times E}) \to 0.$$ To check flatness of
$I_{\Delta _0}\spcheck$ and the remaining statements of (ii), it suffices
to show that, locally,  $Ext^1(\scrO_{\Delta_0},\scrO_{E\times E}) \cong
\scrO_{\Delta_0}$. Clearly $Ext^1(\scrO_{\Delta_0},\scrO_{E\times E})$ is a
sheaf of $\scrO_{\Delta_0}$-modules and thus it is identified with a sheaf
on $E$ via the first projection. If
$\pi_1, \pi_2\: E\times E\to E$ are the projections, we have the relative
Ext sheaves $Ext^i_{\pi_1}(\scrO_{\Delta_0},\scrO_{E\times E})$. (See for
example
\cite{2} for properties of these sheaves.) The curve $E$ is Gorenstein and
thus
$\Ext^1(\Cee_x, \scrO_E) \cong \Cee$ for all $x\in E$. By base change,
$Ext^1_{\pi_1}(\scrO_{\Delta_0},\scrO_{E\times E})$ is a line bundle on
$E$. On the other hand, by the local to global spectral sequence, 
$$Ext^1_{\pi_1}(\scrO_{\Delta_0},\scrO_{E\times E})=
\pi_1{}_*Ext^1(\scrO_{\Delta_0},\scrO_{E\times E}).$$ Thus
$Ext^1(\scrO_{\Delta_0},\scrO_{E\times E})$ can be identified with a line
bundle on $\Delta _0$, and so it is locally isomorphic to
$\scrO_{\Delta_0}$. Dualizing this argument gives an exact sequence
(locally)
$$0 \to \Cal P_0\spcheck \to \scrO_{E\times E} \to
Ext^1(\scrO_{\Delta_0},\scrO_{E\times E}) \to 0,$$ and so (locally again)
$\Cal P_0\spcheck = (I_{\Delta _0})\ddual \cong I_{\Delta _0}$. In
particular $\Cal P_0\spcheck$ is also flat over $E$.

To see (iii), suppose that $S$ and $\Cal L$ are as in (iii). By base change,
$\pi_1{}_*(\Cal L\otimes\pi_2^*\scrO_E(p_0)) = M^{-1}$ is a line bundle on
$S$, and the morphism
$$\pi _1^* \pi_1{}_*(\Cal L\otimes\pi_2^*\scrO_E(p_0)) =\pi _1^*M ^{-1}\to
\Cal L\otimes\pi_2^*\scrO_E(p_0)$$ vanishes along a subscheme $\Cal Z$ of
$S\times E$, flat over $S$ and of degree one on every slice. Thus
$\Cal Z$ corresponds to a morphism $f\: S\to E=  \operatorname{Hilb}^1E$,
such that
$\Cal Z$ is the pullback of
$\Delta_0\subset E\times E$ by $(f\times \Id)^*$. This proves (iii).
\endproof

A very similar argument proves the corresponding result for the dual of the
ideal of the diagonal in $Z\times _BZ$, where $\pi\: Z \to B$ is a flat
family of Weierstrass cubics. In this case, we let $\Delta _0$ be the ideal
sheaf of the diagonal in $Z\times _BZ$, and set $\Cal P_0 = I_{\Delta
_0}\spcheck \otimes \pi _2^*\scrO_Z(-\sigma)$, where $\sigma$ is the
section. Then $\Cal P_0$ is flat over both factors $Z$, and has the
properties (i)--(iii) of (0.3). We leave the details of the formulation and
the proof to the reader.

Finally we discuss a local result which will be needed to handle semistable
sheaves on a singular $E$. (In the application, $R$ is the local ring of
$E$ at a singular point.)

\lemma{0.4} Let $R$ be a local Cohen-Macaulay domain of dimension one and
let $Q$ be a finitely generated torsion free $R$-module. Then $\Ext^1_R(Q,
R) = 0$.
\endstatement
\proof By a standard argument, if $Q$ has rank $n$ there exists an inclusion
$Q\subseteq R^n$. Thus necessarily the quotient $R^n/Q$ is a torsion
$R$-module
$T$. Now $\Ext^1_R(Q, R) \cong \Ext^2_R(T, R)$. Since $R$ is
Cohen-Macaulay, if
$\frak m$ is the maximal ideal of $R$, then $\Ext^2_R(R/\frak m, R) = 0$. An
induction on the length of $T$ then shows that $\Ext^2_R(T, R)=0$ for all
$R$-modules $T$ of finite length. Hence $\Ext^1_R(Q, R) = 0$.
\endproof

\ssection{0.3. Semistable bundles and sheaves on singular curves.}

Let $E$ be a Weierstrass cubic and let $\Cal S$ be a torsion free sheaf on
$E$.  The {\sl normalized degree\/} or {\sl slope\/}
$\mu (\Cal S)$ of $\Cal S$ is defined to be $\deg \Cal
S/\operatorname{rank} \Cal S$. A torsion free sheaf 
$\Cal S$ is {\sl semistable\/} if, for every subsheaf $\Cal S'$ of
$\Cal S$ with
$0< \operatorname{rank} \Cal S' <\operatorname{rank}  \Cal S$, then we have
$\mu (\Cal S') \leq \mu  (\Cal S)$,  and  it is {\sl unstable\/} if it is
not semistable. Equivalently, $\Cal S$ is semistable if, for all
surjections $\Cal S
\to \Cal S''$, where $\Cal S''$ is torsion free and nonzero, we have $\mu
(\Cal S'') \geq \mu (\Cal S)$. A torsion free rank one sheaf is semistable.
Given an exact sequence 
$$0 \to \Cal S' \to \Cal S \to \Cal S'' \to 0,$$ with $\mu(\Cal S') =\mu
(\Cal S) = \mu (\Cal S'')$, $\Cal S$ is semistable if and only if both
$\Cal S'$ and $\Cal S''$ are semistable. If
$\Cal S$ is a torsion free semistable sheaf of negative degree, then (for
$E$ of arithmetic genus one)
$h^0(\Cal S) = 0$ and hence
$h^1(\Cal S) = -\deg \Cal S$, and if
$\Cal S$ is  a torsion free semistable sheaf of strictly positive degree,
then since $h^1(\Cal S)$ is dual to
$\Hom (\Cal S, \scrO_E)$, it follows that $h^1(\Cal S) = 0$ and that
$h^0(\Cal S) = \deg \Cal S$. Every torsion free sheaf
$\Cal S$ has a canonical Harder-Narasimhan filtration, in other words a 
filtration by subsheaves
$F^0 \subset F^1\subset
\cdots$ such that $F^{i+1}/F^i$ is torsion free and semistable and $\mu
(F^i/F^{i-1}) > \mu (F^{i+1}/F^i)$ for all $i\geq 1$.

\definition{Definition 0.5} Let $V$ and $V'$ be two semistable torsion free
sheaves on $E$. We say that $V$ and $V'$ are {\sl $S$-equivalent\/} if
there exists a connected  scheme $S$ and a coherent sheaf $\Cal V$ on
$S\times E$, flat over $S$, and a point $s'\in S$ such that $V\cong \Cal
V|\{s\}\times E$ if $s\neq s'$ and
$V'\cong
\Cal V|\{s'\}\times E$.  We define {\sl
$S$-equivalence\/} to be the equivalence relation on semistable torsion free
sheaves generated by the above relation.  Suppose that $V$ and $V'$ are two
semistable bundles on $E$. We say that $V$ and $V'$ are {\sl restricted
$S$-equivalent\/} if there exists a connected  scheme $S$, a vector bundle
$\Cal V$ on $S\times E$, and a point $s'\in S$ such that
$V\cong \Cal V|\{s\}\times E$ if
$s\neq s'$ and
$V'\cong
\Cal V|\{s'\}\times E$.  We define {\sl restricted 
$S$-equivalence\/} to be the equivalence relation on semistable bundles
generated by the above relation.
\enddefinition

\section{1. A coarse moduli space for semistable bundles over a Weierstrass
cubic.}

Fix a Weierstrass cubic $E$ with an origin
$p_0$ and consider semistable vector bundles of rank $n$ and trivial
determiniant over $E$.  Our goal in this section will be to construct a
coarse moduli space of such bundles, which we will identify with the linear
system
$|np_0|$. Given a vector bundle $V$, we associate to $V$ a point
$\zeta(V)$ in the projective space $|np_0|$ associated to the linear system
$\scrO_E(np_0)$ on
$E$. In case $E$ is smooth, $\zeta(V)$ records the unordered set of degree
zero line bundles that occur as Jordan-H\"older quotients of any maximal
filtration of
$V$. More generally, if
${\Cal V}\to S\times E$ is an algebraic (or holomorphic) family of bundles
of the above type on $E$, then the function $\Phi\colon S\to |np_0|$
defined by $\Phi(s)=\zeta\left({\Cal V}|\{s\}\times E\right)$  is  a
morphism. If $E$ is smooth, two semistable bundles $V$ and $V'$ are are
$S$-equivalent if and only $\zeta (V) =\zeta (V')$. This identifies 
$|np_0|$ as a (coarse) moduli space of $S$-equivalence classes of
semistable rank $n$ bundles with trivial determinant on $E$. A similar
result holds if $E$ is cuspidal. In case
$E$ is nodal, however, there exist $S$-equivalent bundles $V$ and $V'$ such
that
$\zeta (V)
\neq\zeta (V')$. It seems likely that, in case $E$ is nodal,  $\zeta (V)
=\zeta (V')$ if and only if $V$ and $V'$ are restricted $S$-equivalent
(0.5).

The moduli space $|np_0|$ is not a fine moduli space, for two reasons. One
problem is the issue of $S$-equivalence versus isomorphism. To deal with
this problem, we will attempt to choose a ``best" representative for each
$S$-equivalence class, the regular representative. In case $E$ is smooth,
a  regular bundle $V$ is one whose automorphism group has dimension equal
to its rank, the minimum possible dimension. Even after choosing the
regular representative, however, $|np_0|$ fails to be a fine moduli space
because the bundles $V$ are never simple. This allows us to twist universal
bundles by line bundles on an $n$-sheeted cover of $|np_0|$, the {\sl
spectral cover}. This construction will be described in Section 2.

\ssection{1.1.   The Jordan-H\"older constituents of a semistable bundle.} 

The two main results of this section are the following:

\theorem{1.1} Let $V$ be a semistable torsion free sheaf of rank $n$ and
degree zero over $E$. Then  $V$ has a  Jordan-H\"older filtration 
$$0\subset F^0\subset F^1\subset \cdots\subset F^n=V$$ so that each
quotient  $F^i/F^{i-1}$ is a rank one torsion free sheaf of degree zero. For
$\lambda$ a rank one torsion free sheaf of degree zero, define $V(\lambda)$
to be the sum of all the subsheaves of $V$ which have a filtration such
that all of the successive quotients are isomorphic to $\lambda$. Then $V
=\bigoplus _\lambda V(\lambda)$. In particular, if $V$ is locally free,
then $V(\lambda)$ is locally free for every $\lambda$.
\endstatement

\theorem{1.2} Let $V$ be a semistable torsion free sheaf of rank $n$ and 
degree zero on $E$. Then 
\roster
\item"{(i)}" $h^0(V\otimes 
\scrO_E(p_0))=n$ and   the natural evaluation map
$$ev\colon H^0(V\otimes  \scrO_E(p_0))\otimes_\Cee \scrO_E\to V\otimes
\scrO_E(p_0)$$  is  an isomorphism over the generic point of $E$. 
\item"{(ii)}" Suppose that $V$ is locally free with $\det V
=\scrO_E(e-p_0)$. The induced map on determinants defines a map
$$\wedge ^nev\: \det H^0(V\otimes \scrO_E(p_0))\otimes_\Cee \scrO_E \cong
\scrO_E\to \det \left(V\otimes \scrO_E(p_0)\right) \cong
\scrO_E((n-1)p_0+e).$$  Thus $\wedge ^nev$ defines  a non-zero section of
$\scrO_E((n-1)p_0+e)$ up to a nonzero scalar multiple, i\.e\. a point of
$|(n-1)p_0+e|$. We denote this element by  $\zeta(V)$. In particular, if
$e=p_0$, then $\zeta (V)\in |np_0|$.
\endroster
\endstatement
\demo{Proof of Theorem \rom{1.2}} Let $V$ be a semistable sheaf of degree
zero and rank $n$ on $E$. The degree of  $V\otimes 
\scrO_E(p_0)$ is $n$. By definition, $h^0(V\otimes  \scrO_E(p_0)) -
h^1(V\otimes  \scrO_E(p_0)) = n$. By Serre duality,
$h^1(V\otimes  \scrO_E(p_0)) =\dim \Hom(V, \scrO_E(-p_0))$. Since $V$ is
semistable, $\Hom(V, \scrO_E(-p_0)) =0$. Thus
$h^0(V\otimes  \scrO_E(p_0)) =  n$. Next we claim that the induced map
$ev\colon H^0(V\otimes  \scrO_E(p_0))\otimes_\Cee \scrO_E
\to V\otimes \scrO_E(p_0)$ is an isomorphism over the generic point of $E$;
equivalently, its image $I\subset V\otimes \scrO_E(p_0)$ has rank $n$.  To
prove this, we use the following lemma.

\lemma{1.3} Let $E$ be a Weierstrass cubic,  let $I$ be a torsion free
sheaf on
$E$ and let
$\mu _0(I)$ be the maximal value of $\mu (J)$ as $J$ runs over all torsion
free subsheaves of $I$. Then
$$h^0(I) \leq \max(\mu _0(I), 1)\operatorname{rank}I.$$
\endstatement
\proof If $0 \subset F^0 \subset \cdots \subset F^k = I$ is the
Harder-Narasimhan filtration of $I$, then $\mu _0(I) = \mu (F^0)$,
$F^{i+1}/F^i$ is semistable, and
$\mu (F^{i+1}/F^i) < \mu _0(I)$ for all $i \geq 1$. Furthermore
$$h^0(I) \leq \sum _ih^0(F^{i+1}/F^i).$$ Now if $\mu (F^{i+1}/F^i) > 0$,
then since $h^1(F^{i+1}/F^i)= \dim \Hom (F^{i+1}/F^i, \scrO_E) = 0$, it
follows that $h^0(F^{i+1}/F^i) = \deg (F^{i+1}/F^i) = \mu
(F^{i+1}/F^i)\cdot \operatorname{rank}(F^{i+1}/F^i) \leq \mu_0(I)
\operatorname{rank}(F^{i+1}/F^i)$. If $\mu (F^{i+1}/F^i) < 0$, then
$$h^0(F^{i+1}/F^i) = 0 \leq
\operatorname{rank}(F^{i+1}/F^i).$$ There remains the case that $\mu
(F^{i+1}/F^i) = 0$. In this case, we claim that
$$h^0(F^{i+1}/F^i) \leq \operatorname{rank}(F^{i+1}/F^i).$$ In fact since
$F^{i+1}/F^i$ is semistable, this follows from the next claim.

\lemma{1.4} If $V$ is a semistable torsion free sheaf on $E$ with $\mu (V) =
0$, then
$h^0(V) \leq \operatorname{rank} V.$
\endstatement

\proof Argue by induction on $\operatorname{rank}V$. If
$\operatorname{rank} V = 1$ and $h^0(V) \geq 1$, then there exists a
nonzero map $\scrO_E\to V$, and  since
$\mu (\scrO_E) = \mu (V)$, this map must be an isomorphism. Thus
$h^0(V) = 1$.  In general, if $\operatorname{rank}V = n+1$ and
$h^0(V) \neq 0$, choose a nonzero map
$\scrO_E \to V$. Since $V$ is semistable, the cokernel $Q$ of this map is 
torsion free and thus is also semistable, with $\mu (Q) =0$. Since the rank
of $Q$ is $n$, by induction we have $h^0(V) \leq 1 + h^0(Q) \leq n+1$.
\endproof

Returning to the proof of (1.3), we see that in all cases
$$h^0(F^{i+1}/F^i) \leq \max(\mu _0(I),
1)\operatorname{rank}(F^{i+1}/F^i).$$  Summing over $i$ gives the statement
of (1.3).
\endproof

We continue with the proof of Theorem 1.2. There is the map
$$ev\colon H^0(V\otimes \scrO_E(p_0))\otimes_\Cee\scrO_E \to V\otimes
\scrO_E(p_0).$$ Let $I$ be its image. By construction
$I$ is  a subsheaf of a locally free sheaf and hence is torsion free. Also,
by construction the map $H^0(I) \to H^0(V\otimes \scrO_E(p_0))$ is an
isomorphism, and thus 
$h^0(I) = n$. Since $V\otimes  
\scrO_E(p_0)$ is semistable and
$\mu (V\otimes  \scrO_E(p_0)) = 1$, we have $\mu _0(I) \leq 1$. Thus, by
(1.3), $n = h^0(I) \leq \operatorname{rank} I \leq n$, and so
$\operatorname{rank} I = n$. Equivalently, the image of  $ev$ is equal to
$V\otimes  \scrO_E(p_0)$ at the generic point.  From this, the remaining
statements in Theorem 1.2 are clear.
\endproof

\demo{Proof of Theorem \rom{1.1}} Let us first show that $V$ has a
Jordan-H\"older filtration as described. The proof is by induction on the
rank
$n$ of
$V$. If $n=1$, there is nothing to prove. For arbitrary $n$, we shall show
that there exists a nonzero map $\lambda \to V$, where $\lambda$ is a rank
one torsion free sheaf of degree at least zero. By semistability, the
degree of
$\lambda$ is exactly zero and $V/\lambda$ is torsion free. We can then apply
induction to $V/\lambda$. 

The proof of Theorem 1.2 above shows that, if
$V$ is a semistable torsion free sheaf of rank $n$ and degree zero, then
there is an injective map
$\scrO_E^{\oplus n}
\to V\otimes \scrO_E(p_0)$ whose image has rank $n$. Thus there is a map
$\scrO_E(-p_0)^{\oplus n} \to V$ whose image has rank $n$. The cokernel of
this map must be a torsion sheaf $\tau$. Note that, in case $V$ is locally
free,
$\tau$ is supported exactly at the points in the support of $\zeta (V)$.
Since
$\deg V=0$, $\tau \neq 0$. Choose a point $x$ in the support of $\tau$. If
$R$ is the local ring of $E$ at $x$ and $\frak m$ is the maximal ideal of
$x$, then $\tau_x$ is annihilated by some power of $\frak m$. Let $k$ be
such that $\frak m^k\tau \neq 0$ but $\frak m^{k+1}\tau =0$. Choosing a
section of
$\frak m^k\tau$ produces  a subsheaf $\tau_0$ of $\tau$ which is isomorphic
to
$\Cee_x$, in other words is isomorphic to $R/\frak m$ as an $R$-module. 

Let $V_0\subseteq V$ be the inverse image of $\tau _0$. Then $V_0$
corresponds to an extension of $\Cee_x$ by $\scrO_E(-p_0)^{\oplus n}$, and
hence to an extension class in 
$$\Ext^1(\Cee_x, \scrO_E(-p_0)^{\oplus n}) \cong H^0(Ext^1(\Cee_x,
\scrO_E(-p_0)^{\oplus n}) \cong \Ext^1_R(R/\frak m, R^n).$$ The ring $R$ is
a Gorenstein local ring of dimension one, and so $\Ext^1_R(R/\frak m, R)
\cong \Cee$. (Of course, this could be verified directly for the local
rings $R$ under consideration.) In fact, if $x$ is a smooth point of $E$
and $t$ is a local parameter at $x$, then the unique nontrivial extension
of $R/\frak m$ by $R$ corresponds to the exact sequence
$$0 \to R @>{\times t}>> R \to R/\frak m \to 0,$$ whereas if $x$ is a
singular point then the nontrivial extension is given by
$$0 \to R \to \tilde R \to R/\frak m \to 0.$$ Let $\xi$ be the extension
class corresponding to $V_0$ in 
$$\Ext^1(\Cee_x,
\scrO_E(-p_0)^{\oplus n}) \cong \Ext^1_R(R/\frak m, R^n) \cong \Cee^n.$$ In
the local setting, let $M$ be the $R$-module corresponding to $V_0$, and
suppose that we are given an extension
$$0\to R \to N \to R/\frak m \to 0,$$ with a corresponding extension class
$\eta \in \Ext^1_R(R/\frak m, R)$ and a homomorphism $f\: R\to R^n$ such
that $f_*(\eta) =\xi$. By a standard result, there is a homomorphism $N\to
M$ lifting $f$, viewed as a homomorphism
$R\to M$. In particular, this says that the image of $R$ in $M$ is
contained in a strictly larger rank one torsion free $R$-module. 

Returning to the global situation, let $\lambda$ be the unique nontrivial
extension of $\Cee_x$ by $\scrO_E(-p_0)$, and let $\eta$ be the
corresponding  extension class, well-defined up to a nonzero scalar. Thus
$\lambda$ is a rank one torsion free sheaf of degree zero. Since
$Hom(\scrO_E(-p_0),
\scrO_E(-p_0)^{\oplus n})$ is generated by its global sections, there
exists a homomorphism 
$f\: \scrO_E(-p_0)\to \scrO_E(-p_0)^{\oplus n}$ such that the image of
$\eta$ under $f_*$ in 
$\Ext^1(\Cee_x, \scrO_E(-p_0)^{\oplus n})$ is $\xi$. Then the inclusion
$\scrO_E(-p_0)\to \scrO_E(-p_0)^{\oplus n} \to V_0 \to V$ factors through a
nonzero map $\lambda \to V$, necessarily an inclusion with torsion free
cokernel. Thus we have proved the existence of the Jordan-H\"older
filtration by induction.

By using the  fact that $\Ext^1(\lambda, \lambda') =0$ if
$\lambda\neq \lambda'$, an easy argument left to the reader shows that
$V(\lambda)\neq 0$ if and only if $\Hom (V,\lambda) \neq 0$ if and only if
$\Hom (\lambda, V)\neq 0$. Thus we can always arrange that, if $\lambda$ is
a sheaf appearing as one of the quotients in Theorem 1.1, then there exists
a filtration for which $\lambda =F^0$ is the first such sheaf which
appears, and also one for which
$\lambda = F^n/F^{n-1}$ is the last such sheaf which appears.

Fix a rank one torsion free sheaf $\lambda$ of degree zero, and let
$V'(\lambda)$ be the sum of all subsheaves of $V$ which have a filtration
by rank one torsion free sheaves of degree zero which are not isomorphic to
$\lambda$. Let $V(\lambda) = V/V'(\lambda)$. Clearly $V(\lambda)$ is a
torsion free semistable sheaf, such that all of the quotients in a
Jordan-H\"older filtration of $V$ are isomorphic to $\lambda$. Again using
$\Ext^1(\lambda, \lambda') =0$ if
$\lambda\neq \lambda'$, one checks that $\Ext^1(V(\lambda), V'(\lambda))
=0$. Thus, by induction on the rank, $V$ is isomorphic to the direct sum of
the $V(\lambda)$.  This concludes the proof of Theorem 1.1.
\endproof

The construction of Theorem 1.2 works well in families.

\theorem{1.5}  Let $E$ be a Weierstrass cubic, and let $S$ be a scheme or
analytic space. Let $\Cal V$ be a rank $n$ vector bundle over $S\times E$
such that on each slice $\{s\}\times E$,
$\Cal V$ restricts to a semistable vector bundle $V_s$ of trivial
determinant. Then there exists a morphism $\Phi\: S \to |np_0| = \Pee
^{n-1}$ such that, for all
$s\in S$, we have $\Phi(s)=\zeta(V_s)$.  In particular, if $V$ and $V'$ are
restricted $S$-equivalent, then $\zeta (V) =\zeta (V')$.
\endstatement

\proof Let $p_1, p_2$ be the projections from $S\times E$ to $S$ and
$E$. To construct a morphism from $S$ to $|np_0|$ we shall construct a
homomorphism $\Psi\colon p_1^*L_0 \to p_1^*L_1\otimes p_2^*\scrO_E(np_0)$,
where $L_0,L_1$ are line  bundles  on $S$, with the property that the 
restriction of
$\Psi$ to  each slice $\{s\}\times E$ determines a  nonzero section of
$\scrO_E(np_0)$ (which is thus well-defined mod scalars), agreeing with
$\wedge^nev$. The map $\Psi$ is defined in the next  lemma.

\lemma{1.6}  The sheaf
$p_1{}_*(\Cal V \otimes p_2^*\scrO_E(p_0))$ is a locally free sheaf of rank
$n$ on $S$. Let $L_0$ be its determinant line bundle.  If $\hat\Psi\colon
p_1^*p_1{}_*(\Cal V
\otimes p_2^*\scrO_E(p_0))\to\Cal V\otimes p_2^*\scrO_E(p_0)$ is the
natural evaluation map, then its restriction to each slice
$\{s\}\times E$ is generically an isomorphism, agreeing with $ev$. Thus
$$\Psi=\det\hat\Psi\colon p_1^*L_0\to \det{\Cal V}\otimes
p_2^*\scrO_E(np_0)$$ has the property that its restriction to each slice
$\{s\}\times E$ is is nonzero and agrees with $\wedge^nev$.
\endstatement

\proof It follows from Theorem  1.2 that, if $\Cal V_s$ is the restriction
of $\Cal V$ to the slice $\{s\}\times E$, then 
$h^0({\Cal V}_s\otimes  \scrO_E(p_0))=n$ is independent of $s$.  Standard
base change arguments \cite{10, Theorem 12.11, pp\. 290--291} show that,
even if
$S$ is nonreduced,
$p_1{}_*(\Cal V
\otimes  p_2^*\scrO_E(p_0))$ is a locally free sheaf of rank $n$ on $S$,
and the natural map $p_1{}_*(\Cal V \otimes p_2^*\scrO_E(p_0))_s\to
H^0(V_s\otimes \scrO_E(p_0))$ is an isomorphism for every $s\in S$. Thus the
induced morphism $\hat\Psi\: p_1^*p_1{}_*(\Cal V \otimes p_2^*\scrO_E(p_0))
\to
\Cal V \otimes p_2^*\scrO_E(p_0)$ is a morphism between two vector bundles
of  rank $n$. Let $V_s$ be the restriction of $\Cal V$ to the slice
$\{s\}\times E$. Again by base change, the natural map $p_1{}_*(\Cal V
\otimes p_2^*\scrO_E(p_0))_s \to H^0(V_s\otimes 
\scrO_E(p_0))\otimes_\Cee\scrO_E$  is surjective. The result is now
immediate from Theorem 1.2.
\endproof

Next notice that since for every $s\in S$, $\det{\Cal V}|(\{s\}\times E)$
is trivial, it follows that
$\det {\Cal V}$ is isomorphic to $p_1^*L_1$ for some line bundle $L_1$ on
$S$. To complete the proof of Theorem 1.5 we need to check that the section
$\Psi(s)=\zeta(V_s)$ for all $s\in S$.  This is immediate from the
corresponding statement in Theorem 1.2.
\endproof

In fancier terms, Theorem 1.5 says that there is a morphism of functors
from the deformation functor of semistable vector bundles of rank $n$ and
trivial determinant on $E$ to the functor represented by the scheme
$|np_0|$. In general, this morphism is far from smooth; for example, at the
trivial bundle $\scrO_E^n$, the derivative of the morphism is identically
zero. However, if we restrict to regular semistable bundles (to be defined
in \S 1.2 below), then it will follow from (v) in Theorem 3.2 that the
derivative is always an isomorphism.

The sheaves $\lambda$ which appear as successive quotients of $V$ in
Theorem 1.1 are the {\sl Jordan-H\"older quotients\/} or {\sl
Jordan-H\"older constituents\/} of $V$. They appear with multiplicities and
the multiplicity of $\lambda$ in $V$ is independent of the choice of the
filtration. The summands $V(\lambda)$ of $V$ are canonically defined. It is
easy to see from the construction that $\zeta (V) = \sum _\lambda \zeta
(V(\lambda))$. More generally, $\zeta$ is additive over exact sequences of
semistable vector bundles of degree zero. Also, if
$\det V = \scrO_E(e-p_0)$ and $e'$ is a smooth point of $E$, then $e'$ lies
in the support of $\zeta (V)$ as a divisor in
$|(n-1)p_0 +e|$ if and only if
$\lambda =\scrO_E(e'-p_0)$ is a Jordan-H\"older constituent of $V$. Thus,
if the rank of $V(\lambda)$ is $d_\lambda$, then 
$$\zeta (V) = \sum _{\lambda \neq \Cal F}d_\lambda e_\lambda + e_{\Cal
F},$$ where $\lambda \cong \scrO_E(e_\lambda -p_0)$ and $e_{\Cal F}$ is a
divisor of degree
$d_{\Cal F}$ supported at the singular point of $E$. In this way we can
associate a point of the
$n^{\text{th}}$ symmetric product of
$E$ with such a $V$: namely
$$\sum_{\lambda\neq \Cal F}
\operatorname{rank}(V(\lambda))e_\lambda+d_{\Cal F}\cdot s$$ where
$\lambda\cong\scrO_E(e_\lambda-p_0)$ and $s\in E$ is the singular point.
Note that there is a morphism $|np_0| \to \Sym ^nE$, which is a closed
embedding if $E$ is smooth, or more generally away from the elements of
$|np_0|$ whose support meets the singular point of $E$.

Suppose that $E$ is smooth. Since  a degree zero line bundle on $E$ is
identified with a point of
$E$ via the correspondence $\lambda \mapsto q$ if $\lambda\cong
\scrO_E(q-p_0)$, the map which assigns to a semistable bundle
$V$ the unordered $n$-tuple of its Jordan-H\"older quotients, including
multiplicities, is the same as the map assigning to 
$V$ an unordered $n$-tuple $\zeta(V)$ of points of $E$, i.e., a point 
$$\zeta(V)\in\underbrace{(E\times \cdots \times E)}_{\text{$n$
times}}/{\frak S}_n,$$  where ${\frak S}_n$ is the symmetric group on $n$
letters. If $\zeta(V)=(e_1,\ldots,e_n)$, then the condition that the
determinant of $V$ is trivial means that
$\sum_{i=1}^ne_i=0$ in the  group law of $E$, or equivalently that the
divsior $\sum_{i=1}^ne_i$ is linearly equivalent to $np_0$. Thus, the
unordered $n$-tuple
$(e_1,\ldots,e_n)$ associated to
$V$ can be identified with a point in the complete linear system $|np_0|$,
and this point is exactly $\zeta (V)$.

An important difference in case $E$ is singular is that, while a point of
$|np_0|$ determines a  point on the symmetric $n$-fold product of $E$, in
general it contains more information at the singular point than just its
multiplicity. Thus, the function $\Phi$ should be viewed not as a point in
the $n$-fold symmetric product but as a point in the linear system
$|np_0|$. For example, if $E$ is nodal and $n>2$, then an element of
$|np_0|$ supported entirely at the singular point corresponds to a
hyperplane in $\Pee^{n-1}$ meeting the image of $E$ embedded by the complete
linear system $|np_0|$ just at the singular point. As such, it is specified
by two positive integers $a$ and $b$ with $a+b = n$, the orders of contact
of the hyperplane with the two branches of $E$ at the node.

\ssection{1.2. Regular bundles over a Weierstrass cubic.}

Let $E$ be a Weierstrass cubic. Every semistable bundle is of the form
$$V\cong \bigoplus_\lambda V(\lambda)$$ where $\lambda$ ranges over the
isomorphism classes of rank one torsion free sheaves on
$E$ of degree zero. Let us first analyze
$V(\lambda)$ in case
$\lambda$ is a line bundle. If $V(\lambda)$ is a semistable bundle with the
property that all Jordan-H\"older quotients of $V$ are isomorphic to
$\lambda$, or in other words
$H^0(\lambda'\otimes V(\lambda))=0$ for all
$\lambda'\not=\lambda^{-1}$, then the associated graded to every
Jordan-H\"older filtration of $V(\lambda)$ is a direct sum of line bundles
isomorphic to
$\lambda$.  Of course,  one possibility for $V(\lambda)$ is the split one:
$$V(\lambda)\cong\lambda ^{\oplus r}.$$ At the other extreme we have the
maximally non-split case:

\lemma{1.7} Let $E$ be a Weierstrass cubic, possibly singular. For each
natural number
$r>0$ and each line bundle
$\lambda$ of degree zero there is a unique bundle $I_r(\lambda)$ up to
isomorphism with the following properties:
\roster
\item"{(i)}" the rank of $I_r(\lambda)$ is $r$.
\item"{(ii)}" all the  Jordan-H\"older quotients of $I_r(\lambda)$ are
isomorphic to $\lambda$.
\item"{(iii)}" $I_r(\lambda)$ is indecomposable under direct sum.
\endroster Furthermore, for all $r>0$ and all line bundles $\lambda$,
$I_r(\lambda)$ is semistable,
$I_r(\lambda)\spcheck = I_r(\lambda^{-1})$, $\det I_r(\lambda)=\lambda ^r$,
and
$\dim \Hom (I_r(\lambda), \lambda)=\dim \Hom (\lambda, I_r(\lambda))=1$. 
\endstatement

\proof We first construct the bundle $I_r=I_r(\scrO_E)$ by induction on
$r$. For $r=1$ we set $I_r=\scrO_E$. Suppose inductively that we have
constructed $I_{r-1}$ with the properties given in the lemma.  Suppose in
addition that 
$H^0(I_{r-1})\cong\Cee$.  Since the degree of $I_{r-1}$ is zero, it follows
that
$H^1(I_{r-1})\cong\Cee$ and hence there is exactly one non-trivial
extension, up to scalar multiples, of the form
$$0\to I_{r-1}\to I_r\to \scrO_E\to 0.$$ One checks easily all the
inductive hypotheses for the total space of this extension.  This proves
the existence of $I_r$ for all $r>0$. Uniqueness is easy and is left to the
reader.

We define $I_r(\lambda)=I_r\otimes\lambda$. The statements of (1.7) are
then clear.
\endproof

 The bundle $I_r(\lambda)$ has an increasing filtration by subbundles
isomorphic to $I_k(\lambda)$, $k\leq r$.  We denote this filtration by
$$\{0\} \subset F_1I_r(\lambda) \subset \cdots \subset F_rI_r(\lambda) =
I_r(\lambda),$$ and refer to $F_iI_r(\lambda)$ as the {\sl $i^{\text{th}}$
filtrant\/} of
$I_r(\lambda)$. When the bundle is clear from the context, we denote the
subbundles in this filtration by $F_i$. Notice that $F_t\cong I_t(\lambda)$
and that
$I_r(\lambda)/F_{r-t}\cong I_t(\lambda)$.

Let us note some of the basic properties of the bundles
$I_r(\lambda)$.

\lemma{1.8} Let $J$ be a proper degree zero subsheaf of $I_r(\lambda)$.
Then $J$ is contained in $F_{r-1}$. In fact, $J = F_t$ for some $t<r$.
\endstatement
\proof By the semistability of $I_r(\lambda)$, $I_r(\lambda)/J$ is a nonzero
semistable torsion free sheaf of degree zero. Clearly all of its
Jordan-H\"older quotients are isomorphic to $\lambda$. In particular there
is a nonzero map
$I_r(\lambda)/J \to \lambda$. The composition $I_r(\lambda)\to
I_r(\lambda)/J \to
\lambda$ defines a nonzero map from $I_r(\lambda)$ to $\lambda$ containing
$J$ in its kernel. By (1.7), there is a unique such nonzero map mod
scalars, and its kernel is $F_{r-1}$. Thus $J\subseteq F_{r-1}$. Applying
induction to the inclusion $J\subseteq F_{r-1}\cong I_{r-1}(\lambda)$, we
see that $J = F_t$ for some
$t<r$.
\endproof

Our next result is that the filtration is canonical, i.e.,  invariant under
any automorphism of $I_r(\lambda)$.

\corollary{1.9} If $\varphi\colon I_r(\lambda)\to I_r(\lambda)$ is a
homomorphism, then, for all
$i\le r$,    
$$\varphi(F_i)\subseteq F_i.$$ It follows that if $\varphi$ is an
automorphism, then for all $i$ we have $\varphi(F_i)=F_i$. More generally,
if $\varphi\colon I_r(\lambda)\to I_t(\lambda)$ is a homomorphism, then
$$\varphi(F_s(I_r(\lambda)))\subseteq F_s(I_t(\lambda)).$$
\endstatement

\proof It suffices to prove the last statement. By the semistability of
$I_r(\lambda)$ and $I_t(\lambda)$, $\varphi(F_s(I_r(\lambda)))$ is a degree
zero subsheaf of $I_t(\lambda)$ of rank at most $s$. Thus it is contained in
$F_s(I_t(\lambda))$.
\endproof

\lemma{1.10} Fix $0\le t\le r$ and let $q_{r,t}\colon I_r(\lambda)\to
I_t(\lambda)$ be the natural quotient map. Then $q_{r,t}$ induces a
surjective homomorphism from the endomorphism algebra of
$I_r(\lambda)$ to that of $I_t(\lambda)$. A similar statement holds for the
automorphism groups. Finally, as a $\Cee$-algebra, $\Hom(I_r(\lambda),
I_r(\lambda)) \cong \Cee[t]/(t^r)$.
\endstatement

\proof That $q_{r,t}$ induces a map on endomorphism algebras is immediate
from (1.9). Let us show that it is surjective. We might as well assume that
$t=r-1>0$ since the other cases will then follow by induction.  Let
$A\colon I_{r-1}(\lambda)\to I_{r-1}(\lambda)$ be an endomorphism. Since
the map $q_{r,r-1}\colon I_r(\lambda)\to I_{r-1}(\lambda)$ induces the zero
map on $\Hom(\lambda,\cdot)$, it follows by duality that the  map
$q_{r,r-1}^* \colon \Ext^1(I_{r-1}(\lambda), \lambda)\to
\Ext^1(I_r(\lambda), \lambda)$ is zero. Hence the composition $A\circ
q_{r,r-1}\colon I_r(\lambda)\to I_{r-1}(\lambda)$ lifts to a map $\hat
A\colon I_r(\lambda)\to I_r(\lambda)$. Thus the restriction map on
endomorphism algebras is surjective. To see the statement on automorphism
groups, suppose that
$A$ is an isomorphism.  We wish to show that $\hat A$ is an isomorphism.  
To see this, perform the construction for $A^{-1}$ as well, obtaining a map
$\widehat{A^{-1}}\colon I_r(\lambda)\to I_r(\lambda)$.  The composition
$B=\widehat{A^{-1}}\circ \hat A\colon I_r(\lambda)\to I_r(\lambda)$
projects to the identity on $I_{r-1}(\lambda)$. This means that
$B-\operatorname{Id}\colon I_r(\lambda)\to F_{1}(I_\lambda))$.  Since
$r>1$, this map is nilpotent, and hence
$B=\operatorname{Id}+(B-\operatorname{Id})$ is an isomorphism. 

Finally we prove the last statement. Let $A_r\: I_r(\lambda) \to
I_r(\lambda)$ be any endomorphism defined by a composition
$$I_r(\lambda) \twoheadrightarrow I_{r-1}(\lambda) \hookrightarrow
I_r(\lambda).$$ Note that $A_r^r =0$ and that the restriction of $A_r$ to
$I_r(\lambda)/F_1\cong I_{r-1}(\lambda)$ is of the form $A_{r-1}$. Suppose
by induction that $\Hom (I_{r-1}(\lambda), I_{r-1}(\lambda))
=\Cee[A_{r-1}]$. Then every endomorphism $T$ of $I_r(\lambda)$ is of the
form $T = p(A_r) + T'$, where $p$ is a polynomial of degree at most $r-2$
in $A_r$ and $T'$ induces the zero map on $I_r(\lambda)/F_1$. In this case
$T'$ is given by a map from $I_r(\lambda)$ to $F_1$, necessarily zero on
$F_{r-1}$, and it is easy to check that $T'$ must in fact be a multiple of
$A_r^{r-1}$. Thus $\Hom (I_r(\lambda), I_r(\lambda)) =\Cee[A_r] \cong
\Cee[t]/(t^r)$.
\endproof

We need to define an analogue of $I_r(\lambda)$ in case the Jordan-H\"older
quotients are all isomorphic to the non-locally free sheaf $\Cal F$.  We
say that a semistable degree zero bundle
$I(\Cal F)$ concentrated at the singular point of $E$ is {\sl strongly
indecomposable} if $\Hom(I(\Cal F),{\Cal F})\cong\Cee$. Notice that since
$\Hom(V(\Cal F),\Cal F)\not= 0$ for any non-trivial semistable bundle
$V(\Cal F)$ concentrated at the singular point, it follows that if
$I(\Cal F)$ is strongly indecomposable, then it is indecomposable as a
vector bundle in the usual sense. However, the converse is not true: there
exist indecomposable vector bundles which are not strongly indecomposable.
It is natural to ask if every vector bundle supported at $\Cal F$ is an
extension of strongly indecomposable bundles. Unlike the smooth case, it is
also  not true that $I(\Cal F)$ is determined up to isomorphism by its rank
and the fact that it is strongly indecomposable. Nor is it true that
$I(\Cal F)$ always has a unique filtration with successive quotients
isomorphic to $\Cal F$. As we shall show in Section 3, $I(\Cal F)$ is
determined up to isomorphism by its rank and the point $\zeta(I(\Cal F))$.

There is the following analogue for $I(\Cal F)$ of (1.8):

\lemma{1.11} Suppose that $I(\Cal F)$ is strongly indecomposable. Let
$\rho\: I(\Cal F)\to \Cal F$ be a nonzero homomorphism, unique up to scalar
multiples, and  let $X=\Ker\rho$. If $J$ is a  proper  degree zero subsheaf
of $I(\Cal F)$, then $J$ is contained in $X$.
\endstatement
\proof Let $J\subset I(\Cal F)$ be a subsheaf of degree zero.  The quotient
$Q=I(\Cal F)/J$ must be torsion-free, for otherwise $J$ would be contained
in a larger subsheaf $\hat J$ of the same rank and bigger degree,
contradicting the semistability of $I(\Cal F)$. This means that
$Q$ is semistable    of degree zero.  Clearly, it is concentrated at the
singular point.  Thus, there is a nontrivial map $Q\to {\Cal F}$. By the
strong indecomposability of $I(\Cal F)$, the composition $I(\Cal F)\to Q\to
{\Cal F}$  is some nonzero multiple of $\rho$. In particular, the kernel of
this composition is $X$. This proves that $J\subset X$.
\endproof

\definition{Definition 1.12} Let $V$ be a semistable bundle with trivial
determinant over a Weierstrass cubic. We say that $V$ is {\sl regular} or 
{\sl maximally  nonsplit\/} if,
$$V\cong \bigoplus_iI_{r_i}(\lambda_i)\oplus I(\Cal F)$$ where the
$\lambda_i$ are pairwise distinct line bundles and
$I(\Cal F)$ is a strongly indecomposable bundle concentrated at the
singular point.
\enddefinition

For $E$ smooth, Atiyah proved \cite{1} that every vector bundle $V$, all of
whose Jordan-H\"older quotients are isomorphic to $\lambda$, can be written
as a direct sum $\bigoplus _iI_{r_i}(\lambda)$. The argument carries over
to the case where
$E$ is singular, provided that $\lambda$ is a line bundle $\scrO_E(e-p_0)$.
Thus, in this case there is a unique regular bundle $V$ of rank $r$ such
that the support of $\zeta(V)$ is $e$. More generally, given a divisor $e_1+
\cdots + e_n\in |np_0|$ supported on the smooth points, there is a  unique
regular semistable rank $n$ vector bundle $V$ of trivial determinant over
$E$ such that $\zeta(V)=(e_1,\ldots,e_n)$. An  analogue of Atiyah's theorem
for the singular points has been established by T. Teodorescu \cite{12}. In
this paper, we shall show in Section 3 that, given a Cartier divisor $D$ in
$|np_0|$ whose support is the singular point, then there is a unique
regular semistable rank $n$ vector bundle
$V$ of trivial determinant such that $\zeta (V) = D$.

Regular bundles have an extremely nice property:  Their automorphism groups
have minimal possible dimension. We shall show this for smooth $E$ in the
next lemma.  To put this property in context, let us consider first the
centralizers of elements in $GL_n(\Cee)$. The centralizer of any element
has dimension at least $n$. Elements in $GL_n(\Cee)$ whose centralizers 
have dimension exactly $n$ are said to be {\sl regular} elements.  Every
element in $GL_n(\Cee)$ is
$S$-equivalent to a unique regular element up to conjugation. Here two
elements $A,B\in GL_n(\Cee)$ are said to be $S$-equivalent if every
algebraic function on $GL_n(\Cee)$ which is invariant under conjugation
takes the same value on $A$ and $B$. Said another way, $A$ and $B$ are
$S$-equivalent if there is an element
$C\in GL_n(\Cee)$ which is in the closure of the orbits of both $A$ and $B$
under the conjugation action of $GL_n(\Cee)$ on itself. From our point of
view regular bundles are the analogue of regular elements.  In fact, for a
smooth elliptic curve $E$, one way to construct a holomorphic vector bundle
over
$E$  is to fix an element $u$ in the Lie algebra of $SL_n(\Cee)$. Define a 
holomorphic connection on the trivial bundle 
$$\overline\partial_u=\overline\partial+ud\overline z$$ where
$\overline\partial$ is the usual operator on the trivial bundle. If  $u$ is
close to the origin in the Lie algebra, then the automorphism group of this
new holomorphic bundle will be the centralizer of $u$ in $GL_n(\Cee)$. In
particular, this bundle will be regular  and have trivial determinant if
and only if
$U=\operatorname{exp}(u)$ is a regular element in $SL_n(\Cee)$. For
example, if $U$ is a regular semisimple element of $SL_n(\Cee)$  then the
corresponding vector bundle over $E$ will be a sum of distinct line bundles
of degree zero.  More generally, the decomposition of $U$ into its
generalized eigenspaces will correspond to the decomposition of $V$ into
its components $V(\lambda)$. Clearly, $S$-equivalent elements of
$GL_n(\Cee)$ yield $S$-equivalent bundles. 

Here is the analogue of the dimension statements  for vector bundles over a
smooth elliptic curve.

\lemma{1.13} Let $V$ be a semistable rank $n$ vector bundle over a smooth
elliptic curve $E$.
\roster
\item"{(i)}" $\dim \Hom (V,V) \geq n$.
\item"{(ii)}" $V$ is regular if and only if $\dim \Hom (V,V) = n$. In this
case, if $V = \bigoplus _iI_{d_i}(\lambda _i)$, then the
$\Cee$-algebra $\Hom (V,V)$ is isomorphic to $\bigoplus
_i\Cee[t]/(t^{d_i})$. In particular, $\Hom (V,V)$ is an abelian
$\Cee$-algebra.
\item"{(iii)}" $V$ is regular if and only if, for all line bundles
$\lambda$ of degree zero on $E$, $h^0(V\otimes \lambda ^{-1})\leq 1$. 
\endroster
\endstatement
\proof It is easy to check that $\Hom(V(\lambda), V(\lambda ')) \neq 0$ if
and only if $\lambda = \lambda '$, and (using Corollary 1.9 and Lemma 1.10)
that $\Hom(I_d, I_d) \cong
\Cee[t]/(t^d)$. The statements (i) and (ii) follow easily from this and from
Atiyah's theorem. To see (iii), note that, for a line bundle $\mu$ of
degree zero,
$V(\mu)$ is regular if and only if
$h^0(V(\mu) \otimes \mu ^{-1}) = 1$, which implies that $h^0(V\otimes
\lambda ^{-1})\leq 1$ for all $\lambda$ since $h^0(V(\mu) \otimes \lambda
^{-1}) = 0$ if $\lambda \neq \mu$. 
\endproof

We will prove a partial analogue of  (ii) in  Lemma 1.13 for singular
curves in Section 3.

Very similar arguments show:

\lemma{1.14} Let $E$ be a Weierstrass cubic and let $V$ be a semistable
rank $n$ vector bundle over $E$. Then $V$ is regular if and only if, for
every rank one torsion free sheaf $\lambda$ of degree zero on $E$, $\dim
\Hom (V, \lambda)\leq 1$. Moreover, suppose that $V$ is regular and that
$\Cal S$ is a semistable torsion free sheaf of degree zero on $E$. Then
$\dim\Hom (V, \Cal S) \leq
\operatorname{rank}\Cal S$.
\qed
\endstatement

\section{2. The spectral cover construction.}   

In this section we shall construct families of regular semistable bundles
over a smooth elliptic curve $E$.  The main result is Theorem 2.1, which
gives the basic construction of a universal bundle over
$|np_0|\times E$, where $|np_0|\cong \Pee^{n-1}$ is the coarse moduli space
of the last section. We prove that the restriction of the universal bundle
to every slice is in fact regular, and that every regular bundle occurs in
this way. By twisting by a line bundle on the spectral cover, we construct
all possible families of universal bundles (Theorem 2.4) and show how they
are all related by elementary modifications. In Theorem 2.8, we generalize
this result to families of regular semistable bundles parametrized by an
arbitrary base scheme. In case $E$ is singular, we  establish slightly
weaker versions of these results. Most of this material will be redone from
a different perspective in the next section. Finally, we return to the
smooth case and give the formulas for the Chern classes of the universal
bundles.

{\bf Throughout this section, unless otherwise noted,  $E$ denotes a smooth
elliptic curve with origin $p_0$.}

\ssection{2.1. The spectral cover of $|np_0|$.}

Let $E^{n-1}$ be embedded in $E^n$ as the set of $n$-tuples $(e_1,
\dots, e_n)$ such that $\sum _ie_i = 0$ in the group law on $E$, or
equivalently, such that the divisor $\sum _ie_i$ on $E$ is linearly
equivalent to $np_0$. The natural action of the symmetric group
$\frak S_n$ on $E^n$ thus induces an action of $\frak S_n$ on
$E^{n-1}$. As we have seen, the quotient $E^{n-1}/\frak S_n$ is naturally
the projective space $|np_0| \cong
\Pee ^{n-1}$. View $\frak S_{n-1}$ as the subgroup of $\frak S_n$ fixing
$n$, and let $T= E^{n-1}/\frak S_{n-1}$. Corresponding to the inclusion
$\frak S_{n-1}\subset \frak S_n$ there is a morphism $\nu \: T
\to \Pee ^{n-1}$ which realizes $T$ as an $n$-sheeted cover of $\Pee
^{n-1}$. Here $\nu$ is unbranched over
$e_1+ \dots + e_n\in |np_0|$ if and only if the $e_i$ are distinct. The
branch locus of $\nu$ in $\Pee  ^{n-1}$ is naturally the dual hypersurface
to the elliptic normal curve defined by the embedding of $E$ in the dual
projective space (except in case $n=2$, where it corresponds to the four
branch points of the map from $E$ to $\Pee ^1$). The map $\nu\: T\to
|np_0|$ is called {\sl the spectral cover of
$|np_0|$}. We will discuss the reason for this name later.

The sum map $(e_1, \dots, e_n) \mapsto -\sum _{i=1}^{n-1}e_i$ is a
surjective homomorphism from $E^n$ to $E$, and its restriction to $E^{n-1}$
is again surjective, with fibers invariant under $\frak S_{n-1}$. Thus
there is an induced morphism
$r\: T \to E$. In fact, $r(e_1,
\dots, e_n) = e_n$ and
$$F_e=r^{-1}(e) = \{\, (e_1, \dots, e_{n-1}, e): \sum _{i=1}^{n-1}e_i + e =
np_0\,\},$$  modulo the obvious $\frak S_{n-1}$-action. Thus the fiber of
$r$ over $e$ is the projective space $|np_0-e|$, of dimension $n-2$.
Globally, $T$ is the projectivization of the rank
$n-1$ bundle $\Cal E$ over $E$ defined by the exact sequence
$$0 \to \Cal E \to H^0(E; \scrO_E(np_0)) \otimes _\Cee \scrO_E \to
\scrO_E(np_0) \to 0,$$ where the last map is evaluation and is surjective
since
$\scrO_E(np_0)$ is generated by its global sections.  The fiber of
$r$ over a point $e\in E$  consists of those sections of
$\scrO_E(np_0)$ vanishing at $e$, and the corresponding projective space is
just $|np_0 - e|$.  We see that there is an induced morphism on projective
bundles
$$g\: \Pee \Cal E  \to \Pee \left(H^0(E;\scrO_E( np_0)) \otimes _\Cee
\scrO_E\right) = |np_0| \times E \cong
\Pee ^{n-1} \times E,$$ such that $g$ is a closed embedding of $T$  onto
the incidence divisor in $|np_0| \times E$, and that $r$ is just the
composition of this morphism with the projection $|np_0|\times E \to E$.
Clearly $\nu$ is the composition of the morphism
$g\: \Pee \Cal E  \to |np_0| \times E$ with projection to the first factor,
or equivalently $g=(\nu, r)$. Given
$e\in E$, let $F_e = r^*e$ be the fiber over $e$ and let $\zeta$ be the
divisor class corresponding to $c_1(\scrO_T(1))$, viewing $T$ as $\Pee \Cal
E $. Since
$\Cal E$ sits inside the trivial bundle, it follows that $\zeta = g^*\pi
_1^*h$, where
$h = c_1(\scrO_{\Pee ^{n-1}}(1))$, and thus $\zeta = \nu ^*h$. Note also
that each fiber $F_e$ of $T= \Pee \Cal E  \to E$ is mapped linearly into
the corresponding hyperplane $H_e = |np_0 - e|$ of
$\Pee ^{n-1} = |np_0|$ consisting of divisors containing $e$ in their
support. Thus as divisor classes
$\nu _*[F_e] = h$.

There is a special point $\bold o = \bold o_E= np_0 \in |np_0|$. (In terms
of regular semistable bundles, $\bold o$ corresponds to $I_n$.)  It is one
of the
$n^2$ points of ramification of order $n$ for the map $T\to |np_0|$,
corresponding to the $n$-torsion points of $E$. 

\ssection{2.2. A universal family of regular semistable bundles.}

Next we turn to the construction of a universal family of regular
semistable bundles $E$.  It will be given by a bundle  $U_0$ over
$|np_0|\times E$. Over
$E^{n-1}\times E$, we have the diagonal divisor
$$\{\, (e_1, \dots, e_n, e): e = e_n, \sum _{i=1}^ne_i = 0\,\},$$ which is
invariant under the $\frak S_{n-1}$-action and so descends to a divisor
$\Delta$ on $T\times E$, which is the graph of the map $r\colon T\to E$. .
Note that $\Delta \cong T$ and that
$$\Delta = (r\times \Id)^*\Delta _0,$$ where
$\Delta _0$ is the diagonal $\{\, (e, e): e\in E\,\}$. Let $G = T\times
\{p_0\}$. Then the divisor $\Delta - G$ has the property that its
restriction to a slice
$\{(e_1, \dots, e_{n-1}, e_n)\}\times E$ can be identified with the line
bundle $\scrO_E(e_n - p_0)$. We define ${\Cal L}_0\to T\times E$ to be the
line bundle $\scrO_{T\times E}(\Delta - G)$, and we set
$$U_0=\left(\nu\times\operatorname{Id}\right)_*{\Cal L}_0.$$

\theorem{2.1} Let $E$ be a smooth elliptic curve. The sheaf $U_0$ over
$|np_0|\times E$ constructed above is a vector bundle of rank $n$.  For
each $x\in |np_0|$ the restriction of
$U_0$ to $\{x\}\times E$ is a regular semistable bundle $V_x$ with trivial
determinant and with the property that $\zeta(V_x)=x$.
\endstatement

\proof
 Since $\nu\times\operatorname{Id}$ is an $n$-sheeted covering of smooth
varieties, it is a finite flat morphism and hence $U_0$ is a vector bundle
of rank $n$. If
$\{(e_1, \dots, e_{n-1}, e_n)\}$ is not a branch point of $\nu$, or in other
words if the $e_i$ are pairwise distinct, then
$U_0 = (\nu\times \Id)_*\scrO_{T\times E}(\Delta - G)$ restricts over the
slice
$\{(e_1, \dots, e_{n-1}, e_n)\}\times E$ to a bundle isomorphic to the
direct sum 
$$\bigoplus_{e_i}{\Cal L}_0|\{e_i\}\times E $$ which is clearly isomorphic
to 
$$\scrO_E(e_1 - p_0) \oplus \cdots \oplus \scrO_E(e_n - p_0).$$ This shows
that for a  generic point $s\in |np_0|$ the restriction
$U_0|\{s\}\times E$ is as claimed: it is the unique regular semistable
bundle with the given Jordan-H\"older quotients.

In general, consider a point $x\in |np_0|$ of the form $\sum _{i=1}^\ell
r_ie_i$, where the $e_i\in E$ and the $r_i$ are positive integers with
$\sum _ir_i = n$. We claim that the Jordan-H\"older quotients of the
corresponding bundle are
$\scrO_E(e_i-p_0)$, with multiplicity $r_i$. The preimage of $x$ in $T$
consists of $\ell$ points $y_1, \dots, y_\ell$, each of multiplicity $r_i$.
Viewing $T$ as the incidence correspondence in $\Pee^{n-1}\times E$, the
point $y_i$ corresponds to $\left(\sum _{i=1}^\ell r_ie_i, e_i\right)$. If
$R = \scrO_T/\frak m_x\scrO_T$ is the coordinate ring of the fiber over
$x$, then
$R$ is the product of $\ell$ local rings $R_i$ of lengths $r_1, \dots,
r_\ell$. It clearly  suffices to prove the following claim.

\claim{2.2}In the above notation, $\Cal L_0\otimes R_i$  has a filtration
all of whose successive quotients are isomorphic to
$\lambda_i$  where $\lambda_i\cong\scrO_E(e_i-p_0)$. In particular, the
restriction of $U_0$ to this slice is semistable and has  determinant
$\lambda ^{r_i}$.
\endstatement

\proof The ring $R_i$ has dimension $r_i$ and is filtered by ideals whose
successive quotients are isomorphic to $\Cee_{y_i}$. Thus $\Cal L_0\otimes
R_i$ is filtered by subbundles whose quotients are all isomorphic to the
line bundle $\Cal L_0|\{y_i\}\times E$. But by construction this
restriction is $\scrO_E(e_i-p_0)$.
\endproof

At this point, we have seen that $U_0$ is a family of semistable bundles on
$E$ whose restriction to every fiber has trivial determinant and with the
``correct" Jordan-H\"older quotients.  It remains to  show that $U_0$ is a
family of regular bundles over $E$.

\claim{2.3} The restriction of $U_0$ to every slice $\{e\}\times E$ is 
regular. 
\endstatement
\demo{Proof of the claim} To see that the restriction to each slice is
regular, note that a semistable  $V$ of degree $0$ on $E$ is regular if and
only if, for all line bundles $\lambda$ on $E$ of degree zero,
$h^0(V\otimes \lambda ^{-1})\leq 1$. By Riemann-Roch on $E$,
$h^0(V\otimes
\lambda ^{-1}) = h^1(V\otimes \lambda ^{-1})$. Thus we must show that
$h^1(V\otimes \lambda ^{-1})\leq 1$.

First we calculate $R^1\pi _1{}_*(U_0\otimes \pi _2^*\lambda ^{-1})$, where
$\pi _1\:
\Pee ^{n-1} \times E\to \Pee^{n-1}$  is the projection to the first factor.
Let $q_1\: T\times E \to T$ be the first projection. Consider the diagram
$$\CD T\times E @>{\nu \times \Id}>> \Pee ^{N-1} \times E\\ @V{q_1}VV
@VV{\pi _1}V\\ T @>{\nu}>> \Pee ^{n-1}.
\endCD$$ Since $\nu$ and $\nu \times \Id$ are affine, we obtain
$$R^1\pi _1{}_*\left[(\nu \times \Id)_*\scrO_{T\times E }(\Delta - G)
\otimes
\pi _2^*\lambda ^{-1}\right] =
\nu _*R^1q_1{}_*\left[\scrO_{T\times E}(\Delta - G ) \otimes q_2^*
\lambda ^{-1}\right].$$ Now apply flat base change to the Cartesian diagram
$$\CD T\times E @>{r \times \Id}>> E \times E\\ @V{q_1}VV @VV{p _1}V\\ T
@>{r}>> E.
\endCD$$ We have $\scrO_{T\times E }(\Delta - G) = (r\times
\Id)^*\scrO_{E\times E}(\Delta _0 - (E \times \{p_0\})  )$, and thus the
sheaf
$R^1q_1{}_*\left[\scrO_{T\times E }(\Delta - G) \otimes q_2^* \lambda
^{-1}\right]$ is isomorphic to
$$r^*R^1p_1{}_*\left[\scrO_{E\times E }(\Delta _0 -( E \times
\{p_0\}) )\otimes p_2^* \lambda ^{-1}\right].$$  Rrestricting to the slice
$\{e\}\times E$, we see that
$$R^1p_1{}_*\left(\scrO_{E\times E }(\Delta _0 - (E \times
\{p_0\}))\otimes p_2^* \lambda ^{-1}\right)$$ is supported at the point
$e$ of $E$ corresponding to the line bundle $\lambda$ (i\.e\.
$\lambda =
\scrO_E(e-p_0)$), and the calculation of
\cite{6}, Lemma 1.19 of Chapter 7, shows that the length at this point is
one.  Thus taking
$r^*$ gives the sheaf $\scrO_{F_e}$, and $\nu_*\scrO_{F_e}=
\scrO_{H_e}$, where
$H_e$ is a reduced hyperplane in $\Pee ^{n-1}$.  Thus we have seen that
$R^1\pi _1{}_*(U\otimes \pi _2^*\lambda ^{-1})$ is (up to twisting by a
line bundle)
$\scrO_{H_e}$, where $H_e$ is the hyperplane in $\Pee ^{n-1}$ corresponding
to
$|np_0 - e|$. Since $\pi_1$ has relative dimension one, $R^2\pi _1
{}_*(U\otimes
\pi _2^*\lambda ^{-1})=0$. It follows by the theorem on cohomology and base
change \cite{10} Theorem 12.11(b) that the map
$R^1\pi _1 {}_*(U\otimes \pi _2^*\lambda ^{-1})\to H^1(V\otimes
\lambda ^{-1})$ is surjective, and thus $h^1(V\otimes \lambda ^{-1})\leq 1$
as desired.
\endproof
\enddemo

\ssection{2.3. All universal families of regular semistable bundles.}

We have constructed a bundle $U_0$ over $|np_0|\times E$  with given
restriction to each slice.  Our next goal is to understand all such bundles.

\theorem{2.4} Let $E$ be a smooth elliptic curve. Let $\pi_1\colon
|np_0|\times E\to |np_0|$ be the projection onto the first factor, and let
$U_0$ be the bundle constructed in Theorem \rom{2.1}. Then:
 \roster
\item"{(i)}" The sheaf   $\pi _1{}_*Hom(U_0,U_0)$ is a locally free sheaf
of algebras of rank
$n$ over $|np_0|$ which is isomorphic to $\nu _*\scrO_T$.
\item"{(ii)}" Let $U'$ be a rank $n$ vector bundle over $|np_0|\times E$
with the following property.  For each $x\in |np_0|$ the restriction of
$U'$ to $\{x\}\times E$ is isomorphic to the restriction of $U_0$ to
$\{x\}\times E$. Then $U'=(\nu \times
\Id)_*\left[\scrO_{T\times E }(\Delta-G)\otimes q _1^*L\right]$ for a
unique line bundle  $L$  on $T$.
\endroster
\endstatement

\proof In view of  Claim 2.3 and the definition of a regular bundle, $\pi
_1{}_*Hom(U_0,U_0)$ is a locally free sheaf of algebras of  rank $n$ over
$|np_0|$. To see that it is isomorphic to $\nu _*\scrO_T$, note that
multiplication by functions defines a homomorphism
$\nu _*\scrO_T \to \pi _1{}_*Hom(U_0,U_0)$ which is clearly an inclusion of
algebras. Since both sheaves of algebras are rank $n$ vector bundles over
$\Pee^{n-1}$, they agree at the generic point of
$|np_0|$. Thus, over every affine open susbet of $|np_0|$ the rings
corresponding to $\pi _1{}_*Hom(U_0,U_0)$ and $\nu _*\scrO_T$ are two
integral domains with the same quotient fields. Since
$T$ is normal and
$\pi _1{}_*Hom(U_0,U_0)$ is finite over $\nu _*\scrO_T$ (since it is finite
over
$\scrO_{\Pee^{n-1}}$), the two sheaves of algebras must coincide. This
proves (i).

Now suppose that $U'$ satisfies the hypotheses (ii) of (2.4).  By base
change $\pi _1{}_*Hom(U',U_0)$ is a locally free rank
$n$ sheaf over $|np_0|$. Composition of homomorphisms induces the structure
of a
$\pi _1{}_*Hom(U_0,U_0)$-module on $\pi _1{}_*Hom(U',U_0)$. Thus
$\pi _1{}_*Hom(U',U_0)$ corresponds to a $\nu_*\scrO_T$-module. We claim
that, as an
$\scrO_T$-module, $\pi _1{}_*Hom(U',U_0)$ is locally free rank of rank one.
To see this, fix a point $x$ in $|np_0|$ and let $V', V$ be the vector
bundles corresponding to the restrictions of $U', U_0$ to the slice
$\{x\}\times E$. Of course, by hypothesis $V'$ and $V$ are isomorphic.
Choose an isomorphism
$s\: V'\to V$ and extend it to a local section of $\pi _1{}_*Hom(U',U_0)$
in a neighborhood of $x$, also denoted $s$. The map $\pi _1{}_*Hom(U_0,U_0)
\to \pi _1{}_*Hom(U',U_0)$ defined by multiplying against the section $s$
is then surjective at $s$, and hence in a neighborhood. Viewing both sides
as locally free rank $n$ sheaves over $|np_0|$, the map is then a local
isomorphism. But this exactly says that $\pi _1{}_*Hom(U',U_0)$ is a
locally free $\pi _1{}_*Hom(U_0,U_0)$-module of rank one. Thus $\pi
_1{}_*Hom(U',U_0)$ corresponds to a line bundle on
$T$, which we denote by $L^{-1}$. Of course, for any line bundle $M$ on
$T$, 
 setting
$$U_0[M] =  (\nu \times
\Id)_*\left[\scrO_{T\times E }(\Delta-G)\otimes q _1^*M\right].$$ we have
$$\align
\pi _1{}_*Hom (U', U_0[M]) &= \pi  _1{}_*(\nu \times \Id)_*Hom((\nu
\times \Id)^*U', \scrO_{T\times E}(\Delta-G)\otimes q _1^*M))\\
 &=\nu _*q_1{}_*\left[q_1^*M\otimes Hom((\nu \times \Id)^*U',
\scrO_{T\times E}(\Delta-G))\right]\\ &= \nu _*\left[M\otimes q_1{}_*
Hom((\nu \times \Id)^*U', \scrO_{T\times E}(\Delta-G))\right], 
\endalign$$ The case $M=\scrO_T$ tells us that, as $\nu_*\scrO_T$-modules,
$$\nu_*L^{-1}=\pi_1{}_*Hom (U', U_0)= \nu _*\left[q_1{}_* Hom((\nu \times
\Id)^*U', \scrO_{T\times E}(\Delta-G))\right].$$ Hence, we have
$$\pi _1{}_*Hom (U', U_0[M]) =\nu_*(M\otimes L^{-1}).$$ Taking $M=L$ we have
$$\pi _1{}_*Hom (U', U_0[L]) \cong\nu_*\scrO_T.$$ Via this identification,
the section  $1\in H^0(\scrO_T)$ then defines an isomorphism from $U'$ to
$U_0(L)$, as claimed.
\endproof

In view of the previous result, we need to describe all line bundles on
$T$. Since $T$ is a $\Pee^{n-2}$ bundle over $E$, we have:

\lemma{2.5} The projection mapping $r\colon T\to E$ induces an injection 
$$r^*\Pic E\to \Pic T.$$ If $n=2$, $r$ is an isomorphism and thus $\Pic
T\cong \Pic E$.  For $n> 2$, since
$T$ is included in $\Pee \Cal E \subset \Pee^{n-1}\times E$, we can define
by restriction the line bundle 
$\scrO_{\Pee^{n-1}}(1)|T=\scrO_T(1)$ on $T$. Then
$$\Pic T =r^*\Pic E\oplus \Zee[\scrO_T(1)].
\qquad \qed$$
\endstatement

In view of Lemma 2.5, we make the following definition. For $p\in E$, let 
$F_p\subset T$ be the divisor which is the preimage of $p$. 

\definition{Definition 2.6}  For every integer $a$, let $U_a =(\nu\times
\Id)_*\scrO_{T\times E }(\Delta - G - a(F_{p_0}\times E))$. More generally,
given $e\in E$, define
$$U_a[e] = (\nu\times
\Id)_*\scrO_{T\times E }(\Delta - G - (a+1)(F_{p_0}\times E)+(F_e\times
E)).$$ Thus  $U_a[p_0] = U_a$. By Lemma 2.5, every vector bundle obtained
from $U_0$ by twisting by a line bundle on the spectral cover is of the
form $U_a[e]\otimes \pi _1^*\scrO_{\Pee^{n-1}}(b)$ for some $b\in \Zee$ and
$e\in E$. (For $n=2$, we have the relation $\nu ^*\scrO_{\Pee^1}(1) =
\scrO_E(2p_0)$, and thus $U_a\otimes  \pi _1^*\scrO_{\Pee ^{1}}(b)\cong
U_{a-2b}$.)
\enddefinition

The next lemma says that the $U_a$ are all elementary modifications of each
other: 

\lemma{2.7} Let $H =\nu(F_{p_0})$ be the hyperplane in $\Pee^{n-1} =|np_0|$
of divisors whose support contains $p_0$, and let $i\: H \to \Pee^{n-1}$ be
the inclusion. Then there is an exact sequence
$$0 \to U_a \to U_{a-1} \to (i\times \Id)_*\scrO_{H\times E} \to 0.$$
Moreover $\dim \Hom (U_{a-1}|H\times E, \scrO_{H\times E}) = 1$, so that
the above exact sequence is the unique elementary modification of this
type. Likewise,
$U_a[e]$ is given as an elementary modification:
$$0 \to U_{a+1} \to U_a[e] \to (i\times \Id)_*\scrO_{H_e\times E}\otimes
\pi_2^*\scrO_E(e-p_0)
\to 0.$$
\endstatement
\proof Consider the exact sequence
$$\gather 0 \to \scrO_{T\times E }(\Delta - G - a(F_{p_0}\times E)) \to \\
\to \scrO_{T\times E }(\Delta - G - (a-1)(F_{p_0}\times E))\to
\scrO_{F_{p_0}\times E }(\Delta - G - (a-1)(F_{p_0}\times E))\to 0.
\endgather$$ Clearly the restriction of the line bundle $\scrO_{T\times E
}(F_{p_0}\times E)$ to $F_{p_0}\times E$ is trivial, and $G$ and $\Delta$
both restrict to the divisor
$F_{p_0}\times \{p_0\}\subset F_{p_0}\times E$. Hence the last term in the
above sequence is
$\scrO_{F_{p_0}\times E }$. Applying  $(\nu\times \Id)_*$ to the sequence
gives the exact sequence of (2.7). For  $V$ a bundle corresponding to a
point of $H$, $\dim \Hom(V, \scrO_E) = 1$. Thus $\pi _1{}_*Hom(U_a|H\times
E, \scrO_{H\times E})$ is a line bundle on $H$. The given map $U_a|H\times
E \to \scrO_{H\times E}$ constructed above is an everywhere generating
section of this line bundle, so that $\pi _1{}_*Hom(U_a|H\times E,
\scrO_{H\times E})$ is trivial and $\dim \Hom (U_a|H\times E,
\scrO_{H\times E}) = 1$.

The proof of the exact sequence relating $U_{a+1}$ and $U_a[e]$ is similar.
\endproof

In fact, suppose that we have an elementary modification
$$0 \to U' \to U_a \to \scrO_{D\times E}\otimes \pi _2^*\lambda \to 0,$$
where $D$ is a hypersurface in $|np_0|$ and $\lambda$ is a line bundle of
degree zero on $E$. Then it is easy to check that necessarily $D=H_e$ for
some $e$ and
$\lambda =\scrO_E(e-p_0)$. Of course, it is also possible to make elementary
modifications along certain hyperplanes corresponding to taking higher rank
quotients of $U_a$.

\ssection{2.4. Families of bundles over more general parameter spaces.} 

Now let us examine in what sense the bundles $U\to |np_0|\times E$ that we
have constructed are universal.

\theorem{2.8} Let $E$ be a smooth elliptic curve and let $S$ be a scheme or
analytic space. Suppose that
${\Cal U}\to S\times E$ is a rank $n$ holomorphic vector bundle whose 
restriction to each slice $\{s\}\times E$ is a regular semistable bundle
with trivial determinant.  Let $\Phi\colon S\to |np_0|$ be the morphism
constructed in Theorem \rom{1.5}. Let $\nu_S\colon\tilde S\to S$ be the
pullback via
$\Phi$ of the spectral covering $T\to |np_0|$\rom:
$$\tilde S=S\times_{|np_0|}T,$$ and let $\tilde\Phi\colon \tilde S\to T$ be
the map covering
$\Phi$. Let $q_1\colon \tilde S\times E\to \tilde S$ be the projection onto
the first factor. 
 Then there is a line bundle ${\Cal M}\to \tilde S$ and an isomorphism of
$\Cal U$ with 
$$(\nu_S\times\operatorname{Id})_*
\left((\tilde\Phi\times \Id)^*(\scrO_{T\times E}(\Delta-G))\otimes
q_1^*{\Cal M}\right).$$ 
\endstatement

\proof By construction the bundle 
$$(\nu_S\times\operatorname{Id})_* (\tilde\Phi\times \Id)^*(\scrO_{T\times
E}(\Delta-G))$$  is a family of regular semistable bundles with trivial
determinant $E$, which fiber by fiber have the same Jordan-H\"older
quotients as the family ${\Cal U}$. But regular semistable bundles are
determined up to isomorphism by their Jordan-H\"older quotients. This means
that the two families are isomorphic on each slice
$\{s\}\times E$.  Now the argument in the proof of Theorem 2.4 applies to
establish the existence of the line bundle $\Cal M$ on the spectral
covering $\tilde S$ as required.
\endproof

 We can also construct the spectral cover
$\tilde S$ of $S$ directly.  This construction will also the explain the
origin of the name {\sl spectral cover\/}. If $p_1, p_2$ are the
projections of $S\times E$ to the first and second factors, then by
standard base change results  $p_1{}_*Hom (\Cal U, \Cal U)$ is a locally
free sheaf of coherent
$S$-algebras. Moreover, by the classification of regular semistable
bundles, it is commutative. Thus there is a well-defined space
$\tilde S=\bold{Spec}\,p_1{}_*Hom (\Cal U, \Cal U)$ and a morphism $\nu\:
\tilde S \to S$ such that $\scrO_{\tilde S} = p_1{}_*Hom (\Cal U, 
\Cal U)$. It is easy to check directly that $\tilde S = S\times 
_{|np_0|}T$. By construction, there is an action of
$\scrO_{\tilde S}$ on
$\Cal U$ that  commutes with the action of $\scrO_E$, and thus $\Cal U$
corresponds to a coherent sheaf
$\Cal L$ on $\tilde S \times E$. Again by the classification of regular 
semistable  bundles, it is straightforward to check directly that $\Cal L$
is locally free of rank one. Clearly,
$(\nu\times \Id)_*\Cal L = \Cal U$.

We can view Theorem 2.8 as allowing us to replace a family of possibly
non-regular, semistable bundles with trivial determinant on $E$  with a
family of regular semistable bundles without changing the Jordan-H\"older
quotients on any slice.
  Suppose that ${\Cal V}\to S\times E$ is any family of semistable bundles
with trivial determinant over $E$.  We have the map $\Phi\colon S\to
|np_0|$ of Theorem 1.5, and
$(\Phi\times \operatorname{Id})^*U_0\to S\times E$ is a family of regular
semistable  bundles with the same Jordan-H\"older quotients as ${\Cal V}$
along each slice $\{s\}\times E$.  Of course, the new bundle will not be
isomorphic to ${\Cal V}$ (even after twisting with a line bundle on the
spectral cover) unless the original family is a family of regular bundles.

\ssection{2.5. The case of singular curves.}

There is an analogue of these constructions  for singular curves.  Let $E$
be a Weierstrass cubic. The constuction given at the beginning of this
section is valid in this context and produces a $\Pee^{n-2}$-bundle
$T=\Pee\Cal E$ over $E$ and an $n$-fold covering map $\nu\colon T\to
|np_0|$. By the description of $T$ as $\Pee\Cal E$, the projection $T\to
\Pee^{n-1}$ is a finite flat morphism. 

Let $\Omega\subset |np_0|$ be the Zariski open subset of all divisors whose 
support does not contain the singular point of $E$, and let
$T_\Omega\subset T$   be $\nu^{-1}(\Omega)$. We denote by $\nu_\Omega$ the
restriction of
$\nu$ to $T_\Omega$. It is a finite surjective morphism of degree $n$
between smooth varieties. As before, we have the divisor $\Delta\subset
T\times E$.  We denote by $\Delta_\Omega\subset T_\Omega\times E$ the
restriction of $\Delta$ to this open subset. We form the line bundle
$${\Cal L}_0^\Omega=\scrO_{T_\Omega\times E}(\Delta_\Omega-G),$$ where $G$
is the divisor $T_\Omega\times \{p_0\}$. Let $U_0^\Omega$ be the sheaf
$(\nu_\Omega\times
\operatorname{Id})_*({\Cal L}_0^\Omega)$ over $\Omega\times E$.  It is a
vector bundle of rank $n$. The arguments in the proof of Claim 2.3 apply in
this context and show that $U_0^\Omega$ is a family of regular bundles on
$E$ parametrized by $\Omega$.

The arguments in the proof of Theorem 2.4 apply  to yield the following
result.

\proposition{2.9} Let $E$ be a Weierstrass cubic. Let $\Omega\subset
|np_0|$ be the Zariski open subset defined in the previous paragraph. 
 Let $\pi^\Omega_1\colon \Omega\times E\to \Omega$ be the projection onto
the first factor, and let
$U_0^\Omega=(\nu_\Omega\times \operatorname{Id})_*{\Cal L}_0^\Omega$ be the
bundle over $\Omega\times E$ constructed in the previous paragraph. Then:
\roster
\item"{(i)}" The sheaf $(\pi_1^\Omega)_* Hom(U_0^\Omega,U_0^\Omega)$ is a 
locally free sheaf of algebras of rank
$n$ over $\Omega\subset |np_0|$ which is isomorphic to $\nu
_*\scrO_{T_\Omega}$. 
\item"{(ii)}" Let $U'$ be a rank $n$ vector bundle over $\Omega\times E$
with the following property.  For each $x\in \Omega$ the restriction of
$U'$ to $\{x\}\times E$ is isomorphic to the restriction of $U_0$ to
$\{x\}\times E$. Then $U'=(\nu_\Omega \times
\Id)_*\left[\scrO_{T_\Omega\times E }(\Delta_\Omega-G)\otimes q
_1^*L\right]$ for a unique line bundle  $L$  on $T_\Omega$.
\endroster
\endstatement

In Section 3, we shall show how to extend this construction over the
singular points of $E$.

\ssection{2.6. Chern classes.}

Finally, we return to the case where $E$ is smooth and  give the Chern
classes of the various bundles over $|np_0|\times E$  in  case $E$ is
smooth. The proof of (2.10) will be given  in the next section, and we will
prove the remaining results assuming (2.10).

\proposition{2.10} Identify $h \in H^2(\Pee ^{n-1})$ with its pullback to
$\Pee ^{n-1} \times E$. Then the total Chern class $c(U_0)$ and the Chern
character $\ch(U_0)$ of $U_0$ are given by the formulas\rom:
$$\align c(U_0) &= (1-h+ \pi _2^*[p_0]\cdot h)(1-h)^{n-2}\\
\ch U_0 &= ne^{-h} + (1- \pi _2^*[p_0])(1- e^{-h}).
\endalign$$
\endstatement

Once we have (2.10), we can calculate the Chern classes of all the universal
bundles.

\proposition{2.11} Let $U_a$ and $U_a[e]$ be defined as in \rom{(2.6)}. Let
$h$ be the class of a hyperplane in $\Pic\Pee^{n-1}$, which we also view by
pullback as an element of $\Pic (\Pee^{n-1}\times E)$.
\roster
\item"{(i)}" $c(U_a) = (1-h+ \pi _2^*[p_0]\cdot h)(1-h)^{a+n-2}$.
\item"{(ii)}" $\ch (U_a\otimes \pi _1^*\scrO_{\Pee ^{n-1}}(b)) =
ne^{(b-1)h} + (1-a- \pi _2^*[p_0])(e^{bh}- e^{(b-1)h})$.
\item"{(iii)}" $\det U_a[e] = -(a+n-1)h$.
\item"{(iv)}" Let $\tilde c_2$ denote the refined Chern class of a vector
bundle  in the Chow group
$A^2(\Pee^{n-1}\times E)$. Let $A^2_0(\Pee^{n-1}\times E)$ be the subgroup
of
$A^2(\Pee^{n-1}\times E)$ of all cycles homologous to zero, so that
$A^2_0(\Pee^{n-1}\times E) \cong E$. Then
$$\tilde c_2(U_a[e]\otimes \pi _1^*\scrO_{\Pee ^{n-1}}(b)) - \tilde
c_2(U_a\otimes \pi _1^*\scrO_{\Pee ^{n-1}}(b)) = e$$ as an element of
$A^2_0(\Pee^{n-1}\times E) \cong E$.
\endroster
\endstatement
\proof By (2.7), 
$$c(U_a) = c(U_{a-1})c((i\times \Id)_*\scrO_{H\times E})^{-1}$$ and likewise
$$\ch U_a = \ch U_{a-1} - \ch ((i\times \Id)_*\scrO_{H\times E}).$$ Using
the exact sequence
$$0 \to \scrO_{\Pee^{n-1}\times E}(-H\times E) \to \scrO_{\Pee^{n-1}\times
E} \to  (i\times \Id)_*\scrO_{H\times E} \to 0,$$ we have 
$$\align c((i\times \Id)_*\scrO_{H\times E}) &= (1-h)^{-1};\\
\ch ((i\times \Id)_*\scrO_{H\times E}) &= 1-e^{-h}.
\endalign$$ A little manipulation, starting with (2.10), gives (i) and
(ii). To see (iii), note that  by construction $\det U_a[e]$ is the
pullback of a class in $\Pic\Pee^{n-1}$. Moreover, it is independent of the
choice of $e\in E$. Thus we may as well take
$e=p_0$, in which case $U_a[p_0]=U_a$. In this case, the result is
immediate from (i).  (iv)  follows by using the elementary modification
relating $U_a[e]$ and
$U_{a+1}$.
\endproof

Note that
$$c_1\left(U_a\otimes \pi _1^*\scrO_{\Pee ^{n-1}}(b))\right) = 0$$ if and
only if $a-1 = n(b-1)$. A natural solution to this equation is $a=b=1$.The
bundle $U = U_1\otimes \pi _1^*\scrO_{\Pee ^{n-1}}(b) = (\nu\times
\Id)_*\scrO_{T\times E}(\Delta - G -(F_{p_0}\times E))$ is singled out in
this way as
$(\nu \times \Id)_*\Cal P$, where $\Cal P$ is the pullback to
$T\times E$ of the symmetric line bundle $\scrO_{E\times E }(\Delta _0 -
\{p_0\}\times E -f\times \{p_0\})$, which is a Poincar\'e line bundle for
$E\times E$. In this case $\ch U = n + \pi _2^*[p_0](1- e^h)$. Moreover,
one can check that $c_1(U) = 0$ and $c_k(U) =(-1)^kh^{k-1}\pi _2^*[p_0]$
for $k\geq 2$.

It is easy to check that, for $n>2$,  $U_a\otimes \pi _1^*\scrO_{\Pee
^{n-1}}(b) = U_{a'}\otimes
\pi _1^*\scrO_{\Pee ^{n-1}}(b')$ if and only if $a=a'$ and $b=b'$. It is
also possible to vary
$aF$ within its algebraic equivalence class, which is a family isomorphic
to $E$, and  this difference is detectable by looking at
$c_2\left(U_a\otimes \pi _1^*\scrO_{\Pee ^{n-1}}(b)\right)$ in the Chow
group
$A^2(\Pee ^{n-1}\times E)$. More precisely, we have the following:

\proposition{2.12} Given two vector bundles  $U' =(\nu\times
\Id)_*\left(\scrO_{T\times E }(\Delta - G )\otimes M'\right)$ and $U''
=(\nu\times
\Id)_*\left(\scrO_{T\times E }(\Delta - G )\otimes M''\right)$, where $M'$
and
$M''$ are line bundles on
$E$, then $U'$ and $U''$ are isomorphic if and only if they have the same
Chern classes as elements of $A^*(\Pee ^{n-1}\times E)$.
\endstatement 
\proof For simplicity, we shall just consider the case $n>2$. Using the
notation of (2.6) and the description of $\Pic T$, it suffices to show
that, if the Chern classes of
$U_a[e]\otimes \pi _1^*\scrO_{\Pee ^{n-1}}(b)$ and of $U_{a'}[e']\otimes \pi
_1^*\scrO_{\Pee ^{n-1}}(b')$ are equal in the Chow ring, then $a=a', b=b'$,
and $e=e'$. Following the above remarks, the Chern classes of
$U_a[e]\otimes \pi _1^*\scrO_{\Pee ^{n-1}}(b)$ in rational cohomology,
which are of course the same as those of $U_a\otimes \pi _1^*\scrO_{\Pee
^{n-1}}(b)$, determine $a$ and $b$ (use $c_3$ to find $b$ and $c_1$ to find
$a$). By (iii) of (2.11), the class $\tilde c_2$ then determines $e$.
\endproof

\section{3. Moduli spaces via extensions.}

In this section we shall describe a completely different approach to
constructing universal bundles over $\Pee ^{n-1}\times E$. The idea here is
to consider the space of extensions of fixed (and carefully chosen) bundles
over $E$.  From this point of view the projective space is the  projective
space of the relevant extension group, which is {\it a priori\/} a very
different animal from $|np_0|$. We shall show however (Proposition 3.13)
that this projective space is naturally identified with $|np_0|$. There are
several reasons  for considering this alternative approach. First of all it
works as well for singular curves as for smooth ones, so that the
restrictions of the last section to smooth curves or to bundles 
concentrated away from the singularities of a singular curve can be
removed. Also, this method works well for a family of elliptic curves, not
just a single elliptic curve.  Lastly, this approach has a  natural
generalization to all holomorphic principal bundles with structure group 
an arbitrary complex simple group $G$, something which so far is not clear
for the spectral cover approach. The generalization to $G$-bundles is
discussed in \cite{8}. The disadvantange of the approach of this section is
that it constructs some but not all of the families that the spectral cover
approach gives.  The reason is that from this point of view one cannot see
directly the analogue of twisting by a general line bundle on the spectral
cover to produce the general family of regular bundles.

The main results of this section are as follows. In Theorem 3.2, we
consider the set of relevant  extensions and show that every such extension
is a regular semistable bundle with trivial determinant. Conversely, every
regular semistable bundle with trivial determinant arises as such an
extension. In the construction of bundles of rank $n$ over $E$ we must
choose an integer  $d$ with
$1\le d<n$.  We show that constructions for different $d$ are related to one
another (Proposition 3.11 and Theorem 3.12). Next, we compare the extension
moduli space, which is a
$\Pee^{n-1}$, to the coarse moduli space which is $|np_0|$. We find a
natural cohomological identification of these two projective spaces
(Theorem 3.13) and check that it corresponds to the morphism $\Phi$ of
Section 1 (Proposition 3.16). Next we show how the universal bundles
defined via the extension approach lead to the spectral covers of Section 2
(Theorem 3.21). In this way, we can both identify the universal bundles
constructed here with those constructed via spectral covers (Theorem 3.23
and Corollary 3.24), and extend the spectral cover construction to the case
of a singular $E$.

{\bf Throughout this section,  $E$ denotes a Weierstrass cubic with origin
$p_0$.}

\ssection{3.1. The basic extensions.} We begin by recalling a result,
essentially due to Atiyah, which produces the basic bundles for our
extensions:

\lemma{3.1} For each $d\geq 1$, there is a stable bundle
$W_d$ of rank $d$ on $E$ whose determinant is isomorphic to
$\scrO_E(p_0)$. It is unique up to isomorphism. For every rank one torsion
free sheaf $\lambda$ of degree zero, $h^0(W_d\otimes \lambda )=1$ and
$h^1(W_d\otimes
\lambda )=0$.
\endstatement

\proof We briefly outline the proof. An inductive construction of
$W_d$ is as follows: set $W_1 =\scrO_E(p_0)$. Assume inductively that
$W_{d-1}$ has been constructed and that
$h^0(W_{d-1}) = 1$. It then follows by Riemann-Roch that
$h^1(W_{d-1}) = 0$, and thus that $h^0(W_{d-1}\spcheck) = 0$,
$h^1(W_{d-1}\spcheck) = 1$. We then define $W_d$ by taking the unique
nonsplit extension
$$0 \to \scrO_E \to W_{d} \to W_{d-1} \to 0.$$ By construction $W_d$ has a
filtration whose successive quotients, in increasing order, are
$\scrO_E, \dots, 
\scrO_E$, $\scrO_E(p_0)$, and such that all of the intermediate extensions
are not split.  It is the unique bundle with this property.  An easy
induction shows that
$W_d$ is stable. To see this, note that $W_d$ is stable if and only if
every proper subsheaf $J$ of $W_d$ has degree at most zero. But if $J$ is a
proper subsheaf of
$W_d$ of positive degree, then the image of $J$ in $W_{d-1}$ also has
positive degree, and hence
$J\to W_{d-1}$ is surjective. But since the rank of $J$ is at most $d-1$,
the projection of $J$ to
$W_{d-1}$ is an isomorphism. This says that $W_d$ is a split extension of
$W_{d-1}$ by $\scrO_E$, a contradiction. Thus $W_d$ is stable.

The uniqueness statement is clear in the case of rank one. Now assume
inductively that we have showed that, for $d<n$, every stable bundle of
rank $d$ whose determinant is isomorphic to $\scrO_E(p_0)$ is isomorphic to
$W_d$. Let $W$ be a stable bundle of rank $n$ such that $\det W =
\scrO_E(p_0)$. By stability, $h^1(W) = \dim \Hom (W, \scrO_E) = 0$, and so
$h^0(W) = 1$. If $\scrO_E \to W$ is the map corresponding to a nonzero
section, then by stability the cokernel $Q$ is torsion free. An argument as
in the proof that $W_d$ is stable shows that
$Q$ is stable. If $E$ is smooth, then $Q$ is automatically locally free.
When
$E$ is singular, Lemma 0.4 implies that $W$ is locally isomorphic to
$Q\oplus
\scrO_E$. Thus, if $W$ is locally free, then $Q$ is locally free as well.
Once we know that
$Q$ is locally free, we are done by induction. 

To see the final statement, first note that, since $\deg (W_d\otimes
\lambda )=1$, we have by definition that
$$h^0(W_d\otimes \lambda ) - h^1(W_d\otimes \lambda ) = 1.$$ It will thus
suffice to show that $h^1(W_d\otimes \lambda ) = 0$. By Serre duality, 
$$h^1(W_d\otimes \lambda )=\dim \Hom (W_d\otimes \lambda , \scrO_E) =\dim
\Hom(W_d, \lambda\spcheck).$$ Since $\lambda\spcheck$ is also a rank one
torsion free sheaf of degree zero, $\Hom(W_d, \lambda\spcheck) =0$ by
stability.
\endproof

\remark{Exercise} We have defined $\Cal E$ in the previous section as the
rank
$n-1$ vector bundle which is the kernel of the evaluation map
$H^0(\scrO_E(np_0))\otimes \scrO_E \to \scrO_E(np_0)$. Show that
$$\Cal E \cong W_{n-1}\spcheck \otimes \scrO_E(-p_0).$$
\endremark

Now we are ready to see how extensions of the $W_d$ can be used to make
regular semistable bundles.

\theorem{3.2} Let $V$ be an extension of the form
$$0 \to W_d\spcheck \to V \to W_{n-d} \to 0.$$ Then:
\roster
\item"{(i)}"  $V$ has trivial determinant
\item"{(ii)}" $V$ is semistable if and only if the above extension is not
split. In this case $V$ is regular.
\item"{(iii)}" Suppose that  $V$ is semistable, i\.e\. that the above
extension is not split. Then $\dim
\Hom (V,V) = n$ and $\Hom(V,V)$ is an abelian $\Cee$-algebra. Moreover,
every homomorphism
$W_d\spcheck \to V$ is of the form $\phi
\circ \iota$, where $\phi \in \Hom (V,V)$ and $\iota$ is the given inclusion
$W_d\spcheck \to V$. If $V$ and $V'$ are given as extensions as above, then
$V$ and $V'$ are isomorphic if and only if their extension classes in
$\Ext^1(W_{n-d},W_d\spcheck)$ are multiples of each other.
\item"{(iv)}" If $V$ is  a regular semistable  vector bundle of rank $n>1$
with trivial determinant, then $V$ can be written as an extension as above. 
\item"{(v)}" If $V$ is a nontrivial extension of $W_{n-d}$ by
$W_d\spcheck$, and
$ad(V)$ is the sheaf of trace free endomorphisms of $V$, then $H^0(ad(V))
\cong
\Ker\{\, \Hom (W_d\spcheck, W_{n-d})\to H^1(Hom (W_d\spcheck,
W_d\spcheck))\cong
\Cee\,\}$ and $H^1(ad(V)) \cong \Ext^1(W_{n-d},W_d\spcheck)/\Cee\xi$, where
$\xi$ is the extension class corresponding to $V$.
\endroster
\endstatement
\proof (i) This is clear since
$\det W_d \cong\det W_{n-d}$.
\smallskip

\noindent (ii) If $V$ is unstable, let $W$ be the maximal destabilizing
subsheaf. Then $W$ is stable of positive degree and rank $r$ for some
$r< n$. Since $\Hom (W, W_d\spcheck) = 0$, the induced map
$W \to W_{n-d}$ is nonzero. Now it is easy to see by the stability of
$W_s$ that if there is a nonzero map
$W \to W_s$, where $W$ has positive degree and rank $r$, then $r\geq s$,
and every nonzero such map is surjective. (From this it follows in
particular that, for $r\geq s$, $\Hom(W_r, W_s) \cong \Hom (W_s, W_s) =
\Cee$.) If $r> n-d$, the kernel of the map $W
\to W_{n-d}$ is a subsheaf of degree at least zero of $W_d\spcheck$, and
since
$W_d\spcheck$ is a stable bundle of degree $-1$, the kernel is zero. Hence
$W\cong W_{n-d}$, which means that the extension is split. Conversely, if
the extension is split then $V$ is unstable.

Next we show that $V$ is regular. Since
$W_{n-d}$ is a stable bundle of degree
$1$, $\Hom (W_{n-d}, \lambda) =0$ for every rank one torsion free sheaf
$\lambda$ of degree zero. Moreover, with $\lambda$ as above, 
$h^0((W_d\otimes \lambda) =
\dim \Hom (W_d\spcheck, \lambda) = 1$ by the last sentence in Lemma 3.1. 
Thus
$\dim \Hom (V,\lambda)\leq 1$ for every $\lambda$ of degree zero, so that,
by (1.14),  $V$ is regular.

\smallskip

\noindent (iii) Consider the exact sequence
$$\Hom (W_{n-d}, W_d\spcheck) \to \Hom (W_{n-d}, V) \to \Hom (W_{n-d},
W_{n-d})
\to \Ext ^1(W_{n-d}, W_d\spcheck).$$ Note that $\Hom (W_{n-d}, W_d\spcheck)
= 0$ by stability. Since $W_{n-d}$ is stable, it is simple, and so $\Hom
(W_{n-d}, W_{n-d}) =
\Cee \cdot \Id$. But the image of $\Id$ in $\Ext ^1(W_{n-d}, W_d\spcheck)$
is the extension class. Since this class is nonzero,
$\Hom (W_{n-d}, V) = 0$ as well.

Next consider the exact sequence
$$0= \Hom (W_{n-d}, V) \to \Hom (V, V) \to \Hom (W_d\spcheck, V) \to
H^1(W_{n-d}\spcheck \otimes V).$$ Since $\Hom (W_{n-d}, V) = 0$, the map
$\Hom (V, V) \to \Hom (W_d\spcheck, V)$ is an injection. Moreover, we have a
commutative diagram
$$\minCDarrowwidth{.175 in}
\CD
\Cee\cdot \Id @>>> \Hom (V,V) @>>> H^0(ad(V)) @>>> 0\\
 @V{\cong}VV @VVV @VVV @.\\
 \Hom(W_d\spcheck, W_d\spcheck) @>>> \Hom (W_d\spcheck, V) @>>> \Hom
(W_d\spcheck, W_{n-d}) @>>> H^1(Hom(W_d\spcheck, W_d\spcheck)).
\endCD$$ Thus we see that $H^0(ad(V)) \cong
\Ker\{\, \Hom (W_d\spcheck, W_{n-d})\to H^1(Hom (W_d\spcheck,
W_d\spcheck))\cong
\Cee\,\}$ and by duality $H^1(ad(V)) \cong
\Ext^1(W_{n-d},W_d\spcheck)/\Cee\xi$, where
$\xi$ is the extension class corresponding to $V$. This proves (v).

Let us assume that $\dim \Hom (V,V) \geq n$, which we have already checked
in case $E$ is smooth. (We will establish this for singular curves after
proving Part (iv), as well as checking the fact that $\Hom (V,V)$ is
abelian. These results are not used in the proof of Part (iv) of the
theorem.) If we can show that $\dim \Hom (W_d\spcheck, V) = n$, then $\Hom
(V, V) \to
\Hom (W_d\spcheck, V)$ is an isomorphism, and in particular $\dim
\Hom (V,V) = n$ as well.

We compute the dimension of $\Hom(W_d\spcheck,V)$. Consider the exact
sequence
$$0 \to \Hom (W_d\spcheck, W_d\spcheck) \to \Hom (W_d\spcheck, V) \to
\Hom (W_d\spcheck, W_{n-d}) 
\to H^1(Hom (W_d\spcheck, W_d\spcheck)) $$ Since $W_d\spcheck$ is stable,
$\dim
\Hom (W_d\spcheck, W_d\spcheck) = 1$. Next we claim that
$$\dim \Hom (W_d\spcheck, W_{n-d}) = h^0(W_d\otimes W_{n-d}) = n.$$

\claim{3.3} If $E$ is a Weierstrass cubic, then
$h^0(W_d\otimes W_{n-d}) = n$ and $h^1(W_d\otimes W_{n-d}) = 0$. Dually,
$h^0(W_d\spcheck\otimes W_{n-d}\spcheck) = 0$ and
$h^1(W_d\spcheck\otimes W_{n-d}\spcheck) = n$.
\endstatement
\proof If
$d=1$, this follows from the exact sequence
$$0 \to W_{n-1} \to W_{n-1}\otimes \scrO_E(p_0) \to (\Cee _{p_0}) ^{n-1}
\to 0,$$ together with the fact that $h^1(W_{n-1}) = 0$ by stability. The
general case follows by induction on $n$, by tensoring the exact sequence
$$0 \to \scrO_E \to W_d \to W_{d-1}\to 0$$ by $W_{n-d}$.
\endproof

By Riemann-Roch,
$h^1(Hom (W_d\spcheck, W_d\spcheck))=1$. Thus by counting dimensions, to
show that $\dim \Hom (W_d\spcheck, V) = n$ it will suffice to show that
$\Hom (W_d\spcheck, W_{n-d}) \to H^1(Hom (W_d\spcheck, W_d\spcheck))$ is
surjective. Equivalently we must show that the map from  $H^1(Hom
(W_d\spcheck, W_d\spcheck))
\to H^1(Hom (W_d\spcheck, V))$ is zero. But this map is dual to the map
$\Hom (V, W_d\spcheck) \to \Hom (W_d\spcheck, W_d\spcheck) = \Cee
\cdot \Id$. A lifting of
$\Id$ to a homomorphism $V\to W_d\spcheck$ would split the exact sequence,
contrary to assumption. This completes the proof of all of Part (iii)
except for the last sentence. 

We turn to the last statement in Part (iii). If $V$ is a split extension,
then it is unstable and so $V'$ is unstable and therefore a split extension
as well. Thus we may suppose that
$V$ and $V'$ are nontrivial extensions of the given type and that
$\psi\colon V'\to V$ is an isomorphism.  Using $\psi$ to identify
$V$ and
$V'$, suppose that we are given two inclusions $\iota _1, \iota _2\:
W_d\spcheck
\to V$ such that both quotients are isomorphic to $W_{n-d}$. By the first
part of (iii), there is an endomorphism
$A$ of $V$ such  that $A\circ \iota _1 = \iota_2$.   Since $W_{n-d}$ is
simple, the induced map on the quotient $W_{n-d}$ factors must be a
multiple $\alpha\in \Cee$ of the identity.  This multiple $\alpha$ cannot
be  zero, since  otherwise $A$ would define a splitting of the extension
corresponding to $\iota _2$.  In particular,
$A$ is an automorphism of $V$. Furthermore, we see that the extension class
for
$V'$ is $\alpha$ times the extension class for $V$.
 This completes the proof of (iii).

\smallskip

To prove Part (iv) of the theorem, given a semistable $V$, we seek
subbundles of
$V$ isomorphic to $W_d\spcheck$ such that the quotient is isomorphic to
$W_{n-d}$.

\lemma{3.4}  Fix $d>0$. For any $r>0$ and any line bundle $\lambda$ of
degree zero there is a map 
$$W_d\spcheck \to I_r(\lambda)$$ whose image is not contained in a proper
degree zero subsheaf of $I_r(\lambda)$. Likewise, for any  strongly
indecomposable, degree zero, semistable bundle $I(\Cal F)$ concentrated at
the singular point  of a singular curve, there is a map $W_d\spcheck\to
I(\Cal F)$ whose image is not contained in a  proper degree zero subsheaf.
\endstatement

\proof We consider  case of $I_r(\lambda)$ first. It suffices by (1.8) to
show that there is a map $W_d\spcheck\to I_r(\lambda)$ whose image is not
contained in
$F_{r-1}\cong I_{r-1}(\lambda)$. Tensoring the exact sequence
$$0\to I_{r-1}(\lambda)\to I_r(\lambda)\to \lambda\to 0$$ with $W_d$, we
see that there is an exact sequence
$$0 \to \Hom (W_d\spcheck,I_{r-1}(\lambda))\to \Hom
(W_d\spcheck,I_r(\lambda))\to
\Hom (W_d\spcheck,\lambda)\to 0$$  By the last statement in (3.1), there is
a nonzero element of $\Hom (W_d\spcheck,\lambda)$, and by induction on $r$,
$H^1(W_d\otimes I_{r-1}(\lambda)) =0$. Thus there is a map
$W_d\spcheck\to I_r(\lambda)$ not in the image of a homomorphism into
$I_{r-1}(\lambda)$.

Now let us consider the case of a strongly indecomposable bundle
$I(\Cal F)$. Since every semistable bundle concentrated at the singular
point is filtered with associated gradeds isomorphic to ${\Cal F}$, we have
a short exact sequence 
$$0\to X\to I(\Cal F)\to {\Cal F}\to 0.$$ Direct cohomology computations as
above show that there is a map
$W_d\spcheck\to I(\Cal F)$ which has nontrivial image in the quotient
${\Cal F}$. Clearly, the image of this map is not contained in $X$. But, by
(1.11), every proper degree zero subsheaf of $I(\Cal F)$ is contained in
$X$, proving the result in this case as well.
\endproof

We can generalize (3.4) to every regular semistable bundle
$V$.

\corollary{3.5}  Let $V$ be a regular semistable bundle and let $d$ be a
positive integer. Then there is a map $W_d\spcheck\to V$ whose image is not
contained in any proper degree zero subsheaf of $V$. 
\endstatement
\proof This is immediate from the previous result and the fact that $V$
decomposes uniquely as a direct sum 
$\bigoplus_iI_{r_i}(\lambda_i)\oplus I(\Cal F)$, where the $\lambda_i$ are
pairwise distinct line bundles of degree zero and $I(\Cal F)$ is a strongly
indecomposable bundle concentrated at the node.  Since the $\lambda_i$ are
pairwise distinct, any degree zero subsheaf of $V$ is a direct sum of
subsheaves of the factors. Thus,  for each summand $I_{r_i}(\lambda_i)$  or
$I(\Cal F)$, choose  a map $W_d\spcheck$ to the corresponding summand whose
image is not contained in any proper degree zero subbundle of the summand.
The induced map of
$W_d\spcheck$ into the direct sum is as desired. 
\endproof

Note that, if instead $\lambda _i = \lambda _j$ for some $i\neq j$, then
there would exist degree zero subsheaves of the direct sum which were not a
direct sum of subsheaves of the summands, and in fact (3.5) always fails to
hold in this case.

Now let us show that the quotient  of a map satisfying the conclusions of
(3.5) is
$W_{n-d}$.

\proposition{3.6} Let $V$ be a semistable regular bundle of rank $n$ with
trivial determinant and let $\iota\: W_d\spcheck\to V$ be a map whose image
is not contained in any proper degree zero subsheaf of $V$. If the rank of
$V$ is strictly greater than $d$, then $\iota$ is an inclusion and the
quotient $V/W_d\spcheck$ is isomorphic to $W_{n-d}$. Conversely, if $\iota$
is the inclusion of $W_d\spcheck$ in $V$ so that the quotient
$V/W_d\spcheck$ is isomorphic to $W_{n-d}$, then the image of $\iota$ is
not contained in a proper subsheaf of degree zero.
\endstatement
\proof Let $V$ have rank $n\ge d+1$,  and  suppose that we have a map
$\iota\:W_d\spcheck\to V$  whose image is not contained in a proper degree
zero subsheaf of $V$. In particular,
$\iota$ is nontrivial. If $\iota$ is not injective, then by the stability of
$W_d\spcheck$, the image of $\iota$ is a subsheaf of $V$ of rank
$\leq d-1$ and degree $\geq 0$, and hence is a proper subsheaf of
$V$ of degree zero, contrary to assumption. Likewise, if the cokernel of
$\iota$ is not torsion free, then the image of $\iota$ is contained in a
proper subsheaf of $V$ whose degree is strictly larger than 
$-1$, and thus the degree is at least zero.  This again contradicts our
assumption about the map and the  fact that the rank of $V$ is at least
$d+1$.  Thus $\iota$ is injective and its cokernel $W$ is torsion free.
Using (0.4), $W$ is locally a direct summand of $V$, and thus
$W$ is locally free.  It follows that $W$ is a rank
$n-d$ vector bundle whose determinant is
$\scrO_E(p_0)$. 

To conclude that the quotient $W$ is isomorphic to $W_{n-d}$, it suffices
to show that $W$ is stable.  If $W$ is not stable, then there is a proper
subsheaf $U$ of
$W$ with degree at least one. Let
$U''\subset V$ be the preimage of $U$.  The degree of $U''$ is at least
zero, and hence, by the semistablility of $V$  is of degree zero. Clearly,
$U''$  contains the image of $\iota$.  Hence by our hypothesis on
$\iota$, $U''=V$, and consequently, $U=W$. This is a contradiction, so that
$W$ is stable. 

Finally we must show that, if $V$ is written as an extension of $W_{n-d}$ by
$W_d\spcheck$, then the subbundle $W_d\spcheck$ cannot be contained in a
proper subsheaf $U$ of $V$ of degree zero. If $U$ is a such a subsheaf, then
$U/W_d\spcheck$ would be a proper subsheaf of $W_{n-d}$ of degree at least
one, contradicting the stability of $W_{n-d}$.
\endproof

 Corollary 3.5 and Proposition 3.6 show that any regular semistable bundle
over $E$ can be written as an extension of $W_{n-d}$ by $W_d\spcheck$. 
This completes the proof of  Part (iv).

Now let us return to the  point in the proof of (iii) where it is  claimed
that
$\dim \Hom (V,V) \geq n$ for all $V$ which are given as a nonsplit
extension of
$W_d$ by $W_{n-d}\spcheck$. In order to establish this result, we first 
describe the space of all such extensions, which is an immediate
consequence of (3.3):

\lemma{3.7} The space $\Ext^1(W_{n-d}, W_d\spcheck) = H^1(W_{n-d}\spcheck
\otimes W_d\spcheck)$ has dimension $n$, and thus the associated projective
space is a $\Pee ^{n-1}$. \qed
\endstatement

By general properties, there is a universal extension
 $\bold U(d;n)$ over $\Pee_d^{n-1}\times E$ of the form
$$0 \to \pi _2^*W_d\spcheck \otimes \pi _1^*\scrO_{\Pee_d ^{n-1}}(1)
\to \bold U(d;n) \to \pi _2^*W_{n-d} \to 0,$$ with the restriction of
$\bold U(d;n)$ restricted to any slice
$\{x\}\times E$ being isomorphic the bundle $V$ given by the line
$\Cee\cdot x\subset \Ext^1(W_{n-d},W_d\spcheck)$. When $n$ is clear from the
context, we shall abbreviate $\bold U(d;n)$ by $\bold U(d)$.

Next we claim that there is a nonempty open subset of $\Pee ^{n-1}_d$ such
that for $V$ a vector bundle corresponding to a point of this subset,
$\dim \Hom (V,V) = n$. In fact, suppose that $V = \bigoplus _{i=1}^n\lambda
_i$, where the
$\lambda _i$ are distinct line bundles of degree zero whose product is
trivial. By Part (iv) of the theorem, $V$ can be written as a nonsplit
extension of
$W_d$ by $W_{n-d}\spcheck$, and we have seen that $\dim \Hom (V,V) = n$. A
straightforward argument by counting dimensions shows that the set of such
$V$ is an open subset of $\Pee^{n-1}_d$; indeed, we will identify this set
more precisely in (3.17) below as corresponding to the set of all sections
in $|np_0|$ consisting of $n$ smooth points on $E$. Thus there is a
nonempty open subset of bundles $V$ such that $\dim \Hom (V,V) = n$. By
upper semicontinuity applied to the bundle $Hom(\bold U(d;n), \bold
U(d;n))$ over
$\Pee_d^{n-1}\times E$, it follows that $\dim \Hom (V,V) \geq n$ for all
bundles
$V$ corresponding to a point of $\Pee_d^{n-1}$.

Finally, we must show that $\Hom (V,V)$ is abelian. Using the universal
extension
$\bold U(d)$ as above, we can fit together the $\Hom (V,V)$ to a rank $n$
vector bundle $\pi_1{}_*Hom (\bold U(d), \bold U(d))$, which is a coherent
sheaf of algebras over $\Pee^{n-1}$. Consider the map
$$\pi_1{}_*Hom (\bold U(d), \bold U(d))\otimes \pi_1{}_*Hom (\bold U(d),
\bold U(d)) \to \pi_1{}_*Hom (\bold U(d), \bold U(d))$$ defined by
$(A,B)\mapsto AB-BA$. Since $\Hom(V,V)$ is abelian for $V$ in a Zariski
open subset of
$\Pee^{n-1}_d$, namely those $V$ which are a direct sum of $n$ distinct line
bundles of degree zero, this map is identically zero. By base change, the
fiber of
$\pi_1{}_*Hom (\bold U(d), \bold U(d))$ at a point corresponding to $V$ is
$\Hom(V,V)$. Thus
$\Hom (V,V)$ is abelian.
\endproof

The following was checked directly in Lemma 1.13 if $E$ is smooth, but is
by no means obvious in the singular case:

\corollary{3.8} Let $V$ be a regular semistable bundle of rank $n$ over a
Weierstrass cubic. Then: 
\roster
\item"{(i)}" $\Hom(V,V)$ is an abelian $\Cee$-algebra of rank $n$.
\item"{(ii)}" The dual bundle $V\spcheck$ is a regular semistable bundle.
\item"{(iii)}" For all rank one torsion free sheaves $\lambda$ of rank zero
on
$E$, $\dim \Hom (\lambda, V) \leq 1$.
\endroster
\endstatement

\proof The first part is immediate from Parts (iv) and (iii) of Theorem
3.2. (ii) is clear since if $V$ is a nonsplit extension of $W_{n-d}$ by
$W_d\spcheck$, then
$V\spcheck$ is a nonsplit extension of $W_d$ by $W_{n-d}\spcheck$. (iii)
follows from (ii), since $\Cal F\spcheck \cong \Cal F$ and  $\Hom (\lambda,
V)\cong \Hom ( V\spcheck,\lambda\spcheck)$ for all rank one torsion free
sheaves $\lambda$.
\endproof

\remark{Question} For if $V$ is a semistable bundle of degree zero whose
support is concentrated at a smooth point of $E$, then $V$ is regular if and
only if $\dim \Hom (V,V)= \operatorname{rank}V$. Does this continue to hold
at the singular point of a singular curve? For $V$ strongly indecomposable,
what is the structure of the algebra $\Hom(V,V)$?
\endremark

\ssection{3.2. Relationship between the constructions for various $d$.}

For each $d$ with $1\le d<n$ we have a family of regular semistable bundles
parametrized by the projective space
$\Pee_d^{n-1}=\Pee(\Ext^1(W_{n-d},W_d\spcheck) )$, and given as a universal
extension
$$0 \to \pi _2^*W_d\spcheck \otimes \pi _1^*\scrO_{\Pee_d ^{n-1}}(1)
\to \bold U(d) \to \pi _2^*W_{n-d} \to 0.$$ In this section we shall
identify the $\Pee_d^{n-1}$ for the various
$d$, although under this identification the bundles  $\bold U(d)$ are
different for different $d$. Using the universal bundle $\bold U(d)$ and
Theorem 1.5, there is a morphism $\Pee_d^{n-1} \to |np_0|$, which is easily
checked to be of degree one and thus an isomorphism. Thus all of the
$\Pee_d^{n-1}$ are identified with $|np_0|$ and hence with each other, but
we want to find a direct identification here so as to be able to compare
universal bundles.

\lemma{3.9} Let $d,n-d\ge 1$. The natural injection
$$W_d\spcheck\otimes W_{n-d}\spcheck\to W_{d+1}\spcheck\otimes
W_{n-d}\spcheck $$ induces an injective map on $H^1$.  The image of this
map on $H^1$ is the kernel of the  map induced by the tensor products of the
projections 
$$H^1(W_{d+1}\spcheck\otimes W_{n-d}\spcheck)\to
H^1(\scrO_E\otimes\scrO_E)=H^1(\scrO_E).$$  The extensions $X$ of $W_{n-d}$
by $W_{d+1}\spcheck$ which are in the image of the above map are exactly
the extensions $X$ such that $\Hom (X, \scrO_E)\neq 0$. 
\endstatement

\proof We have a short exact sequence
$$0\to  W_d\spcheck\otimes W_{n-d}\spcheck \to W_{d+1}\spcheck\otimes
W_{n-d}\spcheck \to
\scrO_E\otimes W_{n-d}\spcheck\to 0.$$ By (3.3), all the $H^0$ terms
vanish. Thus, the injectivity of the map on $H^1$ is immediate.
Furthermore, the image is identified with the kernel of the map 
$$H^1(W_{d+1}\spcheck\otimes W_{n-d}\spcheck)\to H^1(\scrO_E\otimes
W_{n-d}\spcheck).$$  The last term is one-dimensional and the projection 
$\scrO_E\otimes W_{n-d}\spcheck\to \scrO_E\otimes\scrO_E$ induces an
isomorphism on $H^1$. 

Finally, a bundle $X$ corresponds to an extension in the image of the map on
$H^1$'s if and only if
$X$ is the pushout of an extension of $W_{n-d}$ by $W_d\spcheck$ under the
inclusion
$W_d\spcheck \to W_{d+1}\spcheck$. Thus, if $X$ is the image of an
extension $V$, the quotient of the inclusion $V\to X$ is $\scrO_E$.
Conversely, if there is a nontrivial map $X\to \scrO_E$, then the induced
map $W_{d+1}\spcheck \to \scrO_E$ is nonzero and thus surjective, and the
kernel of $X\to \scrO_E$ is then an extension $V$ of $W_{n-d}$ by
$W_d\spcheck$ such that $X$ is the pushout of $V$.
\endproof

The symmetry of the situation with respect to the two factors leads
immediately to the following.

\corollary{3.10} If $n-d\ge 2$, then  the natural inclusions of bundles
induce the maps 
$$H^1(W_d\spcheck\otimes W_{n-d}\spcheck)\to H^1(W_{d+1}\spcheck\otimes
W_{n-d}\spcheck)$$ and
$$H^1(W_{d+1}\spcheck\otimes W_{n-d-1}\spcheck)\to
H^1(W_{d+1}\spcheck\otimes W_{n-d}\spcheck)$$ which are injections with the
same images. In particular, this produces a natural identification of
$\Ext^1(W_{n-d},W_d\spcheck)$ with
$\Ext^1(W_{n-d-1},W_{d+1}\spcheck)$, and hence of the projective spaces
$\Pee_d^{n-1}\cong \Pee_{d+1}^{n-1}$. Finally, an extension $X$ of
$W_{n-d}$ by $W_{d+1}\spcheck$ is obtained from an extension $V$ of
$W_{n-d-1}$ by
$W_{d+1}\spcheck$ via pullback if and only if $\Hom (X, \scrO_E)\neq 0$.
\qed
\endstatement

Now let us see how the bundles described by extensions which are identified
under this isomorphism are related.

\proposition{3.11} Let $\epsilon_d\in \Ext^1(W_{n-d},W_d\spcheck)$ and
$\epsilon_{d+1}\in \Ext^1(W_{n-d-1},W_{d+1}\spcheck)$ be nonzero classes
that correspond under the identification given in Corollary \rom{3.10}. 
Let $V$ and $V'$, respectively, be the total spaces of these extensions. 
Then
$V$ and $V'$ are isomorphic bundles.
\endstatement

\proof Let $X$ be the bundle of rank $n+1$ obtained by pushing out the
extension $V$ by the map $W_d\spcheck\to W_{d+1}\spcheck$. Clearly, we have
a short exact sequence
$$0\to V\to X\to \scrO_E\to 0.$$ Similarly, let $X'$ be the rank
$n+1$-bundle obtained by pulling back the extension $V'$ along the map
$W_{n-d}\to W_{n-d-1}$. Dually, we have an exact sequence
$$0\to \scrO_E\to X'\to V'\to 0.$$ It follows easily from Theorem 3.2 that
writing $V=Y\oplus I_r(\scrO_E)$ with $H^0(Y)=0$, we have $X\cong Y\oplus
I_{r+1}(\scrO_E)$. Similarly, writing $V'=Y'\oplus I_s(\scrO_E)$ we have
$X'\cong Y'\oplus I_{s+1}(\scrO_E)$.

The fact that $\epsilon_d$ and $\epsilon_{d+1}$ are identified means that
the extensions for $X$ and $X'$ are isomorphic.  In particular, $X$ and
$X'$ are isomorphic bundles. This implies that $r=s$ and that $Y$ and $Y'$
are isomorphic. But then $V$ and
$V'$ are isomorphic as well.
\endproof

Notice that the isomorphism produced by the previous result is canonical on
$Y\subseteq V$ but is not canonical on the $I_r(\scrO_E)$ factor. We shall
see later that the families of bundles $\bold U(d)$ and $\bold U(d+1)$ are
not isomorphic, which means that there cannot be a canonical isomorphism in
general between corresponding bundles. In practice, this means the
following: suppose that $V$ is given as an extension of $W_{n-d}$ by
$W_d\spcheck$, with $n-d >1$. Then $W_{n-d}$ has the distinguished subbundle
isomorphic to $\scrO_E$. Let $W'$ be the preimage in $V$ of this bundle, so
that
$W'$ is an extension of $\scrO_E$ by $W_d\spcheck$. Then $W'\cong
W_{d+1}\spcheck$ if and only if $h^0(V) = 0$, if and only if the support of
$V$ does not contain
$\scrO_E$, but otherwise $W' \cong W_d\spcheck\oplus \scrO_E$.

The direct comparison of the extension classes given above leads to a
comparison of universal bundles. 

\theorem{3.12} Let $H$ be the divisor in $\Pee^{n-1}_d$ such that, if $x\in
H$ and
$V$ is the corresponding extension, then $h^0(V) = 1$. Let $i\:H\to
\Pee^{n-1}_d$ be the inclusion. Then there is an exact sequence 
$$0 \to \bold U(d-1) \to \bold U(d) \to (i\times \Id)_*\scrO_{H\times E}(1)
\to 0,$$ which expresses $\bold U(d-1;n)$ as an elementary modification of
$\bold U(d;n)$. Moreover, this elementary modification is unique in an
appropriate sense.
\endstatement
\proof Let $H'$ be the hyperplane in $\Pee^n_d$ corresponding to the set of
extensions $X$ of $W_{n-d+1}$ by $W_d\spcheck$ such that $\Hom (X, \scrO_E)
\neq 0$. By the last statements of (3.9) and (3.10),
$H'$ is the image of
$\Pee^{n-1}_d$ in
$\Pee^n_d$ as well as the image of $\Pee^{n-1}_{d+1}$. By base change,
$\pi_1{}_*Hom (\bold U(d;n+1)|H'\times E,
\scrO_{H'\times E})$ is a line bundle on $H'$. By looking at the exact
sequence
$$\gather 0 \to Hom (\pi _2^*W_{n-d+1}, \scrO_{\Pee^n_d\times E}) \to Hom
(\bold U(d;n+1)|H'\times E,
\scrO_{\Pee^n_d\times E}) \to \\
\to Hom (\pi _2^*W_d\spcheck\otimes \pi _1^*\scrO_{\Pee^n_d}(1),
\scrO_{\Pee^n_d\times E})\to 0,
\endgather$$ and restricting to $H'\times E$, we see that this line bundle
is $\scrO_{H'\times E}(-1)$. Thus there is a surjection
$$\bold U(d;n+1)|H'\times E \to \scrO_{H'\times E}(1).$$ We claim that the
kernel of this surjection is identified with $\bold U(d-1;n)$. In fact, if
$\bold U'$ denotes the kernel, we have a commutative diagram with exact
rows and columns
$$\minCDarrowwidth{.3 in}
\CD @. 0 @. 0 @. @.\\ @. @VVV @VVV @. @.\\ 0@>>>
\pi_2^*W_{d-1}\spcheck\otimes \scrO_{H'\times E}(1) @>>>\bold U'
@>>>\pi_2^*W_{n-d+1} @>>> 0\\ @. @VVV @VVV @| @.\\ 0@>>>
\pi_2^*W_d\spcheck\otimes \scrO_{H'\times E}(1) @>>>\bold U(d;n+1)|H'\times
E @>>>\pi_2^*W_{n-d+1} @>>> 0\\ @. @VVV @VVV @. @.\\ @. \scrO_{H'\times
E}(1) @=  \scrO_{H'\times E}(1) @. @.\\ @. @VVV @VVV @. @.\\ @. 0 @. 0 @.
@.\\
\endCD$$ and tracing through the diagram identifies $\bold U'$ with $\bold
U(d-1;n)$, compatibly with the identification of $H'$ with $\Pee^{n-1}_d$.

Now we can also consider the line bundle $\pi _1{}_*(\bold
U(d;n+1)|H'\times E)$. A very similar argument shows that this line bundle
is isomorphic to $\scrO_{H'}$, and that the quotient $\bold U''$ of $\bold
U(d;n+1)|H'\times E$ via the natural map 
$$\scrO_{H'\times E} = \pi^*\pi _1{}_*(\bold U(d;n+1)|H'\times E) \to \bold
U(d;n+1)|H'\times E$$ is isomorphic to $\bold U(d;n)$. Putting these two
constructions together, we see that we have found
$$\bold U(d-1;n) \to \bold U(d;n+1)|H'\times E \to \bold U(d;n).$$ Away
from $H$, which is the image of $\Pee^{n-2}_d$ in $\Pee^{n-1}_d$, the
inclusion $\bold U(d-1;n) \to \bold U(d;n)$ is an equality. To summarize,
then, there is a commutative diagram
$$\CD @.  @. 0 @. @.\\ @. @. @VVV @. @.\\ @. @. \bold U(d-1;n) @= \bold
U(d-1;n) @.\\  @. @. @VVV @VVV @.\\ 0@>>>  \scrO_{H'\times E} @>>>\bold
U(d;n+1)|H'\times E @>>>\bold U(d;n) @>>> 0\\
 @. @| @VVV @. @.\\ @. \scrO_{H'\times E} @>>> \scrO_{H'\times E}(1) @. @.\\
@. @. @VVV @. @.\\ @.  @. 0 @. @.\\
\endCD$$ The map $\scrO_{H'\times E} \to \scrO_{H'\times E}(1)$ can only
vanish along $H$. Thus it vanishes exactly along $H$, and the quotient of
$\bold U(d;n)$ by $\bold U(d-1;n)$ is a line bundle supported on $H\times
E$. By what we showed above for
$\bold U(d;n+1)$, this line bundle is necessarily 
$\scrO_{H\times E}(1)$. Hence we have found an exact sequence
$$0 \to \bold U(d-1;n) \to \bold U(d;n) \to (i\times \Id)_* \scrO_{H\times
E}(1) \to 0,$$ realizing $\bold U(d-1;n)$ as an elementary modification of
$\bold U(d;n)$. The uniqueness is straightforward and left to the reader.
This completes the proof of Theorem 3.12.
\endproof

\ssection{3.3. Comparison of coarse moduli spaces.}

We have succeeded in identifying the $\Pee_d^{n-1}$ for the various
$d$, $1\le d<n$, in a purely cohomological way and in showing that
extension classes in different groups which are identified produce
isomorphic vector bundles. Next we wish to identify these projective spaces
with   the projective space $|np_0|$ which is the parameter space of regular
semistable rank $n$ bundles with trivial determinant in the spectral cover
construction of these bundles. Of course, the existence of the bundle $\bold
U(d)$ and Theorem 1.5 give us one such identification. However, although we
shall not need this in what follows, we want to find a direct cohomological
comparison between
$\Ext^1(W_{n-d}, W_d\spcheck)$ and $H^0(\scrO_E(np_0))$.  We have two
identifications: one purely cohomological and the other using the bundles to
identify the spaces.  We shall show that these identifications agree.

Let us  begin with the purely cohomological identification. (We will be
pedantic here about identifying various one-dimensional vector spaces with
$\Cee$ in order to carry out the discussion in families in the next
section.)

\proposition{3.13} Let $H^0_{n-1}=H^0(\scrO_E(p_0)\otimes W_{n-1})$. It is
an $n$-dimensional vector space. Let $D=H^1(\scrO_E)$ be the dualizing line.
Let
$$I\colon H^1(\scrO_E(-p_0)\otimes W_{n-1}\spcheck)\to
H^0(\scrO_E(np_0))\otimes
\det(H^0_{n-1})^{-1}\otimes D$$ be the composition
$$\gather H^1(\scrO_E(-p_0)\otimes W_{n-1}\spcheck) @>{S}>>
(H^0_{n-1})^*\otimes D @>{A}>> \bigwedge^{n-1}H^0_{n-1}\otimes
\operatorname{det}(H^0_{n-1})^{-1}\otimes D  \to \\
@>{ev\otimes\Id\otimes\Id}>> H^0(\det(\scrO_E(p_0)\otimes W_{n-1}))\otimes 
\det(H^0_{n-1})^{-1}\otimes D \\ = H^0(\scrO_E(np_0))\otimes
\det(H^0_{n-1})^{-1}\otimes D, 
\endgather$$ where $S$ is Serre duality,  $A$ is the map induced by taking
adjoints from the natural pairing
$$H^0_{n-1}\otimes \bigwedge^{n-1}H^0_{n-1}\to \det(H^0_{n-1}),$$
$ev$ is the map 
$$ev\colon \bigwedge^{n-1}H^0(\scrO_E(p_0)\otimes W_{n-1})\to
H^0(\bigwedge^{n-1}(\scrO_E(p_0)\otimes W_{n-1})).$$  Then $I$ is an
isomorphism.
\endstatement

\proof On general principles $S$ is an isomorphism.  Since $H^0_{n-1}$ is
$n$-dimensional, the adjoint map $A$ is clearly an isomorphism.  What is
less obvious that the map $ev$ is an isomorphism, which follows from the
next claim.

\claim{3.14} The evaluation map
$$ev\colon \bigwedge^{n-1}H^0(\scrO_E(p_0)\otimes W_{n-1})\to
H^0(\bigwedge^{n-1}(\scrO_E(p_0)\otimes W_{n-1}))$$ is an isomorphism.
\endstatement

\proof We prove this by induction. The case $n=2$ is clear since $ev$ is
the identity.

Assume the result for $n-1\ge 2$. There is a short exact sequence
$$0\to \scrO_E(p_0)\to \scrO_E(p_0)\otimes W_{n-1}\to
\scrO_E(p_0)\otimes W_{n-2}\to 0\tag{$*$}$$  leading to a short exact
sequence
$$0\to H^0(\scrO_E(p_0))\to H^0(\scrO_E(p_0)\otimes W_{n-1})\to
H^0(\scrO_E(p_0)\otimes W_{n-2})\to 0,$$ since by (3.3)  the $H^1$ terms
vanish. Since the first term has dimension one, we have a short exact
sequence
$$0\to H^0(\scrO_E(p_0))\otimes \bigwedge^{n-2}H^0_{n-2}\to
\bigwedge^{n-1}H^0_{n-1}\to \bigwedge^{n-1}H^0_{n-2}\to 0.$$ 

Taking determinants in ($*$) yields an isomorphism 
$$\bigwedge^{n-1}(\scrO_E(p_0)\otimes W_{n-1})\cong  \scrO_E(p_0)\otimes
\bigwedge^{n-2}(\scrO_E(p_0)\otimes W_{n-2}).$$ Tensoring the inclusion
$\scrO_E\to \scrO_E(p_0)$ with
$\bigwedge^{n-2}(\scrO_E(p_0)\otimes W_{n-2})$ and using the above
isomorphism leads to a short exact sequence
$$0\to \bigwedge^{n-2}(\scrO_E(p_0)\otimes W_{n-2})\to
\bigwedge^{n-1}(\scrO_E(p_0)\otimes W_{n-1})\to
\bigwedge^{n-1}(\scrO_E(p_0)\otimes W_{n-1})|_{p_0}\to 0.$$ Unraveling the
definitions one sees that the map
$ev$ induces a commutative diagram, with exact columns,
$$\CD 0 @. 0\\
 @VVV   @VVV\\
 H^0(\scrO_E(p_0))\otimes \bigwedge^{n-2}H^0_{n-2} @>>>
 H^0(\scrO_E(p_0))\otimes H^0(\bigwedge^{n-2}(\scrO_E(p_0)\otimes
W_{n-2}))\\
 @VVV   @VVV \\
\bigwedge^{n-1}H^0_{n-1} @>{ev}>> H^0(\bigwedge^{n-1}(\scrO_E(p_0)\otimes
W_{n-1})) \\
 @VVV   @VVV \\
\bigwedge^{n-1}H^0_{n-2}  @>{e}>>  (\bigwedge^{n-1}(\scrO_E(p_0)\otimes
W_{n-1}))|_{p_0}\\
 @VVV   @VVV\\ 0 @. 0 ,
\endCD$$ and that the restriction of $ev$ to the first term is simply the
tensor product of  the identity on 
$H^0(\scrO_E(p_0))$ and the evaluation map, with $n-2$ replacing $n-1$. By
induction, the top horizontal map is an isomorphism.

To finish the proof of (3.14), and thus of (3.13), it suffices by the
$5$-lemma to show that $e$ is an isomorphism. Now $e$ is the $(n-1)$-fold
wedge product of a map
$$\overline{e}\colon H^0(\scrO_E(p_0)\otimes W_{n-2})\to
(\scrO_E(p_0)\otimes W_{n-1})|_{p_0} $$ defined as follows.  For any
section $\psi$ of $\scrO_E(p_0)\otimes W_{n-2}$, lift to a section $\tilde
\psi$ of $\scrO_E(p_0)\otimes W_{n-1}$, and then restrict $\tilde \psi$ to
$p_0$. Thus, it suffices to prove:

\claim{3.15} The map $\overline{e}$ described above is an isomorphism.
\endstatement

\proof First notice that if $\tilde \psi$ and $\tilde \psi'$ are lifts of
$\psi$ to sections of $\scrO_E(p_0)\otimes W_{n-1}$, then they differ by a
section of $\scrO_E(p_0)\subset \scrO_E(p_0)\otimes W_{n-1}$.  But any
section of $\scrO_E(p_0)$ vanishes at $p_0$ so that $\tilde \psi$ and
$\tilde \psi'$ have the same restriction to $p_0$.  This shows that
$\overline{e}$ is well-defined.   From the diagram
$$\CD @.   @.  0 @.   @.  \\ @.  @.  @VVV @. @. \\ @.   @.   W_{n-1} @. 
@.  \\ @.  @.  @VVV @. @. \\ 0 @>>> \scrO_E(p_0) @>>> \scrO_E(p_0)\otimes
W_{n-1}  @>>>
\scrO_E(p_0)\otimes W_{n-2} @>>> 0, \\ @.  @.  @VVV @. @. \\ @.   @. 
\scrO_E(p_0)\otimes W_{n-1}|_{p_0} @.  @.  \\ @.  @.  @VVV @. @. \\ @.  
@.  0 @.   @.  \\
\endCD$$ the fact that all the $H^1$ terms vanish, and the fact that both
$H^0(\scrO_E(p_0))$ and $H^0(W_{n-1})$ are one dimensional, the claim comes
down to the statement that the images in
$H^0(\scrO_E(p_0)\otimes W_{n-1})$ of $H^0(\scrO_E(p_0))$ and of
$H^0(W_{n-1})$ are equal. But we also have a commutative square
$$\CD
\scrO_E @>>> \scrO_E\otimes W_{n-1} \\ @VVV   @VVV \\
\scrO_E(p_0) @>>> \scrO_E(p_0)\otimes W_{n-1}
\endCD$$ with the top arrow and the left arrow inducing isomorphisms on
$H^0$. Claim 3.14 and hence (3.13) now follow.
\endproof
\enddemo\enddemo

Proposition 3.13 and Corollary 3.10 have produced cohomological
isomorphisms between the extension groups
$\Ext^1(W_{n-d},W_d\spcheck)$ and $H^0(\scrO_E(np_0))$. On the other hand,
as remarked previously, from the existence of the bundle
$\bold U(d)\to
\Pee_d^{n-1}\times E$, Theorem 1.5 produces isomorphisms
$$\Phi_d\colon \Pee_d^{n-1}\to |np_0|$$ sending $x\in \Pee_d^{n-1}$ to the
point $\zeta(V_x)$, where $V$ is the the extension determined by the point
$x$. By Proposition 3.11 the identification of  $\Pee_d^{n-1}$ with
$\Pee_{d+1}^{n-1}$  given in Corollary 3.10 identifies $\Phi_d$ with
$\Phi_{d+1}$.  Still, it remains to compare the map $\Phi_1$ with the
projectivization of the map produced by Proposition 3.13.

\proposition{3.16} The map
$\Phi_1\colon \Pee_1^{n-1}\to |np_0|$ is the projectivization of the
identification 
$$I\colon H^1(\scrO_E(-p_0)\otimes W_{n-1}\spcheck)\to
H^0(\scrO_E(np_0))\otimes M,$$ where $M$ is the line
$\det(H^0(\scrO_E(p_0)\otimes W_{n-1}))^{-1}\otimes D$.  In other words, if
$V$ is a nontrivial extension corresponding to $\alpha \in
\Ext^1(W_{n-1},\scrO_E(-p_0))$, then the point $\zeta(V)\in |np_0|$
corresponds to  the line
$$\Cee\cdot I(\alpha)\subset H^0(\scrO_E(np_0))\otimes M.$$ In particular,
$\Phi_1$ is an isomorphism, and hence so is $\Phi_d$ for every $1\le d<n$. 
\endstatement

\proof Let $\bar I\: \Pee^{n-1}_1\to |np_0|$ be the projectivization of $I$.
We begin by determining when a line bundle $\lambda$ of degree zero is in
the support of $V$.

\claim{3.17} Let $V$ be given by an extension class $\alpha\in
H^1(\scrO_E(-p_0)\otimes W_{n-1}\spcheck)$.  Then
$\Hom(V,\lambda)\not=0$ if and only if the image of $\alpha$ in
$H^1(\lambda \otimes W_{n-1}\spcheck)$ under the map induced by the
inclusion $\scrO_E(-p_0)\to \lambda$ is zero.
\endstatement

\proof There is a nonzero map $\scrO_E(-p_0) \to \lambda$, and it is unique
up to a scalar. Let $V'$ be the pushout of the extension $V$  by the map
$\scrO_E(-p_0)\to \lambda$.  Then the pushout extension is trivial, i\.e\.
the image of $\alpha$ in $H^1(\lambda \otimes W_{n-1}\spcheck)$ is zero, if
and only if there is a map $V'\to \lambda$ splitting the inclusion of
$\lambda$ into $V'$. Such a map is equivalent to a map $V\to \lambda$ so
that the composition $\scrO_E(-p_0)\to V\to \lambda$ is the inclusion.  By
(3.6), if $V$ has a nontrivial map to $\lambda$, then this map restricts to
$\scrO_E(-p_0)$ to be the inclusion (since otherwise the image of
$\scrO_E(-p_0)$ would be contained in a proper subbundle of degree zero).
Thus, there is such a map if and only if the
$\lambda$ component of $V$ is nonzero, which is equivalent to the existence
of a nonzero map
$V\to \lambda$. 
\endproof

\claim{3.18} Suppose that $\lambda\cong\scrO_E(q-p_0)$ for some $q\in E$.  
Let $V$ be given by an extension class $\alpha\in H^1(\scrO_E(-p_0)\otimes
W_{n-1}\spcheck)$.  Then
$\lambda$ is in the support of $V$ if and only if $q$ is in the support of
$\bar I(\alpha) \in |np_0|$.
\endstatement

\proof Applying Serre duality $S$ and the adjoint map $A$ to the previous
claim, and tracing through the identifications, we see that
$\lambda=\scrO_E(q-p_0)$ is in the support of
$V$ if and only if the corresponding map
$$\bigwedge^{n-1}H^0_{n-1}\otimes L\to \bigwedge^{n-1}(\scrO_E(p)\otimes
W_{n-1}\otimes L)|_q$$ vanishes, if and only if the  section giving
$I(\alpha)$ in
$\det(\scrO_E(p_0)\otimes W_{n-1})\otimes L$ vanishes at $q$, if and only
if $q$ is  in the support of $\bar I(\alpha)$.
\endproof

Now we can prove Proposition 3.16.  We have two maps
$\bar I$, the projectivization of the linear map $I$, and
$\Phi_1$, mapping $\Pee_1^{n-1}\to |np_0|$. We wish to show $\bar I$ and
$\Phi_1$ are equal.  We know that $\bar I$ is an isomorphism.  Thus, it
suffices to show that, for  an open dense subset $U$ of $|np_0|$,
$\Phi_1(x)=\bar I(x)$ for all $x\in \Phi_1^{-1}(U)$. Choose $U$ to be the
open subset of divisors in $|np_0|$ whose support is $n$ distinct smooth
points of $E$.  Let $x=\sum _ie_i\in U$. The extension determined by
$\alpha =\Phi_1^{-1}(x)$ is a semistable bundle $V$ which is written as
$\bigoplus_i\lambda_i$ for $n$  distinct line bundles $\lambda_i$, where
$\lambda_i=\scrO_E(e_i-p_0)$. According to Claim 3.18, $\bar I(\alpha)$
contains  $e_i,\ 1\le i\le n$ in its support, and hence $\bar I\circ
\Phi_1^{-1}(x)=x$. This completes the proof of Proposition 3.16.
\endproof

\ssection{3.4. From universal bundles to spectral covers.}

We now take another look at the spectral cover construction,and generalize
it to singular curves. Fix an integer $d$, $1\leq d\leq n-1$, and consider
the sheaf
$\bold  A = \pi_1{}_*Hom (\bold U(d), \bold U(d))$. It is a locally free
rank $n$ sheaf of commutative algebras over $\Pee^{n-1}$, and in case $E$
is smooth we have identified this sheaf with $\nu_*\scrO_T$ in (2.4).
(There is nothing special about taking $\bold U(d)$; we could replace
$\bold U(d)$ by any ``universal" bundle, once we know how to construct
one.)  We propose to reverse this procedure:  starting with
$\bold A$, define $T$ to be the corresponding space
$\bold{Spec}\,\bold A$. In particular, this gives a definition of $T$ in
case $E$ is singular.

\lemma{3.19} Let $E$ be a Weierstrass cubic. With $T$ as defined above,
there is a finite flat morphism
$\nu\:T
\to \Pee^{n-1}$. Moreover, $T$ is reduced.
\endstatement
\proof By construction, there is a finite morphism $\nu\: T\to
\Pee^{n-1}$. In fact, since $\bold A$ is locally free, $T$ is flat over
$\Pee^{n-1}$ of degree $n$. 

The fact that $T$ is reduced follows from the fact that $\Pee^{n-1}$ is
reduced and that $T$ is generically reduced, and as such is a general fact
concerning finite flat morphisms. Cover $\Pee^{n-1}$ by affine open sets
$\Spec R$ such that
$\nu^{-1}(\Spec R) = \Spec R'$, where $R'$ is a free rank $n$ $R$-module.
It will suffice to show that $R'$ is reduced for every such $R$. If $f\in
R'$, then $f$ does not vanish on a Zariski open set, since $R'$ is locally
free and $R$ is reduced. Thus the restriction of $f$ to a general fiber of
$\nu$, consisting of $n$ distinct (reduced) points, is nonzero. It follows
that $f^k\neq 0$ for every $k>0$. Thus $R'$ is reduced, and so $T$ is
reduced.
\endproof

In case $E$ is smooth, Theorem 2.4 shows that the $T$ defined above is the
same as the spectral cover $T$ defined in Section 2, although even in this
case it will be useful to define
$T$ as we have above in order to compare $\bold U(d)$ with the bundles
$U_a$. 

The points of $T$ are by definition in one-to-one correspondence with pairs
$(V, \frak m)$, where $V$ is a regular semistable bundle with trivial
determinant and
$\frak m$ is a maximal ideal in $\Hom(V,V)$. Let us describe such maximal
ideals:

\lemma{3.20} If $V$ is a regular semistable bundle of rank $n$, then the
maximal ideals $\frak m$ of $\Hom(V,V)$ are in  one-to-one correspondence
with nonzero homomorphisms $\rho\: V\to \lambda$ mod scalars, where
$\lambda$ is a torsion free rank one sheaf of degree zero. The
correspondence is as follows: given $\rho$, we set 
$$\frak m= \{\, A\in \Hom (V,V):\rho\circ A = 0\,\},$$ and given a maximal
ideal
$\frak m$, we define $\lambda = V/\frak m \cdot V$ and take $\rho$ to be the
obvious projection.
\endstatement
\proof If $V = \left(\bigoplus _iI_{r_i}(\lambda _i)\right) \oplus I(\Cal
F)$, then
$\Hom (V,V)$ is a direct sum 
$$\left(\bigoplus _i\Hom (I_{r_i}(\lambda _i),I_{r_i}(\lambda _i))\right)
\oplus \Hom (I(\Cal F),I(\Cal F)),$$ and it will clearly suffice to
consider the case where $V$ is either
$I_{r_i}(\lambda _i)$ or $I(\Cal F)$. For simplicity, we assume that
$V=I(\Cal F)$. Thus there is a unique $\rho$ mod scalars by definition. If
we set $\frak m= \{\, A\in \Hom (V,V):\rho\circ A = 0\,\}$, then $\frak m$
is an ideal of $\Hom (V,V)$. In fact, there is an induced homomorphism
$\Hom (V,V) \to
\Hom (\Cal F, \Cal F) =H^0(\scrO_{\tilde E}) = \Cee$, and $\frak m$ is the
kernel of this homomorphism. Thus $\frak m$ is a maximal ideal.

Next we claim that $\frak m$ is the unique maximal ideal in $\Hom(V,V)$. It
suffices to show that $\frak m$ contains every non-invertible element of
$\Hom(V,V)$. If
$A\in \Hom(V,V)$ is not invertible, then $\operatorname{Im}A$ is a proper
torsion free subsheaf of $V$ of rank smaller than $n$ and degree at least
zero. It follows that $\deg \operatorname{Im}A = 0$. But then, by Lemma
1.11,
$\operatorname{Im}A\subseteq \Ker \rho$. It follows that $\rho\circ A =0$,
so that by definition $A\in \frak m$.  Thus $\frak m$ is the unique maximal
ideal of
$\Hom(V,V)$.

Finally we claim that $V/\frak m\cdot V \cong \Cal F$. By definition, the
surjection $\rho\: V \to \Cal F$ factors through the quotient $V/\frak
m\cdot V$, so that $\frak m\cdot V\subseteq \Ker \rho $. Choosing a basis
$A_1,
\cdots, A_{n-1}$ for $\frak m$, we see that $\frak m\cdot V$ is of the form
$A_1(V) +\cdots + A_{n-1}(V)$, and thus it is a subsheaf of $V$ of degree
at least zero. Hence it has degree exactly zero, and thus it is filtered by
subsheaves whose quotients are isomorphic to $\Cal F$. If $\frak m\cdot V$
has rank $r$, it follows by Lemma 1.14 that $\dim \Hom (V, \frak m\cdot V)
\leq r$. But clearly
$\Hom (V,V) = \Hom (V,\frak m\cdot V)\oplus \Cee\Id$, and since $\dim \Hom
(V,V)= n$, we must have
$r=n-1$. Since both $\frak m\cdot V$ and $\Ker \rho$ have degree zero and
rank
$n-1$, and $\frak m\cdot V\subseteq \Ker \rho $,
$\frak m\cdot V = \Ker \rho$. Thus $V/\frak m\cdot V \cong \Cal F$.
\endproof

Next, given the spectral cover $T$, by construction $\bold U(d)$ is a
module over
$\scrO_T= \pi _1{}_*Hom(\bold U(d), \bold U(d))$, and thus $\bold U(d)$
corresponds to a coherent sheaf $\Cal L_d$ over $T\times E$. By
construction,
$(\nu\times
\Id)_*\Cal L_d =\bold U(d)$. In case $E$ is smooth, or more generally in
case $(V,
\frak m)$ is a point of $T$ such that the support of $V$ does not contain
the singular point of
$E$, then it is easy to check directly that $\Cal L_d$ is a line bundle
near $(V,
\frak m)\times E$. 

We can now summarize our description of $T$ as follows:

\theorem{3.21} There  is an isomorphism $(\nu,  r)$ of
$T$ onto the incidence correspondence in $|np_0|\times E$ with the following
property: Let
$\Delta_0$ be the diagonal in $E\times E$, with ideal sheaf $I_{\Delta
_0}$, and let $\Cal P_0$ be the sheaf on
$E\times E$ defined by $Hom(I_{\Delta _0}, \scrO_{E\times E}(-E\times
\{p_0\})$. Then there exists a line bundle $M$ on $T$ such that $\Cal L_d =
( r\times
\Id)^*\Cal P_0\otimes \pi _1^*M$.
\endstatement
\proof We have shown in (0.3) that $\Cal P_0$ is flat over the first factor
in the product $E\times E$ and identifies the first factor with $\bar
J(E)$, the compactified generalized Jacobian of $E$. Let
$T'$ be the incidence correspondence in
$|np_0|\times E$. Note that
$T'$ is irreducible; in fact, projection onto the second factor makes $T'$ a
$\Pee^{n-2}$-bundle over $E$, namely $T'=\Pee\Cal E$ as in Section 2. We
will first find a morphism from $T'$ to $T$ which is a bijection as a
set-valued function. Let $\nu'\: T'\to \Pee^{n-1}$ and $r'\: T' \to E$ be
the projections to the first and second factors. 

By construction, for a point $(x, e)\in T'$, if $V$ is the vector bundle
over $E$ corresponding to $x$ and $\lambda$ is the rank one torsion free
sheaf of degree zero corresponding to $e$, $\dim \Hom(V, \lambda) = 1$.
Thus, by base change, with
$\pi _1\: T'\times E\to T'$ the projection,
$$\pi _1{}_*\left[(\nu'\times \Id)^*\bold  U(d)\spcheck \otimes  (r'\times
\Id)^*\Cal P_0\right] =M$$ is a line bundle on $T'$. After replacing
$(r'\times \Id)^*\Cal P_0$ by
$(r'\times \Id)^*\Cal P_0\otimes M^{-1}=\Cal P'$, we can assume that there
is a surjection 
$$\rho\: (\nu'\times \Id)^*\bold  U(d) \to \Cal P'.$$ On every fiber, the
homomorphism $V \to \lambda$ is preserved up to scalars by every
endomorphism of $V$. Thus, there is an induced homomorphism
$$\pi_1{}_*Hom ((\nu'\times \Id)^*\bold  U(d), (\nu'\times \Id)^*\bold 
U(d)) \to 
\pi_1{}_*Hom (\Cal P', \Cal P').$$ Now by the flatness of $\Cal P_0$, it is
easy to check that $Hom (\Cal P', \Cal P')$ is flat over $T'$ and that the
natural multiplication map $\scrO_{T'} \to
\pi_1{}_*Hom (\Cal P', \Cal P')$ is an isomorphism. By base change, the
first term
$\pi_1{}_*Hom ((\nu'\times \Id)^*\bold  U(d), (\nu'\times \Id)^*\bold 
U(d)) $ in the above homomorphism is the pullback to $T'$ of the sheaf of
algebras $\bold A =
\pi _1{}_*Hom (\bold U(d), \bold U(d))$ over $\Pee^{n-1}$, and hence it is
just the structure sheaf
$\scrO_{T'\times _{\Pee^{n-1}}T}$ of the fiber product $T'\times
_{\Pee^{n-1}}T$. The homomorphism $\scrO_{T'\times _{\Pee^{n-1}}T}\to
\scrO_{T'}$ corresponds to a morphism $T'\to T'\times _{\Pee^{n-1}}T$,
which is a section of the natural projection $T'\times _{\Pee^{n-1}}T \to
T'$. Such a section is the same thing as a morphism $T' \to T$ (covering
the given maps to $\Pee^{n-1}$). On the level of points, this morphism is
as follows: take an element $D$ in
$|np_0|$ and a point $e$ in the support of $D$. Pass to the corresponding
vector bundle $V$ and the morphism $V\to \lambda$, where $\lambda$ is the
rank one degree zero torsion free sheaf corresponding to $e$, and then set
$\frak m$ to be the maximal ideal corresponding to $V\to \lambda$. It is
then clear that $T'\to T$, as constructed above, is a bijection of sets. In
particular, $T$ is irreducible.

Now we want to construct a morphism which is the inverse of the morphism
$T'\to T$. It suffices to find the morphism $r\: T \to E$. Viewing $E$ as
isomorphic to the compactified generalized Jacobian of $E$, we can find
such a morphism once we know that the sheaf $\Cal L_d$ is flat over $E$:

\lemma{3.22} The sheaf $\Cal L_d$ is flat over $T$. If $t\in T$ corresponds
to the  the pair $(V, \frak m)$, and $\lambda$ is the rank one torsion free
sheaf of degree zero given by $V/\frak m\cdot V$, then the restriction of
$\Cal L_d$ to the slice
$\{t\}\times E$ is $\lambda$. Thus $\Cal L_d$ is a flat family of rank one
torsion free sheaves on $T\times E$.
\endstatement
\proof First let us show that, in the above notation, the restriction of
$\Cal L_d$ to the slice $\{t\}\times E$ is $\lambda$. In fact, suppose that
$t$ corresponds to the pair $(V, \frak m)$ and view
$V$ as a rank one module over $\Hom(V,V)$. Then the restriction of $V$ to
$\{t\}\times E$ is given by $V/\frak m \cdot V =\lambda$. Now the Hilbert
polynomial
$P_\lambda (n) = \chi(E;
\lambda\otimes \scrO_E(np_0))$ is independent of the choice of $\lambda$.
As we have proved above, $T$ is irreducible since it is the image of $T'$,
and thus, since it is reduced, it is integral. The proof of Theorem 9.9 on
p\. 261 of
\cite{10} then shows that $\Cal L_d$ is flat over $T$. The last statement
is then clear.
\endproof

By (0.3), as $\Cal L_d$ is flat over
$T$, it defines  a morphism $ r\: T\to E$ (viewing $E$ as $\bar J(E)$).
Thus we also have the product morphism $(\nu,  r)\: T \to
\Pee^{n-1}\times E$, whose image is $T'$. Clearly, on the level of sets, the
morphism $T\to T'$ is the inverse of the morphism $T'\to T$ constructed
above. Since both $T$ and $T'$ are reduced, in fact the two maps are
inverses as morphisms.  By the functorial properties of the compactified
Jacobian (0.3), $\Cal L_d = ( r\times
\Id)^*\Cal P\otimes \pi _1^*M$. This then concludes the proof of (3.21).
\endproof

We have now lined up the spectral covers, and proceed to identify the
bundles
$\bold U(d)$ in terms of $T$. It suffices to identify the bundle $\pi
_1^*M$ in (3.21). We do this first for $d=1$. In order to do so, we first
make the following preliminary remarks. Let
$\scrO_{T\times E}(\Delta)$ denote the rank one torsion free sheaf $(
r\times
\Id)^*Hom(I_{\Delta _0}, \scrO_{E\times E})$. Suppose that  $\Cal L$ is any
flat family of rank one torsion free sheaves on $T\times E$ such that there
exists an injection $\scrO_{T\times E}
\to
\Cal L$, with the cokernel exactly supported along $\Delta$ and with
multiplicity one at a nonempty Zariski open subset of the smooth points. We
claim that in this case $\Cal L = \scrO_{T\times E}(\Delta)$. First, the
universal property of the compactified Jacobian implies that
$\Cal L = (\alpha \times \Id)^*Hom(I_{\Delta _0}, \scrO_{E\times E}) \otimes
\pi_1^*M$ for some morphism $\alpha \: T \to E$ and line bundle
$N$ on $T$. By hypothesis, $\alpha = r$ on a Zariski open dense subset of
$T$, and thus everywhere. Next, since $\scrO_T \to
\pi_1{}_*\scrO_{T\times E}(\Delta)$ is an isomorphism,
$H^0(N) =H^0(\pi_1{}_*\scrO_{T\times E}(\Delta) \otimes N ) = H^0(\Cal L)$,
and every section of $\Cal L$ is given by multiplying the natural section
of 
$\scrO_{T\times E}(\Delta)$ by a section $s$ of $N$. In this case, the
cokernel is supported at $\Delta \cup \pi_1^{-1}(D)$, where $D$ is the
divisor of zeroes of $s$. Thus, if the support of the cokernel is $\Delta$,
then 
$N$ must have a nowhere vanishing section, and so is trivial. We may thus
conclude that
$\Cal L =
\scrO_{T\times E}(\Delta)$.

\theorem{3.23} $\bold U(1) = (\nu\times \Id)_*\scrO_{T\times E}(\Delta -
G)\otimes \pi _1^*\scrO_{\Pee^{n-1}}(1)$.
\endstatement
\proof An equivalent formulation is: 
$$(\nu\times \Id)_*\scrO_{T\times E}(\Delta) \cong \bold U(1) \otimes
\pi_2^*\scrO_E(p_0) \otimes \pi _1^*\scrO_{\Pee^{n-1}}(-1).$$ We will find
a section of the torsion free rank one sheaf $\Cal L$ on $T\times E$
corresponding to $\bold U(1) \otimes
\pi_2^*\scrO_E(p_0) \otimes \pi _1^*\scrO_{\Pee^{n-1}}(-1)$ which vanishes
to order one along $\Delta$. By the remarks before the proof, this will
imply that
$\Cal L = \scrO_{T\times E}(\Delta)$. Now there is an inclusion
$$\scrO_{\Pee^{n-1}\times E} \to \bold U(1) \otimes
\pi_2^*\scrO_E(p_0) \otimes \pi _1^*\scrO_{\Pee^{n-1}}(-1)$$ whose cokernel
is 
$\pi_2^*(\scrO_E(p_0)\otimes W_{n-1})\otimes \pi
_1^*\scrO_{\Pee^{n-1}}(-1)$. Thus $h^0(\Cal L) = h^0(\bold U(1) \otimes
\pi_2^*\scrO_E(p_0) \otimes \pi _1^*\scrO_{\Pee^{n-1}}(-1)) = 1$. To see
where the unique section of $\Cal L$ vanishes, fix a point $x\in
\Pee^{n-1}$ corresponding to a regular semistable $V$, and consider where
the corresponding section of
$V\otimes \scrO_E(p_0)$ vanishes. This section arises from a homomorphism
$\scrO_E(-p_0) \to V$ constructed in (3.4) and (3.5). For example, if
$V=\bigoplus _i\scrO_E(e_i-p_0)$, then, up to the action of $(\Cee^*)^n$,
the map is the direct sum of the natural inclusions $\scrO_E(-p_0)\to
\scrO_E(e_i-p_0)$. At each fiber $\{(V, e_i)\}\times E$ of $T\times E$
lying over $\{V\}$, the section therefore vanishes simply at $((V, e_i),
e_i)$. For a general point $t=(V, \frak m)$ of $T$, the restriction of the
section of
$\Cal L$ to the fiber $\{t\}\times E$ vanishes at the point of $E$ where the
corresponding section of the composite map
$$\scrO_E(-p_0) \to V \to V/\frak m\cdot V \cong \lambda$$ vanishes. By the
construction of (3.4), the composite map $\scrO_E(-p_0) \to V \to
\lambda$  is not identically zero, and hence vanishes exactly at the point
$e$ of
$E$ corresponding to
$\lambda$.   Thus the section of $\Cal L$ vanishes exactly along
$\Delta$, with multiplicity one on a Zariski open and dense subset, proving
(3.23).
\endproof

\corollary{3.24} Suppose that $E$ is smooth. For all $d\in \Zee$ with
$1\leq d\leq n-1$,
$\bold U(d) =  U_{1-d}\otimes
\pi _1^*\scrO_{\Pee^{n-1}}(1)$.
\endstatement
\proof This follows by writing both sides as successive elementary
modifications of
$\bold U(1)$, resp. $U_0\otimes
\pi _1^*\scrO_{\Pee^{n-1}}(1)$, and applying (3.23).
\endproof

\ssection{3.5. The general spectral cover construction.}

For every Weierstrass cubic $E$, we have now constructed a finite cover
$T\to
\Pee^{n-1}$ and a torsion free rank one sheaf $\Cal L_0=\scrO_{T\times
E}(\Delta)\otimes \pi _2^*\scrO_E(-p_0)$. The proof of Theorem 2.4  goes
over word-for-word to show:

\theorem{3.25} Let $E$ be a Weierstrass cubic and let $U'$ be a rank $n$
vector bundle over $|np_0|\times E$ with the following property.  For each
$x\in |np_0|$ the restriction of $U'$ to $\{x\}\times E$ is isomorphic to
the restriction of $U_0$ to $\{x\}\times E$. Then $U'=(\nu \times
\Id)_*\left[\scrO_{T\times E }(\Delta-G)\otimes q _1^*L\right]$ for a
unique line bundle  $L$  on $T$.
\qed
\endstatement

We may  define $U_a$ and, for a smooth point $e\in E$, $U_a[e]$ exactly as
in Definition 2.6, and the proof of Lemma 2.7 shows that $U_a$ is an
elementary modification of $U_{a-1}$, and similarly for $U_a[e]$. Since
$\Pic T\cong r^*\Pic E\oplus \Zee$, every bundle
$U'$ as described in Theorem 3.25 is of the form $U_a[e]\otimes
\scrO_{\Pee^{n-1}}(b)$ for integers $a,b$ and a smooth point $e\in E$.

\remark{Question} In case $E$ is singular, $T$ is singular as well. Is
there an analogue of twisting by Weil divisors on $T$ which are not
Cartier, which produces bundles which are not regular, or perhaps sheaves
which are not locally free, over points of $\Pee^{n-1}$ corresponding to
the singular points of $T$? See Section 6 for a related construction in the
smooth case.
\endremark
\medskip

The following is proved as in the proof of Theorem 2.8.

\theorem{3.26} Let $E$ be a Weierstrass cubic and let $S$ be a scheme or
analytic space. Suppose that
${\Cal U}\to S\times E$ is a rank $n$ holomorphic vector bundle whose 
restriction to each slice $\{s\}\times E$ is a regular semistable bundle
with trivial determinant.  Let $\Phi\colon S\to |np_0|$ be the morphism
constructed in Theorem \rom{1.5}. Let $\nu_S\colon\tilde S\to S$ be the
pullback via
$\Phi$ of the spectral covering $T\to |np_0|$:
$$\tilde S=S\times_{|np_0|}T,$$ and let $\tilde\Phi\colon \tilde S\to T$ be
the map covering
$\Phi$. Let $q_1\colon \tilde S\times E\to \tilde S$ be the projection onto
the first factor. 
 Then there is a line bundle ${\Cal M}\to \tilde S$ and an isomorphism of
${\Cal U}$ with 
$$(\nu_S\times\operatorname{Id})_*
\left((\tilde\Phi\times \Id)^*(\scrO_{T\times E}(\Delta-G))\otimes
q_1^*{\Cal M}\right).\qquad\qed$$ 
\endstatement

\ssection{3.6. Chern classes.}

\theorem{3.27} For all $d$ with $1\leq d \leq n-1$, the total Chern class
and the Chern character of $\bold U(d)$ are given by\rom:
$$\align  c(\bold U(d)) &= (1+h+\pi_2^*[p_0])(1+h)^{d-1};\\
\ch \bold U(d) &= (d -\pi _2^*[p_0])e^h + (n-d) + [p_0].
\endalign$$ Thus $c(U_0) = (1-h+ \pi _2^*[p_0]\cdot h)(1-h)^{n-2}$ and $\ch
U_0 = ne^{-h} + (1- \pi _2^*[p_0])(1- e^{-h})$.
\qed
\endstatement
\proof In $K$-theory, $W_d$ is the same as the direct sum of $d-1$ trivial
bundles  and the line bundle $\scrO_E(p_0)$. Thus 
$$\align c(\bold U(d)) &= (1-\pi_2^*[p_0] + h)(1+h)^{d-1}(1+\pi_2^*[p_0])\\
&=(1+h+\pi_2^*[p_0])(1+h)^{d-1},
\endalign$$ and
$$\align
\ch \bold U(d) &= \pi _2^*\ch (W_d\spcheck)\cdot \pi _1^*\ch (\scrO_{\Pee
^{n-1}}(1)) + \pi _2^*\ch W_{n-d}\\ & = (d-1+ \pi _2^*e^{-[p_0]})e^h +
(n-d-1 + \pi _2^*e^{[p_0]}) \\ &=(d -\pi _2^*[p_0])e^h + (n-d) + [p_0],
\endalign$$ since $e^{-[p_0]}= 1-[p_0]$ and similarly for $e^{[p_0]}$.  The
formulas for $c(U_0)$ and 
$\ch U_0$ then follow from (3.23).
\endproof

\section{4. A relative moduli space for elliptic fibrations.}

Our goal in this section is to do the constructions of the last three
sections in the relative setting of a  family $\pi\: Z \to B$ of elliptic
curves (possibly with singular fibers), in order to produce families of 
bundles whose restriction to every fiber of $\pi$ is regular semistable and
with trivial determinant. First we identify the relative coarse moduli
space as a projective bundle over
$Z$. The extension picture generalizes in a straightforward way to give
$n-1$ ``universal" bundles $\bold U(d)$, $1\leq d\leq n-1$, and they are
related via elementary modifications. Using these bundles, we can
generalize the spectral covers picture as well. Finally, we compute the
Chern classes of the universal bundles.

Let $\pi \: Z \to B$ be an elliptic fibration with a section $\sigma$.
Following the notational conventions of the introduction, we shall always
let
$L^{-1} = R^1\pi _*\scrO_Z$, which we can also identify with the normal
bundle
$\scrO_Z(\sigma)|\sigma$.

\ssection{4.1. A relative coarse moduli space.}

Our first task is to find a relative version of $|np_0|$ for a single
elliptic curve. The relative version of the vector space $H^0(E;
\scrO_E(np_0))$ is just the rank $n$ vector bundle $\pi
_*\scrO_Z(n\sigma)=\Cal V_n$, and the relative moduli space will then be
the associated projective bundle.  From the exact sequence
$$0 \to \scrO_Z((n-1)\sigma) \to \scrO_Z(n\sigma) \to \pi ^*L^{-n}|\sigma
\to 0,$$ we obtain for $n\geq 2$ an exact sequence
$$0 \to \Cal V_{n-1} \to \Cal V_n \to L^{-n} \to 0.$$ (For $n=1$ the
corresponding sequence identifies $\pi_*\scrO_Z(\sigma)$ with
$\scrO_B$ and shows that there is an isomorphism $\scrO_Z(\sigma)|\sigma\to
R^1\pi _*\scrO_Z=L^{-1}$.) Thus $\Cal V_n$ is naturally filtered by
subbundles such that the successive quotients are decreasing powers of $L$.
The following well-known lemma shows that this filtration is split:

\lemma{4.1}
$\pi _*\scrO_Z(\sigma) = \scrO_B$, and, for
$n\geq 2$,
$$\pi _*\scrO_Z(n\sigma) = \Cal V_n = \scrO_B\oplus L^{-2} \oplus
\cdots \oplus L^{-n}.$$
\endstatement

\proof Since $h^0(E;\scrO_E(np_0))=n$ for all the fibers $f$ of $\pi$, it
follows from base change  that $\pi_*\scrO_Z(n\sigma)$ is a vector bundle of
rank $n$.  Furthermore, the local sections of this bundle over an open
subset $U\subset B$ are simply the meromorphic functions on
$\pi^{-1}(U)$ which have poles of order at most $n$ along $\sigma\cap U$.
For $U$ sufficiently small, there are functions $X$ with a pole of order $2$
along $\sigma$ and $Y$ with a pole of order $3$. Moreover, if we require
that
$X$ and $Y$ satisfy a Weierstrass equation, then $X$ and
$Y$ are unique up to  nowhere vanishing functions in $U$ and transform as
sections of
$L^{-2}, L^{-3}$ respectively. We can also use the defining equation of
$Z$ to write $Y^2$ as a cubic polynomial in $X$. Now every section of
$\pi_*\scrO_Z(n\sigma)$ can be written uniquely as
$$(\alpha_0+\alpha_1X+\cdots+\alpha_kX^k)+
Y(\beta_0+\beta_1X+\cdots+\beta_\ell X^\ell)$$  where the $\alpha_i$ are
holomorphic sections of $L^{-2i}$ and the
$\beta_j$ are holomorphic sections of $L^{-2j-3}$ and   $2k\le n$ and
$2\ell+3\le n$. The $\alpha_i,\beta_j$ determine the isomorphism claimed in
the statement of the lemma.
\endproof

Notice that the inclusion  $\Cal V_{n-1}\subset \Cal V_n$ corresponds to the
natural inclusion
$$ \scrO_B\oplus L^{-2} \oplus \cdots \oplus L^{-(n-1)}\subset
\scrO_B\oplus L^{-2} \oplus \cdots \oplus L^{-n}.$$ In particular, the
distinguished points $\bold o_E=np_0\in |np_0|$ corresponding to the
bundles with all Jordan-H\"older quotients trivial fit together to  make a
section $\bold o_Z$ of $\Pee\Cal V_n$.  This section is  the
projectivization $\Pee \scrO_B$ of the first factor $\scrO_B$ in the above
decomposition.

We call the above splitting the {\sl $X$-$Y$ splitting} of
$\pi_*\scrO_Z(n\sigma)$.  While this decomposition of
$\pi_*\scrO_Z(n\sigma)$ is natural it is  not the only possible
decomposition, even having the property described in the previous
paragraph.  For example, another splitting was suggested to us by P.
Deligne. There is a global  holomorphic differential $\omega$ on $E$ which
is given on a Zariski open subset of $E$ by $dX/Y$. There is a local
complex coordinate $\zeta$ for $E$ centered at $p_0$ with the property that
on the open set on which this local coordinate is defined we have
$\omega=d\zeta$. Of course, there is a homomorphism $\Cee\to E$ which pulls
$\zeta$ back to the usual coordinate on $\Cee$. Every meromorphic function
on $E$ with a pole of order at most $n$ at $p_0$ can be expanded as a
Laurent series in $\zeta$:
$$f=\sum_{i=-n}^\infty b_i\zeta^i.$$ The coefficient $b_i$ in this
expansion is a section of $L^i$. We can then use the coefficients
$b_{-n},\ldots,b_{-2},b_0$ to define a splitting of $\pi_*\scrO_Z(n\sigma)$.
(If $f$ is  a meromorphic function on $E$ whose only pole is at $p_0$ then
$b_{-1}$ is determined by the
$b_{-i}$ for $-n\le -i\le -2$.) This splitting is different from the $X$-$Y$
splitting, but both splittings induce the same filtration on
$\pi_*\scrO_Z(n\sigma)$.

In Theorem 1.2 we showed how a semistable bundle of rank $n$ and trivial
determinant on a smooth elliptic curve $E$ determines a point of the linear
series $\scrO_E(np_0)$. This works well for bundles over families of
elliptic curves.

\lemma{4.2} Let  $p\:\Pee \Cal V_n  \to B$ be the projection. Thus, the
fiber of $p$ over $b\in B$ is the complete linear system $|np_0|$, where
$E_b=\pi ^{-1}(b)$ and $p_0$ is the smooth point $\sigma \cap E_b$. If
$V\to Z$ is a rank $n$ vector bundle whose restriction to each fiber of
$Z\to B$ is a semistable bundle with trivial determinant, then $V$
determines a section 
$$A(V)\colon B\to \Pee \Cal V_n ,$$ with the property that, for each $b\in
B$,
$$A(V)(b)=\zeta(V| E_b ).$$
\endstatement

\proof Arguing as in (1.6),  there is an induced morphism
$$\Psi\: \pi^*\pi_*(V \otimes \scrO_Z(\sigma)) \to V \otimes
\scrO_Z(\sigma).$$ The determinant of this morphism is a section of $\pi^*
M \otimes
\scrO_Z(n\sigma)$, for some line bundle $M$ on $B$, and it gives a
well-defined section $A(V)$ of $\Pee \Cal V_n $ over $B$.
\endproof

We note that the proof of (4.2) does not require that $B$ be smooth, or even
reduced.

There is also an analogue for families of elliptic curves of Theorem 1.5.

\lemma{4.3} Let $S$ be a scheme or analytic space over $B$ and let
$\Cal V$ be a rank $n$ vector bundle over $S\times _BZ$, such that the
restriction of $\Cal V$ to every fiber $p_1^{-1}(s) \cong \pi ^{-1}(b)$ is
semistable with trivial determinant, where $p_1, p_2$ are the projections
of $S\times _BZ$ to the first and second factors and $b$ is the point of
$B$ lying under $s$. Then there is an induced morphism $\Phi \: S\to \Pee
\Cal V_n $ of spaces over $B$, which agrees over each
$b\in B$ with the morphism defined in \rom{(1.5)}.
\endstatement

\proof Let $\hat Z=S\times _BZ$, with $\hat\pi \: \hat Z \to S$ the first
projection.  Then $\hat Z$ is an elliptic scheme over
$S$ which maps naturally to $Z$ covering the map of $S\to B$. Let
$\hat\sigma$ be the  induced section.  Set $\widehat {\Cal V}_n = \hat\pi
_*\scrO_{\hat Z}(n\hat
\sigma)$. Clearly $\Pee \widehat {\Cal V}_n$ is identified with the
pullback of
$\Pee \Cal V_n $. Now  apply the above result to this elliptic scheme to
produce a section $S\to \Pee \widehat {\Cal V}_n$ which when composed with
the natural map $\Pee \widehat {\Cal V}_n \to \Pee\Cal V_n$ is the morphism
$\Phi$ of the proposition.
\endproof

\ssection{4.2. Construction of bundles via extensions.}

Our goal for the remainder of this section is to construct various
``universal" bundles over $\Pee \Cal V_n \times _BZ$. The first and easiest
construction of the universal moduli space is via the extension approach,
generalizing what we did in Section 3 for a single elliptic curve.

In order to make the extension construction in families, we first need to
extend the basic bundle $W_k$ over $E$ to bundles over the elliptic scheme
$Z$. 

\proposition{4.4} There is a vector bundle $\Cal W_d$ on $Z$ such that
$\Cal W_d$ is filtered, with successive quotients $\pi ^*L^{d-1}, \pi
^*L^{d-2}, \dots, \scrO_Z(\sigma)$, and such that on every fiber $\Cal W_d$
restricts to $W_d$. Moreover, $\Cal W_d$ is uniquely specified by the above
properties. In fact, if $\Cal W$ is a vector bundle on $Z$ such that $\Cal
W$  restricts to $W_d$ on every fiber, then there exists a line bundle $M$
on $B$ such that $\Cal W = \Cal W_d\otimes \pi ^*M$. Finally, 
$R^0\pi _*\Cal W_d = L^{(d-1)}$ and
$R^1\pi _*(\Cal W_d\spcheck) = L^{-d}$.
\endstatement
\proof In case $d=1$,  take $\Cal W_1= \scrO_Z(\sigma)$. Now suppose
inductively that $\Cal W_{d-1}$ has been defined, and that
$R^1\pi _*(\Cal W_{d-1}\spcheck) = L^{-(d-1)}$. We seek an extension of
$\Cal W_{d-1}$ by a line bundle trivial on ever fiber of
$\pi$, and thus of the form $\pi^*M$ for some line bundle $M$ on $B$, and
such that
$H^0(R^1\pi _*(\Cal W_{d-1}\spcheck\otimes \pi^*M))$ has an everywhere
generating section. Now
$$R^1\pi _*(\Cal W_{d-1}\spcheck\otimes \pi^*M) = R^1\pi _*(\Cal
W_{d-1}\spcheck)
\otimes M = L^{-(d-1)}\otimes M.$$ Thus we must have
$M=L^{d-1}$. With this choice, noting that $R^0\pi _*(\Cal 
W_{d-1}\spcheck\otimes \pi ^*L^{d-1}) = 0$ since $W_d\spcheck$ has no
sections, the Leray spectral sequence gives an isomorphism
$$H^1(\Cal W_{d-1}\spcheck\otimes \pi ^*L^{d-1}) \cong H^0(R^1\pi _*(\Cal
W_{d-1}\spcheck\otimes \pi ^*L^{d-1}))=H^0(\scrO_B)$$ and thus a global
extension of $\Cal W_{d-1}$ by $\pi ^*L^{d-1}$ restricting to
$W_d$ on every fiber. Since the unique section of $W_d$ is given by the
inclusion of the canonical subbundle $\scrO_f \to W_d$, we must have
$R^0\pi _*\Cal W_d = L^{(d-1)}$, and a similar argument (or relative
duality) evaluates
$R^1\pi _*(\Cal W_d\spcheck)$.

Finally suppose that $\Cal W$ is another bundle on $Z$  restricting to
$W_d$ on every fiber. Then since $W_d$ is simple, $\pi _*\Hom (\Cal W_d,
\Cal W)$ is a line bundle $M$ on $B$, and thus  $\pi _*\Hom (\Cal
W_d\otimes \pi ^*M, \Cal W)\cong \scrO_B$. The element $1\in H^0(\scrO_B)$
then defines an isomorphism from $\Cal W_d\otimes \pi ^*M$ to $\Cal W$.
\endproof

Note that the formation of $\Cal W_d$ is compatible with base change, in the
following sense. Given a morphism $g\: B'\to B$, let $Z'=Z\times _BB'$,
with $f\: Z'\to Z$ the induced morphism, and let  $\sigma '$ be the induced
section of
$\pi'\: Z'\to B'$. Then the bundle $\Cal W_d'$ constructed for $\pi'\:Z'\to
B'$ and the section $\sigma'$ is $f^*\Cal W_d$.

Next we construct a universal bundle via extensions. First we identify the
relevant bundles to use as the parameter space of the extension:

\lemma{4.5} For $1\leq d \leq n-1$, the sheaves $R^1\pi _*(\Cal
W_{n-d}\spcheck
\otimes \Cal W_d\spcheck)=\Cal V_{n,d}$ are locally free of rank $n$ over
$B$, and are all canonically identified.
\endstatement
\proof The local freeness and the rank statement follow from  Claim 3.3 and 
base change. The canonical identifications follow from Corollary 3.10.
\endproof

Let
$\Cal V_{n,d} = R^1\pi _*(\Cal W_{n-d}\spcheck \otimes \Cal W_d\spcheck)$ as
above, and let 
${\Cal P}_{n-1,d}$ be the associated projective space bundle $\Pee(\Cal
V_{n,d})
\to B$. By the general properties of extensions, there is a universal
extension over ${\Cal P}_{n-1,d}\times_BZ$ of the form
$$0\to \pi_2^*{\Cal W}_d\spcheck\otimes \pi_2^*\scrO_{{\Cal
P}_{n-1,d}}(1)\to {\bold U}(d)\to \pi_2^*{\Cal W}_{n-d}\to 0.$$ Applying
Lemma 4.3 to these bundles produces bundle maps over $B$
$$\Phi_d\colon {\Cal P}_{n-1,d}\to \Pee \Cal V_n .$$

The projective space bundles ${\Cal P}_{n-1,d}$ over $B$ are all
canonically isomorphic. Under these isomorphisms, the universal bundles
${\bold U}(d)$ are all distinct. Nevertheless, the result in Proposition
3.10 shows that there is an isomorphism $I$ which identifies 
$R^1\pi_*(\scrO_Z(-\sigma)\otimes {\Cal W}\spcheck_{n-1})$ with 
$$R^0\pi_*(\det(\scrO_Z(\sigma)\otimes {\Cal W}_{n-1}))\otimes
\det(R^0\pi_*(\scrO_Z(\sigma)\otimes {\Cal W}_{n-1}))^{-1}\otimes
R^1\pi_*\scrO_Z.$$

Let us identify the various factors on the right-hand-side of this
expression.   First of all, it is straightforward given the inductive
definition of the
${\Cal W}_{n-1}$ to show that. 
$$\det(\scrO_Z(\sigma)\otimes {\Cal W}_{n-1})\cong
\scrO_Z(n\sigma)\otimes \pi^*L^{(n-1)(n-2)/2}.$$ It follows that
$R^0\pi_*(\det(\scrO_Z(\sigma)\otimes {\Cal W}_{n-1}))\cong
R^0\pi_*\scrO_Z(n\sigma)\otimes L^{(n-1)(n-2)/2}$. Next, we have exact
sequences
$$0\to R^0\pi_*(\scrO_Z(\sigma)\otimes \pi^*L^{n-2})\to
R^0\pi_*(\scrO_Z(\sigma)\otimes {\Cal W}_{n-1})\to
R^0\pi_*(\scrO_Z(\sigma)\otimes {\Cal W}_{n-2})\to 0.$$ Since by
Proposition 4.4 we have $R^0\pi_*(\pi^*L^a\otimes
\scrO_Z(\sigma))\cong  L^a$, and since 
$$R^0\pi_*(\scrO_Z(2\sigma))\cong L^{-2}\oplus \scrO_B,$$ an easy inductive
argument shows that 
$$\det(R^0\pi_*(\scrO_Z(\sigma)\otimes {\Cal W}_{n-1}))\cong
L^{((n-1)(n-2)/2)-2}.$$ Lastly, 
$R^1\pi_*\scrO_Z\cong L^{-1}$.

Putting all this together, we get:

\theorem{4.6} There is  an isomorphism of vector bundles over $B$
$$I\colon R^1\pi_*(\scrO_Z(-\sigma)\otimes {\Cal W}_{n-1})\cong
R^0\pi_*\scrO_Z(n\sigma)\otimes L,$$ which fiber by fiber agrees with the
map $I$ of Proposition \rom{3.13}. In other words, 
$$\Cal V_{n,1} \cong \Cal V_n\otimes L.$$ Furthermore, the map induced by
projectivizing $I$ agrees with the map
$\Phi_1$ produced by applying  Lemma \rom{4.3}  to the family $\bold U(1)$
over ${\Cal P}_{n-1,1}\times_BZ$. Let $\Phi_d$ be the map ${\Cal
P}_{n-1,d}\to \Pee \Cal V_n $ obtained by applying Lemma \rom{4.3} to the
family $\bold U(d)$. Then, the maps $\Phi_d$ for $1\le d<n$ are compatible
with the identifications coming from Corollary \rom{3.10}, and hence each
of these maps is an isomorphism of projective bundles over $B$. 
\endstatement

Note that, while the $\Pee^{n-1}$-bundles ${\Cal P}_{n-1,d}$ and $\Pee \Cal
V_n $ are isomorphic, the tautological bundles $\scrO_{{\Cal
P}_{n-1,d}}(1)$ and
$\scrO_{\Pee
\Cal V_n}(1)$ differ by a twist by $p^*L$. We shall use $\Cal P_{n-1}$ to
denote the bundle $\Pee \Cal V_n $ together with its tautological line
bundle. If $\zeta =c_1(\scrO_{\Pee \Cal V_n}(1))$ and $\zeta' =
c_1(\scrO_{{\Cal P}_{n-1,d}}(1))$, then $\zeta = \zeta ' + L$.

\corollary{4.7} Via the isomorphism of \rom{(4.6)} and \rom{(4.6)}, 
$$\Cal V_{n,d} = R^1\pi _*(\Cal W_{n-d}\spcheck
\otimes
\Cal W_d\spcheck) \cong L\oplus L^{-1} \oplus \cdots \oplus L^{-(n-1)}.$$
This splitting is compatible with the inclusion of $\Cal V_{n-1, d}$ in
$\Cal V_{n,d}$ as well as that of $\Cal V_{n-1, d-1}$ in $\Cal V_{n,d}$.
\qed
\endstatement

\corollary{4.8} Under the isomorphism $\pi |\sigma\: \sigma \cong B$, there
is a natural splitting
$$\Cal W_n|\sigma \cong L^{n-1} \oplus L^{n-2}\oplus \cdots \oplus L \oplus
L^{-1}.$$ In fact, the extension
$$0 \to \pi ^*L^{n-1} \to \Cal W_n \to \Cal W_{n-1} \to 0$$ restricts to
the split extension over $\sigma$.
\endstatement
\proof Let us first show that the restriction of $\Cal W_n\spcheck$ to
$\sigma$ is split. Begin  with the exact sequence
$$0 \to \scrO_Z(-\sigma) \otimes \Cal W_n\spcheck \to \Cal W_n\spcheck \to
\Cal W_n\spcheck|\sigma \to 0,$$ and apply $R^i\pi _*$. We get an exact
sequence
$$\CD 0 @>>> \pi _*(\Cal W_n\spcheck|\sigma) @>>> R^1\pi
_*(\scrO_Z(-\sigma) \otimes \Cal W_n\spcheck) @>>>  R^1\pi _*\Cal
W_n\spcheck @>>> 0\\ @.  @|   @|   @| @.\\ 0 @>>> \pi _*\Cal
W_n\spcheck|\sigma @>>> L\oplus L^{-1}\oplus \cdots \oplus L^{-n} @>>>
L^{-n} @>>> 0.
\endCD$$ Tracing through the identifications shows that the map  $R^1\pi
_*(\scrO_Z(-\sigma)
\otimes \Cal W_n\spcheck) \to R^1\pi _*\Cal W_n\spcheck$ is the same as the
map 
$$R^1\pi _*(\scrO_Z(-\sigma) \otimes \Cal W_n\spcheck) \to R^1\pi
_*\scrO_Z(-\sigma) \otimes L^{-n+1}= L^{-1}\otimes L^{-n+1} = L^{-n}$$
coming from the long exact sequence associated to
$$0 \to \scrO_Z(-\sigma) \otimes \Cal W_{n-1}\spcheck \to \scrO_Z(-\sigma)
\otimes
\Cal W_n\spcheck \to \scrO_Z(-\sigma)\otimes \pi^*L^{-n+1} \to 0.$$ This
identifies the map $L\oplus L^{-1}\oplus \cdots \oplus L^{-n} \to L^{-n}$
with projection onto the last factor. Hence $\Cal W_n\spcheck|\sigma$ is
identified with $L\oplus L^{-1}\oplus \cdots \oplus L^{-n+1}$. Dualizing
gives the splitting of $\Cal W_n|\sigma$. The splitting of the extension 
$$0 \to \pi ^*L^{n-1} \to \Cal W_n \to \Cal W_{n-1} \to 0$$ is similar.
\endproof

Now let us relate the bundles $\bold U(d)$ via elementary modifications. 

\proposition{4.9} Let $\Cal H$ be the smooth divisor which is the image of
$\Pee
\Cal V_{n-1} =\Cal P_{n-2}$ in $\Cal P_{n-1}$ under the natural inclusion
$\pi _*\scrO_Z((n-1)\sigma) \subset\pi _*\scrO_Z(n\sigma)$, and let $i\:
\Cal H \to
\Cal P_{n-1}$ be the inclusion. Then there is an exact sequence
$$0 \to \bold U(d) \to \bold U(d+1) \to (i\times \Id)_*\scrO_{\Cal H\times
_BZ}(1)\otimes
\pi^*L^{-d}
\to 0,$$ where $\scrO_{\Cal H\times _BZ}(1)$ denotes the restriction of
$\scrO_{{\Cal P}_{n-1,d}}(1)= \scrO_{\Pee\Cal V_{n,d}}(1)$ to $\Cal H\times
_BZ$. Thus $\bold U(d)$ is an elementary modification of $\bold U(d+1)$,
and it is the only possible such modification along $\Cal H\times_BZ$.
\endstatement
\proof The construction of the proof of Theorem 3.12 gives an inclusion
$\bold U(d) \to \bold U(d+1)$ whose cokernel is the direct image of a line
bundle supported along $\Cal H\times _BZ$. As in the first paragraph of the
proof of (3.12), this line bundle is the inverse of $\pi _1{}_*Hom(\bold
U(d+1)|\Cal H\times _BZ, \scrO_{\Cal H\times _BZ})$ (where for the rest of
the proof we let
$\pi_1$ be the first projection $\Cal H\times _BZ \to \Cal H$). From the
defining exact sequence for $\bold U(d+1)$,
$$\gather
\pi _1{}_*Hom(\bold U(d+1)|\Cal H\times _BZ, \scrO_{\Cal H\times _BZ})\cong 
\pi _1{}_*Hom(\pi _2^*\Cal W_{d+1}\spcheck \otimes \scrO_{\Cal H\times
_BZ}(1), 
\scrO_{\Cal H\times _BZ})\\ = \pi _1{}_*\pi _2^*\Cal W_{d+1}\otimes
\scrO_{\Cal H\times _BZ}(-1).
\endgather$$ Here by base change $\pi _1{}_*\pi _2^*\Cal W_{d+1}$ is a line
bundle on $\Cal H$ whose restriction to every fiber is the nonzero section
of $W_{d+1}$ on that fiber.  Now
$\Cal W_{d+1}$ is filtered by subbundles with successive quotients
$\pi^*L^d,
\pi ^*L^{d-1}, \dots, \scrO_Z(\sigma)$, and the inclusion of $\pi ^*L^d$ in
$\Cal W_{d+1}$ defines a map $L^d \to \pi _1{}_*\pi _2^*\Cal W_{d+1}$ which
restricts to the nonzero section on every fiber. Thus $\pi _1{}_*\pi
_2^*\Cal W_{d+1} \cong L^d$. Hence 
$$\pi _1{}_*Hom(\bold U(d+1)|\Cal H\times _BZ, \scrO_{\Cal H\times _BZ})
\cong L^d\otimes \scrO_{\Cal H\times _BZ}(-1),$$ and thus the cokernel of
the map $\bold U(d)\to \bold U(d+1)$ is as claimed. The uniqueness is clear.
\endproof

This completes the construction of ``universal" bundles over
$\Pee\Cal V_n\times_BZ$. However, we have constructed only $n-1$ bundles
${\bold U}(d)$ for $1\le d<n$. 

Note that the formation of the  universal bundles $\bold U(d)$ over $\Cal
P_{n-1,d}\times _BZ$ is also compatible with base change $B'\to B$ in the
obvious sense.

\ssection{4.3. The spectral cover construction.}

Now we turn to the generalization of the spectral covering construction.
First let us define the analogues of $E^{n-1}, T$, and $\nu$. By repeating
the construction on each smooth fiber, we could take the
$(n-1)$-fold fiber product
$Z\times _BZ\times _B \cdots \times _BZ$ and its quotient under
$\frak S_n$ and
$\frak S_{n-1}$. However this construction runs into trouble at the
singular fibers, reflecting  the difference between the $n$-fold symmetric
product of $E$ and the linear system $|np_0|$ for a singular fiber. 
Instead, we construct the spectral cover in families as follows: Let $\Cal
E$ be defined by the exact sequence of vector bundles over $Z$
$$0 \to \Cal E \to \pi ^*\pi _*\scrO_Z(n\sigma) \to \scrO_Z(n\sigma)
\to 0,$$ where the last map is the natural evaluation map and is
surjective. We set
$\Cal T = \Pee \Cal E$, with $r\: \Cal T \to Z$ the projection. By
construction
$\Cal T$ is a $\Pee^{n-2}$-bundle over $Z$. There is an inclusion
$$\Cal T \to \Pee (\pi ^*\pi _*\scrO_Z(n\sigma)) = \Pee \Cal V_n \times
_BZ,$$ and  we let $\nu$ be the composition of this morphism with the
projection $q_1\:
\Pee \Cal V_n \times _BZ\to \Pee \Cal V_n $. It is easy to see that $\nu\:
{\Cal T}\to \Pee \Cal V_n$ is an
$n$-sheeted covering, which restricts to the spectral cover described in
Section 2 on each smooth fiber of $Z\to B$.  By analogy with the case of a
single elliptic curve, we would like to consider the sheaf
$${\Cal U}_0=(\nu\times 
\Id)_*\scrO_{\Cal T\times _BZ}(\Delta - \Cal G),$$ where $\Delta = (r\times
\Id)^*(\Delta_0)$, for $\Delta_0$  the diagonal in
$Z\times _BZ$, and $\Cal G=(r\times \Id)^*p_2^*\sigma$ for $p_1, p_2$ the
projections of $Z\times _BZ$ to the first and second factors. Here we can
define
$\scrO_{Z\times _BZ}(\Delta_0)$ to be the dual of the ideal sheaf of
$\Delta_0$ in $Z\times _BZ$. It is an invertible sheaf away from the
singularities of
$Z\times _BZ$. The proof of (0.4) shows that $\scrO_{Z\times
_BZ}(\Delta_0)$ is flat over both factors of $Z\times_BZ$, and identifies
the first factor, say, with the relative compactified generalized Jacobian.
As we shall see,
$\Cal U_0$ is indeed a vector bundle of rank
$n$, and that its restriction over each smooth fiber is the bundle $U_0$
described in Proposition 2.9. In particular, ${\Cal U}_0$ is a family of
regular, semistable bundles with trivial determinant over the  family $Z\to
B$ of elliptic curves.

Although we shall not need this in what follows, for concreteness sake let
us describe the singularities of
$Z\times_BZ$ and
$\Cal T$ explicitly the case where the divisors associated to
$G_2$ and $G_3$ are smooth and meet transversally as in the introduction.
In this case $Z$ has the local equation
$y^2 = x^3 + sx+t$.  The morphism to
$B$ is given locally by $(s,t)$, where $t= y^2-x^3-sx$. Using  $x,y,s$ as
part of a set of local coordinates for
$Z$, the fiber product has local coordinates
$x,y,s,x',y', \dots$ and a local equation
$$y^2-x^3-sx = (y')^2-(x')^3-s(x').$$ Rewrite this equation as
$$y^2- (y')^2 = (x-x')(s+ x^2 + xx' + (x')^2) = h_1h_2, $$ say, where a
local calculation shows that $h_1(s,x,x')$ and $h_2(s,x,x')$ define two
smooth hypersurfaces meeting transversally along $\Gamma \times _B\Gamma$.
It follows that the total singularity in $Z\times _BZ$ is a locally trivial
fibration of ordinary threefold double points, and $\Delta$ is a smooth
divisor which fails to be Cartier at the singularities.

We return now to the case of a general $Z$. Fix a $d$ with $1\leq d\leq
n-1$, and consider the sheaf of algebras $\pi _1{}_*Hom (\bold U(d), \bold
U(d)) =\bold A$ over $\Cal P_{n-1}$. Arguing as in Lemma 3.19, the space
$\bold {Spec}\, \bold A$ is reduced and there is a finite flat morphism
$\nu\: \bold {Spec}\, \bold A \to
\Cal P_{n-1}$ restricts over each fiber to give $\nu\: T\to \Pee^{n-1}$.
Moreover,
$\bold U(d) = (\nu\times \Id)_*\Cal L_d$ for some sheaf $\Cal L_d$ on $\bold
{Spec}\,
\bold A\times _BZ$. The method of proof of Theorem 3.21 then shows:

\theorem{4.10} There is an isomorphism from $\bold {Spec}\, \bold A$ to
$\Cal T$. Under this isomorphism, there is a line bundle $\Cal M$ over
$\Cal T$ such that
$\bold U(d) =(\nu\times \Id)_*\scrO_{\Cal T\times _BZ}(\Delta -\Cal
G)\otimes \pi _1^*\Cal M$.
\qed
\endstatement

Since $\bold U(1)\otimes \scrO_{\Pee\Cal V_{n,d}}(-1) \otimes
\pi_2^*\scrO_Z(\sigma)$ has a section vanishing exactly along $\Delta$, the
proof of 3.23 identifies this line bundle in case $d=1$: 

\theorem{4.11} In the above notation,
$$\align
\bold U(1) &\cong (\nu\times \Id)_*\scrO_{\Cal T\times _BZ}(\Delta -\Cal G)
\otimes \scrO_{\Pee\Cal V_{n,d}}(1)\\ &=(\nu\times
\Id)_*\scrO_{\Cal  T\times _BZ}(\Delta -\Cal G) \otimes \scrO_{\Pee\Cal
V_n}(1)\otimes L^{-1}.\qed
\endalign$$
\endstatement

For every $a\in \Zee$, we can then define $\Cal U_a = (\nu\times
\Id)_*\scrO_{\Cal T\times _BZ}(\Delta -\Cal G -a(r^*\sigma\times _BZ))$. It
follows that $\Cal U_a$ is a vector bundle for every $a\in \Zee$.

\theorem{4.12} With $\Cal U_a$ defined as above, there is an exact sequence
$$0 \to \Cal U_a \to \Cal U_{a-1} \to (i\times \Id)_*\scrO_{\Cal H\times
_BZ}\otimes
\pi ^*L^{a-1} \to 0,$$ which realizes $\Cal U_a$ as an elementary
modification of $\Cal U_{a-1}$. Thus, for $1\leq d\leq n-1$,
$$\bold U(d) \cong\Cal U_{1-d}\otimes \pi_1^*\scrO_{\Pee\Cal V_n}(1)\otimes
L^{-1}.$$
\endstatement
\proof From the definition of $\Cal U_a$, there is an exact sequence
$$\align 0 &\to (\nu\times \Id)_*\scrO_{\Cal T\times _BZ}(\Delta -\Cal G
-a(r^*\sigma\times _BZ)) \\ &\to (\nu\times
\Id)_*\scrO_{\Cal T\times _BZ}(\Delta -\Cal G -(a-1)(r^*\sigma\times _BZ))
\\ &\to (\nu\times
\Id)_*\scrO_{r^*\sigma\times _BZ}(\Delta -\Cal G -(a-1)(r^*\sigma\times
_BZ))\to 0.
\endalign$$ The divisor $r^*\sigma$ is a $\Pee^{n-2}$-bundle over $B$ which
intersects each fiber of $\Cal T\to B$ in the $\Pee^{n-2}$ fiber
$r^{-1}(p_0)$. This fiber is mapped linearly via $\nu$ to the hyperplane
$H_{p_0}$ in $|np_0|$. Thus,
$\nu_*r^*\sigma = \Cal H$. Now both $\Delta$ and $\Cal G$ have the same
restriction to $r^*\sigma\times _BZ$, namely $r^*\sigma\times _B\sigma$.
Also, 
$\scrO_{r^*\sigma\times _BZ}(r^*\sigma\times _BZ)$ is the pullback of the
line bundle $\scrO_Z(\sigma)|\sigma = L^{-1}$. It follows that the quotient
of $\Cal U_{a-1}$ by the image of $\Cal U_a$ is exactly the direct image of
$\scrO_{\Cal H\times _BZ}\otimes \pi ^*L^{a-1}$, as claimed. The final
statement in (4.12) then follows by comparing elementary modifications. 
\endproof

Finally we shall need the analogue of Proposition 2.4 for a single elliptic
curve. It is proved exactly as in (2.4).

\theorem{4.13} Let $\Cal U'$ be a rank $n$ vector bundle over $\Pee \Cal
V_n\times _BZ$ such that, for all $x\in \Pee \Cal V_n$, $\Cal
U'|q_1^{-1}(x) \cong
\Cal U_0|q_1^{-1}(x)$. Then there is a unique line bundle $\Cal M$ on $\Cal
T$ such that, if $\pi _1\: \Cal T\times_BZ \to \Cal T$ is  projection
onto the first factor, then $\Cal U'
\cong (\nu\times
\Id)_*\left(\scrO_{\Cal T\times _BZ}(\Delta -
\Cal G)\otimes \pi _1^*\Cal M\right)$.
\qed
\endstatement

\ssection{4.4. Chern class calculations.}

Recall that we let $\zeta = c_1(\scrO_{\Cal P_{n-1}}(1))$, viewed as  a
class in
$H^2(\Cal P_{n-1})$. By pullback, we can also view $\zeta$ as an element of
$H^2(\Cal P_{n-1}\times _BZ)$. We also have the line bundle $\scrO_{\Cal
P_{n-1,d}}(1))$, and its first Chern class $\zeta'$ is given by $\zeta '
=\zeta - L$ (where we identify $L$ with its first Chern class in $H^2(B)$
and then by pullback in any of the relevant spaces).

\theorem{4.14} The Chern characters of the bundles $\bold U(d)$ and $\Cal
U_a$ are given by:
$$\gather
\ch \bold U(d)  =  (e^{-\sigma} + e^{-L} + \dots + e^{-(d-1)L})e^{\zeta -
L} + (e^{\sigma} + e^L + \dots + e^{(n-d-1)L}); \\
\ch \Cal U_a  =  e^{-\zeta}\fracwithdelims(){1 - e^{(a+n)L}}{1-e^L}-\frac{1
- e^{aL}}{1-e^L} +e^{-\sigma}(1-e^{-\zeta}).
\endgather$$
\endstatement
\proof The first statement is clear from the filtration on the $\Cal W_k$
and the definition of $\zeta$. To see the second, we use (4.12) for $1\leq
d\leq n-1$ and calculate 
$$\gather
\ch(\bold U (d)\otimes \scrO_{\Pee\Cal V_n}(-1)\otimes L) = \ch\bold U
(d)\cdot e^{-\zeta + L}=\\ (e^{-\sigma} + e^{-L} + \dots + e^{-(d-1)L}) +
(e^{\sigma +L} + e^{2L} + \dots + e^{(n-d)L})e^{-\zeta}\\ =(e^{-\sigma}-1) +
(1+e^{-L} + \dots + e^{-(d-1)L}) + \\ +(e^{\sigma +L-\zeta}-e^{L-\zeta}) +
(e^L + e^{2L} + \dots + e^{(n-d)L})e^{-\zeta}.
\endgather$$  Let $a = 1-d$. A little manipulation shows that we can write:
$$\align 1+e^{-L} + \dots + e^{-(d-1)L} &= - \frac{e^L - e^{aL}}{1-e^L};\\
e^L + e^{2L} + \dots + e^{(n-d)L} &= \frac{e^L - e^{(a+n)L}}{1-e^L};\\
(e^{-\sigma}-1) + (e^{\sigma +L-\zeta}-e^{L-\zeta}) &=
-(1-e^{-\sigma})(1-e^{\sigma + L-\zeta}).
\endalign$$ In the last term, note that $1-e^{-\sigma}$ is a power series
without constant term in $\sigma$ and thus annihilates every power series
without constant term in $\sigma +L$, since $\sigma ^2 = -L\cdot
\sigma$. Thus we can replace the last term by
$-(1-e^{-\sigma})(1-e^{-\zeta})$. It follows that
$$\align
\ch\bold U(d)\cdot e^{-\zeta+L}& =e^{-\zeta}\fracwithdelims(){e^L - 
e^{(a+n)L}}{1-e^L}-\frac{e^L - e^{aL}}{1-e^L} -
(1-e^{-\sigma})(1-e^{-\zeta})\\ & =e^{-\zeta}\fracwithdelims(){1 -
e^{(a+n)L}}{1-e^L}-\frac{1 - e^{aL}}{1-e^L} +e^{-\sigma}(1-e^{-\zeta}).
\endalign$$ In particular, we have established the formula in (4.14) for
$\ch \Cal U_a$ provided $a = 1-d$ with $1\leq d\leq n-1$. On the other
hand, the formula for
$\Cal U_a$ as an elementary modification shows that
$$\ch \Cal U_a = \ch \Cal U_{a-1} - \ch (\scrO_{\Cal H\times _BZ}\otimes
L^{a-1}).$$ Now from the exact sequence
$$0 \to \scrO_{\Cal P_{n-1}\times _BZ}(-\Cal H\times _BZ) \to \scrO_{\Cal
P_{n-1}\times _BZ} \to \scrO_{\Cal H\times _BZ}\to 0,$$ we see that 
$$\align
\ch (\scrO_{\Cal H\times _BZ}\otimes L^{a-1}) &= \ch (\scrO_{\Cal H\times
_BZ})\cdot e^{(a-1)L}\\ &= e^{(a-1)L}(1- e^{-\Cal H}).
\endalign$$ Next we claim:

\lemma{4.15}  $[\Cal H] = \zeta - nL$.
\endstatement
\proof We have identified $\Cal H$ with the image of $\Pee\Cal V_{n-1}$ in
$\Pee\Cal V_n$. The lemma now follows from the more general statement
below, whose proof is left to the reader:
\enddemo

\lemma{4.16} Let $\Cal V$ be a vector bundle over a scheme $B$, and suppose
that there is an exact sequence
$$0 \to \Cal V' \to \Cal V \to M \to 0,$$ where $M$ is a line bundle on
$B$. Let $\Cal H$ be the Cartier divisor
$\Pee(\Cal V') \subset \Pee (\Cal V)$. Then, if $p\: \Pee(\Cal V) \to B$ is
the projection,
$$\scrO_{\Pee (\Cal V)}(\Cal H) = \scrO_{\Pee (\Cal V)}(1) \otimes p^*M.
\qed$$
\endstatement

Plugging in the expression for $[\Cal H]$, we see that
$$\ch \Cal U_a - \ch \Cal U_{a-1} = - e^{(a-1)L}(1- e^{-(\zeta - nL)}).$$
Comparing this difference with the formula of (4.14) shows that (4.14)
holds for one value of $a$ if and only if it holds for all values of $a$.
Since we have already checked it for $a=0$, we are done.
\endproof

Similar computations give the Chern class of $\bold U(d)$ and $\Cal U_a$.
We leave the calculations to the reader.

\theorem{4.17} The total Chern class of $\bold U(d)$ is given by the
formula:
$$c(\bold U(d)) = (1+\zeta -L +\zeta\cdot
\sigma)\prod_{r=1}^{d-1}(1-(r+1)L+\zeta)\prod_{s=1}^{n-d-1}(1+sL).$$ If
$a\geq 0$, then
$$c(\Cal U_a) = (1-\zeta +L +\zeta\cdot
\sigma)\prod_{s=1}^{n+a-2}(1+(s+1)L-\zeta)\prod_{r=1}^{a-1}(1+ rL)^{-1}.$$
If $-(n-1)\leq a < 0$, then 
$$c(\Cal U_a) = (1-\zeta +L +\zeta\cdot
\sigma)\prod_{s=1}^{n+a-2}(1+(s+1)L-\zeta)\prod_{r=1}^{-a}(1- rL).$$ If $ a
< -(n-1)$, then 
$$c(\Cal U_a) = (1-\zeta +L +\zeta\cdot
\sigma)\prod_{s=0}^{1-n-a}(1-(s-1)L-\zeta)^{-1}\prod_{r=1}^{-a}(1-
rL).\qed$$
\endstatement

Let us work out explicitly the first two Chern classes of ${\Cal U}_a$.
First,
$$c_1(\Cal U_a) = \left[an + \fracwithdelims(){n^2-n}{2}\right]L -
(n+a-1)\zeta.$$ To give $c_2(\Cal U_a)$, write
$$\frac{1-e^{cx}}{1-e^x} = c + \fracwithdelims(){c^2-c}2x + P(c)x^2 +
\cdots,$$ where
$$P(c) = \frac{c(2c-1)(c-1)}{12} = \frac{2c^3 -3c^2 +c}{12}$$ (if $c$ is a
positive integer then $P(c) = \frac12\sum _{i=1}^{c-1}i^2$). A little
manipulation shows that $c_2(\Cal U_a)$ is equal to
$$\gather
\frac{(a+n-1)(a+n-2)}2\zeta ^2 -(n^2+2an-2n -a)\left(\frac{a+n-1}{2}\right)
\zeta\cdot L+\\
\left[\frac{1}{2}\left(an+\frac{n^2-n}{2}\right)^2-P(a+n)+P(a)\right]
L^2+(\sigma\cdot \zeta).
\endgather$$

Finally, we remark that it is possible to work out the first two terms in
$\ch
\Cal U_a$ by applying the Grothendieck-Riemann-Roch theorem directly to the
description of $\Cal U_a$ as $(\nu\times \Id)_*\scrO_{\Cal T\times
_BZ}(\Delta -\Cal G -a(r^*\sigma\times _BZ))$. This calculation is somewhat
long and painful, and does not give the full calculation of $\ch
\Cal U_a$ because $\Delta$ is not a Cartier divisor.

\section{5. Bundles which are regular and semistable on every fiber.}

So far in this paper we have been working universally with the moduli space
of all regular semistable bundles with trivial determinant over an elliptic
curve or an elliptic fibration. In this section we wish to study bundles
$V$ over an elliptic fibration $\pi\: Z\to B$ with the property that the
restriction of $V$ to every fiber is a regular semistable bundle with
trivial determinant.

\ssection{5.1. Sections and spectral covers.} 

Suppose that $V\to Z$ is a vector bundle of rank $n$ whose restriction to
each fiber is a regular semistable bundle with trivial determinant. Then
for each $b\in B$  the bundle $V|_{E_b}$ determines a point in the fiber of
${\Cal P}_{n-1}$ over
$b$.  This means that $V$ determines a section
$A(V) = A\colon B\to {\Cal P}_{n-1}$, as follows from (4.2). We shall
usually identify $A$ with the image $A(B)$ of $A$ in ${\Cal P}_{n-1}$.
Conversely, given a section $A$ of ${\Cal P}_{n-1}$ we can construct a
bundle $V$ over $Z$ which is regular semistable with trivial determinant on
each fiber and such that the section determined by
$V$  is $A$. There are many bundles with this property and we shall analyze
all such. 

We first begin by describing all sections of $\Cal P_{n-1}$.

\lemma{5.1}  A section $A\colon B\to {\Cal P}_{n-1}$ is equivalent to a
line bundle $M\to B$ and an inclusion of $M^{-1}$ into $\Cal V_n$, or
equivalently to sections of $M\otimes L^{-i}$ for $i=0,2,3,\ldots,n$ which
do not all vanish at any point of $B$, modulo the diagonal action of
$\Cee^*$. Under this correspondence, the normal bundle of $A$ in $\Cal
P_{n-1}$ is isomorphic to $(\Cal V_n\otimes M)/\scrO_B$, where the
inclusion of $\scrO_B$ in $\Cal V_n\otimes M$ corresponds to the inclusion
of $M^{-1}$ into $\Cal V_n$. Finally, if either $h^1(\scrO_B) =0$ or
$h^1(\Cal V_n\otimes M) = 0$, then the deformations of $A$ in $\Cal
P_{n-1}$ are unobstructed.
\endstatement
\proof
Let $A$ be a section, which we identify with its image in $\Cal P_{n-1}$. Of
course, $A\cong B$ via the projection
$p\: \Cal P_{n-1} \to B$.
 We have the inclusion of $\scrO_{\Pee\Cal V_n}(-1)$ in
$p^*\Cal V_n$. Pulling back via $A$, we set $M = \scrO_{\Pee\Cal
V_n}(1)|A$, which is a line bundle such that
$M^{-1}$ is a subbundle of
$p^*\Cal V_n|A = \Cal V_n$. An inclusion
$$M^{-1}\to \Cal V_n =
\scrO_B \oplus L^{-2} \oplus \cdots \oplus L^{-n}$$ is given by a nowhere
vanishing section of $(M \otimes \scrO_B) \oplus (M \otimes L^{-2}) \oplus
\cdots
\oplus (M \otimes L^{-n})$, or equivalently by  sections of the bundles
$(M\otimes \scrO_B)$, $(M\otimes L^{-2})$, \dots,
$(M\otimes L^{-n})$ which do not all vanish simultaneously, and these
sections are well-defined modulo the diagonal $\Cee^*$ action. Conversely,
a nowhere vanishing section of
$\Cal V_n\otimes M$ defines an inclusion $M^{-1} \to
\Cal V_n$ and thus a section of $\Cal P_{n-1}$, and the two constructions
are inverse to each other.

The normal bundle $N_{A/\Cal P_{n-1}}$ to $A$ in $\Cal P_{n-1}$ is just the
restriction to $A$ of the relative tangent bundle $T_{\Cal P_{n-1}/B}$, and
thus it is isomorphic to  $(\Cal V_n\otimes M)/\scrO_B$. The deformations
of the subvariety
$A$ are unobstructed if every element of $H^0(N_{A/\Cal P_{n-1}})$
corresponds to an actual deformation of $A$. From the exact sequence
$$0\to H^0(\scrO_B) \to H^0(\Cal V_n\otimes M) \to H^0(N_{A/\Cal P_{n-1}})
\to H^1(\scrO_B) \to H^1(\Cal V_n\otimes M),$$ we see that, if
$H^1(\scrO_B)=0$, then every section of the normal bundle lifts to a
section of $\Cal V_n\otimes M$, unique mod the image of $H^0(\scrO_B) =
\Cee$, and thus gives an actual deformation of $A$. If $H^1(\Cal V_n\otimes
M) = 0$, then viewing the deformations of $M$ as parametrized by $\Pic B$,
if $M'$ is sufficiently close to $M$ in $\Pic B$, then $H^1(\Cal V_n\otimes
M') = 0$ as well and by standard base change results the groups $H^0(\Cal
V_n\otimes M')$ fit together to give a vector bundle over a neighborhood of
$M$ in $\Pic B$. The associated projective space bundle then gives a smooth
family of deformations of
$A$ such that the associated Kodaira-Spencer map is an isomorphism onto
$H^0(N_{A/\Cal P_{n-1}})$. Thus $A$ is unobstructed in this case as well.
\endproof

\definition{Definition 5.2} Let $A\: B\to {\Cal P}_{n-1}$ be a section, and
let $(A,\Id)$ be the corresponding section of ${\Cal P}_{n-1} \times _BZ
\to  Z$. For all $a\in \Zee$, let
$$V_{A,a} = (A,\Id)^*\Cal U_a.$$

For every pair $(A,a)$, the bundle $V_{A,a}$ is of rank $n$ and the
restriction of
$V_{A,a}$ to every fiber of $\pi$ is regular and  semistable with trivial
determinant. Furthermore, for all $a\in \Zee$, the section determined by
$V_{A,a}$ is $A$.
\enddefinition

More generally, we could take any bundle ${\Cal U}$ over ${\Cal
P}_{n-1}\times_BZ$ obtained by twisting ${\Cal U}_a$ by a line bundle on
the universal spectral cover ${\Cal T}$ over ${\Cal P}_{n-1}$, and form
$V_{A,{\Cal U}}=(A,\text{Id})^*{\Cal U}$ to produce a bundle with these
properties. However, these will not exhaust all the possibilities in
general. To describe all possible bundles $V$ corresponding to $A$, we
shall need to define the spectral cover associated to $A$.

\definition{Definition 5.3} Let $A\subseteq {\Cal P}_{n-1}$ be a section.
The scheme-theoretic inverse  image $\nu ^*A$ of $A$ in $\Cal T$ is a
subscheme
$C_A$ of $\Cal T$, not necessarily reduced or irreducible. The  morphism
$g_A = \nu|A\: C_A\to A\cong B$ is finite and flat of degree
$n$. We call
$C_A$ the {\sl   spectral cover\/} associated to the section $A$.
\enddefinition

In the notation of (5.1), we shall show below that $C_A$ is smooth for $M$
sufficiently ample and for a general section corresponding to $M$. In
general, however, no matter how bad the singularities of $C_A$, we have the
following:

\lemma{5.4} The restriction of $r$ to $C_A$ embeds $C_A$ as a subscheme of
$Z$ which is a Cartier divisor. In fact, if $V$ is a vector bundle with
semistable restriction to every fiber and $A$ is the associated section,
then $C_A$ is the scheme of zeroes of $\det \Psi$, where 
$$\Psi\: \pi^*\pi_*(V\otimes \scrO_Z(\sigma)) \to V\otimes \scrO_Z(\sigma)$$
is the natural map. The line bundle $\scrO_Z(C_A)$ corresponding to $C_A$ is
isomorphic to $\scrO_Z(n\sigma) \otimes \pi ^*M$, where
$M$ is the line bundle corresponding to the section $A$. Moreover, the
image of
$C_A$ in $Z$ determines $A$. Finally, if $C\subset Z$ is the zero locus of a
section of $\scrO_Z(n\sigma) \otimes \pi ^*M$ and the induced morphism from
$C$ to $B$ is finite, then $C = C_A$ for a unique section $A$ of $\Cal
P_{n-1}$.
\endstatement 
\proof  Let $i\colon C_A\to \Cal T$ be the natural embedding. We claim that
$r\circ i\colon C_A\to Z$ is a scheme-theoretic embedding. To see this,
recall that we have $\Cal T \subset {\Cal P}_{n-1}\times _BZ$ via
$(\nu, r)$. In fact, from the defining exact sequence
$$0 \to \Cal E \to \pi^*\pi_*\scrO_Z(n\sigma) \to \scrO_Z(n\sigma) \to 0,$$
we see that $\Cal T=\Pee\Cal E$ is a Cartier divisor in
$\Pee(\pi^*\pi_*\scrO_Z(n\sigma)) ={\Cal P}_{n-1}\times _BZ$ defined by the
vanishing of a section of $\pi_2^*\scrO_Z(n\sigma)\otimes \pi_1^*\scrO_{\Cal
P_{n-1}}(1)$. Clearly, the image of
$i(C_A)$ under the map $C_A\to \Cal T \to {\Cal P}_{n-1}\times _BZ$ is an
embedding of $C_A$  in 
$A\times _BZ\cong Z$. Thus
$r\circ i$ is an embedding of $C_A$ into $Z$. Moreover, $C_A$ is the
restriction of $\Cal T\subset {\Cal P}_{n-1}\times _BZ$ to $A\times _BZ$,
and thus
$C_A$ is a Cartier divsior in $Z$. Essentially by definition, $C_A$ is
defined by the vanishing of $\det \Psi$ (since this holds on every fiber
$E_b$). Moreover,
$\scrO_Z(C_A)$ is the restriction to $A\times _BZ$ of 
$\pi_2^*\scrO_Z(n\sigma)\otimes \pi_1^*\scrO_{\Cal P_{n-1}}(1)$, namely
$\scrO_Z(n\sigma) \otimes \pi ^*M$.

Since the hypersurface $\Cal T\subset {\Cal P}_{n-1}\times _BZ$ is the
incidence correspondence, the line bundle $\scrO_{{\Cal P}_{n-1}\times
_BZ}(\Cal T)$ restricts on every fiber $E$ of $\pi$ to
$\scrO_E(np_0)$, and the effective divisor $\Cal T$ restricts to the
tautological divisor in $|np_0|\times E$ corresponding to the inclusion
$\Cal T\subset {\Cal P}_{n-1}\times _BZ$. Thus, by restriction, if
$\scrO_Z(C_A)$ is the line bundle in
$Z$ corresponding to the Cartier divisor
$C_A$, then for every fiber $E=E_b$ of $\pi$,
$\scrO_Z(C_A)|E\cong
\scrO_E(np_0)$, and the section of 
$\scrO_E(np_0)$ defined by $C_A$ is $A(b)$. Thus the image of $C_A$ in $Z$
determines $A$.

Finally, let $C$ be the zero locus of a section of $\scrO_Z(n\sigma)
\otimes \pi ^*M$. Note that 
$$H^0(Z; \scrO_Z(n\sigma) \otimes \pi ^*M) = H^0(B; \pi_*(\scrO_Z(n\sigma)
\otimes
\pi ^*M)) = H^0(B; \Cal V_n \otimes M),$$ so that sections $s$ of
$\scrO_Z(n\sigma)
\otimes \pi ^*M$ mod $\Cee^*$ correspond to sections $s'$ of $\Cal
V_n\otimes M$. Under this correspondence,
$s'$ vanishes at a point of $B$ if and only if $s$ vanishes along the
complete fiber  $\pi^{-1}(b)$. Thus we see that the subschemes $C$ mapping
finitely onto
$B$ are in $1-1$ correspondence with sections $A$ of $\Cal P_{n-1}$ whose
associated line bundle is $M$.
\endproof

We define
$T_A= C_A\times _BZ\subseteq
\Cal T\times _BZ$, and let $\rho _A\: T_A \to C_A$ be the natural map.
There is an induced map
$\nu _A\: T_A \to Z$ such that the following diagram is Cartesian:
$$\CD T_A @>{\nu _A}>> Z\\ @V{\rho _A}VV @VV{\pi}V\\ C_A@>{g _A}>> B.
\endCD$$ Thus, $T_A$ is an elliptic scheme over $C_A$ pulled back from the
elliptic scheme $Z\to B$ via the natural projection mapping  $C_A\to B$.
Even if $C_A$ is smooth, however, $T_A$ is singular along the intersection
of
$C_A\times _BZ$ with
$\Gamma\times _B\Gamma\subset Z\times _BZ$, at points corresponding to
$\Gamma
\cap C_A\subset Z$. If $\dim B =1$, the generic section $A$ will be such
that
$C_A\cap \Gamma =\emptyset$. However, if $\dim B \geq 2$ and $A$ is
sufficiently ample, $C_A\cap \Gamma$ is nonempty. In the generic situation
described in the last section, where $G_2$ and $G_3$ are smooth and meet
transversally, the singularities of $T_A$ are locally trivial families of
threefold double points. In general, if no component of $\Gamma$ is
contained in $C_A$, the codimension of
$C_A\cap \Gamma$ in $C_A$ is two and the codimension of the corresponding
subset of  $T_A$ is three. If a component of $\Gamma$ is contained in
$C_A$, then the codimension of
$C_A\cap \Gamma$ in $C_A$ is one and the codimension of the corresponding
subset of  $T_A$ is two. Note that  $\Delta$ is a Cartier divisor in the
complement of the subset of $T_A$ consisting of singular points of singular
fibers lying over $C_A\cap \Gamma$.

Let us examine the pullback to $T_A= C_A\times _BZ$ of the divisors in
${\Cal T}$. The section $\sigma \subset Z$ pulls back via $\nu_A^*$ to a
section
$\Sigma _A$ of the elliptic fibration $\nu _A\: T_A\to C_A$. Clearly
$\Sigma _A =\nu_A^*\sigma =
\Cal G|T_A$, where as in the last section $\Cal G$ is the pullback to $\Cal
T\times _BZ$ of $\sigma \subset Z$ by the second projection. The  diagonal
$\Delta_0$ in
$Z\times _BZ$ pulls back to a hypersurface in $T_A$, which is the
restriction of
$\Delta\subset \Cal T\times _BZ$ to $C_A\times _BZ = T_A$. We shall continue
to denote this subvariety by $\Delta$. However
$\Delta$ is not a Cartier divisor along the singular set of $T_A$. On the
other hand, the restriction of $\rho_A$ to $\Delta$ is an isomorphism from
$\Delta$ to
$C_A$, so that in a formal sense $\Delta$ is a section. There is also the
class
$\zeta$, which is obtained as follows: take the class $\zeta$ on $\Cal
P_{n-1}$, pull it back to $\Cal T$, and then restrict to $C_A$. In the
notation of (5.1), this class is just $\alpha = c_1(M)$, pulled back from
$B$.  The remaining ``extra" class $r^*\sigma\times _BZ|T_A$ corresponds to
$\sigma \cdot C_A=F$ in $Z$, and in particular it is pulled back from a
class on $C_A$. Note that $F$ maps isomorphically to its image in $B$.
Using $\nu_*r^*\sigma = \Cal H$, we see that the image of $F$ in $B$
corresponds to $A\cap \Cal H$. If  $D$ is  the divisor in $B$ corresponding
to $A\cap \Cal H$ and $V$ is a bundle with  semistable restriction to every
fiber whose associated section $A(V)$ is $A$, then $V|E_b$ has
$\scrO_E$ as a Jordan-H\"older quotient if and only if $b\in D$. The above
classes, together with the pullbacks of classes from
$B$,  are the only divisor classes that exist ``universally" on $C_A\times
_BZ = T_A$ for all sections $A$.

Using these classes, let us realize the bundles $V_{A,a}$ as pushforwards
from
$T_A$. Note that, from the definition, it is not {\it a priori\/} clear that
$(\nu_A)_*\scrO_{T_A}(\Delta-\Sigma _A)$ is locally free, since $\Delta$
need not be Cartier.

\lemma{5.5} For every section $A$ of ${\Cal P}_{n-1}$ and for every $a\in
\Zee$, we have
$$V_{A,a}=(\nu_A)_*\scrO_{T_A}(\Delta-\Sigma _A-aF).$$
\endstatement

\proof There is a commutative diagram, which is in fact a Cartesian square:
$$\CD T_A @>>> {\Cal T}\times_BZ \\ @V{\nu_A}VV   @VV{\nu\times \Id}V \\ Z
@>{(A,\Id)}>> {\Cal P}_{n-1}\times_BZ.
\endCD$$ Moreover, by definition $V_{A,a} = (A,\Id)^*(\nu\times
\Id)_*\scrO_{\Cal T\times _BZ}(\Delta - \Cal G-aF)$. The morphism
$\nu\times \Id$ is finite. Pulling back by the top horizontal arrow, the
sheaf $\scrO_{\Cal T\times _BZ}(\Delta - \Cal G-aF)$ restricts to
$\scrO_{T_A}(\Delta-\Sigma _A-aF)$. Thus, (5.5) is a consequence of the
following general result:
\enddemo

\lemma{5.6} Let 
$$\CD X'@>{f}>> X\\ @V{\pi'}VV @VV{\pi}V\\ Y'@>{g}>> Y
\endCD$$ be a Cartesian diagram of schemes, with $\pi$ a finite morphism.
Let $\Cal S$ be a sheaf on $X$. Then the natural map $g^*\pi_*\Cal S \to
(\pi')_*f^*\Cal S$ is an isomorphism.
\endstatement
\proof The question is local in $Y$ and $Y'$, so that we may assume that $Y
=\Spec R$ and $Y' =\Spec R'$ are affine. Since $\pi$ and $\pi'$ are finite,
and thus affine, we may thus assume that $X=\Spec S$ and $X' =\Spec S'$,
with $S'= S\otimes _RR'$. Suppose that $\Cal S$ corresponds to the
$S$-module $M$. Let $M_R$ be the
$S$-module $M$, viewed as an $R$-module. The assertion of the lemma is the
statement that
$$(M_R)\otimes _RR' \cong (M\otimes _SS')_{R'}.$$ But $M\otimes _SS' =
M\otimes _S(S\otimes _RR')$, and a standard argument now identifies
$(M\otimes _S(S\otimes _RR'))_{R'}$ with $(M_R)\otimes _RR'$. This proves
the lemma.
\endproof

Once we know that the sheaf $\scrO_{T_A}(\Delta-{\Sigma _A}-aF)$ pushes
down to a  vector bundle on $Z$, the same will be true for the twist of
this sheaf by any line bundle on $C_A$. Conversely, we have the following:

\proposition{5.7} Let $V$ be a vector bundle of rank $n$ on $Z$ such that
$V|E_b$ is a regular semistable bundle with trivial determinant for every
fiber $E_b$. Let $A=A(V)$ be the section determined by $V$ and let $C_A\to
A$ be the induced spectral cover. Then there is a unique bundle $N$ on 
$C_A$, such that
$V\cong (\nu_A)_*\left[\scrO_{T_A}(\Delta - \Sigma _A)\otimes
\rho_A^*N\right]$. \qed
\endstatement

The  proof of this result is similar to the proof of Part (ii) of Theorem
2.4 and will be omitted. 

Next we look at the  deformation theory of $V$.

\proposition{5.8} applying the Leray spectral sequence for $\pi\: Z \to B$
to compute $H^1(Z; Hom(V,V))$, there is an exact sequence
$$0\to H^1(B; \pi_*Hom (V,V)) \to H^1(Z; Hom(V,V)) \to H^0(B; R^1\pi_*Hom
(V,V)).$$
\roster
\item"{(i)}" The first term is $H^1(\scrO_{C_A})$ and corresponds to first
order deformations of a line bundle on the spectral cover $C_A$;
\item"{(ii)}"  If $L$ is not trivial, then $H^0(B; R^1\pi_*Hom (V,V))$ is
the tangent space to $A$ in the space of all sections of $\Cal P_{n-1}$,
and the restriction map 
$$H^1(Z; Hom(V,V)) \to H^0(B; R^1\pi_*Hom (V,V))$$ is the natural one which
associates to a first order deformation of $V$ a first order deformation of
the section $A(V)$.
\item"{(iii)}" Suppose that $L$ is nontrivial and that $C_A$ is smooth, or
more generally that
$h^1(\scrO_{C_A})$ is constant in a neighborhood of $A$. Suppose also either
that $h^1(\scrO_B) = 0$ or that $h^1(\Cal V_n\otimes M)=0$, which will hold
as soon as $M$ is sufficiently ample. Then the local moduli space of
deformations of
$V$ is smooth of dimension equal to
$h^1(Z; Hom(V,V))$. In other words, all first order deformations of $V$ are
unobstructed.
\endroster
\endstatement
\proof By construction $\pi_*Hom (V,V) = (g_A)_*\scrO_{C_A}$, and we leave
to the reader the check that the inclusion $H^1(B; \pi_*Hom (V,V)) \to
H^1(Z; Hom(V,V))$ corresponds to deforming the line bundle on $C_A$. Next,
let us fix for a moment a regular semistable bundle $V$ over a single
Weierstrass cubic $E$. Applying (1.5) with $S=\Cee[\epsilon]$, the dual
numbers, for every deformation of $V$ over $S$, there is an induced
morphism $S\to |np_0|$ which restricts over
$S_{\text{red}}$ to $\zeta (V)$. Thus there is an intrinsic homomorphism
from
$H^1(ad(V))$ to the tangent space $H^0(\scrO_E(np_0))/\Cee\cdot \zeta (V)$
of
$|np_0|$ at $\zeta(V)$. By (v) of Theorem 3.2, if
$V$ is a regular semistable bundle, then there is an exact sequence
$$0 \to \Cee \to H^1(W_{n-d}\spcheck\otimes W_d\spcheck) \to H^1(ad(V)) \to
0.$$ which identifies $H^1(ad(V))$ with the tangent space to $|np_0|$ at
$\zeta(V)$. Using the parametrized version of this construction (Lemma 4.3,
with $S$ equal to 
$\Cee[\epsilon]\times B$), there is an induced morphism from
$H^0(R^1\pi_*ad (V))$ to $\Hom(\Cee[\epsilon]\times B, \Cal P_{n-1};A)$,
the space of morphisms from $\Cee[\epsilon]\times B$ to $\Cal P_{n-1}$
extending the section
$A$. This gives an isomorphism from
$R^1\pi_*ad (V)$ to the relative tangent bundle $T_{\Cal P_{n-1}/B}$
restricted to
$A$. As we have seen in Lemma 5.1, this restriction is just the normal
bundle
$N_{A/\Cal P_{n-1}}$ to
$A$ in $\Cal P_{n-1}$. Clearly the map $H^0(R^1\pi _*ad(V)) \to
H^0(N_{A/\Cal P_{n-1}})$ is the natural map from the tangent space of
deformations of
$V$ to the tangent space to deformations of the section
$A$ in $\Cal P_{n-1}$. Now $Hom (V,V) = ad(V) \oplus \scrO_Z$, and so
$R^1\pi_*Hom (V,V) = R^1\pi_*ad (V)\oplus L^{-1}$.  Either $L^4$ or $L^6$
has a nonzero section, so that $L^{-1}$ has a nonzero section if and only
if $L$ is trivial. Thus, if $L$ is not trivial, then $H^0(L^{-1}) =0$, and
so
$$H^0(B; R^1\pi_*Hom (V,V)) = H^0(R^1\pi_*ad (V))$$  as claimed in (ii). 

To prove (iii), begin by using Lemma 5.1 to find a smooth space $Y$
parametrizing small deformations of the section $A$, of dimension
$h^0(N_{A/\Cal P_{n-1}})$.  If $\Cal A \to Y$ is the total space of this
family, there is an induced family of spectral covers $\Cal C \to Y$. By
assumption, the relative Picard scheme $\Pic(\Cal C/Y)$ is smooth in a
neighborhood of the fiber over $A$. Use this smooth space of dimension
$h^1(\scrO_{C_A})+ h^0(N_{A/\Cal P_{n-1}})$ to find a family of bundles
parametrized by a smooth scheme $S$, which is an open subset of $\Pic(\Cal
C/Y)$ and thus is fibered over the open subset $Y$ of sections of
$\Cal P_{n-1}$. This implies that the Kodaira-Spencer map of this family,
followed by the map from $H^1(Z;  Hom(V,V))$ to $H^0(B; R^1\pi_*Hom(V,V))$
is onto, and then that the Kodaira-Spencer map is an isomorphism onto
$H^1(Z;  Hom(V,V))$. Thus, the first order deformations of
$V$ are unobstructed.
\endproof

\ssection{5.2. Relationship to the extension point of view.}

Next we relate the description of bundles constructed out of  sections
$A$ of ${\Cal P}_{n-1}$ with the point of view of extensions. As usual,
this will enable us to construct some of the bundles previously constructed
via spectral covers, but not all.

We have already constructed the bundles
$\Cal W_k$ over $Z$ as well as the universal extension $\bold U(d)$,
$1\le d<n$, which sits in an exact sequence
$$0 \to \pi _2^*\Cal W_d\spcheck \otimes \pi _1^*\scrO_{\Cal P_{n-1,d}}(1)
\to
\bold U(d) \to \pi _2^*\Cal W_{n-d} \to 0.$$ Here the projective space
$\Cal P_{n-1,d}$ of the vector space of extensions is identified with
${\Cal P}_{n-1}$, but, by Theorem 4.6, under this identification 
$$\scrO_{\Cal P_{n-1,d}}(1) \otimes \pi^*L =\scrO_{{\Cal P}_{n-1}}(1).$$
Finally, we have
$$\bold U(d) = \Cal U_{1-d}\otimes \pi _1^*\scrO_{{\Cal P}_{n-1}}(1)\otimes
L^{-1}.$$ Thus there is an exact sequence
$$0 \to \pi _2^*\Cal W_d\spcheck \to \Cal U_{1-d} \to \pi _2^*\Cal
W_{n-d}\otimes  
\pi _1^*\scrO_{{\Cal P}_{n-1}}(-1)\otimes L \to 0.$$

Given a section $A$ of $\Cal P_{n-1,d}={\Cal P}_{n-1}$ such that
$\scrO_{\Cal P_{n-1,d}}(1)|A = M'$, we can pull back the defining extension
for $\bold U(d)$ to obtain an extension
$$0 \to \Cal W_d\spcheck \otimes \pi^*M' \to U_A \to \Cal W_{n-d} \to 0.$$
(Of course, $M'$ is  $M\otimes L^{-1}$.)  Conversely, suppose that we are
given an extension of
$\Cal W_{n-d}$ by
$\Cal W_d\spcheck \otimes \pi^*M'$, where $M'$ is a line bundle on $B$
which we can write as $M\otimes L^{-1}$. In this case, by the Leray
spectral sequence
$$\gather H^1(\Cal W_{n-d}\spcheck \otimes \Cal W_d\spcheck \otimes
\pi^*M') \cong H^0(R^1\pi _*(\Cal W_{n-d}\spcheck \otimes \Cal W_d\spcheck)
\otimes M')\\= H^0(\Cal V_{n,d}\otimes M\otimes L^{-1})=H^0(\Cal V_n\otimes
M).
\endgather$$ Thus nontrivial extensions of
$\Cal W_{n-d}$ by
$\Cal W_d\spcheck \otimes \pi^*M'$ which restrict to nontrivial extensions
on every fiber can be identified with sections of
$\Cal P_{n-1,d}$ corresponding to the line bundle $M$. Finally, we see
that, for
$1\leq d\leq n-1$, we can write $V_{A,1-d}$ as an extension
$$0 \to \Cal W_d\spcheck \to V_{A, 1-d} \to \Cal W_{n-d}\otimes
\pi^*(M^{-1}\otimes L)\to 0.$$

We can also relate the deformation theory of $U_A$ above to the bundles
$\Cal W_d$ and $\Cal W_{n-d}$. Thus, the tangent space to
$\Ker\{\,(g_a)_*\: \Pic C_A \to
\Pic B\,\}$ is $H^1(B; \pi _*(\Cal W_d\otimes \Cal W_{n-d})\otimes
M^{-1}\otimes L)$, and the tangent space to deformations of the section $A$
is 
$H^0(B; R^1\pi _*(\Cal W_d\spcheck\otimes \Cal W_{n-d}\spcheck)\otimes
M\otimes L^{-1})$, provided that $L$ is not trivial.

\ssection{5.3. Chern classes and determinants.}

Let $A$ be a section of $\Cal P_{n-1}$. Corresponding to $A$, there is the
line bundle $M$ on $B$ which is the restriction to $A$ of $\scrO_{\Cal
P_{n-1}}(1)$.   We denote by $\alpha$ the class $c_1(M) \in H^2(B;\Zee)$.
Our goal is to express the Chern classes of $V_{A,a}$ in terms of $\alpha$
and the standard classes on
$Z$. We will also consider more general bundles arising from twisting by a
line bundle on the spectral cover.

First we shall determine the Chern classes of $V_{A,a}$. We begin with the
following lemma:

\lemma{5.9} Let  $A$ be a section of $\Cal P_{n-1}$ corresponding to the  
inclusion of a line bundle $M^{-1}$ in
$\Cal V_n$. Then, for $k\geq 0$, we have 
$p_*([A]\cdot \zeta ^k) =\alpha ^k\in H^{2k}(B;\Zee)$.
\endstatement
\proof Note that by definition  $\zeta|_A=c_1(M)=\alpha$ when we identify
$A$ and
$B$ in the obvious way.  It follows that $\zeta^k|_A=\alpha^k$. This means
that $p_*([A]\cdot \zeta^k)=\alpha^k$.
\endproof

Using (5.9), we can  compute the Chern classes $c_i(V_{A,a})$ by taking the
formula for $c_i(\Cal U_a)$  and replacing $\zeta ^i$ by $\alpha^i$. Thus

\theorem{5.10} Suppose that  $A$ is a section of
${\Cal P}_{n-1}$ such that the corresponding line bundle $M$ has
$c_1(M)=\alpha \in H^2(B)$ \rom(or $\Pic B$\rom). Then 
$$\ch(V_{A,a}) = e^{-\alpha}\fracwithdelims(){1 -
e^{(a+n)L}}{1-e^L}-\frac{1 - e^{aL}}{1-e^L} +e^{-\sigma}(1-e^{-\alpha}).$$
Moreover, in $\pi^*\Pic B \subset \Pic Z$,
$$\det (V_{A,a}) = -(n+a-1)\alpha + \left[an +
\fracwithdelims(){n^2-n}{2}\right]L.\qed$$
\endstatement

There is also a formula for $c(V_{A,a})$ which follows similarly from the
formula for $c(\Cal U_a)$.

Now let us consider the effect of twisting by a line bundle on the spectral
cover. If $N$ is a line bundle  on the spectral cover $C_A$ associated to 
$A$, let
$$V_{A, 0}[N] = (\nu_A)_*\left[\scrO_{T_A}(\Delta - \Sigma _A )\otimes
\rho_A^* N\right].$$ For example, suppose that $N$ is of the form
$\scrO_{C_A}(-aF)\otimes g_A^*N_0$, where $N_0$ is a line bundle on $B$.
Then 
$$V_{A, 0}[N] =  V_{A, a}\otimes \pi^*N_0.$$ In particular, we see that if 
$N=\scrO_{C_A}(-aF)\otimes g_A^*N_0$, for some  line bundle $N_0$ on $B$
and some integer $a$, then 
$$\ch(V_{A, 0}[N])= \left[e^{-\alpha}\fracwithdelims(){1 -
e^{(a+n)L}}{1-e^L}-\frac{1 - e^{aL}}{1-e^L}
+e^{-\sigma}(1-e^{-\alpha})\right]\cdot e^{c_1(N_0)}.$$

For more general line bundles $N$ on $C_A$, we can calculate the
determinant of
$V_{A, 0}[N]$. In what follows, we identify $\Pic B$ with a subgroup of
$\Pic Z$ via $\pi^*$ and write the group law additively.

\lemma{5.11} With $V_{A, 0}[N]$ as defined above, the following formula
holds in
$\Pic B$:
$$c_1(V_{A, 0}[N]) = -(n-1)\alpha + \fracwithdelims(){n^2-n}{2}L +
(g_A)_*c_1(N).$$ Thus, for a fixed section $A$ of $\Cal P_{n-1}$ and a
fixed line bundle $\Cal N$ on
$B$, the set of bundles
$V$ on $Z$ which are regular semistable on every fiber, with $A(V) = A$ and
$\det V =\pi^*\Cal N$ is a principal homogeneous space over $\Ker
\{g_A{}_*\: \Pic C_A \to \Pic B\}$, which is a generalized abelian variety
times a finitely generated abelian group. 
\endstatement
\proof Since it is enough to compute the determinant in the complement of a
set of codimension two, we may restrict attention to the open subset of
$T_A$ where
$\Delta$ is a Cartier divisor. Now it is a general formula that, for a
Cartier divisor $D$ on
$T_A$, 
$$c_1 \left[(\nu_A)_*\scrO_{T_A}(D)\right] = c_1
\left[(\nu_A)_*\scrO_{T_A}\right]+ (\nu_A)_*D.$$ Thus, applying this
formula to  $\scrO_{T_A}(\Delta -\Sigma _A)$ and $\scrO_{T_A}(\Delta
-\Sigma _A)\otimes
\rho_A^* N$, we see that
$$c_1(V_{A, 0}[N]) = c_1(V_{A, 0}) + (\nu_A)_*\rho_A^* c_1(N).$$ But we
have calculated $c_1(V_{A, 0}) = \dsize  -(n-1)\alpha +
\fracwithdelims(){n^2-n}{2}L$, and $(\nu_A)_*\rho_A^* c_1(N) =
\pi^*(g_A)_*c_1(N)$ since $T_A = C_A\times _BZ$. Putting these together
gives the formula in (5.11).
\endproof

If $\dim B \geq 2$ and $M$ is sufficiently ample, we we will see in the next
subsection that the generalized abelian variety $\Ker
\{g_A{}_*\: \Pic C_A \to \Pic B\}$ is in fact a finitely generated abelian
group, with no component of positive dimension.

Using (5.11), let us consider the following problem: Given the section $A$,
when can we find a line bundle $N$ such that $V_{A, 0}[N]$ actually has
trivial determinant?  We are now in position to answer this question in
this case if we consider twisting only by line bundles which exist
universally for all spectral covers.

\proposition{5.12}  Given a section $A$, suppose that
$N=\scrO_{C_A}(-aF)\otimes g_A^*N_0$ for a line bundle $N_0$ on $B$ and an
integer $a$. Then 
$V_{A, 0}[N]$ has trivial determinant for some choice of an $N$ as above if
at least one of the following conditions holds:
\roster
\item"{(i)}" $n$ is odd,
\item"{(ii)}" $L$ is divisible by $2$ in $\Pic B$, or
\item"{(iii)}" 
$\alpha\equiv L\bmod 2$ in
$\Pic B$.
\endroster
\endstatement
\proof It suffices to show that there exists an $a\in \Zee$ such that
$\det (V_{A,a})$ is divisible by $n$. For then, for an appropriate line
bundle $N_0$ on
$B$, we can arrange that $V=V_{A,a}\otimes N$ has trivial determinant. By
(5.10),  we must have
$$(a-1)\alpha \equiv \frac{n(n-1)}2L \bmod n.$$ In the first two cases we
simply take 
$a\equiv 1 \bmod n$. Lastly, let us suppose that $n$ is even and that
$L$ is  not divisible by $2$. Then the condition $\dsize (a-1)\alpha
\equiv
\frac{n(n-1)}2L \bmod n$  is a nontrivial condition on 
$\alpha$. It is satisfied for the appropriate $a$ if $\alpha\equiv L\bmod
2$ in
$\Pic B$.
\endproof

We leave it to the reader to write out necessary and sufficient conditions
for the equation $\dsize (a-1)\alpha \equiv \frac{n(n-1)}2L \bmod n$ to
have a solution in general.

For a general line bundle  $N$ on $C_A$, we can use the
Grothendieck-Riemann-Roch theorem to calculate the higher Chern classes of
$\ch(V_{A, 0}[N])$, but only in the range where $\Delta$ is a Cartier
divisor. Thus, we are essentially only able to compute $c_2$ by this method
for a general line bundle $N$:

\proposition{5.13} Suppose that no component of $\Gamma$ is contained in
$C_A$. Let
$\ch_2$ be the degee two component of the Chern character. Then
$$\gather
\ch_2(V_{A, 0}[N])- \ch_2(V_{A,0})= \\ (\nu_A)_*\left(\left( \Delta -\Sigma
_A +
\frac12(\nu_A^*K_Z- K_{T_A})\right)\cdot (\rho_A)^*( N)
\right)+(\pi_A)^*(g_A)_*\frac{( N)^2}{2}.
\endgather$$
\endstatement  

\proof Working where $\Delta$ is Cartier, we can apply the
Grothendieck-Riemann-Roch theorem to the local complete intersection
morphism
$\nu_A\: T_A\to Z$ to determine the Chern character of $V_{A,0} =
(\nu_A)_*\scrO_{T_A}(\Delta - \Sigma _A)$:
$$\ch (V_{A,a}) = (\nu_A)_*\left(e^{\Delta -\Sigma}\Todd(T_A/Z)\right),$$
valid under our assumptions through terms of degree two. Applying the same
method to calculate the Chern character of $V_{A, 0}[N]$, we find that, at
least through degree two,
$$\ch (V_{A, 0}[N])- \ch (V_{A,0})= (\nu_A)_*\left((e^N-1)(e^{\Delta
-\Sigma}\Todd(T_A/Z)\right).$$ Expanding this out gives (5.13).
\endproof

\ssection{5.4. Line bundles on the spectral cover.}

In this section, we look at the problem of finding extra line bundles on the
spectral cover $C_A$, under the assumption that $C_A$ is smooth and that
$M$ is sufficiently ample. As we shall see, the discussion falls naturally
into three cases: $\dim B =1$, $\dim B = 2$, $\dim B \geq 3$. 

First let us consider the case that $B$ is a curve, with $M$ arbitrary but
$C_A$ assumed to be smooth, or more generally reduced.  Let
$A$ correspond to the line bundle $M$ on $B$. Given $V=V_{A, 0}[N]$, we seek
$\det V$ and $c_2(V)$.  First, by (5.11), working in $\Pic B$ written
additively,
$$\det V_{A, 0}[N]= -(n -1)M +
\fracwithdelims(){n^2-n}{2}L+ (g_A)_*N.$$   Since $g_A{}_*\: \Pic C _A\to
\Pic B$ is surjective in  case $C_A$ is  reduced, we can arrange that the
determinant is in fact trivial, and then the line bundle
$\scrO_{C_A}(D)$ is determined up to the subgroup $\Ker \{g_A{}_*\: \Pic
C_A \to \Pic B\}$. If $C_A$ is smooth, then this subgroup is the product of
an abelian variety and a finite group

We may summarize this discussion as follows:

\theorem{5.14} Suppose that $\dim B =1$. Given a section $A$ of $\Cal
P_{n-1}$ such that
$C_A$ is reduced, the set of bundles $V$ with trivial determinant such that
$A(V) = A$ is a nonempty principal homogeneous space over $\Ker \{g_A{}_*\:
\Pic C_A \to \Pic B\}$. The same statement holds if we replace the
condition that $V$ has trivial determinant by the condition that the
determinant of $V$ is $\pi ^*\lambda$ for some fixed line bundle $\lambda$
on $B$.
\qed
\endstatement

The remaining Chern class is $c_2(V)$. In this case, in $H^4(Z; \Zee)$,
with no assumptions on $\Gamma$, we have (as computed in
\cite{3} in case $n=2$):

\proposition{5.15} For every line bundle $N$ on $C_A$,
$$c_2(V_{A,0}[N]) = c_2(V) = \sigma \cdot \alpha = \deg M.$$
\endstatement
\proof First assume that $C_A$ is reduced. Write $N \cong \scrO_{C_A}(\sum
_ip_i)$, where the $p_i$ are points in the smooth locus of $C_A$ which lie
under smooth fibers. Thus $\rho_A^{-1}(p_i) = f_i$ is a smooth fiber of
$T_A$. In this case, we can obtain $V_{A,0}[N]$ as a sequence of elementary
modifications of the form 
$$0 \to V_{A,0}[N_j] \to V_{A,0}[N_{j+1}] \to (i_j)_*\lambda _j \to 0,$$
where $E_j$ is the fiber on $Z$ corresponding to $f_j\subset T_A$, $i_j\:
E_j \to Z$ is the inclusion, and $\lambda _j = \scrO_{T_A}(\Delta - \Sigma
_A)|f_j$ is a line bundle of degree zero. By standard calculations, 
$$c_2(V_{A,0}[N_j]) = c_2(V_{A,0}[N_{j+1}])$$ and so $c_2(V_{A,0}[N]) =
c_2(V_{A,0}) = \sigma \cdot \alpha$.

In case $C_A$ is not reduced, a similar argument applies, where we replace
$p_i$ by a Cartier divisor whose support is contained in the smooth locus of
$(C_A)_{\text{red}}$ and
$E_j$ by a thickened fiber.
\endproof

\remark{Remark} On the level of Chow groups, the refined Chern class $\tilde
c_2(V_{A,0}[N])$ essentially records the extra information coming from the
natural map $\Pic C_A\to A^2(Z)$.
\endremark
\medskip

Next we consider the case where $\dim B > 1$. First we have the following
result, with no assumption on $C_A$, concerning the connected component of
$\Pic C_A$.

\lemma{5.16} Suppose that $\dim B \geq 2$ and that $M$ is sufficiently
ample. More precisely, suppose that 
$$H^i(B; L^{-1}\otimes M^{-1}) = H^i(B; L\otimes M^{-1}) = \cdots = H^i(B;
L^{n-1}\otimes M^{-1}) = 0$$ for $i=0,1$. Then the natural map from $H^1(Z;
\scrO_Z)$ to $H^1(C_A; \scrO_{C_A})$ is an isomorphism. Finally, if in
addition $L$ is not trivial, then the norm map from
$\Pic^0C_A$ to $\Pic ^0B$ is surjective with finite kernel. Thus $\Ker
\{g_A{}_*\: 
\Pic C_A \to \Pic B\}$ is a finitely generated abelian group.
\endstatement
\proof From the exact sequence
$$0 \to \scrO_Z(-n\sigma) \otimes \pi^*M^{-1} \to \scrO_Z \to \scrO_{C_A}
\to 0,$$ we see that there is a long exact sequence
$$H^1(\scrO_Z(-n\sigma) \otimes \pi^*M^{-1}) \to H^1(\scrO_Z) \to
H^1(\scrO_{C_A})
\to H^2(\scrO_Z(-n\sigma) \otimes \pi^*M^{-1}).$$ Applying the Leray
spectral sequence to $\scrO_Z(-n\sigma) \otimes \pi^*M^{-1}$, we have that
$$H^i(\scrO_Z(-n\sigma) \otimes \pi^*M^{-1}) = H^{i-1}(R^1\pi
_*\left[\scrO_Z(-n\sigma) \otimes \pi^*M^{-1}\right]).$$ Now, by duality,
$$\gather R^1\pi_*\left[\scrO_Z(-n\sigma) \otimes \pi^*M^{-1}\right] = 
R^1\pi _*\scrO_Z(-n\sigma) \otimes M^{-1}\\
 = \left(L^{-1}\oplus L\oplus \cdots \oplus L^{n-1}\right)  \otimes M^{-1}.
\endgather$$ Thus by our assumptions the map $H^1\scrO_Z) \to
H^1(\scrO_{C_A})$ is an isomorphism.  By  applying the Leray spectral
sequence to $\scrO_Z$, we see that there is an exact sequence
$0\to H^1(\scrO_B) \to H^1\scrO_Z) \to H^0(L^{-1})$. As we saw in the proof
of (ii) of (5.8), if $L$ is not trivial, then $H^0(L^{-1}) = 0$ and the
pullback map
$H^1(\scrO_B)
\to H^1\scrO_Z)$ is an isomorphism. The last statement of the lemma is then 
clear.
\endproof

\lemma{5.17} If $M$ is sufficiently ample on $B$, then
$C_A$ is an ample divisor in $Z$.
\endstatement

\proof Equivalently, we must show that for $M$ sufficiently ample on $B$,
$\pi ^*M\otimes \scrO_Z(n\sigma)$ is  ample. But $\scrO_Z(n\sigma)$ is
relatively ample, and thus by a standard result $\pi ^*M\otimes
\scrO_Z(n\sigma)$ is  ample for $M$ sufficiently ample (compare \cite{10,
p\. 161. (7.10)(b)} for the case where $\scrO_Z(n\sigma)$ is relatively
very ample).
\endproof

\corollary{5.18} If $\dim B \geq 3$, $M$ is sufficiently ample, and $Z$ and
$C_A$ are smooth, then $\Pic Z \cong \Pic C_A$. If $\dim B =2$, $M$ is
sufficiently ample,  and  $Z$ and $C_A$ are smooth, then the restriction
mapping $\Pic Z \to
\Pic C_A$ is injective.
\endstatement
\proof This is immediate from the Lefschetz theorem and (5.17).
\endproof

\remark{Remark} If $\dim B = 2$ and $M$ is sufficiently ample, it is natural
to expect an analogue of the Noether-Lefschetz theorem to hold: for generic
sections $C_A$ of $\pi ^*M\otimes \scrO_Z(n\sigma)$,
$\Pic Z \cong \Pic C_A$. However, in the next section, we will see how to
construct sections $A$ such that the spectral cover $C_A$ is smooth but has
larger Picard number than expected.
\endremark

\ssection{5.5. Symmetric bundles.}

Next we turn to bundles with a special invariance property.

\definition{Definition 5.19} Let $\iota\: Z \to Z$ be the involution which
is
$-1$ in every fiber. A bundle $V$ is {\sl symmetric\/} if $\iota^*V \cong
V\spcheck$.
\enddefinition

We shall now analyze when a bundle $V$  is symmetric. We fix a section $A$,
corresponding to the class $\alpha$ and denote $C_A, \nu _A$, $T_A$, $g_A$
simply by
$C, \nu, T,g$.

\proposition{5.20} For a suitable choice of $N\in
\Pic C$ the bundle $V_{A,0}[N]$ is symmetric if and only if
$g^*(L+\alpha) + nF$ is divisible by $2$ in $\Pic C$. In this case, for a
fixed section $A$, the set of all symmetric bundles whose section is $A$ is
a principal homogeneous space over the $2$-torsion in $\Pic C$. 
\endstatement
\proof  Suppose that $V=V_{A,0}[N] = \nu _*\left[\scrO_T(\Delta -
\Sigma_A )\otimes \rho^*\scrO_{C}(N)\right]$, where $N$ is a divisor on
$C$.   For our purposes,  since both $\iota^*V$ and $V\spcheck$ are
bundles, they are isomorphic if and only if they are isomorphic outside the
complement of a set of codimension two in $Z$.  Thus, we shall work as if
$\Delta$ is a Cartier divisor.

There is an induced involution on $T$, also denoted by $\iota$, for which
$\nu$ is equivariant. Thus 
$$\align
\iota ^*V&= \iota^*\nu _*\left[\scrO_T(\Delta  - \Sigma _A )\otimes
\rho^*\scrO_{C}(N)\right]\\ &=\nu_*\iota^*\left[\scrO_T(\Delta  - \Sigma
_A)\otimes
\rho^*\scrO_{C}(N)\right].
\endalign$$ Now $\iota^*\Sigma _A = \Sigma _A$ and $\iota
^*\rho^*\scrO_{C}(N) =
\rho^*\scrO_{C}(N)$.  One the other hand, $\iota ^*\Delta$ is linearly
equivalent to
$2\Sigma _A -\Delta$ on a generic fiber. This says that
$$\iota ^*\Delta = 2\Sigma _A -\Delta + \rho^*D$$ for some divisor
$D$ on $C$. To determine $D$, restrict both sides above to
$\Sigma _A$ where $\iota$ acts trivially. We find that
$D = 2\Delta \cdot \Sigma _A -2\Sigma _A^2$, viewed in the obvious way as a
divisor class on $C$. Thus
$$\iota ^*\Delta = 2\Sigma _A -\Delta + 2\rho^*D_0$$ where $D_0$ is the
fixed divisor class  $\Delta \cdot \Sigma _A -\Sigma _A^2$, viewed as a
divisor on $C$. Here the main point will be the factor of
$2$. However we  note that
$\Sigma _A^2 = -[L']$, where $L'= g^*L$ is the line bundle for the elliptic
scheme $T$, and
$$\Delta \cdot \Sigma _A  = \Delta \cdot \nu^*\sigma =
\nu^*(\nu_*\Delta )\cdot \sigma =\nu ^*(C\cdot \sigma),$$ which after
pullback corresponds to the divisor class $F$ on $T$. (Here $\Delta =
C\times _BC\subset C\times _BZ$, and so $\nu_*\Delta =C$ since $\nu$ is
just the natural projection of $T=C\times _BZ$ to $Z$.)

Next we calculate $V\spcheck$. Relative duality for the finite flat morphism
$\nu$ says that, for every Cartier divisor $D$ on $T$, $\left[\nu
_*\scrO_T(D)\right]\spcheck = \nu _*\left[\scrO_T(-D)\otimes
K_{T/Z}\right]$, where $K_{T/Z} = K_T\otimes
\nu^*K_Z^{-1}$ is the relative dualizing sheaf of the morphism $\nu$. Thus
we must have
$$\Sigma _A -\Delta   -\rho ^*N + K_T-\nu ^*K_Z  = \Sigma _A -\Delta +
2\rho^*D_0  +\rho ^*N.$$ Equivalently, we must have
$$K_T-\nu ^*K_Z  = 2\rho ^*N + 2\rho^*D_0= 2\rho ^*(N +[L'] + F).$$ 
Conversely, given that the above equality holds, the corresponding vector
bundles will be symmetric. To see if this equality holds for the
appropriate choice of  $N$, we must calculate $K_T-\nu ^*K_Z$. Since $Z$ is
an elliptic fibration, $K_Z =
\pi ^*(K_B+L)$, and likewise $K_T = \rho ^*(K_C+L')$, where $L' = g^*L$.
Thus
$K_T-\nu ^*K_Z = \rho^*(K_C - g^*K_B)$. To calculate $K_C$, we use (5.4),
which says that $K_C = K_Z+C|C = K_Z + \pi^*L + n\sigma |C$. On the other
hand, $K_Z - \pi ^*K_B = \pi ^*L$. Restricting to $C$ gives:
$$K_C - g^*K_B = g^*(L+\alpha) + nF.$$ Putting this together, we see that,
if $V$ is symmetric, then we must have $\rho^*(g^*(L+\alpha) + nF)$
divisible by $2$ in $\rho^*\Pic C$, and conversely. 

Next we claim that $\rho^*\: \Pic C \to \Pic T$ is injective. It suffices to
show that $\rho_*\scrO_T =\scrO_C$, for then $\rho_*\rho^*N = N$ for every
line bundle $N$ on $C$. But by flat base change $g^*\pi_*\scrO_Z = \rho_*\nu
^*\scrO_Z = \rho_*\scrO_T$. Since $\pi_*\scrO_Z =\scrO_B$, we have that
$g^*\pi_*\scrO_Z = g^*\scrO_B = \scrO_C = \rho_*\scrO_T$. Hence
$\rho_*\scrO_T =\scrO_C$, and so $\rho^*$ is injective.

Thus, $V$ is symmetric if and only if $g^*(L+\alpha) + nF$ divisible by $2$
in
$\Pic C$. Moreover the set of possible line bundles $N$ for which
$V_{A,0}[N]$ is symmetric is a principal homogeneous space over the
$2$-torsion in $\Pic C$, as claimed. This concludes the proof of (5.20).
\endproof

If $\dim B \geq 3$, $Z$ and $C$ are smooth, and $M$ is sufficiently ample,
then
$g^*(L+\alpha) + nF$ is divisible by $2$ in $\Pic C$ if and only if
$\pi^*(L+\alpha) + n\sigma$ is divisible by $2$ in $\Pic Z$. This can only
happen if $n$ is even and $\alpha \equiv L \bmod 2$. A similar statement is
likely to hold if $\dim B =2$ and $A$ is also assumed to be general.

We can see the conditions  $n$ is even and $\alpha
\equiv L \bmod 2$ clearly in terms of extensions. In this case $n=2d$, and
we can  write $V_{A,1-d}$ as an extension
$$0 \to \Cal W_d\spcheck \to V_{A, 1-d} \to \Cal W_d\otimes M^{-1}\otimes
L\to 0.$$ Under the assumption that $M^{-1}\otimes L =M_0^{\otimes 2}$ for
some line bundle
$M_0$, we can write $V_{A, 1-d}\otimes M_0^{-1}$ as an extension of $\Cal
W_d\otimes M_0$ by the dual bundle $\Cal W_d\spcheck \otimes M_0^{-1}$, and
then check directly that the corresponding bundles are symmetric.

\ssection{5.6. The case of the trivial section.}

We turn to bundles which have reducible or non-reduced spectral covers. We
begin with the extreme case of the trivial section
 $\bold o = \bold o_Z =\Pee \scrO_B \subset {\Cal P}_{n-1}$. To construct
this section we take $M=\scrO_B$ and take a nowhere vanishing section of
$\scrO_B$ and the zero section of $L^{-a}$ for all $a>0$.  Since
$M=\scrO_B$, the class $\alpha$ is zero. The spectral cover $C=C_{\bold
o}\subset Z$ is simply the nonreduced scheme 
$n\sigma$, and the associated reduced subscheme $C_{\text{red}}$ is
identified with $B$.  The bundles associated to this section have the
property that their restrictions to each fiber of $Z$ are isomorphic to
$I_n(\scrO)$.  Conversely, if we have such a bundle $V$ over $Z$, then the
section it determines is $\bold o$.

By our general existence theorem we immediately conclude:

\corollary{5.21} For every $n\ge 1$ there is a vector bundle $V\to Z$ whose
restriction to each fiber $E_b\subset Z$ is isomorphic to
$I_n(\scrO_{E_b})$. \qed
\endstatement

The structure sheaf $\scrO_C$ is filtered by subsheaves with successive
quotients
$$L^{n-1}, L^{n-2}, \dots, \scrO_B.$$ The restriction of
$\scrO_{C\times _BZ}(\Delta -\Sigma _A -aF)$ to $C_{\text{red}}\times
_BZ\cong Z$, is isomorphic to 
$$\scrO_{C\times _BZ}(\Delta -\Sigma _A -aF)|(C_{\text{red}}\times _BZ)
\cong \scrO_Z(\sigma -\sigma)\otimes L^a) =L^a.$$ From this it follows that
$V_{\bold o, a}$ has a filtration by subbundles  with successive quotients
$L^{a+n-1}, L^{a+n-2}, \dots, L^a$. Consequently,
$$\ch (V_{\bold o, a}) = \frac{e^{aL} - e^{(a+n)L}}{1-e^L},$$ which agrees
with the formula in Theorem 5.10 since  $\alpha = 0$.

We have the inclusion $B=C_{\text{red}}\subset C$ and the projection
$C\to B$ so that $\scrO_C$ splits as a module over $\scrO_B$  into
$\Cal S\oplus \scrO_B$ with $\Cal S$  a locally free sheaf of rank $n-1$
over $\scrO_B$. 
 From the filtration of $\scrO_C$ as an $\scrO_B$-module, we see that
$\Cal S$ has a filtration with successive quotients 
 $L^{n-1}, L^{n-2}, \dots,L$. Thus, $\Pic C \cong\Pic B \oplus H^1(\Cal
S)$, and
$H^1(\Cal S)$ is a vector group. In particular, as far as Chern classes are
concerned, we may as well just twist by line bundles $N$ on $C$ which are
pulled back from
$B$. Even if the line bundle
$N$ on $C$ is not pulled back from $B$, if $N_0$ is the restriction of
$N$ to
$C_{\text{red}} \cong B$, it is still clear that $V_{(\bold o, 0)}[N]$ has
a filtration with successive quotients $L^{n-1}\otimes N_0, L^{n-2}\otimes
N_0,
\dots, L\otimes N_0$. We have
$$\ch(V_{\bold o, 0}[N])=\frac{1 - e^{nL}}{1-e^L}\cdot e^{N_0}.$$

\remark{Remark} (1) Note that, unless $L$ is a torsion line bundle, the
bundles
$V_{\bold o, 0}[N]$ are unstable with respect to every ample divisor.

\noindent (2) By contrast with (5.14), even if $\dim B =1$, we cannot always
arrange trivial determinant for $V_{\bold o, 0}[N]$.
\endremark
\medskip

If instead we try to construct $V_{\bold o, a}$ directly as a sequence of
global extensions on $Z$, we run into the following type of question.
Suppose for  simplicity that $n=2$ and that $a=0$. In this case we try to
find a bundle on $Z$ which restricts over every fiber $f$ of $Z$ to be the
nontrivial extension of $\scrO_f$ by $\scrO_f$, in other words to
$I_2$. We may as well try to write it as an extension of $\scrO_Z$ by the
pullback of a line bundle $N$ on $B$.  To do this we need a class
$H^1(\pi^*N)$ whose restriction to every fiber is non-trivial. That is to
say, we need an element in $H^1(\pi^*N)$ whose image under the natural map
$\psi$ in the Leray spectral sequence  (which is an exact sequence in this
case)
$$H^1(\pi _*\pi ^*N) \to H^1(\pi ^*N)
\buildrel\psi\over\longrightarrow H^0(R^1\pi _*\pi ^*N) \to H^2(\pi _*\pi
^*N)$$  is a nowhere zero section of $R^1\pi _*\pi ^*N$.  Of course, $\pi
_*\pi ^*N \cong N$ and $R^1\pi _*\pi ^*N \cong N\otimes R^1\pi _*\scrO_Z=
N\otimes L^{-1}$. Thus if there is to exist a nowhere vanishing  section of
$H^0(R^1\pi _*\pi ^*N)$, it must be the case that $N=L$. But we also need
the condition that the map $H^1(\pi ^*L) \to H^0(L\otimes L^{-1})=
H^0(\scrO_B)$ is  surjective. This is not immediately obvious    from the
spectral sequence since there is no reason for $H^2(B;L)$ to vanish.
Nevertheless, it follows from our construction of $V_{\bold o, 0}$ that the
map $\psi$ is onto in the case $N=L$. Finally, the set  of possible
extensions is a principal homogeneous space over
$H^1(B;L)$, which is identified with the kernel of the natural map
$\Pic (2\sigma) \to \Pic B$.

\ssection{5.7. Deformation to a reducible spectral cover.}

For every choice of a rank $n > \dim B$ and for all sections $A$ of
$\Cal P_{n-1}$ which correspond to a sufficiently ample line bundle, we have
constructed vector bundles $V_{A,a} = V_{A,a}(n)$. In this subsection, we
try to relate the $V_{A,a}(n)$ for various choices of $n$. To this end, 
let $\Cal H = \Cal P^{n-2} = \Pee (\scrO_B \oplus L^{-2} \oplus
\cdots \oplus L^{-n+1}) \subset {\Cal P}_{n-1}$. We begin by considering
what happens when the section $A$ lies in the subbundle
$\Cal H$, but is otherwise generic. To insure that there are actually
sections  of
$\Cal H$ as opposed to just rational sections, it is reasonable to assume
that $n \geq \dim B +2$. A section $A$ of
$\Cal H$ is given by a line bundle $M$ and by $n$ sections $\sigma_0, \dots,
\sigma_{n-1}$ of $M, M\otimes L^{-2}, \dots,  M\otimes L^{-(n-1)}$ which
have no common zeroes. If $M$ is sufficiently ample, the section $A=A_0$
will then move in a family
$A_t$ of sections of ${\Cal P}_{n-1}$, by choosing a nonzero section
$\sigma_n$ of
$M\otimes L^{-n}$ and considering the family defined by the sections
$A_t=(\sigma_0,
\dots, \sigma_{n-1}, t\sigma_n)$. Roughly speaking,
$V_{A_0,a}(n)$ is obtained from the bundle $V'$ of rank $n-1$ corresponding
to $A_0$, viewed as a section of
$\Cal P_{n-2}$. Along each fiber $f$ we add a trivial $\scrO_f$ factor to
the restriction of $V'$. This statement is correct as long as the
restriction of $V'$ to the fiber  does not itself contain an $\scrO_f$
factor, or more generally a summand of the form $I_d(\scrO_f)$ for some
$d\leq n$.  The simplest possibility would be that $V_{A_0,a}(n)$ is a
deformation of $V_{A,a}(n-1)\oplus \scrO_Z$, but a calculation with Chern
classes rules this out. Likewise, $V_{A_0,a}(n)$ is not a deformation of
$V_{A,a}(n-1)\oplus \pi^*N$ for any line bundle $N$ on $B$. Instead, we
shall see that $V_{A_0,a}(n)$ is a deformation of a suitable elementary
modification of $V_{A,a}(n-1)\oplus \pi^*L^a$. Finally, we shall use the
construction to check the Chern class calculations.
 
To make this construction, it is best to begin by working universally
again. We have the
$n$-to-$1$ map $\nu\: \Cal T \to {\Cal P}_{n-1}$. Inside
${\Cal P}_{n-1}$, there is the smooth divisor $\Cal H  = \Cal P_{n-2}$. Now
in
$\Cal T= \Cal T_{n-1}$ there is a smooth divisor $\Cal T'\cong \Cal
T_{n-2}$ defined by the diagram
$$\CD 0 @>>> \Cal E @>>> \pi ^*\pi _*\scrO_Z(n\sigma) @>>>
\scrO_Z(n\sigma) @>>> 0\\ @. @AAA @AAA @AAA @.\\ 0 @>>> \Cal E' @>>>
\pi ^*\pi _*\scrO_Z((n-1)\sigma) @>>> \scrO_Z((n-1)\sigma) @>>> 0.
\endCD$$ We take $\Cal T' = \Pee (\Cal E') \subset \Pee (\Cal E) =\Cal T$.
The restriction of $\nu$ to $\Cal T'$ defines the corresponding map $\Cal
T_{n-2} \to \Cal P_{n-2}$, and in particular
$\nu |\Cal T'$ has degree $n-1$.  Clearly, we have an equality of smooth
divisors in $\Cal T$:
$$\nu ^*\Cal H = \Cal T' + r^*\sigma.$$ The intersection $\Cal T'
\cap r^*\sigma$ is clearly the smooth divisor $\Cal P_{n-3} \subset
r^*\sigma \cong  \Cal P_{n-2}$; it lies over $\Cal P_{n-3}$. A local
calculation shows that $\Cal T'$ and $ r^*\sigma$ meet transversally at the
generic point of $\Cal P_{n-3}$ and thus everywhere. Note that
$\Cal T_1\cong Z$, $r\: \Cal T_1 \to Z$ is the identity, and the
intersection of
$\Cal T_1$ and $r^*\sigma$ in $\Cal T_2$ is $\sigma \subset \Cal T_1$. This
is compatible with the convention
$\bold o \cong \Cal P_0  \cong B$.

Let $\Cal D = \Cal T'\times _BZ$ and, as usual, let $F = r^*\sigma
\times _BZ$. Then $F$ is a smooth divisor and $\Cal D$ is smooth away from
the singularities of $\Cal T\times _BZ$. The divisors $\Cal D$ and $F$ meet
in a reduced  divisor $\Cal P_{n-3}\times _BZ$. We thus have an exact
sequence:
$$0 \to \scrO_{\Cal D + F} \to \scrO_{\Cal  D}\oplus \scrO_F \to
\scrO_{\Cal D\cap F} \to 0.$$ Tensoring the above exact sequence by the
sheaf
$\scrO_{\Cal T\times _BZ}(\Delta - \Sigma _A -aF)$, using the fact that
$\Delta \cap F = \Sigma _A\cap F$, gives a new  exact sequence
$$0 \to \scrO_{\Cal D + F}(\Delta - \Sigma _A -aF) \to
\scrO_{\Cal  D}(\Delta - \Sigma _A -aF)\oplus \scrO_F(-aF)\to
\scrO_{\Cal D\cap F}( -aF) \to 0.$$  (In a neighborhood of $F$, $\Delta$ is
Cartier, and so the above sequence is still exact.) Of course,
$F|F = -L|F$. Now apply 
$(\nu \times \Id)_*$ to the above exact sequence. To keep track of the
ranks, we shall write $\Cal U_a(n)$ when we want to denote the appropriate
vector bundle of rank
$n$, and similarly for $V_{A,a}(n)$. (However, in the notation,
$V_{A,a}(n-1)$ will be a general rank $(n-1)$-bundle but $V_{A,a}(n)$ will
be the special rank $n$ bundle corresponding to a reducible section. Of
course, this will not affect Chern class calculations.) We have:
$$0 \to \Cal U_a(n)|\Cal P_{n-2}\times _BZ \to \Cal U_a(n-1) \oplus
(L^a|\Cal P_{n-2}\times _BZ) \to  L^a|\Cal P_{n-3}\times _BZ \to 0.$$ Let
$A$ be a section of ${\Cal P}_{n-1}$ lying in $\Cal P_{n-2}$ and otherwise
general. Pulling back the above exact sequence via $A$, we get an exact
sequence relating the special rank $n$ bundle
$V_{A,a}(n)$ with a general rank $(n-1)$-bundle $V_{A,a}(n-1)$ obtained by
viewing
$A$ as a section of $\Cal P_{n-2}$:
$$0\to V_{A,a}(n) \to V_{A,a}(n-1)\oplus \pi ^*L^a \to (\pi ^*L^a)|D
\to 0,$$ where $D$ is the divisor in $Z$ corresponding to $\Cal P_{n-3}
\cap A$. In particular $D$ is pulled back from $B\cong A$. Thus we have
realized the special bundle $V_{A,a}(n)$ as an elementary modification of
$V_{A,a}(n-1)\oplus \pi ^*L^a$ along the divisor $D$.

To calculate the cohomology class of $D$, note that the class of $\Cal
P_{n-3}$ in
$\Cal P_{n-2}$ is given by $\zeta - (n-1)L$ (by applying (4.15) with
$n$ replaced by $n-1$), and  so the class of
$D$ is given by $p_*([A]\cdot (\zeta - (n-1)L))$. By (5.9),
$$[D] = \alpha - (n-1)L. \tag5.22$$ For $M$ sufficiently ample and $A$
general, $D$ is a smooth divisor, and we get $V_{A,a}(n)$ by an elementary
modification of the direct sum $V_{A,a}(n-1)\oplus \pi ^*L^a$ along
$D$. Here, of course, the surjection from $V_{A,a}(n-1)$ to $\pi ^*L^a|D$
arises because on every fiber $f$ over a point of $D$,
$V_{A,a}(n-1)$ has a trivial quotient $\scrO_f$.

Note that, assuming we are the range where the calculations are correct, we
obtain an inductive formula for $\ch V_{A,a}(n)$:
$$\ch V_{A,a}(n) = \ch  V_{A,a}(n-1) + \ch (L^a) -\ch (L^a|D).$$ Now from
the exact sequence
$$0 \to L^a\otimes \scrO_Z(-D) \to   L^a \to  L^a|D \to 0,$$ we see that
$\ch (L^a|D) = \ch (L^a) - \ch (L^a\otimes \scrO_Z(-D))$, and thus using
(5.22)
$$\ch V_{A,a}(n)=\ch V_{A,a}(n-1)+e^{(a+n-1)L+\alpha}.$$  Note that this is
consistent with the formula given in (5.10) for $\ch V_{A,a}$. 

This inductive picture must be modified for small values of $n$. For
example, in case $\dim B =3$,   a general section in $\Cal P_2$ degenerates
to a rational section of $\Cal P_1$ plus some exceptional fibers, and there
is a further problem in the passage from $\Cal P_1$ to $\Cal P_0 = \bold
o$. However, we will not discuss these matters further.

\ssection{5.8. Subsheaves of $V$ and reducible spectral covers.}

\proposition{5.23} Let $V$ be a rank $n$ bundle on $Z$ whose restriction to
every fiber is regular and  semistable with trivial determinant. Then the
spectral cover $C=C_A$ associated to $V$ is reduced and irreducible if and
only if  there is no subsheaf $V'\subset V$ whose restriction to the
generic fiber is a semistable bundle of degree zero and rank $r$ with  
$0<r<n$, if and only if  there is no quotient sheaf
$V''$ of $V$ which is torsion free and whose restriction to the generic
fiber is a semistable bundle of degree zero and rank $r$ with   $0<r<n$.
\endstatement

\proof Clearly, $V$ has a subsheaf $V'$ as in the statement of the
proposition if and only if it has a quotient sheaf $V''$ as described
above.  

If $C$ is not reduced and irreducible, then there is a proper closed
subvariety $C'\subset C$ which maps surjectively onto $B$ and is  finite of
degree $r,\ 0<r<n$ over $B$. We may assume that $C'$ is reduced. Let
$T'=T\times _BC'$ be the corresponding subscheme of
$T=T_A$. The surjection $\scrO_T \to \scrO_{T'}$ and the fact that $\nu
=\nu_A$ is finite leads to a surjection
$$V=(\nu\times \Id)_*\left[\scrO_T(\Delta - \Sigma _A)\otimes \rho^*N\right]
\twoheadrightarrow (\nu\times \Id)_*\left[\scrO_{T'}(\Delta - \Sigma
_A)\otimes
\rho^*N\right]=V''.$$ By construction, $V''$ is a torsion free sheaf on $Z$
of rank $r$ with $0< r < n$. Restrict to a generic smooth fiber $\pi
^{-1}(b)=E_b$ of $\pi$ such that the fiber of the projection
$C'\to B$ has $r$ distinct points $e_1, \dots, e_r \in E_b$ over $b$.  By
Lemma 5.6, the restriction of
$V''$ to $E_b$ is a direct sum of the $r$ line bundles
$\scrO_{E_b}(e_i-p_0)$, and in particular it is semistable (and in fact
regular).

Conversely, suppose that there is an exact sequence
$$0 \to V' \to V \to V'' \to 0,$$ where both $V'$ and $V''$ are nonzero
torsion free sheaves whose restrictions to a generic fiber are semistable.
Let $r'$ be the rank of $V'$ and $r''$ be the rank of $V''$. After
restricting to a nonempty Zariski open subset of $Z$, we may assume that
$V'$ and
$V''$ are locally free. Consider now the commutative diagram
$$\minCDarrowwidth{.2 in}
\CD 0 @>>> \pi^*\pi_*(V'\otimes \scrO_Z(\sigma)) @>>> \pi^*\pi_*(V\otimes
\scrO_Z(\sigma)) @>>> \pi^*\pi_*(V''\otimes \scrO_Z(\sigma)) @>>> 0\\ @.
@V{\Psi'}VV  @V{\Psi}VV   @V{\Psi''}VV @.\\ 0 @>>> V'\otimes
\scrO_Z(\sigma) @>>> V\otimes
\scrO_Z(\sigma) @>>> V''\otimes \scrO_Z(\sigma) @>>> 0.
\endCD$$ By definition, $C$ is the Cartier divisor which is the scheme of
zeroes of
$\det\Psi$. On the other hand, we clearly have $\det \Psi = \det \Psi'\cdot
\det
\Psi''$. If $C'$ is the scheme of zeroes of $\det \Psi '$, and $C''$ is the
scheme of zeroes of $\det \Psi ''$, then $C = C'+C''$ on a nonempty Zariski
open subset of $Z$. Furthermore, $C'$ maps to $B$ with degree $r'$ and
$C''$ maps to $B$ with degree $r''$, so that neither of $C', C''$ is
trivial. It follows that the restriction of $C$ to a nonempty Zariski open
subset of $Z$ is either nonreduced or reducible, and so the same is true
for $C$ as well.
\endproof

Finally, let us remark that if $V$ is merely assumed to be regular and
semistable on a generic fiber, so that $A(V)$ is just a rational section,
the above proof still goes through.

\section{6. Bundles which are not regular and semistable on every fiber.}

Let $\pi\: Z \to B$ be an elliptic fibration with $\dim B =d$, and  let
$E_b =\pi^{-1}(b)$. In this section, we consider some examples of bundles
$V$, such that
$\det V$ has trivial restriction to each fiber, which fail to be regular or
semistable on every fiber
$E_b$. From the general principles mentioned in the introduction, it is
reasonable to consider only those bundles whose restriction to the generic
fiber is semistable. We shall further assume here that the restriction to
the generic fiber is regular (this will exclude, for example, the tangent
bundle of an elliptic fibration whose base $B$ has dimension at least two).
Thus, we shall consider bundles $V$ such that, for a nonempty proper closed
subset $Y$ of $B$ and for all
$y\in Y$, either
$V|E_y$ is unstable or it is semistable but not regular. There is an
important difference between the case $\dim Y =d-1$ and $\dim Y<d-1$. In
the first case, $V$ is not determined by its restriction to $\pi
^{-1}(B-Y)$ and can be obtained via elementary modifications from a
``better" bundle (or reflexive sheaf). In this case, there is a lot of
freedom in creating such $V$ where $V|E_y$ is unstable along a
hypersurface. By contrast, it is more difficult to arrange that $V|E_y$ is
semistable but not regular along a hypersurface. If $\dim Y < d-1$, then,
since
$V$ is a vector bundle, it is determined by its restriction to $\pi
^{-1}(B-Y)$ and the behavior of
$V$ is much more tightly controlled by the rational section $A(V)$ of $\Cal
P_{n-1}$. Here the case where
$V|E_y$ is unstable for $y\in Y$ (as well as the case where $V$ is
reflexive but not locally free) corresponds to the case where $A(V)$ is
just a quasisection, i\.e\. where the projection $A(V) \to B$ has degree
one but is not an isomorphism. The case where $V|E_y$ is semistable but not
regular for $y\in Y$ corresponds to the case where there are singularities
in the spectral cover $C_A$, and $V$ is obtained by twisting by a line
bundle on $C_A|B-Y$ which does not extend to a line bundle on $B$. As will
be clear from the examples, a wide variety of behavior is possible, and we
shall not try to give an exhaustive discussion of all that can occur.

\ssection{6.1. Codimension one phenomena and elementary modifications.}

First we shall discuss the phenomena which occur in codimension one, and
which amount to generalized elementary modifications. As will be clear,
when we make the most general elementary modifications, we lose control in
codimension two on $B$. Thus for example many of the constructions lead to
reflexive sheaves which are not locally free. For this reason, we shall
concentrate to a certain extent on the case $\dim B = 1$, which will
suffice for the generic behavior in codimension one when $\dim B$ is
arbitrary.

The first very general lemma says that, locally, every possible bundle with
a given restriction to the generic fiber arises as an elementary
modification.

\lemma{6.1} Let $V$ be a vector bundle on $Z$ whose restriction to every
fiber $E_b$ is  semistable and whose restriction to the generic fiber is
regular. Suppose that $A(V)=A$ is the section of
$\Cal P_{n-1}$ corresponding to $V$. Let 
$$Y = \{\, b\in B: \text{ $V|E_b$ is not regular}\,\}.$$ Then $Y$ is a
Zariski closed subset of $B$. For every $y\in Y$, there exists a Zariski
neighborhood $\Omega$ of $y$ in $B$ and a morphism $\varphi\: V_{A,0}|\pi
^{-1}(\Omega)\to V|\pi ^{-1}(\Omega)$ which is an isomorphism over a
nonempty Zariski open subset of $\Omega$. Moreover,  we can choose a
$\varphi$ which extends to a homomorphism
$V_{A,0} \otimes \pi ^*M^{-1} \to V$, where $M$ is a sufficiently ample
line bundle on $B$.

More generally, suppose that $V$ is merely assumed to have regular
semistable restriction to the generic fiber, so that $V|E_b$ may be
unstable for some fibers. Then there exists a closed subset $X$ of $B$ of
codimension at least two such that the section $A(V)$ extends over $B-X$
and, with $Y$ as above, for every
$y\in Y-X$, there exists a Zariski neighborhood $\Omega$ of $y$ in $B$ and
a morphism $\varphi\: V_{A,0}|\pi ^{-1}(\Omega)\to V|\pi ^{-1}(\Omega)$
which is an isomorphism over a nonempty Zariski open subset of $\Omega$.
Finally, we can choose a $\varphi$ which extends to a homomorphism
$V_{A,0} \otimes \pi ^*M^{-1} \to V|B-X$, where $M$ is a sufficiently ample
line bundle on $B-X$.
\endstatement
\proof Let us first consider the case where the restriction of $V$ to every
fiber is semistable. In this case the section $A=A(V)$ is defined over all
of $B$. Consider the sheaf
$\pi_*Hom (V_{A,0}, V)$. On
$B-Y$, this sheaf is locally free of rank $n$. On a sufficiently small open
set
$\Omega$, we can thus find a section $\varphi$ of $\pi_*Hom (V_{A,0},
V)|\Omega$ which restricts to an isomorphism on a general fiber. Since this
is an open condition, the set of points $b\in \Omega$ such that $\varphi$
fails to be an isomorphism on
$E_b$ is a proper Zariski closed subset of $\Omega$, as claimed. Finally,
if $M$ is sufficiently ample, then $\pi_*Hom (V_{A,0}, V)\otimes M$ is
generated by its global sections. Choosing such a section which restricts
to an isomorphism from
$V_{A,0} \otimes \pi ^*M^{-1}|E_b$ to $V|E_b$ for a fiber $E_b$ defines a
map
$\varphi$ which extends to a homomorphism $V_{A,0} \otimes \pi ^*M^{-1}\to
V$, as claimed.

In case $V$ has unstable restriction to some fibers, the above proof goes
through as long as we are able to define the section $A(V)$. Now the
rational section of
$\Cal P_{n-1}$ defined by $V$ extends to a closed irreducible subvariety of
$\Cal P_{n-1}$, which we shall also denote by $A(V) =A$. The morphism
$p|A\: A\to B$ is birational, and thus over the complement of a codimension
two set $X$ in
$B$ it is an isomorphism. Thus $A$ is a well-defined section over $B-X$,
and so defines a bundle $V_{A,0}$ over $\pi ^{-1}(B-X)$. We may then apply
the first part of the proof.
\endproof

Let $V$ be a vector bundle on $Z$ whose restriction to the generic fiber
$E_b$ is  semistable.  Let 
$$Y = \{\, b\in B: \text{ $V|E_b$ is not semistable}\,\}.$$ Then $Y$ is a
Zariski closed subset of $B$. Suppose that $W=\pi^{-1}(Y)\subset Z$. We can
restrict $V$ to the elliptic fibration $W\to Y$. For simplicity, we shall
assume that $W$ is irreducible (otherwise we would need to work one
irreducible component at a time). By general theory, there exists a torsion
free sheaf
$\Cal S$ over $W$ and a surjection $V|W
\to \Cal S$, such that at a generic point $w$ of $W$, the map $V|E_w \to
\Cal S|E_w$ is the maximal destabilizing quotient of $V|E_w$. Let $i\: W\to
Z$ be the inclusion and let $V'$ be the kernel of the surjection $V \to
i_*\Cal S$. If
$W$ is a hypersurface in $Z$, i\.e\. if $Y$ is a hypersurface in $B$, then
$V'$ is a reflexive sheaf. However, if $W$ has codimension greater than
one, $V'$ fails to be reflexive, and in fact $(V')\ddual = V$.

For example, if $\dim B=1$, $W$ is a finite set of points. Choosing one
such point
$w$, we have that $V|E_w$ is unstable. Let $Q$ be the maximal destabilizing
quotient sheaf for $V|E_w$, and suppose that $\deg Q = e< 0$. Then $V'$
fits into an exact sequence
$$0 \to V' \to V \to i_*Q \to 0,$$ where $i$ is the inclusion of the fiber
$E_w$ in $Z$. Such elementary modifications of $V$ are {\sl allowable\/} in
the terminology of \cite{4},
\cite{5}. As opposed to the general construction of (6.1), allowable
elementary modifications are  canonical, subject to a choice of an
irreducible component of $W$. For the above allowable elementary
modification over an elliptic surface, we have
$$c_2(V') = c_2(V) + e< c_2(V).$$ Thus an allowable elementary modification
always decreases $c_2$.

\lemma{6.2} A sequence of allowable elementary modifications terminates.
The end result is a torsion free reflexive sheaf $V'$ such that the set 
$$\{\, b\in B: \text{ $V'|E_b$ is not semistable}\,\}$$ has codimension at
least two.
\endstatement
\proof We shall just write out the proof in the case $\dim B = 1$. In this
case, by (6.1),  we can fix a bundle $V_0 = V_{A,0}\otimes \pi^*M^{-1}$ for
some section
$A$, together with a morphism
$\varphi\: V_0 \to V$ which is an isomorphism over a general fiber. Thus
$\det
\varphi$ defines an effective Cartier divisor, not necessarily reduced,
supported on a union of fibers of $\pi$. Denote this divisor by $D$.
Clearly $D$ is the pullback of a divisor $\bold d$ on $B$, and thus has a
well-defined length $\ell$, namely the degree of $\bold d$. We claim that
every sequence of allowable elementary modifications has length at most
$\ell$. This is clearly true if $\ell =0$, since then $V_0\to V$ is an
isomorphism and every fiber of $V$ is already semistable. Since a sequence
of allowable elementary modifications will stop only when the restriction
of $V$ to every fiber is semistable, we will get the desired conclusion.

Let $V'$ be an allowable elementary modification of $V$ at the fiber $E_w$.
We claim that $\varphi$ factors through the map $V'\to V$. In this case, it
follows that $E_w$ is in the support of $D$. Thus, if $\varphi'\: V_0 \to
V'$ is the induced map, then $(\det \varphi ') = D-E_w$, which has length
$\ell -1$, and we will be done by induction on the  length
$\ell$. It suffices to prove that the induced map
$V_0 \to i_*Q$ is zero in the above notation. Equivalently, we must show
that the induced map
$V_0|E_w \to Q$ is zero. But $V_0|E_w$ is semistable and $\deg Q <0$, and
so we are done.
\endproof

As a corollary, we have the following Bogomolov type inequality:

\corollary{6.3} Let $V$ be a vector bundle on $Z$ such that the restriction
of $V$ to a generic fiber $E_b$ is regular and semistable. Suppose that
$\dim B = d$. Then, for every ample divisor $H$ on $B$, $c_2(V)\cdot
\pi^*H^{d-1} \geq 0$. Moreover, equality holds if and only if $V$ is
semistable in codimension one and the line bundle $M$ corresponding to the 
rational section $A(V)$ is a torsion line bundle. Finally, $M$ is a torsion
line bundle if and only if either the rational section $A(V) =\bold o$ or
$L$ is a torsion line bundle and $M$ is a power of $L$.
\endstatement
\proof We may assume that $H$ is very ample. By choosing a general curve
which is a complete intersection of $d-1$ divisors linearly equivalent to
$H$, we can further assume that $\dim B = 1$, and must show that $c_2(V)
\geq 0$. Since an allowable elementary modification strictly decreases
$c_2$, we can further assume that the restriction of $V$ to every fiber is
semistable. Choose a nonzero map
$V_0 \to V$, where $V_0$ is regular semistable on every fiber. Defining $Q$
by the exact sequence
$$0 \to V_0 \to V \to Q \to 0,$$
$Q$ is a torsion sheaf supported on some (possibly nonreduced) fibers whose
restriction to a $b\in B$ has a filtration by degree zero sheaves on $E_b$.
It then follows that $c_2(V) = c_2(V_0)$. Now if $M$ is the line bundle 
corresponding to the section $A(V)$ of
$\Cal P_{n-1}$, then by (5.15) $c_2(V_0) = \deg M$. On the other hand, at
least one of $M, M\otimes L^{-2}, \dots, M\otimes L^{-n}$ has a nonzero
section. Thus, for some $i=0, 2, \dots, n$, $\deg M \geq i\deg L$. Now
$\deg L \geq 0$, and
$\deg L = 0$ only if $L$ is a torsion line bundle. Thus, $\deg M \geq 0$,
and
$\deg M = 0$ only if $i=0$, in which case $M$ is trivial, or $L$ is torsion
and there is a nowhere vanishing section of $M\otimes L^{-i}$. In all cases
$M$ is a torsion line bundle and we have proved the statements of the
lemma. 
\endproof

\remark{Remark} (1) If $c_2(V)\cdot \pi^*H^{d-1} =0$ above, in other words
we have equality, it follows that the rational section $A(V)$ is actually a
section.
\smallskip
\noindent (2) If $A$ is a rational section and $A \neq \bold o$, we get
better inequalities along the lines of $$c_2(V)\cdot\pi ^*H^{d-1} \geq
2L\cdot H^{d-1},$$  since we must have nonzero sections of $M\otimes
L^{-i}$, $i=0,2, \dots, n$ for at least two values of $i$. If $A$ is a
section, then except for a small number of exceptional cases we will
actually have $c_2(V)\cdot\pi ^*H^{d-1} \geq (d+1)L\cdot H^{d-1}$.
\endremark
\medskip

The process of taking allowable elementary modifications is in a certain
sense reversible: we can begin with a bundle $V_0$ such that the
restriction of $V_0$ to every fiber is semistable and introduce instability
by making elementary modifications. Let us first consider the case where
$\dim B =1$. At the first stage, fixing a fiber $E_b$ and a stable sheaf
$Q$ on $E_b$ of positive degree, we seek a surjection $V_0|E_b \to Q$. To
analyze when such surjections exist is beyond the scope of this paper.
However, in case $V_0|E_b$ is regular and
$Q=W_k$, then we have seen in Section 3 that such a surjection always
exists; indeed, the set of all surjections is an open subset in $\Hom (V,
W_k)$ which has dimension $n$. Note however that while allowable elementary
modifications are canonical, their inverses are not. To be able to continue
to make elementary modifications along the same fiber, we would also have
to analyze when there exist surjections from
$V|E_b$ to $Q$, where $V$ is a rank $n$ bundle on $E_b$ of degree zero, $Q$
is a torsion free sheaf of rank $r<n$ on $E_b$, and $\mu(Q)$ is larger than
the maximum of
$\mu(\Cal S)$ as $\Cal S$ ranges over all proper torsion free subsheaves of
$V|E_b$. 

In case $\dim B > 1$, further complications can ensue in codimension two.
For example, suppose that $V_0$ has regular semistable restriction to every
fiber of
$\pi$. Let $D$ be a divisor in $B$ and let $W=\pi^{-1}(D)$, with
$\pi'=\pi|W$. Even though we can find a surjection $V|E_b \to W_k$ for
every  $b\in D$, we can only find a global surjection  $V|W \to \Cal
W_k\otimes (\pi')^*N$, for some line bundle $N$ on $D$, under special
circumstances. We can find a nonzero such map in general, but it will
vanish in general in codimension two, leading to a reflexive but not
locally free sheaf.

We turn next to the issue of bundles which are semistable on every fiber,
but which are not regular in codimension one. It turns out that we do not
have the freedom that we did before in introducing instability on a fiber;
there is a condition on the spectral cover in order to be able to make a
bundle not be regular. (See \cite{3}, \cite{6} for the rank two case.) We
shall just state the result in the case where $\dim B =1$. The result is
that, if the spectral cover is smooth, it is not possible to create a
non-regular but semistable bundle over any fiber.

\proposition{6.4} Let $\dim B =1$ and let $V$ be a vector bundle over $Z$
whose restriction to every fiber is semistable and whose restriction to the
generic fiber is regular. Let $A=A(V)$ be the corresponding section and
$C=C_A$ be the spectral cover. If $b\in B$ and $C$ is smooth at all points
lying over $b\in B$, then $V|E_b$ is regular.
\endstatement
\proof Using Lemma 6.1, write $V$ as a generalized elementary modification
$$0 \to V_0 \to V \to Q \to 0,$$ where $V_0$ is regular and semistable on
every fiber, and $Q$ is a torsion sheaf supported on fibers. Looking just
at the part of $Q$ which is supported on $E_b$, this sheaf (as a sheaf on
$Z$) has a filtration whose successive quotients are direct images of
torsion free rank one sheaves of degree zero on $E_b$. By induction on the
length of $Q$, as in the proof of Lemma 6.2, it will suffice to show the
following: if $V_0$ has regular semistable restriction to $E_b$, if $i\:
E_b \to Z$ is the inclusion, and if $C$ is smooth over all points lying over
$b$, then for every exact sequence
$$0 \to V_0 \to V \to i_*\lambda \to 0,$$ where $\lambda$ is a rank one
torsion free sheaf on $E_b$ of degree zero,
$V|E_b$ is again regular. After shrinking $B$, we can assume that $V_0$ is
regular everywhere and that $V_0 \to V$ is an isomorphism away from $b$. 

It will suffice to show that
$V\spcheck|E_b$ is regular. There is the dual exact sequence
$$0 \to V\spcheck \to V_0\spcheck \to i_*\lambda^{-1} \to 0.$$ By
assumption, $\dim \Hom (V_0\spcheck, i_*\lambda^{-1}) = \dim \Hom
(V_0\spcheck|E_b,
\lambda^{-1}) = 1$. Thus there is a unique possible elementary
modification. On the other hand, there is a unique point $b'\in C$ lying
above $b$ and corresponding to the surjection $V_0\spcheck|E_b \to
\lambda^{-1}$. Since by assumption $b'$ is a smooth point of $C$, it is a
Cartier divisor, and the ideal sheaf of $b'$ is the line bundle
$\scrO_C(-b')$. Now we know that $V_0\spcheck$ is of the form $(\nu\times
\Id)_*\left[\scrO_C(\Delta - \Sigma)\otimes \rho^*N\right]=V_{A,0}[N]$ for
a line bundle
$N$ on
$C$. Let $i'$ be the inclusion of the fiber over $b'$ (which is just $E_b$)
into
$T$. Applying $(\nu\times \Id)_*$ to the exact sequence
$$0 \to \scrO_C(\Delta - \Sigma)\otimes \rho^*(N\otimes \scrO_C(-b')) \to
\scrO_C(\Delta - \Sigma)\otimes \rho^*N \to (i')_*\lambda^{-1} \to 0,$$ we
get an exact sequence
$$0 \to V_{A,0}[N\otimes \scrO_C(-b')] \to V_0\spcheck \to i_*\lambda^{-1}
\to 0.$$ By the uniqueness of the map $V_0\spcheck  \to i_*\lambda^{-1}$,
it then follows that 
$$V \spcheck = V_{A,0}[N\otimes \scrO_C(-b')]$$  and in particular it is
regular. Thus the same is true for $V$.
\endproof

\remark{Remark} (1) Of course, Proposition 6.4 gives conditions in case
$\dim B > 1$ as well.

\noindent (2) The condition that $C$ is singular at the point corresponding
to $b$ and $V_0 \to \lambda$ is not a sufficient condition for there to
exist an elementary modification such that the result is not regular over
$E_b$. 
\endremark

\ssection{6.2. The tangent bundle of an elliptic surface.}

As an example of the preceding discussion, we analyze the tangent bundle of
an elliptic surface. Let $\pi\: Z \to B$ be an elliptic surface over the
smooth curve
$B$, with $g(B) = g$. We suppose that $Z$ is generic in the following
sense: $Z$ is smooth, the line bundle $L$ has positive degree $d$, so that
the Euler characteristic of $Z$ is $12d$, all the singular fibers of $\pi$
are nodal curves (and thus there are $12d$ such curves), and the
$j$-function $B\to \Pee^1$ has generic branching behavior in the sense of
\cite{6, p\. 63}. The assumption of generic branching behavior implies that
the Kodaira-Spencer map associated to the deformation
$Z$ of the fibers of $\pi$ is an isomorphism at the curves with $j=0, 1728,
\infty$ and that the Kodaira-Spencer map vanishes simply where it fails to
be an isomorphism. By the Riemann-Hurwitz formula, if $b$ is equal to the
number of points where the Kodaira-Spencer map is not an isomorphism, then
$b= 10d+2g-2$
\cite{6, p\. 68}.

Quite generally, we have the following lemma:

\lemma{6.5} Let $\pi\: Z\to B$ be a smooth elliptic surface, and suppose
that $V$ is a a vector bundle on  $Z$ whose restriction to a general fiber
is $I_2$. Then there is an exact sequence
$$0 \to \pi^*M_1\to V \to \pi^*M_2\otimes I_X \to 0,$$ where $M_1$ and
$M_2$ are line bundles on $B$ and $X$ is a zero-dimensional local complete
intersection subscheme of $Z$. Here $\det V = \pi ^*(M_1\otimes M_2)$ and
$c_2(V) =\ell (X)$.
\endstatement
\proof By assumption, $\pi_*V =M_1$ is a rank one torsion free sheaf on
$B$, and thus it is a line bundle. We have the natural map $\psi\:
\pi^*\pi_*V =\pi^*M_1\to V$. If this map were to vanish along a divisor,
the divisor would have to be a union of fibers. But this is impossible
since the induced map
$$\pi _*\pi^*M_1 = M_1 \to \pi_*V = M_1$$ is the identity. Thus $\psi$ only
vanishes in codimension two. The remaining statements are clear.
\endproof

Of course, in the case of the tangent bundle, we can identify this sequence
precisely as follows:

\lemma{6.6} With $\pi\: Z\to B$ a smooth elliptic surface as before, there
is an exact sequence
$$0 \to T_{Z/B} \to T_Z \to \pi ^*T_B \otimes I_X \to 0.$$ Here
$T_{Z/B}=L^{-1}$ is the sheaf of relative tangent vectors and $I_X$ is the
ideal sheaf of the $12d$ singular points of the singular fibers.
\endstatement
\proof Begin with the natural map $T_Z \to \pi ^*T_B$. This map is
surjective except at a singular point of a singular fiber, where it has the
local form 
$$h_1\frac{\partial}{\partial z_1} + h_2\frac{\partial}{\partial z_2}
\mapsto (z_2h_1+ z_1h_2)\frac{\partial}{\partial t}.$$ Thus the image of
$T_Z$ in $\pi ^*T_B$ is exactly $\pi ^*T_B \otimes I_X$. The kernel of the
map $T_Z \to \pi ^*T_B$ is by definition $T_{Z/B}$, which can be checked
directly to be a line bundle in local coordinates. Moreover, $T_{Z/B}$ is
dual to
$K_{Z/B}=L$, and thus $T_{Z/B} =L^{-1}$.
\endproof

\corollary{6.7} If $E$ is a singular fiber of $\pi$, there is an exact
sequence
$$0 \to n_*\scrO_{\tilde E} \to T_Z|E \to \frak m_x \to 0,$$ where $n\:
\tilde E \to E$ is the normalization and $x$ is the singular point of
$E$. In particular $T_Z|E$ is unstable. If $E$ is a smooth fiber where the
Kodaira-Spencer map is zero, then
$T_Z|E \cong \scrO_E\oplus \scrO_E$. For all other fibers $E$, $T_Z|E\cong
I_2$.
\endstatement
\proof If $E$ is a singular fiber, then by restriction we have a surjection
$T_Z|E
\to \frak m_x$. The kernel must be a non-locally free rank one torsion free
sheaf of degree one, and thus it is isomorphic to $ n_*\scrO_{\tilde E}$.
For a smooth fiber $E$, restricting the tangent bundle sequence to $E$
gives an exact sequence
$$0\to \scrO_E \to T_Z|E \to \scrO_E \to 0,$$ such that the coboundary map 
$$\theta\: H^0(\scrO_E) = H^0(N_{E/Z})\to H^1(\scrO_E) = H^1(T_E)$$ is the
Kodaira-Spencer map. This map is nonzero, then, if and only if $T_Z|E\cong
I_2$, and it is zero if and only if $T_Z|E\cong \scrO_E\oplus \scrO_E$.
\endproof

To go from $T_Z$ to one of our standard bundles, begin by making the
allowable elementary modifications along the singular fibers, by taking
$V'$ to be the kernel of the induced map $T_Z \to \bigoplus _x(i_x)_*\frak
m_x$. Here the sum is over the the singular points, i\.e\. the $x\in X$,
and $i_x$ is the inclusion of the singular fiber containing
$x$ in $Z$. Note that $c_2(V') =0$, so that no further allowable elementary
modifications are possible, and the restriction of $V'$ to every fiber is
semistable. Let $F$ be the union of the singular fibers. Thus as a divisor 
on $Z$, $F =\pi^*\bold f$, where $\bold f$ is a divisor on $B$ of degree
$12d$ which is a section of $L^{12}$. If
$I_F$ is the ideal of
$F$, then there is an inclusion $I_F\subset I_X$ and thus an inclusion $\pi
^*T_B
\otimes I_F \subset \pi ^*T_B \otimes I_X$. Clearly $V'$ is the result of
pulling back the extension $T_Z$ of $\pi ^*T_B \otimes I_X$ by $L^{-1}$ via
the inclusion 
$\pi ^*T_B \otimes I_F \subset \pi ^*T_B \otimes I_X$. Thus there is an
exact sequence
$$0 \to \pi^*L^{-1} \to V' \to \pi^*(T_B \otimes \scrO_B(-\bold f)) \to 0.$$
Taking the map
$$\gather
\Ext^1(\pi ^*T_B \otimes I_F, L^{-1}) =H^1(\pi ^*T_B^{-1}
\otimes\scrO_Z(F)\otimes L^{-1}) \\
\to H^0(R^1\pi_*(\pi ^*T_B^{-1}
\otimes\scrO_Z(F)\otimes L^{-1})) = H^0(B; K_B\otimes L^{-2}\otimes
\scrO_B(\bold f)),
\endgather$$ and using the fact that $\scrO_B(\bold f)\cong L^{12}$, we see
that the extension restricts to the trivial extension over a section of
$K_B\otimes L^{10}$, and thus at $10d+2g-2$ points, confirming the
numerology above. Note that the passage from $T_Z$ to $V'$ was canonical. 

Next we want to go from $V'$ to a bundle $V_0$ which is regular semistable
on every fiber, and thus is isomorphic to $I_2$ on every fiber. We claim
that a further elementary modification of $V'$ will give us back a bundle
which restricts to $I_2$ on every fiber. Quite generally, suppose that $V'$
is given as an extension 
$$0 \to \pi^*L_1\to V' \to \pi^*L_2 \to 0,$$ where the image of the
extension class in $H^0(B; L_2^{-1}\otimes L_1\otimes L^{-1})$ vanishes
simply at $k$ points $x_1, \dots, x_k$. After twisting $V'$ by the line
bundle $\pi^*L_2^{-1}$, we may assume that
$L_2$ is trivial. Thus in the case where we began with the tangent bundle,
and after relabeling $V'$, we wind up with a bundle $V'$ which fits into an
exact sequence
$$0 \to \pi^*L_1\to V' \to \scrO_Z \to 0,$$ where $L_1 = K_B\otimes
L^{11}$. The extension class for $V'$ defines an element of $H^1(Z;
\pi^*L_1)$. Via the Leray spectral sequence, there is a homomorphism from
$H^1(Z; \pi^*L_1)$ to $H^0(B; R^1\pi_*\scrO_Z\otimes L_1) = H^0(B;
L_1\otimes L^{-1})$. Thus there is a section of
$L_1\otimes L^{-1}$, well-defined up to a nonzero scalar, and it  defines a
homomorphism
$\pi^*L \to \pi^*L_1$ and thus a homomorphism $H^1(\pi^*L) \to
H^1(\pi^*L_1)$. Consider the commutative diagram
$$\CD H^1(B; L) @>>> H^1(\pi^*L) @>>> H^0(R^1\pi_*\pi^*L) = H^0(\scrO_B)\\
@VVV @VVV @VVV\\ H^1(B; L_1) @>>> H^1(\pi^*L_1) @>>> H^0(R^1\pi_*\pi^*L_1)
= H^0(L_1\otimes L^{-1}).
\endCD$$ The induced map $H^0(\scrO_B) \to H^0(L_1\otimes L^{-1})$ is just
the given section of $L_1\otimes L^{-1}$. We have seen in \S 5.6 that there
is a class $\xi_0
\in  H^1(\pi^*L)$ mapping to $1\in H^0(\scrO_B)$. Since the map $H^1(B; L)
\to H^1(B; L_1)$ is surjective, we can modify $\xi_0$ by an element in
$H^1(B; L)$ so that its image in $H^1(\pi^*L_1)$ is the same as the
extension class for $V'$, and the resulting element $\xi$ of $H^1(\pi^*L)$
is unique up to adding  an element of the kernel of the map $H^1(B; L) \to
H^1(B; L_1)$. Let $V_0$ be the extension of $\scrO_Z$ by $\pi^*L$
corresponding to $\xi$. Thus $V_0$ is some bundle of the form $V_{\bold o,
0}[N]$. There is an induced map of extensions
$$\CD 0 @>>> \pi^*L_1 @>>> V' @>>> \scrO_Z @>>> 0\\ @. @| @VVV @VVV @.\\ 0
@>>> \pi^*L_1 @>>> V_0\otimes \pi^*(L_1\otimes L^{-1}) @>>> \pi^*(L_1\otimes
L^{-1}) @>>> 0.
\endCD$$  Thus there is an exact sequence
$$0 \to V' \to V_0\otimes \pi^*(L_1\otimes L^{-1}) \to \bigoplus
_i\scrO_{E_{x_i}} \to 0,$$  and we have realized the tangent bundle as
obtained from $V_0$ by elementary modification and twisting. 

Of course, we can construct many other bundles this way, starting from
$V_0$, not just the tangent bundle. Begin with $V_0$ which has restriction
$I_2$ to every fiber. Normalize so that there is an exact sequence
$$0 \to \pi^*L \to V_0 \to \scrO_Z \to 0$$ as in \S 5.6. Here $L
=\pi_*V_0$. The bundle
$\pi^*L$ is destabilizing. Choose
$r$ fibers
$E_{x_i}$ lying over $x_i\in B$, where we make elementary modifications by
taking the unique quotient
$\scrO_{E_{x_i}}$ of $V_0|E_{x_i}$. The result is a new bundle $V'$. The
subbundle
$\pi^*L$ still maps into $V'$, in fact we continue to have $L=\pi_*V'$, and
the quotient is
$\pi^*\scrO_B(-\bold r)$, where
$\bold r$ is the divisor $\sum _ix_i$ of degree $r$ on $B$. The bundle $V'$
is the pullback of the extension $V_0$ by the morphism $\pi^*\scrO_B(-\bold
r) \to
\pi^*\scrO_B$. In particular, by reversing the arguments above, we see that
the restriction of the extension to
$E_{x_i}$ becomes split. Thus $V'|E_{x_i} \cong \scrO_{E_{x_i}} \oplus
\scrO_{E_{x_i}}$ and the restriction of $V'$ to all other fibers is $I_2$.
Note that $\pi^*L$ continues to destabilize
$V'$. 

Choose $s$ fibers lying over points $y_j\in B$ distinct from the $x_i$, and
let $\bold s$ be the divisor $\sum _jy_j$. Choose rank one torsion free
sheaves
$\mu_j$ on
$E_{y_j}$ of degree
$d_j > 0$ and surjections from $I_2$ to $\mu _j$. (Such surjections always
exist.) Take the bundle $V$ defined to be the kernel of the given
surjection $V' \to
\bigoplus _j\mu _j$. Now $\det V = (d- r-s)f$ and $c_2(V) =\sum _j\deg \mu
_j$. The bundle
$\pi^*L$ no longer maps into $V$, since the composed morphism
$\pi^*L|E_{y_j} \to
\mu_j$ is nontrivial for every $j$. In fact, $\pi^*(L\otimes \scrO_B(-\bold
s))$ maps to $V$, and $\pi_*V = L\otimes \scrO_B(-\bold s)$. Note that this
subbundle fails to be destabilizing exactly when $2(d-s) < d-r-s$, or
equivalently $d+r< s$. In this case, for a suitable ample divisor $H$ as
defined in \cite{6}, $V$ is
$H$-stable.

\ssection{6.3. Quasisections and unstable fibers.}

For the rest of this section, we shall assume that $V$ is regular and
semistable in codimension one and consider the phenomena that arise in
higher codimension. Over a Zariski open subset of $B$, we have defined
$A(V)$, and it extends to a subvariety of $\Cal P_{n-1}$ mapping
birationally to $B$, in other words to a {\sl quasisection\/} of $\Cal
P_{n-1}$. Of course, if the restriction of $V$ to every fiber is
semistable, then $A(V)$ is a section. 

\remark{Question} Suppose that $V$ is a vector bundle over $Z$ and that
there exists a closed subset $Y$ of $B$ of codimension at least two such
that, for all
$b\notin Y$, $V|E_b$ is semistable. Suppose further that $A(V)$ is actually
a section. Does it then follow that  $V|E_b$ is semistable for all $b\in B$?
\endremark
\medskip

For the remainder of this subsection, we shall assume that $A(V)$ is an
honest quasisection, in other words that the morphism $A(V) \to B$ is not an
isomorphism, and see what kind of behavior is forced on
$V$.  For example, if $n \leq \dim B$, then with a few trivial exceptions
there are no honest sections of $\Cal P_{n-1}$ and we are forced to consider
quasisections. We will analyze the case where $\dim B =2$ and see that two
kinds of behavior are possible: either
$V$ has unstable restriction to some fibers or $V$ fails to be locally
free  at finitely many points
$Z$. For example, suppose that $1\leq d\leq n-1$ and consider
$V_{A,1-d}$, defined over the complement of a set of codimension $2$ in
$B$. Then as we have seen in
\S 5.2, $V_{A, 1-d}$ is given as an extension of $\Cal W_{n-d}\otimes
\pi^*(M^{-1}
\otimes L)$ by $\Cal W_d\spcheck$. This extension extends over $B$, but it
induces the split extension of $W_{n-d}$ by $ W_d\spcheck$ wherever the
section of $\Cal V_n\otimes M$ vanishes.

Assume that $\dim B =2$ and let $s$ be a  section of $\Cal V_n
\otimes M$ which vanishes simply at finitely many points, but which is
otherwise generic. The corresponding quasisection $A=A(V)$ will contain a
line inside the full fiber of
$\Cal P_{n-1}$ at these points, which is a $\Pee^{n-1}$, and will simply be
the blowup of $B$ over the corresponding points. Pulling back the
$\Pee^{n-1}$-bundle
$\Cal P_{n-1}$ by the morphism $A\to B$, we get an honest section over $A$.
Let
$\tilde Z = Z\times _BA$. Clearly $\tilde Z$ is the blowup of $Z$ along the
fibers over the exceptional points of $B$, and the exceptional divisors of
$\tilde Z \to Z$ are of the form $\Pee^1\times E_b$, where the $\Pee^1$ is
linearly embedded in the $\Pee^{n-1}$ fiber. The section
$A$ of
$\Cal P_{n-1}\times _BA$ defines  a vector bundle
$\tilde \Cal U_a \to \tilde Z$ for every $a\in \Zee$. To decide what
happens over the exceptional points of $B$, we need to understand the
restriction of $\tilde
\Cal U_a$ to the exceptional fibers $\Pee^1\times E_b$. Of course, this is
just the restriction of the universal bundle $U_a$ defined over
$\Pee^{n-1}\times E_b$ to the subvariety $\Pee^1\times E_b$. Thus we need
to know the restriction of
$U_a$ to  $\Pee^1\times \{e\}$. We shall be able to find this restriction
in case 
$-(n-2) \leq a \leq 1$, but for arbitrary $a$ we shall further need to
assume that the
$\Pee^1$ is a generic line in $\Pee^{n-1}$.

\proposition{6.8} Let $E$ be a smooth elliptic curve and let $e\in E$.
Suppose that
$-(n-2) \leq a \leq 1$. Then
$$U_a|\Pee ^{n-1}\times \{e\} \cong \cases \scrO_{\Pee ^{n-1}}^{1-a}\oplus
\scrO_{\Pee ^{n-1}}(-1)^{n-1+a}, &\text{if $a
\neq 1$ or $e\neq p_0$;}\\
\scrO_{\Pee ^{n-1}} \oplus  \Omega^1_{\Pee^{n-1}},  &\text{if $a = 1$ and
$e= p_0$,}
\endcases$$ where $\Omega ^1_{\Pee^{n-1}}$ is the cotangent bundle of
$\Pee^{n-1}$.
\endstatement
\proof Let $i_e$ be the inclusion of $\Pee ^{n-1}$ in $\Pee ^{n-1}\times E$
via the slice $\Pee ^{n-1}\times \{e\}$. Then
$$i_e^*(\nu \times \Id)_*(\Delta - G-aF) = \nu _*\scrO_T(F_e - aF_{p_0}) =
\nu _*r^*\scrO_E(e-ap_0).$$ Set $d=1-a$, then $U_a =\bold U(d) \otimes
\pi_1^*\scrO_{\Pee^{n-1}}(-1)$. Thus, for $e\in E$, the restriction of
$\pi_2^*W_d\spcheck$ to $\Pee^{n-1}\times \{e\}$ is trivial, and similarly
for
$\pi_2^*W_{n-d}$, and the defining  exact sequence
$$0 @>>> \pi _2^*W_d\spcheck @>>> U_a  @>>> \pi_2^*W_{n-d} \otimes
\pi_1^*\scrO_{\Pee^{n-1}}(-1) @>>> 0$$ restricts to the exact sequence
$$ 0 @>>> \scrO_{\Pee^{n-1}}^d @>>> U_a |\Pee^{n-1}\times \{e\} @>>>
\scrO_{\Pee^{n-1}}(-1)^{n-d} @>>> 0.$$ Since $\Ext
^1(\scrO_{\Pee^{n-1}}(-1)^{n-d},\scrO_{\Pee^{n-1}}^d) =
H^1(\scrO_{\Pee^{n-1}}(1) )^{d(n-d)} =0$, this extension splits and we see
that 
$$\align U_a |\Pee^{n-1}\times \{e\} &\cong \scrO_{\Pee^{n-1}}^d \oplus
\scrO_{\Pee^{n-1}}(-1)^{n-d}\\ &\cong \scrO_{\Pee^{n-1}}^{1-a} \oplus
\scrO_{\Pee^{n-1}}(-1)^{n+a-1}.
\endalign$$ Now suppose that $a=1$. In this case $U_1 |\Pee^{n-1}\times
\{e\} = \nu _*r^*\scrO_E(e-p_0)$, and thus $h^0(U_1 |\Pee^{n-1}\times
\{e\})$ is zero if
$e\neq p_0$ and one if $e=p_0$. We have the elementary modification
$$0 \to U_1|\Pee^{n-1}\times \{e\} \to U_0 |\Pee^{n-1}\times \{e\} \to
\scrO_H \to 0,$$ where $H$ is a hyperplane in $\Pee^{n-1}$. Thus we may
write
$$0 \to U_1 |\Pee^{n-1}\times \{e\}\to  \scrO_{\Pee^{n-1}} \oplus
\scrO_{\Pee^{n-1}}(-1)^{n-1}\to \scrO_H \to 0.$$ Clearly $h^0(U_1
|\Pee^{n-1}\times \{e\})=0$ if and only if the induced map
$\scrO_{\Pee^{n-1}} \to \scrO_H$ is nonzero, or equivalently onto. In this
case, we can choose a summand $\scrO_{\Pee^{n-1}}$ of $\scrO_{\Pee^{n-1}}
\oplus
\scrO_{\Pee^{n-1}}(-1)^{n-1}$ such that the map $\scrO_{\Pee^{n-1}} \oplus
\scrO_{\Pee^{n-1}}(-1)^{n-1} \to \scrO_H$ is zero on the factor
$\scrO_{\Pee^{n-1}}(-1)^{n-1}$ and is the obvious map on the first factor.
Thus the kernel is $\scrO_{\Pee^{n-1}}(-1)^n$. 

In the remaining case, corresponding to $e=p_0$ and $a=1$, the map 
$\scrO_{\Pee^{n-1}} \oplus
\scrO_{\Pee^{n-1}}(-1)^{n-1} \to \scrO_H$ is zero on the first factor. Now
$H\cong
\Pee^{n-2}$, and modulo automorphisms of $\scrO_{\Pee^{n-1}}(-1)^{n-1}$
there is a unique surjection $\scrO_{\Pee^{n-1}}(-1)^{n-1} \to \scrO_H$. We
must therefore identify the kernel of this surjection with $\Omega
^1_{\Pee^{n-1}}$. Begin with the Euler sequence
$$0 \to \Omega ^1_{\Pee^{n-1}} \to \bigoplus _{i=1}^n\scrO_{\Pee^{n-1}}(-1)
\to
\scrO_{\Pee^{n-1}} \to 0.$$ After a change of basis in the direct sum, we
can assume that the right hand map restricted to the
$n^{\text{th}}$ factor vanishes along $H$. Thus there is an induced
surjection 
$$\bigoplus _{i=1}^{n-1}\scrO_{\Pee^{n-1}}(-1) \to
\scrO_{\Pee^{n-1}}/\scrO_{\Pee^{n-1}}(-1) =\scrO_H$$ whose kernel is
$\Omega ^1_{\Pee^{n-1}}$, as claimed.
\endproof

We remark that, in case $\dim B$ is arbitrary, $a=1$ and $A$ is a
quasisection corresponding to a simple blowup of $B$, then  one can show
directly from (6.8) that $V_{A,1}$ does not extend to a vector bundle over
$Z$.

When we are not in the range $-(n-2) \leq a \leq 1$, we do not identify
explicitly the bundle $U_a|\Pee^{n-1}\times \{e\}$, except in case $n=2$.
However, the next result identifies its restriction to a generic line.

\proposition{6.9} Let $\ell \cong \Pee ^1$ be a line in $\Pee ^{n-1}$, and
suppose that $\ell$ is not contained in any of the one-dimensional family
of hyperplanes $H_e$. Write $a = a' + nk$, where $-(n-2) \leq a' \leq 1$.
Then
$$U_a|\ell\times E \cong (U_{a'} \otimes
\pi_1^*\scrO_{\Pee^1}(-k))|\Pee^1\times E.$$ In particular
$$U_a|\ell\times \{e\} \cong \cases \scrO_{\Pee ^1}(-k)^{1-a'}\oplus
\scrO_{\Pee ^1}(-k-1)^{n-1+a'}, &\text{if $a'
\neq 1$ or $e\neq p_0$;}\\
\scrO_{\Pee ^1}^{n-2}(-k-1) \oplus \scrO_{\Pee ^1}(-k) \oplus
\scrO_{\Pee ^1}(-k-2),  &\text{if $a' = 1$ and $e= p_0$.}
\endcases$$
\endstatement 
\proof Let $C$ be the preimage of
$\ell$ in $T$. If  $\ell$ is not contained in any of the hyperplanes
$H_e$, then it will meet each $H_e$ in exactly one point. Thus the map
$r|C\: C \to E$ has degree one, and $F_{p_0}\cdot C = p_0$. We claim that,
under the morphism $\nu \: E \to \Pee ^1$, $\scrO_{\Pee ^1}(1)$ pulls back
to
$\scrO_E(np_0)=nF_{p_0}|C$. To see this, let
$\nu^*\scrO_{\Pee^{n-1}}(1)=\zeta\in \Pic T$. Then the class of $\nu^*\ell$
lies in
$\zeta ^{n-1}$. Now $T= \Pee\Cal E$ with $c_1(\Cal E) = -np_0$. Thus, in
$A^{n-1}(T)$,
$$\zeta ^{n-1} =r^*(np_0)\cdot \zeta^{n-2}.$$ Hence $\zeta |C = r^*(np_0)|C
= nF_{p_0}|C$.

Write $a = a'+nk$ with $-(n-2)\leq a' \leq 1$. Then
$$\align &(\nu \times \Id)_*\scrO_{C\times E}(\Delta -G -aF_{p_0}) =\\ &=
(\nu \times
\Id)_*\left(\scrO_{C\times E}(\Delta -G -a'F_{p_0})
\otimes
\pi _1^*\scrO_E(-nkp_0)\right)\\ &=(\nu \times \Id)_* \left(\scrO_{C\times
E}(\Delta -G -a'F_{p_0})\otimes (\nu \times \Id) ^*\scrO_{\Pee
^1}(-k)\right)\\ &=(\nu \times \Id)_*\scrO_{C\times E}(\Delta -G -a'F_{p_0})
\otimes
\pi_1^*\scrO_{\Pee ^1}(-k),
\endalign$$  proving the first claim. The second statement follows from the
special case
$-(n-2)\leq a \leq 1$ proved in (6.8), and the well-known fact (which
follows from the conormal sequence) that
$\Omega ^1_{\Pee^{n-1}} |\Pee^1 \cong \scrO_{\Pee^1}(-1)^{n-2} \oplus
\scrO_{\Pee^1}(-2)$.
\endproof

Now we can analyze what happens to $V_{A,a}$ when $\dim B =2$ and $A$ is a
quasisection, under a slight genericity condition on $A$, generalizing the
case (for $\dim B$ arbitrary) where $-(n-2) \leq a \leq 0$:

\theorem{6.10} Suppose that $\dim B =2$. Let $A$ be a quasisection of $\Cal
P_{n-1}$, and suppose that
$a\not \equiv 1 \bmod n$. Suppose that $A$ is smooth and is the blowup of
$B$ at a finite number of points $b_1, \dots, b_r$, and that the image of
the exceptional $\Pee^1$ is a generic line in the fiber $\Pee^{n-1}$ as in
\rom{(6.9)}, in other words it is not contained in one of the hyperplanes
$H_e$. Then the rank $n$ bundle
$V_{A,a}$, which is defined on $Z-\bigcup _iE_{b_i}$, extends to a vector
bundle over
$Z$, which we continue to denote by
$V_{A,a}$. The restriction of $V_{A,a}$ to a fiber $E_{b_i}$ is the
unstable bundle
$W_d\spcheck
\oplus W_{n-d}$, where $a = a'+nk$ with $-(n-2)\leq a' \leq 1$, and
$d=1-a'$.
\endstatement
\proof By assumption, $A$ is the blowup of $B$ at a finite number of points
$b_1, \dots, b_r$, where the quasisection $A$ contains a
$\Pee^1$ lying in the $\Pee^{n-1}$-fiber of $p\: \Cal P_{n-1} \to B$. As we
have defined earlier, let $\tilde Z = Z\times _BA$, so that $\tilde Z$ is a
blowup of
$Z$ at the fibers $E_{b_i}$. Let $D_i \cong \Pee^1\times E_{b_i}$ be the
exceptional divisor of the blowup $q\: \tilde Z \to Z$ over $E_{b_i}$. 
There is a section of
$\tilde B\to A$ corresponding to the inclusion of $A$ in
$\Cal P_{n-1}$, and hence by pulling back $\Cal U_a$ there is a bundle
corresponding to $A$, which we shall denote by $\tilde V$. Using (6.8) and
(6.9), the restiction of $\tilde V$ to an exceptional divisor $D_i
=\Pee^1\times E_{b_i}$, which is the same as the restriction of $\Cal U_a$,
namely $U_a$, fits into an exact sequence
$$0 \to \pi_2^*W_d\spcheck \otimes \pi _1^*\scrO_{\Pee ^1}(-k+1) \to 
\tilde V|D_i \to \pi_2^*W_{n-d} \otimes \pi _1^*\scrO_{\Pee ^1}(-k)\to 0.$$
Make the elementary modification along the divisor $D_i$ corresponding to
the surjection $\tilde V|D_i \to \pi_2^*W_{n-d} 
\otimes \pi _1^*\scrO_{\Pee ^1}(-k)$. The result is a new bundle $V'$ over
$\tilde Z$, such that over $D_i$ we have an exact sequence
$$0 \to \pi_2^*W_{n-d}  \otimes \pi _1^*\scrO_{\Pee ^1}(-k+1)\to V' |D_i \to
\pi_2^*W_d\spcheck \otimes \pi _1^*\scrO_{\Pee ^1}(-k+1) \to 0.$$ Now, since
$H^1(E_{b_i}; W_d \otimes W_{n-d}) = H^1(\scrO_{\Pee^1}) =0$, it follows
from the K\"unneth formula that
$$\Ext^1(\pi_2^*W_d\spcheck \otimes \pi _1^*\scrO_{\Pee ^1}(-k+1),
\pi_2^*W_{n-d} 
\otimes \pi _1^*\scrO_{\Pee ^1}(-k+1)) = 0.$$ Thus $V'\otimes \scrO_{\tilde
Z}(-(k+1)D_i)|D_i = \pi_2^*(W_d\spcheck \oplus W_{n-d})$. It follows by
standard blowup results that $q_*V'\otimes \scrO_{\tilde Z}(-(k+1)\sum
_iD_i)$ is locally free on $Z$ and its restriction to each fiber
$E_{b_i}$ is $W_d\spcheck \oplus W_{n-d}$. This completes the proof.
\endproof

Finally we must deal with the case $a\equiv 1 \bmod n$.

\theorem{6.11} Suppose that $\dim B =2$. Let $A$ be a quasisection of $\Cal
P_{n-1}$, and suppose that
$a\equiv 1 \bmod n$. Suppose that $A$ is smooth and is the blowup of $B$ at
a finite number of points $b_1, \dots, b_r$, and that the image of the
exceptional $\Pee^1$ is a generic line in the fiber $\Pee^{n-1}$ as in
\rom{(6.9)}, in other words it is not contained in one of the hyperplanes
$H_e$. Then the rank $n$ bundle
$V_{A,a}$, which is defined on $Z-\bigcup _iE_{b_i}$, extends to a reflexive
non-locally free sheaf on $Z$, which we continue to denote by $V_{A,a}$.
The sheaf 
$V_{A,a}$ is locally free except at the points $\sigma \cap E_{b_i}$. Near
such points, $V_{A,a}$ has the local form
$$R^{n-2}\oplus M,$$ where $R=\Cee\{z_1, z_2, z_3\}$, and $M$ is the
standard rank two reflexive non-locally free sheaf given by the exact
sequence
$$0 \to R\to R^3 \to M \to 0,$$ where the map $R\to R^3$ is given by
$1\mapsto (z_1, z_2, z_3)$.
\endstatement
\proof We shall just work near a single fiber $E_b = E_{b_i}$ for some $i$.
Thus let $\tilde Z$ be the blowup of
$Z$ along $E_b$, with exceptional divisor $D \cong \Pee^1\times E_b$. The
basic birational picture to keep in mind is the following: if we blow up
the subvariety
$\Pee^1\times \{p_0\}\subset D$, we get a new exceptional divisor $D_1$ in
$Z_1 = \operatorname{Bl}_{\Pee^1\times \{p_0\}}\tilde Z$. Here $D_1\cong
\Pee(\scrO_{\Pee^1} \oplus \scrO_{\Pee^1}(-1))$, and so $D_1$ is isomorphic
to the blowup $\Bbb F_1$ of $\Pee^2$ at one point. The proper transform
$D'$ of $D$ in $Z_1$ meets $D_1$ along the exceptional divisor in $D_1$,
and can be contracted in $Z_1$. The result is a new manifold $Z_2$,
isomorphic to the blowup of $Z$ at the point $\sigma \cap E_b$, where
$D_1$ blows down to the exceptional divisor $P$ in $Z_2$. 

The quasisection $A$ defines a section of the pullback of $\Cal P_{n-1}$ to
$B$, and thus a bundle $\tilde V$ over $\tilde Z$, which we can then pull
back to $Z_1$. The next step is to show that, after appropriate elementary
modifications,
$\tilde V$ corresponds to a bundle over $Z_2$ whose restriction to $P$ is
just
$(T_P\otimes \scrO_P(-1))\oplus \scrO_P^{n-2}$, where $T_P$ is the tangent
bundle to
$P$. Finally, a local lemma shows that every such bundle has a direct image
on $Z$ which has the local form $M\oplus R^{n-2}$. Since each of these
steps is somewhat involved, we divide the proof into three parts. First we
describe the basic geometry of the blowups involved.

Let $\tilde Z$ be the blowup of $Z$ along $E_b$, with exceptional divisor
$D \cong \Pee^1\times E_b$. Let $Z_1$ be the blowup of $\tilde Z$ along
$\Pee^1\times \{p_0\}\subset D$, with exceptional divisor $D_1$. Let
$D'$ be the proper transform of $D$ in $Z_1$. The divisor $D_1
=\Pee(\scrO_{\Pee^1} \oplus \scrO_{\Pee^1}(-1))$ is isomorphic to $\Bbb
F_1$. Let
$j\: D_1\to Z_1$ be the inclusion and $q\: D_1 \to
\Pee^1$ be the morphism induced by projection from a point. Let $\ell =
\Pee^1\times \{p_0\} = D'\cap D_1$, so that $\ell$ is the exceptional
divisor in $D_1$ viewed as the blowup of $\Pee^2$. Finally we let $s\:
D_1\to
\Pee^2$ be the blowup map. On a fiber $\Pee^1\times \{e\}$ with $e\neq
p_0$, 
$\tilde V \otimes \scrO_{\tilde Z}(-D')$ restricts to $\scrO_{\Pee^1}^n$,
whereas it restricts on $\Pee^1\times \{p_0\}$ to $\scrO_{\Pee^1}(1) \oplus
\scrO_{\Pee^1}(-1) \oplus \scrO_{\Pee^1}^{n-2}$. Thus, if $V_0$ is the
pullback to
$Z_1$ of $\tilde V \otimes \scrO_{\tilde Z}(-D')$, then $V_0$ restricts on
$D_1$ to $q^*\left[\scrO_{\Pee^1}(1) \oplus
\scrO_{\Pee^1}(-1) \oplus \scrO_{\Pee^1}^{n-2}\right]$. 

\claim{1}  Let $V_0$ be the pullback to
$Z_1$ of $\tilde V \otimes \scrO_{\tilde Z}(-D')$. Make the elementary
modification
$$0 \to V'\to V_0 \to j_*q^*\scrO_{\Pee^1}(-1) \to 0.$$ Then $V'$
restricted to $\ell$ is the trivial bundle $\scrO_{\Pee^1}^n$. It follows
that $V'|D'$ is pulled back from the factor $E_b$.
\endstatement
\proof  We have an exact sequence
$$0 \to V'|D'\to V_0|D' \to j_*\scrO_{\Pee^1}(-1) \to 0,$$ where we write
$j$ also for the inclusion of the fiber $\ell =\Pee^1\times \{p_0\}$ in the
ruled surface $D'\cong \Pee^1\times E_b$.  By standard formulas for
elementary modifications, it is straightforward  to compute that
$c_2(V'|D') = c_2(V_0|D') -1$. But
$c_2(V_0|D') =h\pi_2^*[p_0]=1$ by the formulas of \S 2.6. Thus $c_2(V'|D')
=0$. Now by a sequence of allowable elementary modifications $V_0|D',
V'|D'=V_1, \dots, V_r$, we can reach a vector bundle $V_r$ over $D'$ whose
restriction to every fiber
$\Pee^1\times \{e\}$  is semistable and thus trivial; this happens if and
only if
$V_r$ is pulled back from the base, and so has $c_2=0$. But each allowable
elementary modification along the fiber $\Pee^1\times \{p_0\}$ drops
$c_2$ by a positive integer. Since $V'|D'$ already has
$c_2=0$, no further elementary modifications are possible. Hence $V'|\ell$
is already semistable and therefore trivial, and thus 
$V'|D'$ is  pulled back from  $E_b$ as claimed.
\endproof

By construction, $V'|\ell$ is given as an extension
$$0 \to \scrO_{\Pee^1}(-1) \oplus \scrO_{\Pee^1}^{n-2} \to V'|\ell \to 
\scrO_{\Pee^1}(1) \to 0.$$  Now $\Ext^1( \scrO_{\Pee^1}(1),
\scrO_{\Pee^1}(-1) \oplus \scrO_{\Pee^1}^{n-2})
\cong H^1(\scrO_{\Pee^1}(-2)) \cong \Cee$, so there is a unique nonsplit
extension of this type, which is clearly the trivial bundle
$\scrO_{\Pee^1}^n$.

\claim{2} With $V'$ as in Claim \rom1, the restriction of $V'$ to $D_1$ is
the pullback $s^*(T_{\Pee^2}(-1)\oplus \scrO_{\Pee^2}^{n-2})$.
\endstatement
\proof By definition, there is an exact sequence
$$0 \to  q^*[\scrO_{\Pee^1}(-1) \oplus \scrO_{\Pee^1}^{n-2}] \otimes
\scrO_{D_1}(-D_1) \to V'|D_1 \to q^*
\scrO_{\Pee^1}(1) \to 0.$$ Next, a straightforward calculation shows that
$\scrO_{D_1}(-D_1)=\scrO_{D_1}(\ell)\otimes q^*\scrO_{\Pee^1}(1)$. Thus the
extensions of $q^*\scrO_{\Pee^1}(1)$ by
$q^*[\scrO_{\Pee^1}(-1) \oplus \scrO_{\Pee^1}^{n-2}]\otimes
\scrO_{D_1}(-D_1))$ are classified by  
$$H^1(D_1; q^*[\scrO_{\Pee^1}(-1) \oplus
\scrO_{\Pee^1}^{n-2}]
\otimes\scrO_{D_1}(\ell)).$$ It is easy to check that $H^1(D_1;
\scrO_{D_1}(\ell)) = 0$ and that $h^1(
\scrO_{D_1}(\ell)\otimes q^*\scrO_{\Pee^1}(-1) ) = 1$. Thus the dimension
of the Ext group in question is one, so that there just one nontrivial
extension up to isomorphism. Note that
$V'|D_1$ is itself such an extension: it cannot be the split extension
since the restriction of $V'|D_1$ to
$\ell$ is trivial. Thus, to complete the proof of Claim 2, it will suffice
to show that
$s^*(T_{\Pee^2}(-1)\oplus \scrO_{\Pee^2}^{n-2})$ is also given as an
extension of
$q^*\scrO_{\Pee^1}(1)$ by $q^*[\scrO_{\Pee^1}(-1) \oplus
\scrO_{\Pee^1}^{n-2}]\otimes
\scrO_{D_1}(-D_1)$. It clearly suffices to do the case $n=2$, i\.e\. show
that
$s^*T_{\Pee^2}(-1)$ is an extension of
$q^*\scrO_{\Pee^1}(1)$ by $q^*\scrO_{\Pee^1}(-1)$, necessarily nonsplit
since the restriction to $\ell$ is trivial. To see this, note that
$T_{\Pee^2}(-1)$ has restriction $\scrO_{\Pee^1}\oplus \scrO_{\Pee^1}(1)$
to every line. Thus by the standard construction (cf\. \cite{11}, p\. 60)
there is an exact sequence
$$0 \to \scrO_{\Bbb F_1}(\ell) \otimes q^*\scrO_{\Pee ^1}(t) \to 
s^*T_{\Pee ^2}(-1) \to q^*\scrO_{\Pee ^1}(1-t)\to 0$$ for some integer $t$.
By looking at $c_2$, we must have $t=0$ and thus 
$s^*T_{\Pee ^2}(-1)$ is an extension of
$q^*\scrO_{\Pee ^1}(1)$ by 
$\scrO_{\Bbb F_1}(\ell)$, which is nonsplit because its restriction to
$\ell$ is trivial. Thus  we have identified $V'|D_1$ with
$s^*(T_{\Pee^2}(-1)\oplus \scrO_{\Pee^2}^{n-2})$.
\endproof

Let $Z_2$ be the result of contracting $D'$ in $Z_1$. This has the effect of
contracting $\ell\subset D_1$ to a point, so that the image of $D_1$ in
$Z_2$ is an exceptional $\Pee^2$, which we denote by $P$. Moreover, by the
above claims $V'$ induces a vector bundle on $Z_2$ whose restriction to $P$
is identified with
$T_{\Pee^2}(-1)\oplus \scrO_{\Pee^2}^{n-2}$. Thus, the proof of (6.11) will
be complete once we prove the following:

\claim{3} Let $X$ be a manifold of dimension $3$ and let $\tilde X$ be the
blowup of $X$ at a point $x$, with exceptional divisor $P\cong \Pee^2$.
Suppose that $W$ is a vector bundle on $\tilde X$ such that $W|P \cong
T_{\Pee^2}(-1)\oplus
\scrO_{\Pee^2}^{n-2}$. Let $\rho\: \tilde X \to X$ be the blowup map. Then
$\rho_*W$ is locally isomorphic to $M\oplus R^{n-2}$ in the notation above.
In particular,  $\rho_*W$ is reflexive but not locally free.
\endstatement
\noindent {\it Proof.} We shall just do the case $n=2$, the other cases
being similar. By the formal functions theorem, the completion of the stalk
of the direct image  $\rho_*W$ at $x$ is $M'= \varprojlim H^0(W \otimes
\scrO_{nP})$. Now from the exact sequences
$$0 \to \scrO_{\Pee^2}(-1) \to \scrO_P^3 \to  W|P \to 0$$ and the sequence
$$0 \to  W\otimes \scrO_{\tilde X}(-(n+1)P) \to  W\otimes
\scrO_{(n+1)P}\to
 W\otimes
\scrO_{nP}\to 0,$$ it is easy to check that the three sections of
$ W|P$ lift to give three generators of $M'$ as an $R$-module. Hence there
is a surjection $\scrO_{\tilde X}^3 \to\rho^*\rho_*W\to W$, and by checking
determinants the kernel is $\scrO_{\tilde X}(P)$. Now up to an change of
coordinates in $\Cee^3$ the only injective homomorphism from
$\scrO_{\tilde X}(P)$ to $\scrO_{\tilde X}^3$ is given by the three
generators of the maximal ideal of $\Cee ^3$ at the origin. Taking direct
images of the exact sequence
$$0 \to \scrO_{\tilde X}(P)\to  \scrO_{\tilde X}^3 \to  W \to 0$$ and using
the vanishing for the first direct image of $\scrO_{\tilde X}(P)$ gives 
$M'\cong M$ as previously defined. So we have established Claim 3, and hence
(6.11).
\endproof

We give a brief and inconclusive discussion of how the above constructions
begave in families, assuming $\dim B =2$ for simplicity. Let $D$ be the unit
disk in $\Cee$. Suppose that we are given a general family of nowhere 
vanishing  sections $s_t$ of $\Cal V_n$ which at a special point $t=0$
acquires a simple zero at $b\in B$. We can view the family $s=\{s_t\}$ as a
section of the pullback of
$\Cal V_n$ to
$B\times D$, where it has a simple zero at $(b,0)$. Thus, for an integer
$a$, there is a bundle $\Cal V_{s,a}$ over $Z\times D-\{(b,0)\}$, which
completes uniquely to a reflexive sheaf over $Z\times D$, which we continue
to denote by
$\Cal V_{s,a}$. For example, if $-(n-2)\leq a \leq 0$, then it is easy to
see that
$\Cal V_{s,a}$ is a bundle over $Z\times D$, whose restriction to $Z\times
\{0\}$  is everywhere regular semistable except over $E_b$ where it
restricts to
$W_d\spcheck \oplus W_{n-d}$ for the appropriate $d$. One can ask if this
holds for all  $a \not\equiv 1 \bmod n$. Note that, if we consider the
relative deformation theory of the unstable bundle $W_d\spcheck
\oplus W_{n-d}$ over the base $B$, for $n=2$ the codimension of the locus
of unstable bundles forces every deformation of $V$ to have unstable
restriction to some fibers, whereas for $n> 2$ we expect that in the
general deformation
$V_t$ we can arrange that the restriction of $V_t$ to every fiber is
semistable.

If $a\equiv 1 \bmod n$, then $\Cal V_{s,a}$ is a flat family of coherent
sheaves. However, there is no reason {\it a priori\/} why $\Cal
V_{s,a}|Z\times \{0\}$ is reflexive.  In fact, preliminary calculations
suggest that, for $a=1$, the restriction $\Cal V_{s,a}|Z\times \{0\}$ has
the local form $M\oplus \frak m^{n-2}$, where $\frak m$ is the maximal
ideal of the point $\sigma \cap E_b$. Note that the $R$-module
$M$ is not smoothable, even locally, but that $R^k\oplus M$ is smoothable
to a free
$R$-module for all $k\geq 1$. One can also show that the more complicated
$R$-module $\frak m^k\oplus M$ is smoothable to a free
$R$-module for all $k\geq 1$. This agrees  with the picture for  sections of
the bundle $\Cal V_n$: for
$n=2$, if a section has a simple isolated zero, that zero must remain under
deformation, but for $n>2$ we expect in general that we can deform to an
everywhere nonzero section in general. 

\ssection{6.4. Bundles which are not regular in high codimension.}

In this subsection we consider bundles $V$ such that $V|E_b$ is semistable
for all
$b$, and $Y= \{\, b\in B: \text{ $V|E_b$ is not regular}\,\}$ has
codimension at least $2$ in $B$. The first lemma shows that, if the
spectral cover $C_A$ is smooth, then $V$ is in fact everywhere regular.

\lemma{6.12} Let $V$ be a vector bundle over $Z$ such that $V|E_b$ is
semistable for all $b$, and $Y= \{\, b\in B: \text{ $V|E_b$ is not
regular}\,\}$ has codimension at least $2$ in $B$. Suppose that the
associated spectral cover $C_A$ is smooth. Then $V|E_b$ is regular for all
$b\in B$. More generally, suppose that 
$V$ is a vector bundle over $Z$ such that $Y= \{\, b\in B: \text{ $V|E_b$ is
either not semistable or not regular}\,\}$  has codimension at least $2$ in
$B$, that the section $A$ defined by $V$ over $B-Y$ extends to a section
over all of $B$, and that the associated spectral cover
$C_A$ of $B$ is smooth. Then $V|E_b$ is semistable and regular for all
$b\in B$.
\endstatement
\proof We have seen in (5.7) that there is a line bundle $N$ on
$C_A-g_A^{-1}(Y)$ such that $V|Z-\pi^{-1}(Y) \cong V_{A,0}[N]$. Since $C_A$
is smooth, and
$g_A^{-1}(Y)$ has codimension at least two in $C_A$, the line bundle
$N$ on  $C_A-g_A^{-1}(Y)$ extends to a line bundle over $C_A$, which we
continue to denote by $N$. We now have two vector bundles on $Z$, namely
$V$ and
$V_{A,0}[N]$, which are isomorphic over $Z-\pi^{-1}(Y)$. Since the
codimension of 
$\pi^{-1}(Y)$ in $Z$ is at least two, $V$ and $V_{A,0}[N]$ are isomorphic.
But 
$V_{A,0}[N]$ restricts to a regular bundle on every fiber, and so the same
must be true for $V$.
\endproof

We turn to methods for constructing bundles which are semistable on every
fiber but which are not regular in codimension two. Of course, by the above
lemma, the corresponding spectral covers will not be smooth. The idea is to
find such bundles by using a three step filtration, as opposed to the
two-step extensions which have used from Section 3 onwards in our
constructions. Such constructions correspond to nonmaximal parabolic
subgroups in
$SL_n$.

Consider first the case of a single Weierstrass cubic $E$. We seek bundles
of rank $n+1$ which have a filtration $0\subset F^0\subset F^1 \subset V$,
where
$F^0\cong W_k\spcheck, F^1/F^0\cong \scrO_E$, and $V/F^1\cong W_{n-k}$. Such
extensions can be described by a nonabelian cohomology group as in \cite{8}.
However, it is also easy to describe them directly. Note that a fixed $F^1$
is described by an extension class $\alpha _0$ in $\Ext^1(\scrO_E,
W_k\spcheck) \cong H^1(W_k\spcheck) \cong \Cee$. If $\alpha _0 = 0$, then
$F^1 = W_k\spcheck
\oplus \scrO_E$, and if $\alpha _0 \neq 0$ then $F^1\cong W_{k+1}\spcheck$.
Having determined $F^1$, the extension
$F^2$ corresponds to a class in $\Ext^1(W_{n-k}, F^1)$. Since $\Hom
(W_{n-k},
\scrO_E) =
\Ext^2(W_{n-k}, W_k\spcheck) =0$, there is a short exact sequence
$$0 \to \Ext ^1(W_{n-k}, W_k\spcheck) \to \Ext^1(W_{n-k}, F^1) \to
\Ext^1(W_{n-k},
\scrO_E) \to 0,$$ and so $\dim \Ext^1(W_{n-k}, F^1) = n+1$. Thus roughly
speaking the moduli space of filtrations as above is an affine space
$\Cee^{n+2}$. In fact, by general construction techniques there is a
universal bundle $\Cal F^1$ over
$\Ext^1(\scrO_E, W_k\spcheck)\times E =\Cee\times E$. We can then form the
relative Ext sheaf
$$Ext^1_{\pi_1}(\pi _2^*W_{n-k}, \Cal F^1)= R^1\pi_1{}_*(\pi
_2^*W_{n-k}\spcheck\otimes \Cal F^1).$$ It is a vector bundle of rank $n+1$
over $\Cee$, which is necessarily trivial, and thus the total space of this
vector bundle is $\Cee^{n+2}$. There is a universal extension of
$\pi_2^*W_{n-k}$ by $\Cal F^1$ defined over
$\Cee^{n+2}\times E$. It follows that the set of filtrations is indeed
parametrized by a moduli space isomorphic to
$\Cee^{n+2}$, although there is not a canonical linear structure. What is
canonical is the exact sequence
$$0 \to \Ext^1(W_{n-k}, W_k\spcheck) \to \Cee^{n+2} \to \Ext^1(\scrO_E,
W_k\spcheck) \oplus \Ext^1(W_{n-k}, \scrO_E) \to 0.$$ We understand this
sequence to mean that the first term, which is a vector space, acts on the
middle term, which is just an affine space, via affine translations, and
the quotient is the last term, which is again a vector space. Here the
projection to
$\Ext^1(\scrO_E, W_k\spcheck) \oplus \Ext^1(W_{n-k}, \scrO_E)$ measures the
extensions $F^1$ of
$\scrO_E$ by $W_k\spcheck$ and $V/F^0$ of $W_{n-k}$ by $\scrO_E$.  We
denote the image of $\xi\in \Cee^{n+2}$ in $\Ext^1(\scrO_E, W_k\spcheck)
\oplus
\Ext^1(W_{n-k}, \scrO_E) \cong \Cee\oplus
\Cee$ by $(\alpha _0, \alpha_1)$. Here $\alpha _0 \neq 0$ if and only if
$F^1\cong W_{k+1}\spcheck$ and $\alpha_1 \neq 0$ if and only if $V/F^0\cong
W_{n-k+1}$. In  case $\alpha _0 = 0$, say, $F^1\cong W_k\spcheck \oplus
\scrO_E$, and
$\Ext^1(W_{n-k}, F^1)$ naturally splits as $\Ext ^1(W_{n-k},
W_k\spcheck)\oplus
\Ext^1(W_{n-k}, \scrO_E)$. In this case, both the class $\alpha_1$ and the
class
$e\in \Ext ^1(W_{n-k}, W_k\spcheck)$ are well-defined. A similar statement
holds if $\alpha _1 =0$. Note that the affine space $\Cee^{n+2}$
parametrizes filtrations $F^i$ together with {\sl fixed\/} isomorphisms
$F^0\to W_k\spcheck$,
$F^1/F^0\to \scrO_E$, $V/F^1 \to W_{n-k}$.

The  subspace
$\Ext^1(W_{n-k}, W_k\spcheck)$, namely where both $\alpha_0$ and $\alpha _1$
vanish, corresponds to those $V$ of the form $V'\oplus
\scrO_E$, where $V'$ is an extension of $W_{n-k}$ by $W_k\spcheck$. There
is a hyperplane $H$ in $\Ext^1(W_{n-k}, W_k\spcheck)$ where such $V'$
contain a Jordan-H\"older quotient isomorphic to $\scrO_E$, and thus over
the locus $\alpha _0 = \alpha _1 = 0, e\in H$, $V=V'\oplus
\scrO_E$ has a subbundle of the form $\scrO_E\oplus \scrO_E$. Hence, over a
affine subspace of
$\Cee^{n+2}$ of codimension three, the $V$ we have constructed are not
regular.

\lemma{6.13} Suppose that $V$ corresponds to a class $\xi \in \Cee^{n+2}$,
and that $\alpha_0, \alpha _1$ are as above.
\roster
\item"{(i)}" $V$ is unstable if and only $\alpha_0 \alpha _1=0$ and $e=0$
\rom(this statement is well-defined by the above remarks\rom).
\item"{(ii)}" If $\alpha_0 \alpha _1 \neq 0$, then $h^0(V) = 0$ and
conversely.
\endroster
\endstatement
\proof (i) Let us assume for example that $\alpha _0 = e=0$. Then $F^1 =
W_k\spcheck \oplus \scrO_E$ and $V$ is isomorphic either to $W_k\spcheck
\oplus W_{n-k+1}$ or to $W_k\spcheck \oplus \scrO_E \oplus W_{n-k}$, and in
either case it is unstable. Conversely, if $V$ is unstable, then it has a
maximal destabilizing subsheaf $W$ of positive degree, which is stable and
which must map nontrivially onto $W_{n-k}$. Thus $\deg W =1$. Now if $W\cap
W_k\spcheck\neq 0$, then $W\cap W_k\spcheck$ has degree
$\leq -1$ and is contained in the kernel of the map $W\to W_{n-k}$. This
would force the image of $W$ to have degree at most zero, which is
impossible. So $W\cap W_k\spcheck=0$ and thus the map $W\to V/F^0$ is
injective. Now either $V/F^0
\cong W_{n-k+1}$ or $V/F^0
\cong W_{n-k}\oplus \scrO_E$. In the first case, $W\cong W_{n-k+1}$ by the
stability of $ W_{n-k+1}$ and $V\cong W_k\spcheck \oplus W_{n-k+1}$. In
this case
$\alpha _0 = e=0$. In the remaining case, $V/F^0
\cong W_{n-k}\oplus \scrO_E$ and $W\cong W_{n-k}$. In this case $\alpha
_1=e=0$.  In both cases we must have $\alpha_0 \alpha _1=0$ and $e=0$. 

\smallskip
\noindent (ii) First suppose that $\alpha_0 \alpha _1 \neq 0$. Since $\alpha
_0\neq0$, $F^1\cong W_{k+1}\spcheck$. From the  exact sequence
$$0\to F^1\to V \to W_{n-k} \to 0,$$ and the fact that $H^0(F^1) = 0$,
there is an exact sequence $H^0(V) \to H^0(W_{n-k}) \to H^1(F^1)$. If we
compose the map $H^0(W_{n-k}) \to H^1(F^1)$ with the natural map $ H^1(F^1)
\to H^1(\scrO_E)$, the result is $\alpha _1$ up to a nonzero scalar. Thus,
if $\alpha _1\neq 0$, the map $H^0(W_{n-k}) \to H^1(F^1)$ is injective and
so $H^0(V) =0$. Conversely, suppose that either $\alpha _0$ or
$\alpha _1$ is zero. If for example $\alpha _1=0$, then $V/F^0 \cong
W_{n-k}\oplus \scrO_E$, so that $h^0(V/F^0)=2$. Since $H^0(V)$ is the
kernel of the map $H^0(V/F^0) \to H^1(F^0) = H^1(W_k\spcheck)\cong \Cee$,
$H^0(V) \neq 0$. The case $\alpha _0 \neq 0$ is similar and simpler. Thus,
if $h^0(V) =0$, then  
$\alpha_0 \alpha _1 \neq 0$.
\endproof

The group  $\Cee^*\times \Cee^*$ (or more precisely $\Cee^3/\Cee$) acts on
the affine space
$\Cee^{n+2}$ by acting on the identifications of the quotients
$F^{i+1}/F^i$ with the standard bundles. The quotient by this action (which
is not in fact separated) is the set of bundles $V$ of rank $n+1$, together
with filtrations on
$V$ with the appropriate graded object. The action of  $\Cee^*\times
\Cee^*$ is compatible with the projection to
$\Cee^2$. If we normalize the action so that
$(\lambda, \mu)\cdot (\alpha_0, \alpha _1) = (\lambda\alpha_0, \mu\alpha
_1)$, then $(\lambda, \mu)$ acts on the distinguished subspace
$\Ext^1(W_{n-k}, W_k\spcheck) \cong \Cee^n$ by $e\mapsto \lambda\mu e$.
Clearly, the action is free over the set $\alpha_0\alpha _1\neq 0$. The
quotient of the points where
$\alpha _0\neq 0$,
$\alpha _1 = 0$, $e\neq 0$ is a $\Pee^{n-1}$, and this $\Pee^{n-1}$ is
identified with the corresponding $\Pee^{n-1}$ where $\alpha _1\neq 0$,
$\alpha _0 = 0$, $e\neq 0$; in fact, the $\Cee^*\times \Cee^*$-orbits
intersect along the subspace where $\alpha _0=\alpha _1 = 0$, $e\neq 0$.
The points
$\alpha_0\alpha _1=e=0$ are unstable points and do not appear in a GIT
quotient for the action. We also have the map (1.5) 
$\Phi\: \Cee^{n+2}-(\Cee \cup\Cee)$ to the coarse moduli space $\Pee^n$ of
semistable bundles of rank $n+1$ on $E$. By Lemma 6.13, the image of the two
subsets
$\{\,\alpha _0 = 0, e\neq 0\,\}$ and $\{\,\alpha _1 = 0, e\neq 0\,\}$ is
exactly the hyperplane in $\Pee^n$ corresponding to bundles $V$ such that
$h^0(V) \neq 0$, or in other words such that $V$ has $\scrO_E$ as a
Jordan-H\"older quotient.

\lemma{6.14} The map $\Phi\: \Cee^{n+2}-(\Cee \cup\Cee) \to \Pee^n$ is the
geometric invariant theory quotient of $\Cee^{n+2}-(\Cee \cup\Cee)$ by the
action of $\Cee^*\times \Cee^*$.
\endstatement
\proof First suppose that the point $x\in \Cee^{n+2}-(\Cee \cup\Cee)$ lies
in the open dense subset
$\alpha_0\alpha _1\neq 0$ where $\Cee^*\times \Cee^*$ acts freely. Thus if
$V$ is the vector bundle corresponding to
$x$, then  $h^0(V) =0$; equivalently,
$V$ has no Jordan-H\"older quotient equal to $\scrO_E$, and $V$ is a regular
semistable bundle. If $\Phi(x) = \Phi(x')$, then $x'$ also lies in the set
$\alpha_0\alpha _1\neq 0$, and the bundle $V'$ corresponding to $x'$ is also
regular and semistable. Thus $V\cong V'$, and we must determine if the
filtration $F^i$ on $V$ exists is unique up to isomorphism.  First, if $V$
is a regular semistable bundle of rank $n+1$, then it is an extension of
$W_{n-k}$ by
$W_{k+1}\spcheck$, where the subbundle $W_{k+1}\spcheck$ of $V$ is unique
modulo automorphisms of $V$, and taking the further filtration of
$W_{k+1}\spcheck$ by the subbundle $W_k\spcheck $, with quotient $\scrO_E$.
Thus
$V=\Phi(x)$ for some $x$. Conversely, if $V$ has on it a filtration $F^i$
with
$\alpha_0\alpha_1 \neq 0$, then $F^1\cong W_{k+1}\spcheck$. Moreover, if
$H^0(V) = 0$, then every subbundle of $V$ isomorphic to $W_k\spcheck$ is
contained in a subbundle isomorphic to $W_{k+1}\spcheck$ (whereas if
$H^0(V) \neq 0$, this is no longer the case; cf\. \S 3.2). Thus the
filtration $F^i$ is unique up to automorphisms of $V$. The above argument
shows that $\Phi$ induces an isomorphism
$$\left(\Cee^{n+2} - \{\, \alpha_0\alpha _1= 0\,\}\right)/\Cee^*\times
\Cee^* \to
\Pee^n - H.$$

In case $x$ lies in the set $\alpha _0 = 0, , \alpha _1\neq 0, e\neq 0$, a
straightforward argument identifies the quotient by $\Cee^*\times \Cee^*$
with
$H\subset \Pee^n$, and likewise for $\alpha _0 \neq 0, \alpha _1= 0, e\neq
0$, 
$\alpha _0 =\alpha _1= 0, e\neq 0$.
\endproof

The coarse moduli space $\Pee^n$ has its associated spectral cover $T$,
which is an $(n+1)$-sheeted cover of $\Pee^n$. Let $\tilde T \to
\Cee^{n+2}-(\Cee
\cup\Cee)$ be the pulled back cover of $\Cee^{n+2}-(\Cee \cup\Cee)$ via the
morphism $\Phi$. Using Lemma 6.14, we can see directly that $\tilde T$ is
singular, with the generic singularities a locally trivial family of
threefold double points. In fact, the inverse image of $H$ in $T$ is of the
form $H\cup T'$, where $T'$ is the spectral cover of $H\cong \Pee^{n-1}$.
The intersection of $H$ and $T'$ is transverse (see \S 5.7), and $H\cap
T'$, viewed as a subset of
$H\subset \Pee^n$, corresponds to those bundles which have $\scrO_E$ as a
Jordan-H\"older quotient with multiplicity at least two. If $t$ is a local
equation for $H$ in $\Pee^n$ near a generic point of $H\cap T'$, there are
local coordinates on $T$ for which $t=uv$, since $H$ splits into $H\cup
T'$. Thus the local equation for $\tilde T$ is $\alpha_0\alpha _1= uv$,
which is the equation for a family of threefold double points. 
 
We can also do the above constructions in families $\pi\: Z \to B$. We
could take the point of view of \cite{8} and realize the relative
nonabelian cohomology groups as a bundle of affine spaces over $B$.
However, it is also possible to proceed directly as in \S 5.2. We seek
vector bundles $V$ which have a filtration
$0\subset F^0\subset F^1 \subset V$, where $F^0\cong \Cal W_k\spcheck
\otimes
\pi^*M_0$, $F^1/F^0\cong \pi^*M_1$, and $V/F^1\cong \Cal W_{n-k} \otimes
\pi^*M_2$ for line bundles $M_0, M_1, M_2$ on $B$. Of course, we can
normalize by twisting
$V$ so that one of the $M_i$ is trivial. The analysis of such extensions
parallels the analysis for a single $E$. We begin by constructing $F^0$. It
is described by an extension class in
$$\gather H^1(Z; \pi^*M_1^{-1}\otimes \Cal W_k\spcheck \otimes \pi^*M_0)
\cong H^0(B; R^1\pi_*(\Cal W_k\spcheck) \otimes M_1^{-1}\otimes M_0) \\
=H^0(B;L^{-k}\otimes M_1^{-1}\otimes M_0).
\endgather$$ If the difference line bundle $M_1^{-1}\otimes M_0$ is
sufficiently ample, then there will be nonzero sections $\alpha _0$ of
$L^{-k}\otimes M_1^{-1}\otimes M_0$ vanishing along a divisor $D_0$ in $B$.
Next, we seek extensions of $F^1$ by $\Cal W_{n-k} \otimes \pi^*M_2$. Now
$H^0(\Cal W_{n-k}\spcheck \otimes \pi^*M_2^{-1}\otimes \pi^*M_1) = 0$, and
by the Leray spectral sequence
$$H^2(\Cal W_{n-k}\spcheck \otimes \pi^*M_2^{-1}\otimes \Cal W_k\spcheck
\otimes
\pi^*M_0) \cong H^1(B; R^1\pi_*(\Cal W_{n-k}\spcheck \otimes \Cal
W_k\spcheck)\otimes M_2^{-1}\otimes M_0).$$ We assume that $M_2^{-1}\otimes
M_0$ is so ample that the above group is zero. In this case there is an
exact sequence
$$\gather 0 \to H^1(\Cal W_{n-k}\spcheck \otimes \pi^*M_2^{-1}\otimes \Cal
W_k\spcheck 
\otimes \pi^*M_0) \to \Ext^1(\Cal W_{n-k} \otimes \pi^*M_2, F^1) \to \\
H^1(\Cal W_{n-k}\spcheck \otimes \pi^*M_2^{-1}\otimes \pi^*M_1) \to 0.
\endgather$$ The left-hand group is $H^0(R^1\pi_*(\Cal W_{n-k}\spcheck
\otimes \Cal W_k\spcheck)\otimes M_2^{-1}\otimes M_0)$, and the right-hand
group is
$H^0(L^{n-k} \otimes M_2^{-1}\otimes M_1)$. Thus, for $ M_2^{-1}\otimes M_1$
sufficiently ample, there will exist sections $\alpha _1$ of $L^{n-k}
\otimes M_2^{-1}\otimes M_1$, vanishing along a  divisor $D_1$ in $B$, and
we will be able to lift these sections to extension classes in $\Ext^1(\Cal
W_{n-k} \otimes
\pi^*M_2, F^1)$. Moreover, if we restrict, say, to the divisor $D_0 = 0$,
then there is also a well-defined class $e$ in $H^0(R^1\pi_*(\Cal
W_{n-k}\spcheck 
\otimes \Cal W_k\spcheck)\otimes M_2^{-1}\otimes M_0)$. There is a divisor
$D$ on
$B$ corresponding to such extensions which have a factor $\scrO_{E_b}$ for
$b\in D$. In fact, if $\alpha = c_1(L\otimes M_2^{-1}\otimes M_0)$, then it
follows from (4.15) and (5.9) that $[D] = \alpha - nL$. (Compare also
(5.21).) As long as
$M_2^{-1}\otimes M_0$ is also sufficiently ample, we can assume that the
divisors
$D_0, D_1$ and $D$ are smooth and intersect transversally in a subvariety
of $B$ of codimension three. Along this subvariety, $V$ fails to be regular.

Note that the $V$ constructed above are a deformation of $V'\oplus
\scrO_Z$, where
$V'$ is a twist of a bundle of the form $V_{A, 1-k}$; it suffices for
example to take $\alpha _0 = 0$ and $e\neq 0$. 

For generic choices, the spectral cover $C_A$ will acquire  singularities
in codimension three, which will generically be families of  threefold
double points. In particular, there are Weil divisors on $C_A$ which do not
extend to Cartier divisors, as predicted by Lemma 6.12. It is also amusing
to look at the case $\dim B =2$, where for generic choices the spectral
cover will be smooth. The construction then deforms $V'\oplus \scrO_Z$ to a
bundle $V$ which has regular semistable restriction to every fiber.
Starting with a generic $V' = V_{A,a}(n)$ of rank $n$, we cannot in general
deform $V'\oplus \scrO_Z$ to a standard bundle
$V_{A,b}(n+1)\otimes \pi^*N_0$. Instead, the spectral cover $C_A$ has
Picard number larger than expected. In fact, we have the divisor $F =
(r^*\sigma \times _BZ)\cap C_A$, which is mapped isomorphically onto its
image in $B$, and this image is the same as $A\cap \Cal H\subset \Cal
P_{n-1}$. Now $A\cap \Cal H$ corresponds to the bundles $V$ such that
$h^0(V) \neq 0$, and thus by Lemma 6.13 this locus is just
$D_0\cup D_1$. Thus in $C_A$ the divisor $F$ splits into a sum of two
divisors, which we continue to denote by $D_0$ and $D_1$. Using these extra
divisors, we can construct more vector bundles over $Z$, of the form
$V_{A,0}[N]$ for some extra line bundle $N$, which enable us to deform 
$V'\oplus \scrO_Z$ to a bundle which is everywhere regular and semistable.

Let us just give the details in a symmetric case. Let $M$ be a sufficiently
ample line bundle on $B$. There exist bundles
$V$ on $Z$ which have regular semistable restriction to every fiber and
also have a filtration $F^0\subset F^1 \subset  V$, with
$$V/F^1 \cong \Cal W_k\otimes \pi^*M^{-1}; \qquad F^1/F^0\cong \scrO_Z;
\qquad F^0
\cong \Cal W_k\spcheck \otimes \pi^*M.$$ The bundle $V$ is a deformation of
a bundle of the form $V'\oplus \scrO_Z$. Thus, there must exist a line
bundle $N$ on the spectral cover $C_A$ such that
$V_{A,0}[N]$ has the same Chern classes as $V$. Direct calculation with the
Grothendieck-Riemann-Roch theorem shows that this happens for
$$N = M \otimes \scrO_{C_A}(-F + (k+1)D_0 + kD_1)$$ as well as for
$$N = M \otimes \scrO_{C_A}(-F + kD_0 + (k+1)D_1),$$ and that these are the
only two ``universal" choices for $N$.

\remark{Question} Suppose that $\dim B = 3$ and consider spectral covers
which have an ordinary threefold double point singularity. The local Picard
group of the singularity is $\Zee$. Given $a\in \Zee$, we can twist by a
line bundle over the complement of the singularity which correspond to the
element $a\in \Zee$. The result is a vector bundle on $Z$, defined in the
complement of finitely many fibers, and thus the direct image is a coherent
reflexive sheaf on $Z$. What is  the relationship of local behavior of this
sheaf at the finitely many fibers to the integer $a$?
\endremark 

\section{7. Stability.}

Our goal in this final section will be to find sufficient conditions for
$V_{A,a}$, or more general bundles constructed in the previous two
sections, to be stable with respect to a suitable ample divisor. Here
suitable means in general a divisor of the form $H_0+ N\pi^*H$, where $H_0$
is an ample divisor on $Z$ and
$H$ is an ample divisor on $B$, and $N\gg 0$. As we have already see in \S
5.6, for $A =\bold o$, the bundle $V_{\bold o, a}$ is essentially always
unstable with respect to every ample divisor. Likewise, suppose that $A$ is
a section lying in
$\Cal H$ as in \S 5.7, so that $V_{A,a}$ has a surjection to $\pi^*L^a$. If
the line bundle corresponding to $A$ is sufficiently ample, it is easy to
see that for appropriate choices of $a$ we can always arrange
$\mu_H(\pi^*L^a) <
\mu_H(V_{A,a})$, so that $V_{A,a}$ is unstable. Thus, we shall have to make
some assumptions about $A$. More generally, let $V$  be a bundle whose
restriction to the generic fiber is regular and semistable, and let $A$ be
the associated quasisection. It turns out that, if the spectral cover
$C_A$ is irreducible, then $V$ is stable with respect to all divisors of
the form
$H_0+ N\pi^*H$, provided that $N\gg 0$. A similar result holds in families.
However, we are only able to give  an effective estimate for $N$ in case
$\dim B = 1$. In particular, whether there is an effective bound for $N$
which depends only on $Z$, $H_0$, $H$, $c_1(V)$, and $c_2(V)$ is open in
case $\dim B > 1$. We believe that such a bound should exist, and can give
such an explicit bound for a general $B$ in the rank two case for an
irreducible quasisection $A$. (Of course, when $\dim B > 2$, an irreducible
quasisection $A$ will almost never be an actual section.) However, we shall
not give the details in this paper.

\ssection{7.1. The case of a general $Z$.}

Let $\pi\: Z\to B$ be a flat family of Weierstrass cubics with a section. We
suppose in fact that $Z$ is smooth of dimension $d+1$. Fix an ample divisor
$H_0$ on $Z$ and an ample divisor $H$ on $B$, which we will often identify
with $\pi^*H$ on $Z$. 

\theorem{7.1} Let $V$ be a vector bundle of rank $n$ over $Z$ whose
restriction to the generic fiber is regular and semistable, and such that
the spectral cover of the quasisection corresponding to $V$ is irreducible.
Then there exists  an $\epsilon_0 > 0$, depending on $V, H_0, H$, such that
$V$ is  is stable with respect to 
$\epsilon H_0 + H$ for all $0< \epsilon < \epsilon_0$.
\endstatement
\proof Let $W$ be a subsheaf of $V$ with $0 <
\operatorname{rank} W < r$. The semistability  assumption on $V|f$, for a
generic $f$, and the fact that $W|f \to V|f$ is injective for a generic $f$
imply that $c_1(W) \cdot f
\leq 0$. If however $c_1(W) \cdot f = 0$, then $W$ and $V/W$ are also
semistable on the generic fiber. By Proposition 5.22, the spectral cover
corresponding to $V$ would then be reducible (the proof in (5.22) needed
only that $V$ has regular semistable restriction to the generic fiber),
contrary to hypothesis. Thus in fact 
$c_1(W) \cdot f < 0$. Equivalently,
$c_1(W) \cdot H^d  < 0$.

For a torsion free sheaf $W$, define $\mu _H(W) = \dsize \frac{c_1(W)\cdot
H^d}{\operatorname{rank} W}$, by analogy with an ample $H_0$. If $W$ is a
subsheaf of $V$ such that  $0 <
\operatorname{rank} W < n$, then $\mu _H(W)$ is a  strictly negative
rational number with denominator bounded by $n-1$. 

\lemma{7.2} There is a constant $A$, depending only on $V, H_0, H$, such
that
$$\frac{c_1(W)\cdot H^i\cdot H_0^{d-i}}{\operatorname{rank} W} \leq A$$ for
all $i$ with $0\leq i\leq d$ and all nonzero subsheaves $W$ of $V$. 
\endstatement
\proof There exists a filtration 
$$0\subset F^0\subset F^1\subset \cdots \subset F^{n-1}=V$$ such that
$F^j/F^{j-1}$ is a torsion free rank one sheaf, and thus is of the form
$L_j\otimes I_{X_j}$ for $L_j$ a line bundle on $Z$ and $X_j$ a subscheme of
codimension at least two (possibly empty). Suppose that $W$ has rank one.
Then there is a nonzero map from $W$ to $L_j\otimes I_{X_j}$ for some $j$,
and thus $W$ is of the form
$L_j\otimes \scrO_Z(-D)\otimes I_X$ for some effective diviisor $D$ on $Z$
and subscheme
$X$ of codimension at least two (possibly empty). Thus
$$c_1(W)\cdot H^i\cdot H_0^{d-i} \leq c_1(L_j)\cdot H^i\cdot H_0^{d-i}.$$
Thus these numbers are bounded independently of $W$. In case $W$ has
arbitrary rank $r$, $1\leq r\leq n-1$, find a similar filtration of the
bundle
$\bigwedge ^rV$ by subsheaves whose successive quotients are rank one
torsion free sheaves, and use the existence of a nonzero map $\bigwedge
^rW\to \bigwedge ^rV$ to argue as before.
\endproof

Returning to the proof of Theorem 7.1, if $W$ is a subsheaf of $V$ such
that      
$0 < \operatorname{rank} W < n$, it follows that
$$\mu_{\epsilon H_0 + H}(W) = \frac{c_1(W)\cdot (\epsilon H_0 +
H)^d}{\operatorname{rank} W}\leq -\frac{1}{n-1}+O(\epsilon).$$ On the other
hand, since $\det V$ is pulled back from $B$, $c_1(V) \cdot H^d = 0$ and so 
$$\mu_{\epsilon H_0 + H}(V) = \frac{c_1(V)\cdot (\epsilon H_0 + H)^d}{n}=
O(\epsilon).$$  Thus, for $\epsilon$ sufficiently small, for every subsheaf
$W$ of $V$ with $0 <
\operatorname{rank} W < n$, 
$$\mu_{\epsilon H_0 + H}(W) < \mu_{\epsilon H_0 + H}(V).$$ In other words,
$V$ is stable with respect to
$\epsilon H_0 + H$.
\endproof

\corollary{7.3} Let $\Cal V$ be a family of vector bundles over $S\times
Z$, such that, for each $s\in S$, the restriction $V_s=\Cal V|\{s\}\times
Z$ has regular semistable restriction to the generic fiber of $\pi$ and the
corresponding spectral cover is irreducible. Then there exists  an
$\epsilon_0 > 0$, depending on $\Cal V, H_0, H$, such that, for every $s\in
S$,
$V_s$ is  is stable with respect to 
$\epsilon H_0 + H$ for all $0< \epsilon < \epsilon_0$.
\endstatement
\proof We may assume that $S$ is irreducible. The proof of Theorem 7.1 goes
through as before as long as we can uniformly bound the integers
$c_1(W)\cdot H^i\cdot H_0^{d-i}$ as
$W$ ranges over subsheaves of $V_s$ over all $s\in S$. But there exists a
filtration
$$0\subset F^0\subset F^1\subset \cdots \subset F^{n-1}=\Cal V$$ such that
$F^j/F^{j-1}$ is a torsion free rank one sheaf on $S\times Z$, and thus is
of the form
$\Cal L_j\otimes I_{\Cal X_j}$ for $\Cal L_j$ a line bundle on $S\times Z$
and
$\Cal X_j$ a subscheme of $S\times Z$ of codimension at least two, such
that, at a generic point $s$ of $S$, $(\{s\}\times Z)\cap\Cal X_j $ has
codimension at least two in $Z$. On a nonempty Zariski open subset $\Omega$
of
$S$, the filtration restricts to a filtration of $V_s$ of the form used in
the proof of Lemma 7.2, and
$c_1(\Cal L_j| \{s\}\times Z)$ is independent of $s$. Similar filtrations
exist for the exterior powers $\bigwedge ^r\Cal V$. This bounds 
$c_1(W)\cdot H^i\cdot H_0^{d-i}$ as
$W$ ranges over subsheaves of $V_s$ over all $s$ in a nonempty Zariski open
subset of $S$.  By applying the same construction to the components of
$S-\Omega$ and  induction on
$\dim S$, we can find the desired bound for all $s\in S$.
\endproof

\ssection{7.2. The case of an elliptic surface.}

In case $\dim B =1$, there is a more precise result.

\theorem{7.4} Let $\pi\: Z \to B$ be an elliptic surface and let $H_0$ be
an ample divisor on $Z$. Let $f$ be the numerical equivalence class of a
fiber.  Let  
$V$ be a vector bundle of rank $n$ on $Z$ which is regular and semistable
on the generic fiber, with $\det V$ the pullback of a line bundle  on
$B$, and with $c_2(V) = c$, and such that the spectral cover of $V$ is
irreducible. Then for all $\dsize t\geq t_0= \frac{n^3}{4}c_2(V)$, $V$ is
stable with respect to
$H_0+tf=H_t$.
\endstatement
\proof If $V$ is $H_{t_0}$-stable, then as it is $f$-stable (here stability
is defined with respect to the nef divisor $f$) it is stable with respect
to every convex combination of $H_{t_0}$ and $f$ and thus for every divisor
$H_t$ with
$t\geq t_0$. Thus we may assume that $V$ is not $H_{t_0}$-stable for some
$t_0\geq 0$.

\lemma{7.5} Suppose that $V$ is not $H_{t_0}$-stable for some $t_0\geq 0$.
Then there exists a
$t_1 \geq t_0$ and a divisor $D$ such that
$D\cdot H_{t_1} = 0$ and 
$$-\frac{n^3}{2}c_2(V) \leq D^2 < 0.$$
\endstatement 
\proof By  Theorem 7.1, for all $t\gg 0$, $V$ is $H_t$-stable.  Let $t_1$
be the greatest lower bound of the $t$ such that, for all $t'\geq t$, $V$ is
$H_{t'}$-stable. Thus
$t_1\geq 0$. The condition that $V$ is $H_t$-unstable is clearly an open
condition on $t$. It follows that $V$ is strictly $H_{t_1}$-semistable, so
that there is an exact sequence
$$0 \to V' \to V \to V'' \to 0,$$ with both $V'$, $V''$ torsion free and of
strictly smaller rank than $V$, and with
$\mu_{H_{t_1}}(V') = \mu_{H_{t_1}}(V'') = \mu_{H_{t_1}}(V)$. Thus, both
$V'$ and
$V''$ are $H_{t_1}$-semistable. Let $D = r'c_1(V'') - r''c_1(V')$. Then the
equality $\mu_{H_{t_1}}(V') = \mu_{H_{t_1}}(V'')$ is equivalent to $D\cdot
H_{t_1} = 0$. Note that $D$ is not numerically trivial, for otherwise $V$
would be strictly $H_t$-semistable for all $t$, contradicting the fact that
it is
$H_t$-stable for $t\gg 0$. Thus, by the Hodge index theorem, $D^2<0$. Now,
for a torsion free sheaf of rank $r$, define the Bogomolov number (or
discriminant) of
$W$ by
$$B(W) = 2rc_2(W) - (r-1)c_1(W)^2.$$ If $W$ is semistable with respect to
some ample divisor, then $B(W) \geq 0$. Finally, we have the identity
(\cite{5}, Chapter 9, ex\. 4):
$$B(V) = 2nc_2(V) = \frac{n}{r'}B(V') + \frac{n}{r''}B(V'') -
\frac{D^2}{r'r''},$$ and thus, as $B(V')\geq 0$ and $B(V'')\geq 0$ by
Bogomolov's inequality,
$$D^2\geq -(r'r'')2nc_2(V).$$ Since $r'+r'' = n$, $r'r''\leq n^2/4$, and
plugging this in to the above inequality proves the lemma.
\endproof

Returning to the proof of Theorem 7.4, the proof of Lemma 1.2 in Chapter 7
of \cite{6} (see also \cite{5}, Chapter 6, Lemma 3) shows that, if $\dsize
t\geq t_0 =
\frac{n^3}{4}c_2(V)$, then for every divisor $D$ such that
$\dsize D^2 \geq -\frac{n^3}{2}c_2(V)$ and $D\cdot f > 0$, we have $D\cdot
H_t> 0$. Now, if $V$ is not $H_{t_0}$-stable,  we would be able to find a
$t_1 \geq t_0$ and an exact sequence $0 \to V' \to V \to V'' \to 0$ as
above, with
$\mu_{H_{t_1}}(V') = \mu_{H_{t_1}}(V'')$. Now setting $D = r'c_1(V'') -
r''c_1(V')$ as before, we have
$$\gather 0< \mu_f(V) -\mu _f(V') = \frac{c_1(V')\cdot f + c_1(V'')\cdot
f}{n} -
\frac{c_1(V')\cdot f}{r'} \\ = \frac{\left(r'c_1(V'') -
r''c_1(V')\right)\cdot f}{r'n} = \frac{D\cdot f}{r'n},
\endgather$$ so that $D\cdot f > 0$, and likewise $D\cdot H_{t_1}= 0$. Thus
$D\cdot H_{t_0} <0$, contradicting the choice of $t_0$. It follows that,
for all $\dsize t\geq t_0 = \frac{n^3}{4}c_2(V)$, $V$ is $H_t$-stable. This
completes the proof of (7.5).
\endproof
 
As a final comment, the difficulty in finding an effective bound in case
$\dim B > 1$ is the following: For a torsion free sheaf $W$, we can define
$B(W)$ as before, but it is an element of $H^4(Z)$, not an integer. In the
notation of the proof of Lemma 7.5, Bogomolov's inequality can be used to
give a bound for
$B(V')\cdot H_t^{n-2}$ and $B(V'')\cdot H_t^{n-2}$ for some (unknown) value
of
$t$, and thus there is a lower bound for $D^2\cdot H_t^{n-2}$, again for
one  unknown value of $t$. However this does not seem to give enough
information to complete the proof of the theorem, except in the rank two
case where the lower bound can be made explicit for all $t$.

\enddocument